\let\chooseClass2   

\ifx\chooseClass0
\IfFileExists{extarticle.cls}{
\documentclass[14pt]{extarticle}
\emergencystretch 7 pt
\setlength{\oddsidemargin}{-17 mm} 
\setlength{\textwidth}{192 mm}     
\setlength{\textheight}{246 mm}    
\setlength{\topmargin}{-31 mm}     
}{
\documentclass[12pt]{article}
\textheight = 8.5in
\textwidth 6.3in
\setlength{\oddsidemargin}{0mm}
\setlength{\topmargin}{-15 mm}
}{}

\makeatletter
\def\@seccntformat#1{\csname the#1\endcsname.\quad}
\renewcommand\section{\@startsection {section}{1}{\z@}%
                                   {-3.5ex \@plus -1ex \@minus -.2ex}%
                                   {2.3ex \@plus.2ex}%
                                   {\normalfont\large\bfseries}}
\renewcommand\subsection{\@startsection{subsection}{2}{\z@}%
                        {3.25ex plus 1ex minus .2ex}{-.5em}%
                        {\normalfont\normalsize\bfseries}}
\renewcommand\subsubsection{\@startsection{subsubsection}{3}{\z@}%
                        {3.25ex plus 1ex minus .2ex}{-.5em}%
                        {\normalfont\normalsize\it}}
\@addtoreset{equation}{section}
\makeatother

\usepackage{amsthm}
\newtheoremstyle{boldhead}
{\topsep}
{\topsep}
{\slshape}
{}
{\bfseries}
{.}
{ }
{\thmname{#1}\thmnumber{ #2}\thmnote{ (#3)}}

\newtheoremstyle{boldremark}
{\topsep}
{\topsep}
{\upshape}
{}
{\bfseries}
{.}
{ }
{\thmname{#1}\thmnumber{ #2}\thmnote{ (#3)}}

\swapnumbers

\theoremstyle{boldhead}
\newtheorem{theorem}[subsection]{Theorem}
\newtheorem{claim}[subsection]{Claim}
\newtheorem{corollary}[subsection]{Corollary}
\newtheorem{lemma}[subsection]{Lemma}
\newtheorem{proposition}[subsection]{Proposition}
\newtheorem{maintheorem}[subsection]{Main Theorem}

\theoremstyle{boldremark}
\newtheorem{definition}[subsection]{Definition}
\newtheorem{example}[subsection]{Example}
\newtheorem{examples}[subsection]{Examples}
\newtheorem{remark}[subsection]{Remark}
\newtheorem{problem}[subsection]{Problem}
\newtheorem{question}[subsection]{Question}
\newtheorem{warning}[subsection]{Warning}
\newtheorem{acknowledgement}[subsection]{Acknowledgements}

\fi

\ifx\chooseClass1
\documentclass[dvips]{birkmult}

\newtheorem{theorem}[subsection]{Theorem}

\newtheorem{corollary}[subsection]{Corollary}
\newtheorem{lemma}[subsection]{Lemma}
\newtheorem{proposition}[subsection]{Proposition}
\newtheorem{maintheorem}[subsection]{Main Theorem}

\theoremstyle{definition}
\newtheorem{definition}[subsection]{Definition}
\newtheorem{example}[subsection]{Example}
\newtheorem{examples}[subsection]{Examples}
\newtheorem{remark}[subsection]{Remark}
\newtheorem{problem}[subsection]{Problem}
\newtheorem{question}[subsection]{Question}
\newtheorem{warning}[subsection]{Warning}
\newtheorem{acknowledgement}[subsection]{Acknowledgements}

\fi

\ifx\chooseClass2
\documentclass[12pt]{article}
\makeatletter
\def\@seccntformat#1{\csname the#1\endcsname.\quad}
\renewcommand\section{\@startsection {section}{1}{\z@}%
                                   {-3.5ex \@plus -1ex \@minus -.2ex}%
                                   {2.3ex \@plus.2ex}%
                                   {\normalfont\large\bfseries}}
\renewcommand\subsection{\@startsection{subsection}{2}{\z@}%
                        {3.25ex plus 1ex minus .2ex}{-.5em}%
                        {\normalfont\normalsize\bfseries}}
\renewcommand\subsubsection{\@startsection{subsubsection}{3}{\z@}%
                        {3.25ex plus 1ex minus .2ex}{-.5em}%
                        {\normalfont\normalsize\it}}
\@addtoreset{equation}{section}
\makeatother

\usepackage{amsthm}
\newtheoremstyle{boldhead}
{\topsep}
{\topsep}
{\slshape}
{}
{\bfseries}
{.}
{ }
{\thmname{#1}\thmnumber{ #2}\thmnote{ (#3)}}

\newtheoremstyle{boldremark}
{\topsep}
{\topsep}
{\upshape}
{}
{\bfseries}
{.}
{ }
{\thmname{#1}\thmnumber{ #2}\thmnote{ (#3)}}

\swapnumbers

\theoremstyle{boldhead}
\newtheorem{theorem}[subsection]{Theorem}
\newtheorem{corollary}[subsection]{Corollary}
\newtheorem{lemma}[subsection]{Lemma}
\newtheorem{proposition}[subsection]{Proposition}

\theoremstyle{boldremark}
\newtheorem{definition}[subsection]{Definition}
\newtheorem{example}[subsection]{Example}

\newtheorem{remark}[subsection]{Remark}

\fi

\usepackage{    amsmath,
                amsfonts,
                amssymb,
                amsxtra,
                varioref,
}

\numberwithin{equation}{section}

\IfFileExists{url.sty}{\usepackage{url}}{}
\providecommand{\url}[1]{{\tt #1}}

\usepackage{ifpdf}
\ifpdf
    \IfFileExists{hyperref.sty}{\usepackage[pdftex]{hyperref}}{}
\else
    \IfFileExists{hyperref.sty}{\usepackage[hypertex]{hyperref}}{}
\fi

\ifx\chooseClass1
\else
    \ifpdf
        \usepackage[pdftex]{graphics}
    \else
        \usepackage[dvips]{graphics}
    \fi
\fi

\IfFileExists{euscript.sty}{\usepackage[mathcal]{euscript}}{}
\IfFileExists{mathrsfs.sty}{\usepackage{mathrsfs}}{\let\mathscr\mathfrak}

\usepackage{diagrams}
\diagramstyle[height=2em,balance,righteqno,PostScript=dvips,nohug]

\IfFileExists{dsfont.sty}%
{\usepackage[sans]{dsfont}\newcommand\1{{\mathds 1}}}%
{\newcommand\1{{1\mkern-5mu {\mathrm I}}}}

\newlength{\mylabelwidths}
{\end{list}}

\def\rhaha{\raise.24ex\hbox{$\rightharpoonup$}\kern-1em\lower.24ex\hbox{$\rightharpoondown$}}%
\def\lhaha{\raise.24ex\hbox{$\leftharpoonup$}\kern-1em\lower.24ex\hbox{$\leftharpoondown$}}%
\def\dhaha{\downharpoonleft\kern-.22em\downharpoonright\kern.02em}%
\def\uhaha{\upharpoonleft\kern-.22em\upharpoonright\kern.02em}
\newarrowhead{twoharpoons}\rhaha\lhaha\dhaha\uhaha

\newarrow{Id}===={twoharpoons}
\newarrow{Epi}----{triangle}
\newarrow{MapsTo}{mapsto}---{->}
\newarrow{Mono}{boldhook}---{->}
\newarrow{TTo}----{->}
\newarrow{Twoar}===={=>}
\newarrow{Equalsto}=====
\newarrow{DashTo}{}{dash}{}{dash}{->}
\newarrow{DashLine}{}{dash}{}{dash}{}
\newarrow{DotsTo}....{->}

\newcommand\NN{{\mathbb N}}
\newcommand\ZZ{{\mathbb Z}}

\newcommand{\ca}{{\mathcal A}}
\newcommand{\cb}{{\mathcal B}}
\newcommand{\cc}{{\mathcal C}}
\newcommand{\cd}{{\mathcal D}}
\newcommand{\ce}{{\mathcal E}}

\newcommand{\ci}{{\mathcal I}}

\newcommand{\ck}{{\mathcal K}}
\newcommand{\cm}{{\mathcal M}}

\newcommand{\cp}{{\mathcal P}}
\newcommand{\cq}{{\mathcal Q}}
\newcommand{\cs}{{\mathcal S}}
\newcommand{\ct}{{\mathcal T}}

\newcommand{\cv}{{\mathcal V}}
\newcommand{\cw}{{\mathcal W}}

\newcommand{\cz}{{\mathcal Z}}

\newcommand{\Yo}{{\mathscr Y}}

\newcommand{\kf}{{\mathsf k}}

\newcommand{\SSS}{{\mathfrak S}}
\newcommand{\fu}{{\mathscr U}}

\newcommand{\sh}{{[\,]}}

\newcommand{\bull}{{\scriptscriptstyle\bullet}}

\newcommand{\Com}{{{\mathsf C}_\kk}}
\newcommand{\uCom}{{\underline{\mathsf C}_\kk}}

\newcommand{\tdt}{\otimes\dots\otimes}

\newcommand{\uni}{{\mathbf i}}
\newcommand{\unix}{{\sS{_X}\uni^\ca_0}}

\newcommand{\lar}[1]{\hbox{\,\large #1\,}}
\newcommand{\sS}[2]{\vphantom{#2}#1 #2}
\newcommand{\n}[1]{\nobreakdash-\hspace{0pt}}
\newcommand{\ainf}[1]{$A_\infty$\nobreakdash-\hspace{0pt}}
\newcommand{\ainfu}[1]{$A_\infty^u$\nobreakdash-\hspace{0pt}}

\newcommand{\Cat}{{\mathcal C}at}

\newcommand{\kCat}{\kk\text-{\mathcal C}at}
\newcommand{\KCat}{{\mathcal K}\text-{\mathcal C}at}
\newcommand{\Rep}{\mathrm{Rep}}

\let\boxt\boxtimes

\let\emptyset\varnothing
\let\eps\varepsilon
\let\ge\geqslant
\let\kk\Bbbk
\let\le\leqslant

\let\mb\mathbf
\let\msf=\mathsf

\let\tens\otimes

\let\und\underline
\let\wh\widehat

\newcommand\cA{{\mathcal A}}\newcommand\cB{{\mathcal B}}\newcommand\cC{{\mathcal C}}\newcommand\cP{{\mathcal P}}
\newcommand\cR{{\mathcal R}}

\DeclareMathOperator\Ainfty{{\mathsf A}_\infty}
\DeclareMathOperator\Ainftyu{{\mathsf A}_\infty^{\mathrm u}}

\DeclareMathOperator\kMod{\kk\text-\mathbf{Mod}}
\DeclareMathOperator\gr{\mathbf{gr}}

\DeclareMathOperator\dg{\mathbf{dg}}

\newcommand{\dgCat}{{\mathbf{dg}\text-\Cat}}
\newcommand{\grCat}{{\mathbf{gr}\text-\Cat}}

\DeclareMathOperator\Set{\mathcal Set}

\newcommand\cQuiver{{\mathscr Q}}

\newcommand{\dQ}{{{}^d\mkern-5mu\cQuiver}}

\DeclareMathOperator\sspan{span}

\DeclareMathOperator\bicomod{-bicomod}
\DeclareMathOperator\bimod{-bimod}

\DeclareMathOperator\coev{coev}

\DeclareMathOperator\ev{ev}

\DeclareMathOperator\Hom{Hom}

\DeclareMathOperator\id{id}
\DeclareMathOperator\Id{Id}
\DeclareMathOperator\im{Im}

\DeclareMathOperator\inj{in}

\DeclareMathOperator\modul{-mod}

\DeclareMathOperator\Ob{Ob}
\newcommand{\op}{{\operatorname{op}}}
\DeclareMathOperator\Par{Par}
\DeclareMathOperator\pr{pr}

\DeclareMathOperator\src{src}

\DeclareMathOperator\tgt{tgt}

\DeclareMathOperator\vect{-vect}

\newcommand{\appref}[1]{Appendix~\ref{#1}}

\newcommand{\corref}[1]{Corollary~\ref{#1}}
\newcommand{\defref}[1]{Definition~\ref{#1}}

\newcommand{\lemref}[1]{Lemma~\ref{#1}}
\newcommand{\propref}[1]{Proposition~\ref{#1}}
\newcommand{\remref}[1]{Remark~\ref{#1}}
\newcommand{\secref}[1]{Section~\ref{#1}}
\newcommand{\thmref}[1]{Theorem~\ref{#1}}

\newlength{\texthigh}
\newlength{\textwids}
\newcommand{\VCat}{{\mathcal V}\text-{\mathcal C}at}
\newcommand{\WCat}{{\mathcal W}\text-{\mathcal C}at}
\newcommand{\bT}{{\bar{T}}}
\DeclareMathOperator{\perm}{perm}

\newcommand{\straightForward}{This is proven by a straightforward computation}
\newcommand{\proofInArXiv}{The full proof is given in archive version
\cite{LyuMan-AmodSerre} of this article}

\begin{document}
\bibliographystyle{amsalpha}
\title{$A_\infty$-bimodules and Serre $A_\infty$-functors}
\ifx\chooseClass0
\author{Volodymyr Lyubashenko%
\thanks{Institute of Mathematics,
National Academy of Sciences of Ukraine,
3 Tereshchenkivska st.,
Kyiv-4, 01601 MSP,
Ukraine;
lub@imath.kiev.ua}
\ and Oleksandr Manzyuk%
\thanks{Fachbereich Mathematik,
Postfach 3049,
67653 Kaiserslautern,
Germany;
manzyuk@mathematik.uni-kl.de}
}
\fi

\ifx\chooseClass1
\author{Volodymyr Lyubashenko}
\address{Institute of Mathematics NASU\\
3 Tereshchenkivska st.\\
Kyiv-4, 01601 MSP\\
Ukraine}
\email{lub@imath.kiev.ua}

\author{Oleksandr Manzyuk}
\address{Fachbereich Mathematik\\
Postfach 3049\\
67653 Kaiserslautern\\
Germany}
\email{manzyuk@mathematik.uni-kl.de}

\date{December 11, 2006}
\dedicatory{To the memory of Oleksandr Reznikov}
\fi

\ifx\chooseClass2
\author{Volodymyr Lyubashenko%
\thanks{Institute of Mathematics,
National Academy of Sciences of Ukraine,
3 Tereshchenkivska st.,
Kyiv-4, 01601 MSP,
Ukraine;
lub@imath.kiev.ua}
\ and Oleksandr Manzyuk%
\thanks{Fachbereich Mathematik,
Postfach 3049,
67653 Kaiserslautern,
Germany;
manzyuk@mathematik.uni-kl.de}
}
\fi

\ifx\chooseClass1
\else
\maketitle
\fi

\begin{abstract}
We define $A_\infty$\n-bimodules similarly to Tradler and show that
this notion is equivalent to an $A_\infty$\n-functor with two arguments
which takes values in the differential graded category of complexes of
$\Bbbk$\n-modules, where $\Bbbk$ is a ground commutative ring. Serre
$A_\infty$\n-functors are defined via $A_\infty$\n-bimodules likewise
Kontsevich and Soibelman. We prove that a unital closed under shifts
$A_\infty$\n-category $\mathcal A$ over a field $\Bbbk$ admits a Serre
$A_\infty$\n-functor if and only if its homotopy category
$H^0\mathcal{A}$ admits a Serre $\Bbbk$\n-linear functor. The proof
uses categories enriched in $\mathcal K$, the homotopy category of
complexes of $\Bbbk$\n-modules, and Serre $\mathcal K$\n-functors. Also
we use a new $A_\infty$\n-version of the Yoneda Lemma generalizing the
previously obtained result.
\end{abstract}

\ifx\chooseClass1
\maketitle
\else
\fi

\allowdisplaybreaks[1]

Serre--Grothendieck duality for coherent sheaves on a smooth projective
variety was reformulated by Bondal and Kapranov in terms of Serre
functors \cite{MR1039961}. Being an abstract category theory notion
Serre functors were discovered in other contexts as well, for instance,
in Kapranov's studies of constructible sheaves on stratified spaces
\cite{MR1069414}. Reiten and van den Bergh showed that Serre functors
in categories of modules are related to Auslander--Reiten sequences and
triangles \cite{MR1887637}.

Often Serre functors are considered in triangulated categories and it
is reasonable to lift them to their origin -- pretriangulated
$\dg$\n-categories or \ainf-categories. Soibelman defines Serre
\ainf-functors in \cite{MR2095670}, based on Kontsevich and Soibelman
work which is a sequel to \cite{math.RA/0606241}. In the present
article we consider Serre \ainf-functors in detail. We define them via
\ainf-bimodules in \secref{sec-Serre-A8-functors} and use enriched
categories to draw conclusions about existence of Serre \ainf-functors.

\ainf-modules over \ainf-algebras are introduced by
Keller~\cite{math.RA/9910179}. \ainf-bimodules over \ainf-algebras are
defined by Tradler~\cite{xxx0108027,math.QA/0210150}. \ainf-modules and
\ainf-bimodules over \ainf-categories over a field were first
considered by Lef\`evre-Hasegawa~\cite{Lefevre-Ainfty-these} under the
name of polydules and bipolydules. \ainf-modules over
\ainf-categories were developed further by Keller~\cite{math.RT/0510508}.
We study \ainf-bimodules over
\ainf-categories over a ground commutative ring $\kk$ in
\secref{sec-A8-bimodules} and show that this notion is equivalent to an
\ainf-functor with two arguments which takes values in the
$\dg$\n-category $\uCom$ of complexes of $\kk$\n-modules. A similar
notion of \ainf-modules over an \ainf-category $\cc$ from
\secref{sec-A8-modules} is equivalent to an \ainf-functor
\(\cc\to\uCom\). The latter point of view taken by
Seidel~\cite{SeidelP-book-Fukaya} proved useful for ordinary and
differential graded categories as well, see Drinfeld's
article~\cite[Appendix~III]{Drinf:DGquot}.

Any unital \ainf-category $\ca$ determines a $\ck$\n-category $\kf\ca$
\cite{Lyu-AinfCat,BesLyuMan-book}, where $\ck$ is the homotopy category
of complexes of $\kk$\n-modules. Respectively, an \ainf-functor $f$
determines a $\ck$\n-functor $\kf f$. In particular, a Serre
\ainf-functor \(S:\ca\to\ca\) determines a Serre $\ck$\n-functor
 \(\kf S:\kf\ca\to\kf\ca\). We prove also the converse: if $\kf\ca$
admits a Serre $\ck$\n-functor, then $\ca$ admits a Serre \ainf-functor
(\corref{cor-ainf-Serre-exists-iff-K-Serre}). This shows the importance
of enriched categories in the subject.

Besides enrichment in $\ck$ we consider in
\secref{sec-Serre-functor-K-Cat} also categories enriched in the
category $\gr$ of graded $\kk$\n-modules. When $\kk$ is a field, we
prove that a Serre $\ck$\n-functor exists in $\kf\ca$ if and only if
the cohomology $\gr$\n-category
\(H^\bullet\ca\overset{\text{def}}=H^\bullet(\kf\ca)\) admits a Serre
$\gr$\n-functor (\corref{cor-Serre-Hbullet-Serre}). If the
$\gr$\n-category \(H^\bullet\ca\) is closed under shifts, then it
admits a Serre $\gr$\n-functor if and only if the $\kk$\n-linear
category \(H^0\ca\) admits a Serre $\kk$\n-linear functor
(\corref{cor-gr-Serre-k-Serre}, \propref{pro-Serre0-gr-Serre}). Summing
up, a unital closed under shifts \ainf-category \(\ca\) over a field
$\kk$ admits a Serre \ainf-functor if and only if its homotopy category
\(H^0\ca\) admits a Serre $\kk$\n-linear functor
(\thmref{thm-Serre-A8-k-functor}). This applies, in particular, to a
pretriangulated \ainf-enhancement $\ca$ of a triangulated category
\(H^0\ca\) over a field $\kk$.

In the proofs we use a new \ainf-version of the Yoneda Lemma
(\thmref{thm-Yoneda-Lemma}). It generalizes the previous result that
the Yoneda \ainf-functor is homotopy full and faithful
\cite[Theorem~9.1]{Fukaya:FloerMirror-II},
\cite[Theorem~A.11]{math.CT/0306018}, as well as a result of Seidel
\cite[Lemma~2.12]{SeidelP-book-Fukaya} which was proven over a ground
field $\kk$. The proof of the Yoneda Lemma occupies
\appref{sec-The-Yoneda-Lemma}. It is based on the theory of
\ainf-bimodules developed in \secref{sec-A8-bimodules}.

\subsection*{Acknowledgment}
The second author would like to thank the Mathematisches
For\-schungs\-in\-sti\-tut Oberwolfach and its director Prof. Dr.
Gert-Martin Greuel personally for financial support. Commutative
diagrams were typeset with the package \texttt{diagrams.sty} by Paul
Taylor.

\subsection{Notation and conventions.}
Notation follows closely the usage of the book \cite{BesLyuMan-book}.
In particular, $\fu$ is a ground universe containing an element which
is an infinite set, and $\kk$ denotes a $\fu$\n-small commutative
associative ring with unity. A \emph{graded quiver} $\cc$ typically
means a $\fu$\n-small set of objects $\Ob\cc$ together with
$\fu$\n-small $\ZZ$\n-graded $\kk$\n-modules of morphisms \(\cc(X,Y)\),
given for each pair \(X,Y\in\Ob\cc\). For any graded $\kk$\n-module $M$
there is another graded $\kk$\n-module $sM=M[1]$, its
\emph{suspension}, with the shifted grading \((sM)^k=M[1]^k=M^{k+1}\).
The mapping
 \(s:M\to sM\) given by the identity maps \(M^k\rId M[1]^{k-1}\) has
degree $-1$. The composition of maps, morphisms, functors, etc. is
denoted $fg=f\cdot g=\rTTo^f \rTTo^g =g\circ f$. A function (or a
functor) $f:X\to Y$ applied to an element is denoted
$f(x)=xf=x.f=x\bull f$ and occasionally $fx$.

Objects of the (large) Abelian $\fu$\n-category $\Com$ of complexes of
\(\kk\)\n-modules are $\fu$\n-small differential graded
$\kk$\n-modules. Morphisms of $\Com$ are chain maps. The category
$\Com$ is symmetric closed monoidal: the inner hom-complexes
\(\uCom(X,Y)\) are $\fu$\n-small, therefore, objects of $\Com$. This
determines a (large) differential graded $\fu$\n-category $\uCom$. In
particular, $\uCom$ is a non-small graded $\fu$\n-quiver.

Speaking about a symmetric monoidal category \((\cc,\tens,\1,c)\) we
actually mean the equivalent notion of a symmetric Monoidal category
\((\cc,\tens^I,\lambda^f)\)
 \cite[Definitions 1.2.2, 1.2.14]{Lyu-sqHopf},
\cite[Chapter~3]{BesLyuMan-book}. It is equipped with tensor product
functors \(\tens^I:(X_i)_{i\in I}\mapsto\tens^{i\in I}X_i\), where $I$
are finite linearly ordered (index) sets. The isomorphisms
\(\lambda^f:\tens^{i\in I}X_i\to\tens^{j\in J}\tens^{i\in f^{-1}j}X_i\)
given for any map $f:I\to J$ can be thought of as constructed from the
associativity isomorphisms $a$ and commutativity isomorphisms $c$. When
$f$ is non-decreasing, the isomorphisms $\lambda^f$ can be ignored
similarly to associativity isomorphisms. The coherence principle of
\cite[Section~3.25]{BesLyuMan-book} allows to write down canonical
isomorphisms $\omega_c$ (products of $\lambda^f$'s and their inverses)
between iterated tensor products, indicating only the permutation of
arguments $\omega$. One of them,
 \(\sigma_{(12)}:\tens^{i\in I}\tens^{j\in J}X_{ij}
 \to\tens^{j\in J}\tens^{i\in I}X_{ij}\)
is defined explicitly in \cite[(3.28.1)]{BesLyuMan-book}. Sometimes,
when the permutation of arguments reads clearly, we write simply
$\perm$ for the corresponding canonical isomorphism.

A symmetric multicategory $\wh\cc$ is associated with a lax symmetric
Monoidal category \((\cc,\tens^I,\lambda^f)\), where natural
transformations $\lambda^f$ are not necessarily invertible, see
\cite[Section~4.20]{BesLyuMan-book}.

The category of graded $\kk$\n-linear quivers has a natural symmetric
Monoidal structure. For given quivers $\cq_i$ the tensor product quiver
\(\boxt^{i\in I}\cq_i\) has the set of objects
 \(\prod_{i\in I}\Ob\cq_i\) and the graded $\kk$\n-modules of morphisms
 \((\boxt^{i\in I}\cq_i)((X_i)_{i\in I},(Y_i)_{i\in I})
 =\tens^{i\in I}\cq_i(X_i,Y_i)\).

For any graded quiver $\cc$ and a sequence of objects
\((X_0,\dots,X_n)\) of $\cc$ we use in this article the notation
\begin{align*}
\bT^ns\cc(X_0,\dots,X_n) &=s\cc(X_0,X_1)\tdt s\cc(X_{n-1},X_n),
\\
T^ns\cc(X_0,X_n)
&=\oplus_{X_1,\dots,X_{n-1}\in\Ob\cc}\bT^ns\cc(X_0,\dots,X_n).
\end{align*}
When the list of arguments is obvious we abbreviate the notation
\(\bT^ns\cc(X_0,\dots,X_n)\) to \(\bT^ns\cc(X_0,X_n)\). For $n=0$ we
set \(T^0s\cc(X,Y)=\kk\) if $X=Y$, and $0$ otherwise. The tensor quiver
is \(Ts\cc=\oplus_{n\ge0}T^ns\cc=\oplus_{n\ge0}(s\cc)^{\tens n}\).

An \emph{\ainf-category} means a graded quiver $\cc$ with $n$\n-ary
compositions \(b_n:T^ns\cc(X_0,X_n)\to s\cc(X_0,X_n)\) of degree 1
given for all $n\ge1$ (we assume for simplicity that $b_0=0$) such that
$b^2=0$ for the $\kk$\n-linear map $b:Ts\cc\to Ts\cc$ of degree 1
\begin{equation}
b = \sum_{\substack{r+n+t=k\\r+1+t=l}}
1^{\tens r}\tens b_n\tens1^{\tens t}: T^ks\cc \to T^ls\cc.
\label{eq-b-1b1-TkC-TlC}
\end{equation}
The composition \(b\cdot\pr_1=(0,b_1,b_2,\dots):Ts\cc\to s\cc\) is
denoted \(\check{b}\).

The tensor quiver $T\cc$ becomes a counital coalgebra when equipped
with the cut comultiplication
 $\Delta_0:T\cc(X,Y)\to
 \oplus_{Z\in\Ob\cc}T\cc(X,Z)\bigotimes_\kk T\cc(Z,Y)$,
 $h_1\tens h_2\tdt h_n\mapsto \sum_{k=0}^n
 h_1\tdt h_k\bigotimes h_{k+1}\tdt h_n$.
The map $b$ given by \eqref{eq-b-1b1-TkC-TlC} is a coderivation with
respect to this comultiplication. Thus $b$ is a \emph{codifferential}.

We denote by $\mb n$ the linearly ordered index set
\(\{1<2<\dots<n\}\).

An \emph{\ainf-functor} $f:\ca\to\cb$ is a map of objects
 $f=\Ob f:\Ob\ca\to\Ob\cb$, \(X\mapsto Xf\) and $\kk$\n-linear maps
$f:Ts\ca(X,Y)\to Ts\cb(Xf,Yf)$ of degree 0 which agree with the cut
comultiplication and commute with the codifferentials $b$. Such $f$ is
determined in a unique way by its components
\(f_k=f\pr_1:T^ks\ca(X,Y)\to s\cb(Xf,Yf)\), $k\ge1$ (we require that
$f_0=0$). This generalizes to the case of \ainf-functors with many
arguments \(f:(\cA_i)_{i\in\mb{n}}\to\cB\). Such \emph{\ainf-functor}
is a quiver map \(f:\boxt^{i\in\mb n}Ts\ca_i\to Ts\cb\) of degree 0
which agrees with the cut comultiplication and commutes with the
differentials. Denote
 \(\check{f}=f\cdot\pr_1:\boxt^{i\in\mb n}Ts\ca_i\to s\cb\). The
restrictions \(f_{(k_i)_{i\in\mb n}}\) of \(\check{f}\) to
\(\boxt^{i\in\mb n}T^{k_i}s\ca_i\) are called the components of $f$. It
is required that the restriction \(f_{00\dots0}\) of \(\check{f}\) to
\(\boxt^{i\in\mb n}T^0s\ca_i\) vanishes. The components determine
coalgebra homomorphism $f$ in a unique way. Commutation with the
differentials means that the following compositions are equal
\begin{equation*}
\bigl( \boxt^{i\in\mb n}Ts\ca_i \xrightarrow{f} Ts\cb
\xrightarrow{\check{b}} s\cb \bigr)
 =\bigl( \boxt^{i\in\mb n}Ts\ca_i
\rTTo^{\sum_{i=1}^n1^{\boxt(i-1)}\boxt b\boxt1^{\boxt(n-i)}}
\boxt^{i\in\mb n}Ts\ca_i \xrightarrow{\check{f}} s\cb \bigr).
\end{equation*}
The set of \ainf-functors \((\ca_i)_{i\in\mb n}\to\cb\) is denoted by
\(\Ainfty((\ca_i)_{i\in\mb n};\cb)\). There is a natural way to compose
\ainf-functors, the composition is associative, and for an arbitrary
\ainf-category \(\ca\), the identity \ainf-functor
\(\id_\ca:\ca\to\ca\) is the unit with respect to the composition.
Thus, \ainf-categories and \ainf-functors constitute a symmetric
multicategory \(\Ainfty\) \cite[Chapter~12]{BesLyuMan-book}.

With a family \((\ca_i)_{i=1}^n,\cb\) of \ainf-categories we associate
a graded quiver \(\und{\Ainfty}((\ca_i)_{i=1}^n;\cb)\). Its objects are
\ainf-functors with $n$ entries. Morphisms are \ainf-transformations
between such \ainf-functors $f$ and $g$, that is,
$(f,g)$\n-coderivations. Such coderivation $r$ can be identified with
the collection of its components $\check{r}=r\cdot\pr_1$, thus,
\begin{multline*}
s\und{\Ainfty}((\ca_i)_{i=1}^n;\cb)(f,g) \simeq
\prod_{X,Y\in\prod_{i=1}^n\Ob\ca_i}
\uCom\bigl((\boxt^{i\in\mb n}Ts\ca_i)(X,Y),s\cb(Xf,Yg)\bigr)
\\
= \prod_{(X_i)_{i\in\mb n},(Y_i)_{i\in\mb n}\in\prod_{i=1}^n\Ob\ca_i}
\uCom\bigl(\tens^{i\in\mb n}[Ts\ca_i(X_i,Y_i)],
 s\cb((X_i)_{i\in\mb n}f,(Y_i)_{i\in\mb n}g)\bigr).
\end{multline*}
Moreover, \(\und{\Ainfty}((\ca_i)_{i=1}^n;\cb)\) has a distinguished
\ainf-category structure which together with the evaluation
\ainf-functor
\[ \ev^{\Ainfty}:
(\ca_i)_{i=1}^n,\und{\Ainfty}((\ca_i)_{i=1}^n;\cb) \to \cb, \qquad
(X_1,\dots,X_n,f) \mapsto (X_1,\dots,X_n)f
\]
turns the symmetric multicategory $\Ainfty$ into a closed multicategory
\cite{BesLyuMan-book}. Thus, for arbitrary \ainf-categories
\((\ca_i)_{i\in\mb n}\), \((\cb_j)_{j\in\mb m}\), \(\cc\), the map
\begin{align*}
\varphi^{\Ainfty}:
 \Ainfty((\cb_j)_{j\in\mb m};\und\Ainfty((\ca_i)_{i\in\mb n};\cc))
&\longrightarrow \Ainfty((\ca_i)_{i\in\mb n},(\cb_j)_{j\in\mb m};\cc),
\\
f &\longmapsto ((\id_{\ca_i})_{i\in\mb n},f)\ev^{\Ainfty}
\end{align*}
is bijective. It follows from the general properties of closed
multicategories that the bijection \(\varphi^{\Ainfty}\) extends
uniquely to an isomorphism of \ainf-categories
\[ \und{\varphi}^{\Ainfty}:
 \und\Ainfty((\cb_j)_{j\in\mb m};\und\Ainfty((\ca_i)_{i\in\mb n};\cc))
\to\und\Ainfty((\ca_i)_{i\in\mb n},(\cb_j)_{j\in\mb m};\cc)
\]
with \(\Ob\und{\varphi}^{\Ainfty}=\varphi^{\Ainfty}\). In particular,
if \(\cc\) is a unital \ainf-category, \(\varphi^{\Ainfty}\) maps
isomorphic \ainf-functors to isomorphic.

The components $\ev^{\Ainfty}_{k_1,\dots,k_n;m}$ of the evaluation
\ainf-functor vanish if $m>1$ by
formula~\cite[(12.25.4)]{BesLyuMan-book}. For $m=0,1$ they are
\begin{multline}
\ev^{\Ainfty}_{k_1,\dots,k_n;0}:
[\tens^{i\in\mb n}T^{k_i}s\ca_i(X_i,Y_i)]\tens
 T^0s\und{\Ainfty}((\ca_i)_{i=1}^n;\cb)(f,f)
\\
\hfill \rTTo^{f_{k_1,\dots,k_n}}
 s\cb((X_i)_{i\in\mb n}f,(Y_i)_{i\in\mb n}f), \quad
\\
\hskip\multlinegap \ev^{\Ainfty}_{k_1,\dots,k_n;1} =\bigl[
(\tens^{i\in\mb n}T^{k_i}s\ca_i(X_i,Y_i))\tens
 s\und{\Ainfty}((\ca_i)_{i=1}^n;\cb)(f,g)
\rTTo^{1\tens\pr_{k_1,\dots,k_n}} \hfill
\\
[\tens^{i\in\mb n}T^{k_i}s\ca_i(X_i,Y_i)]\tens
\uCom\bigl(\tens^{i\in\mb n}[T^{k_i}s\ca_i(X_i,Y_i)],
 s\cb((X_i)_{i\in\mb n}f,(Y_i)_{i\in\mb n}g)\bigr)
\\
\rTTo^{\ev^\Com} s\cb((X_i)_{i\in\mb n}f,(Y_i)_{i\in\mb n}g) \bigr].
\label{eq-evA80-evA81}
\end{multline}
When \((\ca_i)_{i=1}^n,\cb\) are unital \ainf-categories, we define
\(\und{\Ainftyu}((\ca_i)_{i=1}^n;\cb)\) as a full \ainf-subcategory of
\(\und{\Ainfty}((\ca_i)_{i=1}^n;\cb)\), whose objects are unital
\ainf-functors. Equipped with a similar evaluation \(\ev^{\Ainftyu}\)
the collection $\Ainftyu$ of unital \ainf-categories and unital
\ainf-functors also becomes a closed multicategory. Similarly to the
case of \(\Ainfty\), there is a natural bijection
\begin{align*}
\varphi^{\Ainftyu}:
 \Ainftyu((\cb_j)_{j\in\mb m};\und\Ainftyu((\ca_i)_{i\in\mb n};\cc))
&\longrightarrow \Ainftyu((\ca_i)_{i\in\mb n},(\cb_j)_{j\in\mb m};\cc),
\\
f &\longmapsto ((\id_{\ca_i})_{i\in\mb n},f)\ev^{\Ainftyu}
\end{align*}
for arbitrary unital \ainf-categories \((\ca_i)_{i\in\mb n}\),
\((\cb_j)_{j\in\mb m}\), \(\cc\).

In the simplest version graded spans $\cp$ consist of a $\fu$\n-small
set \(\Ob_s\cp\) of source objects, a $\fu$\n-small set \(\Ob_t\cp\) of
target objects, and $\fu$\n-small graded $\kk$\n-modules \(\cp(X,Y)\)
for all \(X\in\Ob_s\cp\), \(Y\in\Ob_t\cp\). Graded quivers $\ca$ are
particular cases of spans, distinguished by the condition
\(\Ob_s\ca=\Ob_t\ca\). The tensor product $\cp\tens\cq$ of two spans
$\cp$, $\cq$ exists if \(\Ob_t\cp=\Ob_s\cq\) and equals
\[ (\cp\tens\cq)(X,Z)
=\bigoplus_{Y\in\Ob_t\cp} \cp(X,Y)\tens_\kk\cq(Y,Z).
\]
Details can be found in \cite{BesLyuMan-book}.

Next we explain our notation for closed symmetric monoidal categories
which differs slightly from \cite[Chapter~1]{KellyGM:bascec}.

Let \((\cv,\tens,\1,c)\) be a closed symmetric monoidal $\fu$\n-category. For
each pair of objects \(X,Y\in\Ob\cv\), let \(\und\cv(X,Y)\) denote the
inner hom-object. Denote by \(\ev^\cv:X\tens\und\cv(X,Y)\to Y\) and
\(\coev^\cv:Y\to\und\cv(X,X\tens Y)\) the evaluation and coevaluation
morphisms, respectively. Then the mutually inverse adjunction
isomorphisms are explicitly given as follows:
\begin{align*}
\cv(Y,\und\cv(X,Z)) & \longleftrightarrow \cv(X\tens Y,Z), \\
f &\rMapsTo (1_X\tens f)\ev_{X,Z}, \\
\coev_{X,Y}\und\cv(X,g) &\lMapsTo g.
\end{align*}

There is a \(\cv\)\n-category \(\und\cv\) whose objects are those of
\(\cv\), and for each pair of objects \(X\) and \(Y\), the object
\(\und\cv(X,Y)\in\Ob\cv\) is the inner hom-object of \(\cv\). The
composition is found from the following equation:
\begin{multline}
\bigl[ X\tens\und\cv(X,Y)\tens\und\cv(Y,Z)
\rTTo^{1\tens\mu_{\und\cv}\;{}} X\tens\und\cv(X,Z) \rTTo^{\ev^\cv} Z
\bigr]
\\
=\bigl[ X\tens\und\cv(X,Y)\tens\und\cv(Y,Z) \rTTo^{\ev^\cv\tens1}
Y\tens\und\cv(Y,Z) \rTTo^{\ev^\cv} Z \bigr].
 \label{eq-1mevV-evV1evV}
\end{multline}
The identity morphism \(1^{\und\cv}_X:\1\to\und\cv(X,X)\) is found from the
following equation:
\[
\bigl[
X\rTTo^\sim X\tens\1\rTTo^{1\tens
1^{\und\cv}_X}X\tens\und\cv(X,X)\rTTo^{\ev^\cv}X
\bigr]=\id_X.
\]

For our applications we need several categories \(\cv\), for instance
the Abelian category $\Com$ of complexes of \(\kk\)\n-modules and its
quotient category \(\ck\), the homotopy category of complexes of
\(\kk\)\n-modules. The tensor product is the tensor product of
complexes, the unit object is \(\kk\), viewed as a complex concentrated
in degree \(0\), and the symmetry is the standard symmetry
 \(c:X\tens Y\to Y\tens X\), \(x\tens y\mapsto (-)^{xy}y\tens x\).
We shorten up the usual notation \((-1)^{\deg x\cdot\deg y}\) to
\((-)^{xy}\). Similarly, \((-)^x\) means \((-1)^{\deg x}\),
\((-)^{x+y}\) means \((-1)^{\deg x+\deg y}\), etc. For each pair of
complexes \(X\) and \(Y\), the inner hom-object \(\und\ck(X,Y)\) is the
same as the inner hom-complex \(\uCom(X,Y)\) in the symmetric closed
monoidal Abelian category $\Com$. The evaluation morphism
\(\ev^\ck:X\tens\und\ck(X,Y)\to Y\) and the coevaluation morphism
\(\coev^\ck:Y\to\und\ck(X,X\tens Y)\) in \(\ck\) are the homotopy
classes of the evaluation morphism \(\ev^\Com:X\tens\uCom(X,Y)\to Y\)
and the coevaluation morphism \(\coev^\Com:Y\to\uCom(X,X\tens Y)\) in
\(\Com\), respectively.

It is easy to see that \(\mu_{\und\ck}=m^\uCom_2\) and
\(1^{\und\ck}_X=1^\uCom_X\), therefore \(\und\ck=\kf\uCom\).

Also we use as \(\cv\) the category \(\gr=\gr(\kMod)\) of graded
\(\kk\)\n-modules, and the familiar category \(\kMod\) of
\(\kk\)\n-modules.

The following identity holds for any homogeneous $\kk$\n-linear map
$a:X\to A$ of arbitrary degree by the properties of the closed monoidal
category \(\uCom\):
\begin{multline}
\bigl(\uCom(A,B)\tens\uCom(B,C) \rTTo^{m_2} \uCom(A,C)
\rTTo^{\uCom(a,C)} \uCom(X,C)\bigr)
\\
= \bigl(\uCom(A,B)\tens\uCom(B,C) \rTTo^{\uCom(a,B)\tens1}
\uCom(X,B)\tens\uCom(B,C) \rTTo^{m_2} \uCom(X,C)\bigr).
\label{eq-m2C(aC)-C(aB)1m2}
\end{multline}
Let $f:A\tens X\to B$, $g:B\tens Y\to C$ be two homogeneous
$\kk$\n-linear maps of arbitrary degrees, that is,
 \(f\in\uCom(A\tens X,B)^\bull\), \(g\in\uCom(B\tens Y,C)^\bull\). Then
the following identity is proven in \cite{math.CT/0306018} as
equation~(A.1.2):
\begin{multline}\label{eq-identity-m2}
\bigl( X\tens Y \rTTo^{\coev_{A,X}\tens\coev_{B,Y}}
\uCom(A,A\tens X)\tens\uCom(B,B\tens Y)
\\
\hfill \rTTo^{\uCom(A,f)\tens\uCom(B,g)} \uCom(A,B)\tens\uCom(B,C)
\rTTo^{m_2} \uCom(A,C) \bigr) \quad
\\
= \bigl( X\tens Y \rTTo^{\coev_{A,X\tens Y}} \uCom(A,A\tens X\tens Y)
\rTTo^{\uCom(A,f\tens1)} \uCom(A,B\tens Y) \rTTo^{\uCom(A,g)}
\uCom(A,C) \bigr).
\end{multline}

\section{$\cv$-categories}
We refer the reader to \cite[Chapter~1]{KellyGM:bascec} for the basic
theory of enriched categories. The category of unital (resp.
non-unital) $\cv$\n-categories (where $\cv$ is a closed symmetric
monoidal category) is denoted $\VCat$ (resp. \(\VCat^{nu}\)).

\subsection{Opposite $\cv$-categories.}
Let \(\ca\) be a \(\cv\)\n-category, not necessarily unital. Its
opposite \(\ca^\op\) is defined in the standard way. Namely,
\(\Ob\ca^\op=\Ob\ca\), and for each pair of objects \(X,Y\in\Ob\ca\),
\(\ca^\op(X,Y)=\ca(Y,X)\). The composition in \(\ca^\op\) is given by
\begin{multline*}
\mu_{\ca^\op}=\bigl[\ca^\op(X,Y)\tens\ca^\op(Y,Z)=\ca(Y,X)\tens\ca(Z,Y)\rTTo^c\\
\ca(Z,Y)\tens\ca(Y,X)\rTTo^{\mu_\ca}\ca(Z,X)=\ca^\op(X,Z)\bigr].
\end{multline*}
More generally, for each \(n\ge1\), the iterated $n$\n-ary composition
in \(\ca^\op\) is
\begin{multline}
\mu^{\mb n}_{\ca^\op}=\bigl[
\tens^{i\in\mb n}\ca^\op(X_{i-1},X_i)=\tens^{i\in\mb
n}\ca(X_i,X_{i-1})\rTTo^{\omega^0_c}\\
\tens^{i\in\mb n}\ca(X_{n-i+1},X_{n-i})\rTTo^{\mu^{\mb
n}_\ca}\ca(X_n,X_0)=\ca^\op(X_0,X_n)
\bigr],
\label{equ-comp-opposite}
\end{multline}
where the permutation
 $\omega^0=\bigl(
\begin{smallmatrix}
1 & 2 & \dots & n-1 & n \\
n &n-1& \dots &  2  & 1
\end{smallmatrix}
 \bigr)$
is the longest element of $\SSS_n$, and $\omega^0_c$ is the
corresponding signed permutation, the action of $\omega^0$ in tensor
products via standard symmetry. Note that if \(\ca\) is unital, then so
is \(\ca^\op\), with the same identity morphisms.

Let \(f:\ca\to\cb\) be a \(\cv\)\n-functor, not necessarily unital. It
gives rise to a \(\cv\)\n-functor \(f^\op:\ca^\op\to\cb^\op\) with
\(\Ob f^\op=\Ob f:\Ob\ca\to\Ob\cb\), and
\[
f^\op=\bigl[\ca^\op(X,Y)=\ca(Y,X)\rTTo^f\cb(Yf,Xf)=\cb^\op(Xf,Yf)\bigr],\quad X,Y\in\Ob\ca.
\]
Note that \(f^\op\) is a unital \(\cv\)\n-functor if so is \(f\).
Clearly, the correspondences \(\ca\mapsto\ca^\op\), \(f\mapsto f^\op\)
define a functor \(-^\op:\VCat^{nu}\to\VCat^{nu}\) which restricts to a
functor \(-^\op:\VCat\to\VCat\). The functor \(-^\op\) is symmetric
Monoidal. More precisely, for arbitrary \(\cv\)\n-categories \(\ca_i\),
\(i\in\mb n\), the equation
 \(\boxt^{i\in\mb n}\ca^\op_i=(\boxt^{i\in\mb n}\ca_i)^\op\) holds.
Indeed, the underlying \(\cv\)\n-quivers of both categories coincide,
and so do the identity morphisms if the categories \(\ca_i\),
 \(i\in\mb n\), are unital. The composition in
\(\boxt^{i\in\mb n}\ca^\op_i\) is given by
\begin{multline*}
\mu_{\boxt^{i\in\mb n}\ca^\op_i}=\bigl[
 \bigl(\tens^{i\in\mb n}\ca^\op_i(X_i,Y_i)\bigr)\tens
 \bigl(\tens^{i\in\mb n}\ca^\op_i(Y_i,Z_i)\bigr)
\\
\rTTo^{\sigma_{(12)}}
 \tens^{i\in\mb n}(\ca_i(Y_i,X_i)\tens\ca_i(Z_i,Y_i))
\rTTo^{\tens^{i\in\mb n}c\;{}}
 \tens^{i\in\mb n}(\ca_i(Z_i,Y_i)\tens\ca_i(Y_i,X_i))
\\
\rTTo^{\tens^{i\in\mb n}\mu_{\ca_i}\;{}}
 \tens^{i\in\mb n}\ca_i(Z_i,X_i)=\tens^{i\in\mb n}\ca^\op_i(X_i,Z_i)
\bigr].
\end{multline*}
The composition in \((\boxt^{i\in\mb n}\ca_i)^\op\) is given by
\begin{multline*}
\mu_{(\boxt^{i\in\mb n}\ca_i)^\op}=\bigl[
(\boxt^{i\in\mb n}\ca_i)^\op((X_i)_{i\in\mb n},(Y_i)_{i\in\mb n})\tens
(\boxt^{i\in\mb n}\ca_i)^\op((Y_i)_{i\in\mb n},(Z_i)_{i\in\mb
n})\\
=\bigl(\tens^{i\in\mb n}\ca_i(Y_i,X_i)\bigr)\tens\bigl(\tens^{i\in\mb n}\ca_i(Z_i,Y_i)\bigr)
\rTTo^c
\bigl(\tens^{i\in\mb n}\ca_i(Z_i,Y_i)\bigr)\tens\bigl(\tens^{i\in\mb
n}\ca_i(Y_i,X_i)\bigr)\\
\rTTo^{\sigma_{(12)}}\tens^{i\in\mb
n}(\ca_i(Z_i,Y_i)\tens\ca_i(Y_i,X_i))
\rTTo^{\tens^{i\in\mb n}\mu_{\ca_i}}\tens^{i\in\mb n}\ca_i(Z_i,X_i)=\tens^{i\in\mb n}\ca^\op_i(X_i,Z_i)
\bigr].
\end{multline*}
The equation \(\mu_{\boxt^{i\in\mb n}\ca^\op_i}=\mu_{(\boxt^{i\in\mb
n}\ca_i)^\op}\) follows from the following equation in \(\cv\):
\begin{multline*}
\bigl[\bigl(\tens^{i\in\mb n}\ca_i(Y_i,X_i)\bigr)\tens\bigl(\tens^{i\in\mb
n}\ca_i(Z_i,Y_i)\bigr)
\rTTo^{\sigma_{(12)}}\tens^{i\in\mb
n}(\ca_i(Y_i,X_i)\tens\ca_i(Z_i,Y_i))\\
\hfill\rTTo^{\tens^{i\in\mb n}c\;{}}\tens^{i\in\mb
n}(\ca_i(Z_i,Y_i)\tens\ca_i(Y_i,X_i))\bigr]\quad\\
\quad=\bigl[\bigl(\tens^{i\in\mb n}\ca_i(Y_i,X_i)\bigr)\tens\bigl(\tens^{i\in\mb n}\ca_i(Z_i,Y_i)\bigr)
\rTTo^c
\bigl(\tens^{i\in\mb n}\ca_i(Z_i,Y_i)\bigr)\tens\bigl(\tens^{i\in\mb
n}\ca_i(Y_i,X_i)\bigr)\hfill\\
\rTTo^{\sigma_{(12)}}\tens^{i\in\mb
n}(\ca_i(Z_i,Y_i)\tens\ca_i(Y_i,X_i))\bigr],
\end{multline*}
which is a consequence of coherence principle of
\cite[Lemma~3.26, Remark~3.27]{BesLyuMan-book}. Therefore, \(-^\op\)
induces a symmetric multifunctor
\(-^\op:\wh{\VCat^{nu}}\to\wh{\VCat^{nu}}\) which restricts to a
symmetric multifunctor \(-^\op:\wh\VCat\to\wh\VCat\).

\subsection{$\Hom$-functor.}
A \(\cv\)\n-category \(\ca\) gives rise to a
\(\cv\)\n-functor \(\Hom_\ca:\ca^\op\boxt\ca\to\und\cv\) which maps a pair
of objects \((X,Y)\in\Ob\ca\times\Ob\ca\) to \(\ca(X,Y)\in\Ob\und\cv\), and
whose action on morphisms is given by
\begin{multline*}
\Hom_\ca=\bigl[
\ca^\op(X,Y)\tens\ca(U,V)=\ca(Y,X)\tens\ca(U,V)\rTTo^{\coev^\cv}
\\
\hfill\und\cv(\ca(X,U),\ca(X,U)\tens\ca(Y,X)\tens\ca(U,V))
\rTTo^{\und\cv(1,(c\tens1)\mu^{\mb3}_\ca)}
\und\cv(\ca(X,U),\ca(Y,V))
\bigr].
\end{multline*}
Equivalently, \(\Hom_\ca\) is found by closedness of \(\cv\) from the diagram
\begin{diagram}[LaTeXeqno]
\ca(X,U)\tens\ca(Y,X)\tens\ca(U,V) & \rTTo^{1\tens\Hom_\ca} &
\ca(X,U)\tens\und\cv(\ca(X,U),\ca(Y,V))\\
\dTTo<{c\tens1} && \dTTo>{\ev^\cv}\\
\ca(Y,X)\tens\ca(X,U)\tens\ca(U,V) & \rTTo^{\mu^{\mb 3}_\ca} & \ca(Y,V)
\label{equ-Hom-A}
\end{diagram}

\begin{lemma}
Let \(\ca\) be a \(\cv\)\n-category. Then
\[
\Hom_{\ca^\op}=\bigl[\ca\boxt\ca^\op\rTTo^c\ca^\op\boxt\ca\rTTo^{\Hom_\ca}\und\cv\bigr].
\]
\end{lemma}

\begin{proof}
Using~\eqref{equ-comp-opposite}, we obtain:
\begin{multline*}
\Hom_{\ca^\op}=\bigl[\ca(X,Y)\tens\ca^\op(U,V)=\ca^\op(Y,X)\tens\ca^\op(U,V)
\\
\rTTo^{\coev^\cv}
\und\cv(\ca^\op(X,U),\ca^\op(X,U)\tens\ca^\op(Y,X)\tens\ca^\op(U,V))
\\
\hfill \rTTo^{\und\cv(1,(c\tens1)\mu^{\mb3}_{\ca^\op})}
\und\cv(\ca^\op(X,U),\ca^\op(Y,V))\bigr] \quad
\\
\quad =\bigl[
\ca(X,Y)\tens\ca(V,U)\rTTo^{\coev^\cv}\und\cv(\ca(U,X),\ca(U,X)\tens\ca(X,Y)\tens\ca(V,U))
\hfill
\\
\rTTo^{\und\cv(1,(c\tens1)\omega^0_c\mu^{\mb3}_{\ca})}
\und\cv(\ca(U,X),\ca(V,Y)) \bigr],
\end{multline*}
where \(\omega^0=(13)\in\SSS_3\). Clearly,
\((c\tens1)\omega^0_c=(1\tens c)(c\tens1)\), therefore
\begin{multline*}
\Hom_{\ca^\op}=\bigl[
\ca(X,Y)\tens\ca(V,U)\rTTo^{\coev^\cv}\und\cv(\ca(U,X),\ca(U,X)\tens\ca(X,Y)\tens\ca(V,U))
\\
\hfill \rTTo^{\und\cv(1,(1\tens c)(c\tens1)\mu^{\mb3}_\ca)}
\und\cv(\ca(U,X),\ca(V,Y)) \bigr] \quad
\\
=\bigl[ \ca(X,Y)\tens\ca(V,U) \rTTo^c \ca(V,U)\tens\ca(X,Y)
\rTTo^{\coev^\cv} \hfill
\\
\hfill \und\cv(\ca(U,X),\ca(U,X)\tens\ca(V,U)\tens\ca(X,Y))
\rTTo^{\und\cv(1,(c\tens1)\mu^{\mb3}_\ca)} \und\cv(\ca(U,X),\ca(V,Y))
\bigr] \quad
\\
=\bigl[ \ca(X,Y)\tens\ca(V,U)\rTTo^c\ca(V,U)\tens\ca(X,Y)
\rTTo^{\Hom_\ca} \und\cv(\ca(U,X),\ca(V,Y)) \bigr].\hfill
\end{multline*}
The lemma is proven.
\end{proof}

An object \(X\) of \(\ca\) defines a \(\cv\)\n-functor \(X:\1\to\ca\),
\(*\mapsto X\), \(\1(*,*)=\1\rTTo^{1^\ca_X}\ca(X,X)\), whose source
$\1$ is a $\cv$\n-category with one object $*$. This $\cv$\n-category
is a unit of tensor multiplication $\boxt$. The \(\cv\)\n-functors
\(\ca(\_,Y)=\Hom_\ca(\_,Y):\ca^\op\to\und\cv\) and
\(\ca(X,\_)=\Hom_\ca(X,\_):\ca\to\und\cv\) are defined as follows:
\begin{gather*}
\ca(\_,Y)=\bigl[
\ca^\op\rTTo^\sim\ca^\op\boxt\1\rTTo^{1\boxt
Y}\ca^\op\boxt\ca\rTTo^{\Hom_\ca}\und\cv
\bigr],\\
\ca(X,\_)=\bigl[
\ca\rTTo^\sim\1\boxt\ca\rTTo^{X\boxt
1}\ca^\op\boxt\ca\rTTo^{\Hom_\ca}\und\cv
\bigr].
\end{gather*}
Thus, the \(\cv\)\n-functor \(\ca(\_,Y)\) maps an object \(X\) to
\(\ca(X,Y)\), and its action on morphisms is given by
\begin{multline}
\ca(\_,Y)=\bigl[
\ca^\op(W,X)=\ca(X,W)\rTTo^{\coev^\cv}\und\cv(\ca(W,Y),\ca(W,Y)\tens\ca(X,W))\\
\rTTo^{\und\cv(1,c\mu_\ca)}\und\cv(\ca(W,Y),\ca(X,Y))
\bigr].
\label{equ-A-dot-Y}
\end{multline}
Similarly, the \(\cv\)\n-functor \(\ca(X,\_)\) maps an object \(Y\)
to \(\ca(X,Y)\), and its action on morphisms is given by
\begin{multline}
\ca(X,\_)=\bigl[ \ca(Y,Z) \rTTo^{\coev^\cv}
\und\cv(\ca(X,Y),\ca(X,Y)\tens\ca(Y,Z))
\\
\rTTo^{\und\cv(1,\mu_\ca)} \und\cv(\ca(X,Y),\ca(X,Z)) \bigr].
\label{equ-A-X-dot}
\end{multline}

\subsection{Duality functor.}
The unit object \(\1\) of \(\cv\) defines the duality \(\cv\)\n-functor
\(\und\cv(\_,\1)=\Hom_{\und\cv}(\_,\1):\und{\cv}^\op\to\und\cv\).
The functor \(\und\cv(\_,\1)\) maps an object \(M\) to its dual
\(\und\cv(M,\1)\), and its action on morphisms is given by
\begin{multline}
\und\cv(\_,\1)=\bigl[
\und\cv^\op(M,N)=\und\cv(N,M)\rTTo^{\coev^\cv}\und\cv(\und\cv(M,\1),\und\cv(M,\1)\tens\und\cv(N,M))
\rTTo^{\und\cv(1,c)}\\
\hfill\und\cv(\und\cv(M,\1),\und\cv(N,M)\tens\und\cv(M,\1))
\rTTo^{\und\cv(1,\mu_{\und\cv})}
\und\cv(\und\cv(M,\1),\und\cv(N,\1))
\bigr].
\label{eq-V(-1)=coev-V(1cmu)}
\end{multline}

For each object \(M\) there is a natural morphism
\(e:M\to\und\cv(\und\cv(M,\1),\1)\) which is a unique solution of the
following equation in \(\cv\):
\begin{diagram}
\und\cv(M,\1)\tens M & \rTTo^{\quad c\quad} & M\tens\und\cv(M,\1)\\
\dTTo<{1\tens e} && \dTTo>{\ev^\cv}\\
\und\cv(M,\1)\tens\und\cv(\und\cv(M,\1),\1) & \rTTo^{\quad\ev^\cv\quad} &
\1
\end{diagram}
Explicitly,
\begin{multline*}
e=\bigl[
M\rTTo^{\coev^\cv}\und\cv(\und\cv(M,\1),\und\cv(M,\1)\tens
M)\rTTo^{\und\cv(1,c)}\\
\und\cv(\und\cv(M,\1),M\tens\und\cv(M,\1))\rTTo^{\und\cv(1,\ev^\cv)}\und\cv(\und\cv(M,\1),\1)
\bigr].
\end{multline*}
An object \(M\) is \emph{reflexive} if \(e\) is an isomorphism in
\(\cv\).

\subsection{Representability.}
Let us state for the record the following

\begin{proposition}[Weak Yoneda Lemma]\label{prop-repr-K-functors}
Let \(F:\ca\to\und\cv\) be a \(\cv\)\n-functor, \(X\) an object of
\(\ca\). There is a bijection between elements of \(F(X)\),
\emph{i.e.}, morphisms \(t:\1\to F(X)\), and natural transformations
\(\ca(X,\_)\to F:\ca\to\und\cv\) defined as follows: with an element
\(t:\1\to F(X)\) a natural transformation is associated whose
components are given by
\[
\ca(X,Z)\rTTo^{t\tens
F_{X,Z}}F(X)\tens\und\cv(F(X),F(Z))\rTTo^{\ev^\cv}F(Z),\quad Z\in\Ob\ca.
\]
In particular, \(F\) is representable if and only if there is an object
\(X\in\Ob\ca\) and an element \(t:\1\to F(X)\) such that for each
object \(Z\in\Ob\ca\) the above composite is invertible.
\end{proposition}

\begin{proof}
Standard, see \cite[Section~1.9]{KellyGM:bascec}.
\end{proof}

\section{Serre functors for $\cv$-categories}
 \label{sec-Serre-functor-K-Cat}
Serre functors for enriched categories are for us a bridge between
ordinary Serre functors and Serre \ainf-functors.

\subsection{Serre $\cv$-functors.}
Let \(\cc\) be a \(\cv\)\n-category, \(S:\cc\to\cc\) a
\(\cv\)\n-functor. Consider a natural transformation \(\psi\) as in the
diagram below:
\begin{diagram}[LaTeXeqno]
\cc^\op\boxt\cc & \rTTo^{1\boxt S} & \cc^\op\boxt\cc\\
\dTTo<{\Hom^\op_{\cc^\op}} &\ldTwoar^{\psi}& \dTTo>{\Hom_\cc}\\
\und\cv^\op &\rTTo_{\und\cv(\_,\1)} & \und\cv
\label{equ-Serre-natural-transform}
\end{diagram}
The natural transformation \(\psi\) is a collection of morphisms of \(\cv\)
\[\psi_{X,Y}:\1\to\und\cv(\cc(X,YS),\und\cv(\cc(Y,X),\1)),\quad X,Y\in\Ob\cc.\]
Equivalently, \(\psi\) is given by a collection of
morphisms \(\psi_{X,Y}:\cc(X,YS)\to\und\cv(\cc(Y,X),\1)\) of \(\cv\),
for \(X,Y\in\Ob\cc\). Naturality of \(\psi\) may be verified
variable-by-variable.

\begin{definition}
Let \(\cc\) be a \(\cv\)\n-category. A \(\cv\)\n-functor
\(S:\cc\to\cc\) is called a \emph{right Serre functor} if  there exists
a natural isomorphism \(\psi\) as in
\eqref{equ-Serre-natural-transform}. If moreover $S$ is a
self--equivalence, it is called a \emph{Serre functor}.
\end{definition}

This terms agree with the conventions of Mazorchuk and Stroppel
\cite{math.RT/0508119} and up to taking dual spaces with the
terminology of Reiten and van den Bergh \cite{MR1887637}.

\begin{lemma}\label{lem-psiXY-natural}
Let \(S:\cc\to\cc\) be a \(\cv\)\n-functor. Fix an object $Y$ of $\cc$.
A collection of morphisms
\((\psi_{X,Y}:\cc(X,YS)\to\und\cv(\cc(Y,X),\1))_{X\in\Ob\cc}\) of $\cv$
is natural in $X$ if and only if
\begin{multline}
\psi_{X,Y}=\bigl[ \cc(X,YS)
\rTTo^{\coev^\cv}\und\cv(\cc(Y,X),\cc(Y,X)\tens\cc(X,YS))
\rTTo^{\und\cv(1,\mu_\cc)}
\\
\und\cv(\cc(Y,X),\cc(Y,YS)) \rTTo^{\und\cv(1,\tau_Y)}
\und\cv(\cc(Y,X),\1) \bigr],
\label{eq-psiXY-tauY}
\end{multline}
where
\begin{equation}
\tau_Y=\bigl[ \cc(Y,YS) \rTTo^{1^\cc_Y\tens1} \cc(Y,Y)\tens\cc(Y,YS)
\rTTo^{1\tens\psi_{Y,Y}} \cc(Y,Y)\tens\und\cv(\cc(Y,Y),\1)
\rTTo^{\ev^\cv} \1 \bigr].
\label{equ-def-tau}
\end{equation}
\end{lemma}

\begin{proof}
The collection \((\psi_{X,Y})_{X\in\Ob\cc}\) is a natural
$\cv$\n-transformation
\begin{diagram}
\cc^\op &&\rTTo^{\cc(\_,YS)} && \und\cv\\
&\rdTTo_{\cc(Y,\_)^\op} &\dTwoar>{\psi_{-,Y}}
&\ruTTo_{\und\cv(\_,\1)}\\
&&\und\cv^\op &&
\end{diagram}
if the following diagram commutes:
\begin{diagram}
\cc(Z,X) &\rTTo^{\cc(\_,YS)} &\und\cv(\cc(X,YS),\cc(Z,YS))
\\
&&\dTTo^{\und\cv(1,\psi_{Z,Y})}
\\
\dTTo<{\cc(Y,\_)} &= &\und\cv(\cc(X,YS),\und\cv(\cc(Y,Z),\1))
\\
&&\uTTo>{\und\cv(\psi_{X,Y},1)}
\\
\und\cv(\cc(Y,Z),\cc(Y,X)) &\rTTo^{\und\cv(\_,\1)}
&\und\cv(\und\cv(\cc(Y,X),\1),\und\cv(\cc(Y,Z),\1))
\end{diagram}
By closedness, this is equivalent to commutativity of the exterior of
the following diagram:
\begin{diagram}[height=1.5em,nobalance]
\cc(X,YS)\tens\cc(Z,X) &&\rTTo^{1\tens\cc(\_,YS)} &&
\begin{array}{l}
\cc(X,YS)\tens \\
\tens\und\cv(\cc(X,YS),\cc(Z,YS))
\end{array}
\\
&\rdTTo^c &&= & \\
\dTTo<{\psi_{X,Y}\tens1} &&\cc(Z,X)\tens\cc(X,YS) &&\dTTo>{\ev^\cv}
\\
&= &&\rdTTo^{\mu_\cc} &\\
\und\cv(\cc(Y,X),\1)\tens\cc(Z,X) &&\dTTo>{1\tens\psi_{X,Y}} &&\cc(Z,YS)
\\
\dTTo<{1\tens\cc(Y,\_)} &\rdTTo^c &&&
\\
&&\cc(Z,X)\tens\und\cv(\cc(Y,X),\1) &&
\\
\begin{array}{l}
\und\cv(\cc(Y,X),\1)\tens \\
\tens\und\cv(\cc(Y,Z),\cc(Y,X))
\end{array}
&= &\dTTo>{\cc(Y,\_)\tens1} &&\dTTo>{\psi_{Z,Y}}
\\
&\rdTTo^c &&&\\
\dTTo<{1\tens\und\cv(\_,\1)} &&
\begin{array}{r}
\und\cv(\cc(Y,Z),\cc(Y,X))\tens \\
\tens\und\cv(\cc(Y,X),\1)
\end{array}
&&
\\
&= &&\rdTTo^{\mu_{\und\cv}} &\\
\begin{array}{l}
\und\cv(\cc(Y,X),\1)\tens \\
\tens\und\cv(\und\cv(\cc(Y,X),\1),\und\cv(\cc(Y,Z),\1))
\end{array}
&&\rTTo^{\ev^\cv} &&\und\cv(\cc(Y,Z),\1)
\end{diagram}
The right upper quadrilateral and the left lower quadrilateral commute
by definition of \(\cc(\_,YS)\) and \(\und\cv(\_,\1)\)
respectively, see \eqref{equ-A-dot-Y} and its particular
case~\eqref{eq-V(-1)=coev-V(1cmu)}. Since \(c\) is an isomorphism,
commutativity of the exterior is equivalent to commutativity of the
pentagon. Again, by closedness, this is equivalent to commutativity of
the exterior of the following diagram:
\begin{diagram}[nobalance]
\cc(Y,Z)\tens\cc(Z,X)\tens\cc(X,YS) &\rTTo^{1\tens\mu_\cc}
&\cc(Y,Z)\tens\cc(Z,YS)
\\
\dTTo<{1\tens1\tens\psi_{X,Y}} &&\dTTo>{1\tens\psi_{Z,Y}}
\\
\cc(Y,Z)\tens\cc(Z,X)\tens\und\cv(\cc(Y,X),\1)
&&\cc(Y,Z)\tens\und\cv(\cc(Y,Z),\1)
\\
&\rdTTo(2,4)^{\mu_\cc\tens1} &\dTTo>{\ev^\cv}
\\
\dTTo<{1\tens\cc(Y,\_)\tens1} &&\1
\\
&&\uTTo>{\ev^\cv}
\\
\cc(Y,Z)\tens\und\cv(\cc(Y,Z),\cc(Y,X))\tens\und\cv(\cc(Y,X),\1)
&\rTTo^{\ev^\cv\tens1} &\cc(Y,X)\tens\und\cv(\cc(Y,X),\1)
\end{diagram}
The triangle commutes by definition of \(\cc(Y,\_)\), see
\eqref{equ-A-X-dot}. It follows that naturality of \(\psi_{-,Y}\) is
equivalent to commutativity of the hexagon:
\begin{diagram}[LaTeXeqno]
\cc(Y,Z)\tens\cc(Z,X)\tens\cc(X,YS) &\rTTo^{1\tens\mu_\cc\psi_{Z,Y}}
&\cc(Y,Z)\tens\und\cv(\cc(Y,Z),\1)
\\
\dTTo<{\mu_\cc\tens\psi_{X,Y}} &&\dTTo>{\ev^\cv}
\\
\cc(Y,X)\tens\und\cv(\cc(Y,X),\1) &\rTTo^{\ev^\cv} &\1
\label{dia-C(YZ)-C(ZX)-C(XYS)}
\end{diagram}

Assume that \(\psi_{-,Y}\) is natural, so the above diagram commutes,
and consider a particular case, $Z=Y$. Composing both paths of the
diagram with the morphism
\(1^\cc_Y\tens1\tens1:\cc(Y,X)\tens\cc(X,YS)\to\cc(Y,Y)\tens\cc(Y,X)\tens\cc(X,YS)\),
we obtain:
\begin{diagram}[LaTeXeqno]
\cc(Y,X)\tens\cc(X,YS) & \rTTo^{\mu_\cc} & \cc(Y,YS)\\
\dTTo<{1\tens\psi_{X,Y}} && \dTTo>{\tau_Y}\\
\cc(Y,X)\tens\und\cv(\cc(Y,X),\1) & \rTTo^{\ev^\cv} & \1
\label{equ-psi-X-Y-mu-tau-Y}
\end{diagram}
where $\tau_Y$ is given by expression~\eqref{equ-def-tau}. By
closedness, the above equation admits a unique solution \(\psi_{X,Y}\),
namely, \eqref{eq-psiXY-tauY}.

Assume now that \(\psi_{X,Y}\) is given by \eqref{eq-psiXY-tauY}. Then
\eqref{equ-psi-X-Y-mu-tau-Y} holds true. Plugging it into
\eqref{dia-C(YZ)-C(ZX)-C(XYS)}, whose commutativity has to be proven,
we obtain the equation
\begin{diagram}[nobalance]
\cc(Y,Z)\tens\cc(Z,X)\tens\cc(X,YS) &\rTTo^{1\tens\mu_\cc}
&\cc(Y,Z)\tens\cc(Z,YS) &\rTTo^{\mu_\cc} &\cc(Y,YS)
\\
\dTTo<{\mu_\cc\tens1} &&&&\dTTo>{\tau_Y}
\\
\cc(Y,X)\tens\cc(X,YS) &\rTTo^{\mu_\cc} &\cc(Y,YS) &\rTTo^{\tau_Y} &\1
\end{diagram}
which holds true by associativity of composition.
\end{proof}

\begin{lemma}\label{lem-psiXY-natural-inY}
Let \(S:\cc\to\cc\) be a \(\cv\)\n-functor. Fix an object \(X\) of
\(\cc\). A collection of morphisms
\((\psi_{X,Y}:\cc(X,YS)\to\und\cv(\cc(Y,X),\1))_{Y\in\Ob\cc}\) of
\(\cv\) is natural in \(Y\) if and only if for each \(Y\in\Ob\cc\)
\begin{multline}
\psi_{X,Y}=\bigl[
\cc(X,YS)\rTTo^{\coev^\cv}\und\cv(\cc(Y,X),\cc(Y,X)\tens\cc(X,YS))\rTTo^{\und\cv(1,S\tens1)}
\\
\und\cv(\cc(Y,X),\cc(YS,XS)\tens\cc(X,YS))\rTTo^{\und\cv(1,c\mu_\cc)}\und\cv(\cc(Y,X),\cc(X,XS))
\\
\rTTo^{\und\cv(1,\tau_X)}\und\cv(\cc(Y,X),\1)
\bigr],
\label{eq-psiXY-V(1tauX)}
\end{multline}
where \(\tau_X\) is given by \eqref{equ-def-tau}.
\end{lemma}

\begin{proof}
Naturality of \(\psi_{X,-}\) presented by the square
\begin{diagram}
\cc & \rTTo^S & \cc\\
\dTTo<{\cc(\_,X)^\op} &\ldTwoar^{\psi_{X,-}} &
\dTTo>{\cc(X,\_)}\\
\und\cv^\op & \rTTo^{\und\cv(\_,\1)} & \und\cv
\end{diagram}
is expressed by commutativity in \(\cv\) of the following diagram:
\begin{diagram}
\cc(Y,Z) & \rTTo^S & \cc(YS,ZS)
\\
\dTTo<{\cc(\_,X)} && \dTTo>{\cc(X,\_)}
\\
\und\cv(\cc(Z,X),\cc(Y,X)) && \und\cv(\cc(X,YS),\cc(X,ZS))
\\
\dTTo<{\und\cv(\_,\1)} && \dTTo>{\und\cv(1,\psi_{X,Z})}
\\
\und\cv(\und\cv(\cc(Y,X),\1),\und\cv(\cc(Z,X),\1)) &
\rTTo^{\und\cv(\psi_{X,Y},1)} & \und\cv(\cc(X,YS),\und\cv(\cc(Z,X),\1))
\end{diagram}
By closedness, the latter is equivalent to commutativity of the exterior of
the diagram displayed \vpageref{equ-naturality-psi-X-}.
\begin{figure}
\begin{center}
\resizebox{!}{.90\texthigh}{\rotatebox{90}{%
\begin{diagram}[height=3em]
\cc(X,YS)\tens\cc(Y,Z) &\rTTo^{1\tens S} & \cc(X,YS)\tens\cc(YS,ZS) &
\rTTo^{1\tens\cc(X,\_)} &
\begin{array}{l}
\cc(X,YS)\tens \\
\tens\und\cv(\cc(X,YS),\cc(X,ZS))
\end{array}
\\
\dTTo<{\psi_{X,Y}\tens1} &\rdTTo^c &&\rdTTo^{\mu_\cc} &\dTTo>{\ev^\cv}
\\
\und\cv(\cc(Y,X),\1)\tens\cc(Y,Z) && \cc(Y,Z)\tens\cc(X,YS) && \cc(X,ZS)
\\
\dTTo<{1\tens\cc(\_,X)} & \rdTTo^c &\dTTo>{1\tens\psi_{X,Y}}
\\
\begin{array}{l}
\und\cv(\cc(Y,X),\1)\tens \\
\tens\und\cv(\cc(Z,X),\cc(Y,X))
\end{array}
&&\cc(Y,Z)\tens\und\cv(\cc(Y,X),\1) &\framebox{\(*\)}&
\\
&\rdTTo^c &\dTTo>{\cc(\_,X)\tens1} && \dTTo>{\psi_{X,Z}}
\\
\dTTo<{1\tens\und\cv(\_,\1)} &&
\begin{array}{r}
\und\cv(\cc(Z,X),\cc(Y,X))\tens \\
\tens\und\cv(\cc(Y,X),\1)
\end{array}
&&
\\
&&&\rdTTo^{\mu_{\und\cv}} &
\\
\begin{array}{l}
\und\cv(\cc(Y,X),\1)\tens \\
\tens\und\cv(\und\cv(\cc(Y,X),\1),\und\cv(\cc(Z,X),\1))
\end{array}
&&\rTTo^{\ev^\cv} && \und\cv(\cc(Z,X),\1)
\end{diagram}
}}
\end{center}
\caption{\label{equ-naturality-psi-X-}}
\end{figure}
Since \(c\) is an isomorphism, it follows that the polygon marked by
\framebox{\(*\)} is commutative. By closedness, this is equivalent to
commutativity of the exterior of the following diagram:
\begin{diagram}[nobalance]
\cc(Z,X)\tens\cc(X,YS)\tens\cc(Y,Z)
&\rTTo^{\hspace*{-2.5em}1\tens1\tens S}
&\cc(Z,X)\tens\cc(X,YS)\tens\cc(YS,ZS)
\\
\dTTo<{1\tens c} &&\dTTo>{1\tens\mu_\cc}
\\
\cc(Z,X)\tens\cc(Y,Z)\tens\cc(X,YS) &&\cc(Z,X)\tens\cc(X,ZS)
\\
\dTTo<{1\tens1\tens\psi_{X,Y}} &&\dTTo>{1\tens\psi_{X,Z}}
\\
\cc(Z,X)\tens\cc(Y,Z)\tens\und\cv(\cc(Y,X),\1)
&&\cc(Z,X)\tens\und\cv(\cc(Z,X),\1)
\\
&\rdTTo(2,4)^{c\mu_\cc\tens1} &\dTTo>{\ev^\cv}
\\
\dTTo<{1\tens\cc(\_,X)\tens1} &&\1
\\
&&\uTTo>{\ev^\cv}
\\
\cc(Z,X)\tens\und\cv(\cc(Z,X),\cc(Y,X))\tens\und\cv(\cc(Y,X),\1)
&\rTTo^{\ev^\cv\tens1\hspace*{-2em}} &\cc(Y,X)\tens\und\cv(\cc(Y,X),\1)
\end{diagram}
The triangle commutes by \eqref{equ-A-dot-Y}. Therefore, the remaining
polygon is commutative as well:
\begin{diagram}[LaTeXeqno]
\cc(Z,X)\tens\cc(X,YS)\tens\cc(Y,Z) & \rTTo^{1\tens1\tens S} &
\cc(Z,X)\tens\cc(X,YS)\tens\cc(YS,ZS)
\\
\dTTo<{(123)^\sim}&& \dTTo>{1\tens\mu_\cc}\\
\cc(Y,Z)\tens\cc(Z,X)\tens\cc(X,YS) && \cc(Z,X)\tens\cc(X,ZS)
\\
\dTTo<{\mu_\cc\tens1}&& \dTTo>{1\tens\psi_{X,Z}}
\\
\cc(Y,X)\tens\cc(X,YS)&&
\cc(Z,X)\tens\und\cv(\cc(Z,X),\1)
\\
\dTTo<{1\tens\psi_{X,Y}}&&\dTTo>{\ev^\cv}
\\
\cc(Y,X)\tens\und\cv(\cc(Y,X),\1) & \rTTo^{\ev^\cv} & \1
\label{equ-naturality-psi-X-Y-wrt-Y}
\end{diagram}

Suppose that the collection of morphisms
\((\psi_{X,Y}:\cc(X,YS)\to\und\cv(\cc(Y,X),\1))_{Y\in\Ob\cc}\) is
natural in \(Y\). Consider diagram~\eqref{equ-naturality-psi-X-Y-wrt-Y}
with \(Z=X\). Composing both paths with the morphism
\(1^\cc_X\tens1\tens1:\cc(X,YS)\tens\cc(Y,X)\to\cc(X,X)\tens\cc(X,YS)\tens\cc(Y,X)\)
gives an equation:
\begin{diagram}[LaTeXeqno]
\cc(Y,X)\tens\cc(X,YS) &\rTTo^{S\tens1} & \cc(YS,XS)\tens\cc(X,YS) &
\rTTo^{c\mu_\cc} & \cc(X,XS)\\
\dTTo<{1\tens\psi_{X,Y}} &&= && \dTTo>{\tau_X}\\
\cc(Y,X)\tens\und\cv(\cc(Y,X),\1) &&\rTTo^{\ev^\cv} && \1
\label{equ-psi-X-Y-S-mu-tau-X}
\end{diagram}
The only solution to the above equation is given by
\eqref{eq-psiXY-V(1tauX)}.

Conversely, suppose equation~\eqref{equ-psi-X-Y-S-mu-tau-X} holds. It
suffices to prove that diagram~\eqref{equ-naturality-psi-X-Y-wrt-Y}
is commutative. Plugging in the expressions for \((1\tens\psi_{X,Y})\ev^\cv\)
and \((1\tens\psi_{X,Z})\ev^\cv\) into
\eqref{equ-naturality-psi-X-Y-wrt-Y}, we obtain (cancelling a common
permutation of the factors of the source object):
\begin{diagram}
\cc(X,YS)\tens\cc(Y,Z)\tens\cc(Z,X) & \rTTo^{1\tens S\tens S} &
\cc(X,YS)\tens\cc(YS,ZS)\tens\cc(ZS,XS)
\\
\dTTo<{1\tens\mu_\cc} && \dTTo>{\mu_\cc\tens1}
\\
\cc(X,YS)\tens\cc(Y,X) && \cc(X,ZS)\tens\cc(ZS,XS)
\\
\dTTo<{1\tens S} && \dTTo>{\mu_\cc}
\\
\cc(X,YS)\tens\cc(YS,XS) &&\cc(X,XS)
\\
\dTTo<{\mu_\cc} &&\dTTo>{\tau_X}
\\
\cc(X,XS) &\rTTo^{\tau_X} &\1
\end{diagram}
Commutativity of the diagram follows from associativity of \(\mu_\cc\)
and the fact that \(S\) is a \(\cv\)\n-functor. The lemma is proven.
\end{proof}

\begin{proposition}\label{prop-e-K-psi-S-psi}
Assume that \(S:\cc\to\cc\) is a \(\cv\)\n-functor, and \(\psi\) is a
natural transformation as in \eqref{equ-Serre-natural-transform}. Then
the following diagram commutes (in \(\cv\)):
\begin{diagram}
\cc(Y,X) & \rTTo^{e} & \und\cv(\und\cv(\cc(Y,X),\1),\1)\\
\dTTo<{S} && \dTTo>{\und\cv(\psi_{X,Y},\1)}\\
\cc(YS,XS) & \rTTo^{\psi_{YS,X}} & \und\cv(\cc(X,YS),\1)
\end{diagram}
In particular, if for each pair of objects \(X,Y\in\Ob\cc\) the object
\(\cc(Y,X)\) is reflexive, and \(\psi\) is an isomorphism, then \(S\)
is fully faithful.
\end{proposition}

\begin{proof}
By closedness, it suffices to prove commutativity of the following
diagram:
\begin{diagram}
\cc(X,YS)\tens\cc(Y,X) &\rTTo^{1\tens e}
&\cc(X,YS)\tens\und\cv(\und\cv(\cc(Y,X),\1),\1)
\\
\dTTo<{1\tens S} &&\dTTo>{1\tens\und\cv(\psi_{X,Y},\1)}
\\
\cc(X,YS)\tens\cc(YS,XS) &&\cc(X,YS)\tens\und\cv(\cc(X,YS),\1)
\\
\dTTo<{1\tens\psi_{YS,X}} &&\dTTo>{\ev^\cv}
\\
\cc(X,YS)\tens\und\cv(\cc(X,YS),\1) &\rTTo^{\ev^\cv} &\1
\end{diagram}
Using \eqref{equ-psi-X-Y-mu-tau-Y} and the definition of \(e\), the
above diagram can be transformed as follows:
\begin{diagram}
\cc(X,YS)\tens\cc(Y,X) &\rTTo^{c} &\cc(Y,X)\tens\cc(X,YS)
\\
\dTTo<{1\tens S} &&\dTTo>{1\tens\psi_{X,Y}}
\\
\cc(X,YS)\tens\cc(YS,XS) &&\cc(Y,X)\tens\und\cv(\cc(Y,X),\1)
\\
\dTTo<{\mu_\cc} &&\dTTo>{\ev^\cv}
\\
\cc(X,XS) &\rTTo^{\tau_X} &\1
\end{diagram}
It is commutative by \eqref{equ-psi-X-Y-S-mu-tau-X}.
\end{proof}

\propref{prop-e-K-psi-S-psi} implies that a right Serre functor is
fully faithful if and only if \(\cc\) is hom-reflexive, \emph{i.e.}, if
\(\cc(X,Y)\) is a reflexive object of \(\cv\) for each pair of objects
\(X,Y\in\Ob\cc\). If this is the case, a right Serre functor will be a
Serre functor if and only if it is essentially surjective on objects.
The most natural reason for hom-reflexivity is, of course, $\kk$ being
a field. When $\kk$ is a field, an object $C$ of \(\gr(\kk\vect)\) is
reflexive iff all spaces $C^n$ are finite-dimensional. The ring $\kk$
being a field, the homology functor \(H^\bullet:\ck\to\gr(\kk\vect)\)
is an equivalence (see e.g. \cite[Chapter~III, \S~2,
Proposition~4]{GelfMan:Homology}). Hence, an object $C$ of \(\ck\) is
reflexive iff all homology spaces $H^nC$ are finite-dimensional. A
projective module of finite rank over an arbitrary commutative ring
$\kk$ is reflexive as an object of a rigid monoidal category
\cite[Example~1.23]{DelMil}. Thus, an object $C$ of $\gr(\kMod)$ whose
components $C^n$ are projective $\kk$\n-modules of finite rank is
reflexive.

\begin{proposition}\label{prop-existence-Serre-K-functor}
Let \(\cc\) be a \(\cv\)\n-category. There exists a right Serre
\(\cv\)\n-functor \(S:\cc\to\cc\) if and only if for each object
\(Y\in\Ob\cc\) the \(\cv\)\n-functor
\[ \Hom_\cc(Y,\_)^\op\cdot\und\cv(\_,\1)
=\und\cv(\cc(Y,\_)^\op,\1):\cc^\op\to\und\cv
\]
is representable.
\end{proposition}

\begin{proof}
Standard, see \cite[Section~1.10]{KellyGM:bascec}.
\end{proof}

\subsection{Commutation with equivalences.}
Let \(\cc\) and \(\cc'\) be \(\cv\)\n-categories with right Serre
functors \(S:\cc\to\cc\) and \(S':\cc'\to\cc'\), respectively. Let
\(\psi\) and \(\psi'\) be isomorphisms as in
\eqref{equ-Serre-natural-transform}. For objects \(Y\in\Ob\cc\),
\(Z\in\Ob\cc'\), define \(\tau_Y\), \(\tau'_Z\) by \eqref{equ-def-tau}.
Let \(T:\cc\to\cc'\) be a \(\cv\)\n-functor, and suppose that \(T\) is
fully faithful. Then there is a natural transformation
\(\varkappa:ST\to TS'\) such that, for each object \(Y\in\Ob\cc\), the
following equation holds:
\begin{equation}
\bigl[
\cc(Y,YS)\rTTo^T\cc'(YT,YST)\rTTo^{\cc'(YT,\varkappa)}\cc'(YT,YTS')\rTTo^{\tau'_{YT}}\1
\bigr]
=\tau_Y.
\label{equ-tau-Y-T-C1kappa-tau-TY}
\end{equation}
Indeed, the left hand side of equation~\eqref{equ-tau-Y-T-C1kappa-tau-TY}
equals
\begin{multline*}
\bigl[
\cc(Y,YS)\rTTo^T\cc'(YT,YST)\rTTo^{1\tens\varkappa_Y}
\cc'(YT,YST)\tens\cc'(YST,YTS')\\
\rTTo^{\mu_{\cc'}}\cc'(YT,YTS')\rTTo^{\tau'_{YT}}\1
\bigr].
\end{multline*}
Using relation \eqref{equ-psi-X-Y-mu-tau-Y} between \(\tau'_{YT}\) and
\(\psi'_{YST,YT}\), we get:
\begin{multline*}
\bigl[
\cc(Y,YS)\rTTo^T\cc'(YT,YST)\rTTo^{1\tens\varkappa_Y}\cc'(YT,YST)\tens\cc'(YST,YTS')\\
\rTTo^{1\tens\psi'_{YST,YT}}\cc'(YT,YST)\tens\und\cv(\cc'(YT,YST),\1)\rTTo^{\ev^\cv}\1
\bigr].
\end{multline*}
Therefore, equation~\eqref{equ-tau-Y-T-C1kappa-tau-TY} is equivalent to
the following equation:
\begin{multline*}
\bigl[ \cc(Y,YS) \rTTo^{1\tens\varkappa_Y} \cc(Y,YS)\tens\cc'(YST,YTS')
\rTTo^{1\tens\psi'_{YST,YT}}
\\
\cc(Y,YS)\tens\und\cv(\cc'(YT,YST),\1) \rTTo^{1\tens\und\cv(T,\1)}
\cc(Y,YS)\tens\und\cv(\cc(Y,YS),\1) \rTTo^{\ev^\cv} \1 \bigr]=\tau_Y.
\end{multline*}
It implies that the composite
\[
\1\rTTo^{\varkappa_Y}\cc'(YST,YTS')\rTTo^{\psi'_{YST,YT}}\und\cv(\cc'(YT,YST),\1)
\rTTo^{\und\cv(T,\1)}\und\cv(\cc(Y,YS),\1)
\]
is equal to \(\und{\tau_Y}:\1\to\und\cv(\cc(Y,YS),\1)\), the morphism
that corresponds to \(\tau_Y\) by closedness of the category \(\cv\).
Since the morphisms \(\psi'_{YST,YT}\) and \(\und\cv(T,\1)\) are
invertible, the morphism \(\varkappa_Y:\1\to\cc'(YST,YTS')\) is
uniquely determined.

\begin{lemma}\label{lem-T-C1kappa-psi-VT1-psi}
The transformation \(\varkappa\) satisfies the following equation:
\begin{multline*}
\psi_{X,Y}=\bigl[
\cc(X,YS)\rTTo^T\cc'(XT,YST)\rTTo^{\cc'(XT,\varkappa)}\cc(XT,YTS')\\
\rTTo^{\psi'_{XT,YT}}\und\cv(\cc'(YT,XT),\1)\rTTo^{\und\cv(T,\1)}\und\cv(\cc(Y,X),\1)
\bigr],
\end{multline*}
for each pair of objects \(X,Y\in\Ob\cc\).
\end{lemma}

\begin{proof}
The exterior of the following diagram commutes:
\[ \hspace*{-2em}
\begin{diagram}[inline]
&&\cc(Y,YS) &\rTTo^{\tau_Y} & \1\\
&\ruTTo^{\mu_\cc} &\dTTo>T && \uTTo>{\tau'_{YT}}
\\
\cc(Y,X)\tens\cc(X,YS) \hspace*{-1em} &&\cc'(YT,YST)
&\rTTo^{\cc'(YT,\varkappa)} &\cc'(YT,YTS')
\\
&\rdTTo_{T\tens T} &\uTTo>{\mu_{\cc'}} && \uTTo>{\mu_{\cc'}}
\\
&&\hspace*{-3.5em} \cc'(YT,XT)\tens\cc'(XT,YST)
&\rTTo^{1\tens\cc'(XT,\varkappa)} &\cc'(YT,XT)\tens\cc'(XT,YTS')
\end{diagram}
\]
The right upper square commutes by the definition of \(\varkappa\),
commutativity of the lower square is a consequence of associativity of
\(\mu_{\cc'}\). The left quadrilateral is commutative since \(T\) is a
\(\cv\)\n-functor. Transforming both paths with the help of
equation~\eqref{equ-psi-X-Y-mu-tau-Y} yields the following equation:
\begin{multline*}
\bigl[ \cc(Y,X)\tens\cc(X,YS) \rTTo^{1\tens\psi_{X,Y}}
\cc(Y,X)\tens\und\cv(\cc(Y,X),\1) \rTTo^{\ev^\cv} \1 \bigr]
\\
\quad=\bigl[ \cc(Y,X)\tens\cc(X,YS) \rTTo^{T\tens T}
\cc'(YT,XT)\tens\cc'(XT,YST) \rTTo^{1\tens\cc'(XT,\varkappa)} \hfill
\\
\cc'(YT,XT)\tens\cc'(XT,YTS') \rTTo^{1\tens\psi'_{XT,YT}}
\cc'(YT,XT)\tens\und\cv(\cc'(YT,XT),\1) \rTTo^{\ev^\cv} \1 \bigr]
\\
\quad=\bigl[ \cc(Y,X)\tens\cc(X,YS) \rTTo^{1\tens T}
\cc(Y,X)\tens\cc'(XT,YST) \hfill
\\
\rTTo^{1\tens\cc'(XT,\varkappa)} \cc(Y,X)\tens\cc'(XT,YTS')
\rTTo^{1\tens\psi'_{XT,YT}} \cc(Y,X)\tens\und\cv(\cc'(YT,XT),\1)
\\
\rTTo^{1\tens\und\cv(T,\1)} \cc(Y,X)\tens\und\cv(\cc(Y,X),\1)
\rTTo^{\ev^\cv} \1 \bigr].
\end{multline*}
The required equation follows by closedness of \(\cv\).
\end{proof}

\begin{corollary}\label{cor-T-equiv-implies-ST-TS'}
If \(T\) is an equivalence, then the natural transformation
\(\varkappa:ST\to TS'\) is an isomorphism.
\end{corollary}

\begin{proof}
\lemref{lem-T-C1kappa-psi-VT1-psi} implies that
\(\cc'(XT,\varkappa):\cc'(XT,YST)\to\cc'(XT,XTS')\) is an isomorphism,
for each \(X\in\Ob\cc\). Since \(T\) is essentially surjective, it
follows that the morphism
\(\cc'(Z,\varkappa):\cc'(Z,YST)\to\cc'(Z,YTS')\) is invertible, for
each \(Z\in\Ob\cc'\), thus \(\varkappa\) is an isomorphism.
\end{proof}

\begin{corollary}
A right Serre $\cv$\n-functor is unique up to an isomorphism.
\end{corollary}

\begin{proof}
Suppose \(S,S':\cc\to\cc\) are right Serre functors. Applying
\corref{cor-T-equiv-implies-ST-TS'} to the functor
\(T=\Id_\cc:\cc\to\cc\) yields a natural isomorphism \(\varkappa:S\to
S'\).
\end{proof}

\subsection{Trace functionals determine the Serre functor.}
Combining for a natural transformation $\psi$ diagrams
\eqref{equ-psi-X-Y-mu-tau-Y} and \eqref{equ-psi-X-Y-S-mu-tau-X} we get
the equation
\begin{diagram}[LaTeXeqno]
\cc(Y,X)\tens\cc(X,YS) &\rTTo^{\mu_\cc} &\cc(Y,YS) &\rTTo^{\tau_Y} &\1
\\
\dTTo<{S\tens1} &&= &&\uTTo>{\tau_X}
\\
\cc(YS,XS)\tens\cc(X,YS) &\rTTo^c &\cc(X,YS)\tens\cc(YS,XS)
&\rTTo^{\mu_\cc} &\cc(X,XS)
\label{dia-S1-c-mutau-mutau}
\end{diagram}
The above diagram can be written as the equation
\begin{diagram}[LaTeXeqno]
\cc(X,YS)\tens\cc(Y,X) &\rTTo^{1\tens S} &\cc(X,YS)\tens\cc(YS,XS)
\\
\dTTo<c &= &\dTTo>{\phi_{YS,X}}
\\
\cc(Y,X)\tens\cc(X,YS) &\rTTo^{\phi_{X,Y}} &\1
\label{dia-1Sphi-cphi}
\end{diagram}

When $S$ is a fully faithful right Serre functor, the pairing
\begin{equation}
\phi_{X,Y} =\bigl[ \cc(Y,X)\tens\cc(X,YS) \rTTo^{\mu_\cc} \cc(Y,YS)
\rTTo^{\tau_Y} \1\bigr]
\label{equ-def-pairing-phi}
\end{equation}
is perfect. Namely, the induced by it morphism
\(\psi_{X,Y}:\cc(X,YS)\to\und\cv(\cc(Y,X),\1)\) is invertible, and
induced by the pairing
\begin{equation*}
c\cdot\phi_{X,Y} =\bigl[ \cc(X,YS)\tens\cc(Y,X) \rTTo^c
\cc(Y,X)\tens\cc(X,YS) \rTTo^{\phi_{X,Y}} \1\bigr]
\end{equation*}
the morphism \(\psi':\cc(Y,X)\to\und\cv(\cc(X,YS),\1)\) is invertible.
In fact, diagram~\eqref{dia-1Sphi-cphi} implies that
\[ \psi' =\bigl[ \cc(Y,X) \rTTo^S \cc(YS,XS) \rTTo^{\psi_{YS,X}}
\und\cv(\cc(X,YS),\1) \bigr].
\]

Diagram~\eqref{dia-1Sphi-cphi} allows to restore the morphisms
$S:\cc(Y,X)\to\cc(YS,XS)$ unambiguously from \(\Ob S\) and the trace
functionals $\tau$, due to \(\psi_{YS,X}\) being isomorphisms.

\begin{proposition}\label{prop-traces-determine-Serre-functor}
A map \(\Ob S\) and trace functionals $\tau_X$, \(X\in\Ob\cc\),
such that the induced \(\psi_{X,Y}\) from \eqref{eq-psiXY-tauY} are
invertible, define a unique right Serre $\cv$\n-functor
\((S,\psi_{X,Y})\).
\end{proposition}

\begin{proof}
Let us show that the obtained morphisms \(S:\cc(Y,X)\to\cc(YS,XS)\)
preserve the composition in $\cc$. In fact, due to associativity of
composition we have
\begin{multline*}
\bigl[ \cc(X,ZS)\tens\cc(Z,Y)\tens\cc(Y,X) \rTTo^{1\tens S\tens S}
\cc(X,ZS)\tens\cc(ZS,YS)\tens\cc(YS,XS)
\\
\hfill \rTTo^{1\tens\mu_\cc} \cc(X,ZS)\tens\cc(ZS,XS)
\rTTo^{\phi_{ZS,X}} \1 \bigr] \quad
\\
\quad=\bigl[\cc(X,ZS)\tens\cc(Z,Y)\tens\cc(Y,X) \rTTo^{1\tens S\tens S}
\cc(X,ZS)\tens\cc(ZS,YS)\tens\cc(YS,XS) \hfill
\\
\hfill \rTTo^{\mu_\cc\tens1} \cc(X,YS)\tens\cc(YS,XS) \rTTo^{\mu_\cc}
\cc(X,XS) \rTTo^{\tau_X} \1 \bigr] \quad
\\
 \quad = \bigl[\cc(X,ZS)\tens\cc(Z,Y)\tens\cc(Y,X)
\rTTo^{(1\tens S\tens1)(\mu_\cc\tens1)} \cc(X,YS)\tens\cc(Y,X) \hfill
\\
\hfill \rTTo^{1\tens S} \cc(X,YS)\tens\cc(YS,XS) \rTTo^{\mu_\cc}
\cc(X,XS) \rTTo^{\tau_X} \1 \bigr] \quad
\\
 \quad = \bigl[\cc(X,ZS)\tens\cc(Z,Y)\tens\cc(Y,X)
\rTTo^{(1\tens S\tens1)(\mu_\cc\tens1)} \cc(X,YS)\tens\cc(Y,X) \hfill
\\
\hfill \rTTo^c \cc(Y,X)\tens\cc(X,YS) \rTTo^{\mu_\cc} \cc(Y,YS)
\rTTo^{\tau_Y} \1 \bigr] \quad
\\
\quad = \bigl[\cc(X,ZS)\tens\cc(Z,Y)\tens\cc(Y,X)
\rTTo^{(123)_c(1\tens1\tens S)}
\cc(Y,X)\tens\cc(X,ZS)\tens\cc(ZS,YS) \hfill
\\
\hfill \rTTo^{1\tens\mu_\cc} \cc(Y,X)\tens\cc(X,YS) \rTTo^{\mu_\cc}
\cc(Y,YS) \rTTo^{\tau_Y} \1 \bigr] \quad
\\
\quad = \bigl[\cc(X,ZS)\tens\cc(Z,Y)\tens\cc(Y,X)
\rTTo^{(123)_c(1\tens1\tens S)}
\cc(Y,X)\tens\cc(X,ZS)\tens\cc(ZS,YS) \hfill
\\
\hfill \rTTo^{\mu_\cc\tens1} \cc(Y,ZS)\tens\cc(ZS,YS) \rTTo^{\mu_\cc}
\cc(Y,YS) \rTTo^{\tau_Y} \1 \bigr] \quad
\\
\quad = \bigl[\cc(X,ZS)\tens\cc(Z,Y)\tens\cc(Y,X)
\rTTo^{(123)_c(\mu_\cc\tens1)} \cc(Y,ZS)\tens\cc(Z,Y) \hfill
\\
\hfill \rTTo^{1\tens S} \cc(Y,ZS)\tens\cc(ZS,YS) \rTTo^{\mu_\cc}
\cc(Y,YS) \rTTo^{\tau_Y} \1 \bigr] \quad
\\
\quad = \bigl[\cc(X,ZS)\tens\cc(Z,Y)\tens\cc(Y,X)
\rTTo^{(123)_c(\mu_\cc\tens1)} \cc(Y,ZS)\tens\cc(Z,Y) \hfill
\\
\hfill \rTTo^c \cc(Z,Y)\tens\cc(Y,ZS) \rTTo^{\mu_\cc}
\cc(Z,ZS) \rTTo^{\tau_Z} \1 \bigr] \quad
\\
\quad = \bigl[\cc(X,ZS)\tens\cc(Z,Y)\tens\cc(Y,X) \rTTo^{(321)_c}
\cc(Z,Y)\tens\cc(Y,X)\tens\cc(X,ZS) \hfill
\\
\rTTo^{1\tens\mu_\cc} \cc(Z,Y)\tens\cc(Y,ZS) \rTTo^{\mu_\cc}
\cc(Z,ZS) \rTTo^{\tau_Z} \1 \bigr].
\end{multline*}
On the other hand
\begin{multline*}
\bigl[ \cc(X,ZS)\tens\cc(Z,Y)\tens\cc(Y,X) \rTTo^{1\tens\mu_\cc}
\cc(X,ZS)\tens\cc(Z,X)
\\
\hfill \rTTo^{1\tens S} \cc(X,ZS)\tens\cc(ZS,XS)
\rTTo^{\phi_{ZS,X}} \1 \bigr] \quad
\\
\quad = \bigl[ \cc(X,ZS)\tens\cc(Z,Y)\tens\cc(Y,X)
\rTTo^{1\tens\mu_\cc} \cc(X,ZS)\tens\cc(Z,X) \hfill
\\
\hfill \rTTo^c \cc(Z,X)\tens\cc(X,ZS) \rTTo^{\phi_{X,Z}} \1 \bigr]
\quad
\\
\quad = \bigl[\cc(X,ZS)\tens\cc(Z,Y)\tens\cc(Y,X) \rTTo^{(321)_c}
\cc(Z,Y)\tens\cc(Y,X)\tens\cc(X,ZS) \hfill
\\
\hfill \rTTo^{\mu_\cc\tens1} \cc(Z,X)\tens\cc(X,ZS) \rTTo^{\mu_\cc}
\cc(Z,ZS) \rTTo^{\tau_Z} \1 \bigr] \quad
\\
\quad = \bigl[\cc(X,ZS)\tens\cc(Z,Y)\tens\cc(Y,X) \rTTo^{(321)_c}
\cc(Z,Y)\tens\cc(Y,X)\tens\cc(X,ZS) \hfill
\\
\rTTo^{1\tens\mu_\cc} \cc(Z,Y)\tens\cc(Y,ZS) \rTTo^{\mu_\cc}
\cc(Z,ZS) \rTTo^{\tau_Z} \1 \bigr].
\end{multline*}
The last lines of both expressions coincide, hence \((S\tens
S)\mu_\cc=\mu_\cc S\).

Let us prove that the morphisms \(S:\cc(X,X)\to\cc(XS,XS)\) of $\cv$
preserve units. Indeed, the exterior of the following diagram commutes:
\begin{diagram}
&&\cc(X,XS) &\rTTo^{\tau_X} &\1
\\
&\ruId>{\hspace*{1em}=} &\uTTo>{\mu_\cc} &&
\\
\cc(X,XS) &\rTTo_{\lambda^{\msf{\lar.I}}(1_X\tens1)}
&\cc(X,X)\tens\cc(X,XS) &&
\\
\dTTo<{\lambda^{\msf{I\lar.}}}>\wr &= &\uTTo>c &= &\uTTo>{\tau_X}
\\
\cc(X,XS)\tens\1 &\rTTo^{1\tens1_X} &\cc(X,XS)\tens\cc(X,X) &&
\\
&\rdTTo<{1\tens1_{XS}}>{\hspace*{1em}1\tens?} &\dTTo>{1\tens S} &&
\\
&&\cc(X,XS)\tens\cc(XS,XS) &\rTTo^{\mu_\cc} &\cc(X,XS)
\end{diagram}
Therefore, both paths from \(\cc(X,XS)\) to $\1$, going through the
isomorphism \(\lambda^{\msf{I\lar.}}\), sides of triangle marked
`\(1\tens?\)', \(\mu_\cc\) and $\tau_X$, compose to the same morphism
$\tau_X$. Invertibility of \(\psi_{X,X}\) implies that the origin `?'
of the mentioned triangle commutes, that is,
\[ 1_{XS}=\bigl[\1 \rTTo^{1_X} \cc(X,X) \rTTo^S \cc(XS,XS)\bigr].
\]

Summing up, the constructed \(S:\cc\to\cc\) is a $\cv$\n-functor.
Applying \lemref{lem-psiXY-natural} we deduce that \(\psi_{-,Y}\) is
natural in the first argument for all objects $Y$ of $\cc$. Recall that
\(\psi_{X,Y}\) is a unique morphism which makes
diagram~\eqref{equ-psi-X-Y-mu-tau-Y} commutative. Due to
equation~\eqref{dia-S1-c-mutau-mutau} \(\psi_{X,Y}\) makes commutative
also diagram~\eqref{equ-psi-X-Y-S-mu-tau-X}. This means that
\(\psi_{X,Y}\) can be presented also in the
form~\eqref{eq-psiXY-V(1tauX)}. Applying \lemref{lem-psiXY-natural-inY}
we deduce that \(\psi_{X,-}\) is natural in the second argument for all
objects $X$ of $\cc$. Being natural in each variable \(\psi\) is
natural as a whole \cite[Section~1.4]{KellyGM:bascec}.
\end{proof}

\subsection{Base change.}
Let \(\cv=(\cv,\tens^I_\cv,\lambda^f_\cv)\),
\(\cw=(\cw,\tens^I_\cw,\lambda^f_\cw)\) be closed symmetric Monoidal
$\fu$\n-categories. Let
\((B,\beta^I):(\cv,\tens^I_\cv,\lambda^f_\cv)\to(\cw,\tens^I_\cw,\lambda^f_\cw)\)
be a lax symmetric Monoidal functor. Denote by
\(\wh{B}:\wh\cv\to\wh\cw\) the corresponding multifunctor. According to
\cite{Manzyuk-PhD}, \((B,\beta^I)\) gives rise to a lax symmetric
Monoidal \(\Cat\)\n-functor \((B_*,\beta^I_*):\VCat\to\WCat\). Since
the multicategories \(\wh\cv\) and \(\wh\cw\) are closed, the
multifunctor \(\wh B\) determines the closing transformation
 \(\und{\wh B}\). In particular, we have a \(\cw\)\n-functor
\(B_*\und\cv\to\und\cw\), \(X\mapsto BX\), which is denoted by
\(\und{\wh B}\) by abuse of notation, whose action on morphisms is
found from the following equation in \(\cw\):
\begin{equation}
\bigl[ BX\tens B(\und\cv(X,Y))\rTTo^{1\tens\und{\wh{B}}\;\;{}}
BX\tens\und\cw(BX,BY) \rTTo^{\ev^\cw}BY \bigr]=\wh{B}(\ev^\cv).
\label{equ-und-hat-H-0}
\end{equation}
Let \(\wh{B_*}:\wh\VCat\to\wh\WCat\) denote the symmetric
\(\Cat\)\n-multifunctor that corresponds to the lax symmetric Monoidal
\(\Cat\)\n-functor \((B_*,\beta^I_*)\). Clearly, \(\wh{B_*}\) commutes
with taking opposite.

In the sequel, the tensor product in the categories \(\cv\) and
\(\cw\) is denoted by \(\tens\), the unit objects in both categories
are denoted by \(\1\).

Let \(\ca\) be a \(\cv\)\n-category. We claim that the \(\cw\)\n-functor
\[
\wh{B_*}\Hom_\ca\cdot\und{\wh{B}}=
\bigl[
B_*(\ca)^\op\boxt
B_*(\ca)\rTTo^{\wh{B_*}\Hom_\ca\;{}}B_*\und\cv\rTTo^{\und{\wh{B}}}\und\cw
\bigr]
\]
coincides with \(\Hom_{B_*\ca}\). Indeed, both functors send a pair of
objects \((X,Y)\in\Ob\ca\times\Ob\ca\) to the object
\(B(\ca(X,Y))=(B_*\ca)(X,Y)\) of \(\cw\). Applying \(\wh{B}\) to
equation~\eqref{equ-Hom-A} yields a commutative diagram
\begin{diagram}[nobalance,objectstyle=\scriptstyle]
B(\ca(X,U))\tens B(\ca(Y,X))\tens B(\ca(U,V)) &\rTTo^{1\tens\wh{B}\Hom_\ca\;{}}
&B(\ca(X,U))\tens B\und\cv(\ca(X,U),\ca(Y,V))
\\
\dTTo<{c\tens1} && \dTTo>{\wh{B}(\ev^\cv)}
\\
B(\ca(Y,X))\tens B(\ca(X,U))\tens B(\ca(U,V)) &
\rTTo^{\wh{B}(\mu^{\mb3}_\ca)}~\|_{\mu^{\mb3}_{B_*\ca}} & B(\ca(Y,V))
\end{diagram}
Expanding \(\wh{B}(\ev^\cv)\) according to \eqref{equ-und-hat-H-0} we
transform the above diagram as follows:
\begin{diagram}[objectstyle=\scriptstyle]
B(\ca(X,U))\tens B(\ca(Y,X))\tens B(\ca(U,V)) &
\rTTo^{1\tens\wh{B}\Hom_\ca\cdot\und{\wh{B}}\;\;{}} & B(\ca(X,U))\tens
\und{\cw}(B(\ca(X,U)),B(\ca(Y,V)))
\\
\dTTo<{c\tens1} && \dTTo>{\ev^{\cw}}
\\
B(\ca(Y,X))\tens B(\ca(X,U))\tens B(\ca(U,V)) &
\rTTo^{\mu^{\mb3}_{B_*\ca}} & B(\ca(Y,V))
\end{diagram}
It follows that the functors \(\wh{B_*}\Hom_\ca\cdot\und{\wh{B}}\) and
\(\Hom_{B_*\ca}\) satisfy the same equation, therefore they must
coincide by closedness of \(\cw\).

There is a natural transformation of \(\cw\)\n-functors \(\zeta'\)
as in the diagram below:
\begin{diagram}
B_*\und\cv^\op & \rTTo^{B_*\und\cv(\_,\1)} & B_*\und\cv\\
\dTTo<{(\und{\wh{B}})^\op} &\ldTwoar^{\zeta'} & \dTTo>{\und{\wh{B}}}\\
\und\cw^\op & \rTTo^{\und\cw(\_,B\1)} & \und\cw
\end{diagram}
For each object \(X\), the morphism
\(\zeta'_X:B(\und\cv(X,\1))\to\und\cw(BX,B\1)\) in \(\cw\) comes from the map
\(\wh{B}(\ev^\cv):BX\tens B\und\cv(X,\1)\to B\1\) by
closedness of \(\cw\). In other words,
\(\zeta'_X=\und{\wh{B}}_{X,\1}\). Naturality of \(\zeta'\) is expressed
by the following equation in \(\cw\):
\begin{diagram}
B\und\cv(Y,X) & \rTTo^{B\und\cv(\_,\1)} &
B\und\cv(\und\cv(X,\1),\und\cv(Y,\1))\\
\dTTo<{\und{\wh B}} && \dTTo>{\und{\wh B}}\\
\und\cw(BY,BX) && \und\cw(B\und\cv(X,\1),B\und\cv(Y,\1))\\
\dTTo<{\und\cw(\_,B\1)} && \dTTo>{\und\cw(1,\zeta'_Y)}\\
\und\cw(\und\cw(BX,B\1),\und\cw(BY,B\1)) & \rTTo^{\und\cw(\zeta'_X,1)}
& \und\cw(B\und\cv(X,\1),\und\cw(BY,B\1))
\end{diagram}
By closedness of \(\cw\), it is equivalent to the following equation:
\begin{diagram}[nobalance]
B\und\cv(X,\1)\tens B\und\cv(Y,X)
 &\rTTo^{\hspace*{-3.5em}1\tens B\und\cv(\_,\1)}
&B\und\cv(X,\1)\tens B\und\cv(\und\cv(X,\1),\und\cv(Y,\1))
\\
\dTTo<{\zeta'_X\tens\und{\wh B}} &&\dTTo>{\wh{B}(\ev^\cv)}
\\
\und\cw(BX,B\1)\tens\und\cw(BY,BX) &&B\und\cv(Y,\1)
\\
\dTTo<{1\tens\und\cw(\_,B\1)} &&\dTTo>{\zeta'_Y}
\\
\und\cw(BX,B\1)\tens\und\cw(\und\cw(BX,B\1),\und\cw(BY,B\1))
&\rTTo^{\ev^\cw\hspace*{-1.5em}} &\und\cw(BY,B\1)
\end{diagram}
By \eqref{equ-A-dot-Y}, the above equation reduces to the equation
\begin{diagram}
B\und\cv(Y,X)\tens B\und\cv(X,\1) &
\rTTo^{\wh{B}(\mu_{\und\cv})}~\|_{\mu_{B_*\und\cv}} & B\und\cv(Y,\1)
\\
\dTTo<{\und{\wh B}\tens\zeta'_X}~=>{\und{\wh B}_{Y,X}\tens\und{\wh B}_{X,\1}} &&
\dTTo<{\und{\wh B}_{Y,\1}}~=>{\zeta'_Y}
\\
\und\cw(BY,BX)\tens\und\cw(BX,B\1) & \rTTo^{\mu_{\und\cw}} &\und\cw(BY,B\1)
\end{diagram}
which expresses the fact that \(\und{\wh B}:B_*\und\cv\to\und\cw\) is a
\(\cw\)\n-functor.

Suppose that \(\beta^\emptyset:\1\to B\1\) is an isomorphism. Then
there is a natural isomorphism of functors
\[\und\cw(1,(\beta^\emptyset)^{-1}):\und\cw(\_,B\1)\to\und\cw(\_,\1):\und\cw^\op\to\und\cw.\]
Pasting it with \(\zeta'\) gives a natural transformation \(\zeta\) as
in the diagram below:
\begin{diagram}[LaTeXeqno]
B_*\und\cv^\op & \rTTo^{B_*\und\cv(\_,\1)} & B_*\und\cv\\
\dTTo<{(\und{\wh{B}})^\op} &\ldTwoar^{\zeta} & \dTTo>{\und{\wh{B}}}\\
\und\cw^\op & \rTTo^{\und\cw(\_,\1)} & \und\cw
\label{equ-commutation-H-0-duality}
\end{diagram}

\begin{proposition}\label{prop-K-Serre-exists-iff-k-Serre}
Suppose \(\zeta\) is an isomorphism. Let \(\cc\) be a
\(\cv\)\n-category, and suppose \(S:\cc\to\cc\) is a right Serre
\(\cv\)\n-functor. Then \(B_*(S):B_*(\cc)\to B_*(\cc)\) is a right
Serre \(\cw\)\n-functor.
\end{proposition}

\begin{proof}
Let \(\psi\) be a natural isomorphism as in
\eqref{equ-Serre-natural-transform}. Applying the
\(\Cat\)\n-multifunctor \(\wh{B_*}\) and patching the result with
diagram \eqref{equ-commutation-H-0-duality} yields the following
diagram:
\begin{diagram}[LaTeXeqno]
B_*(\cc)^\op\boxt B_*(\cc) & \rTTo^{1\boxt B_*(S)} & B_*(\cc)^\op\boxt
B_*(\cc)
\\
\dTTo<{\wh{B_*}(\Hom_{\cc^\op})^\op} &\ldTwoar^{\wh{B_*}(\psi)}&\dTTo>{\wh{B_*}\Hom_\cc}
\\
B_*\und\cv^\op & \rTTo^{B_*\und\cv(\_,\1)} & B_*\und\cv
\\
\dTTo<{(\und{\wh{B}})^\op} &\ldTwoar^\zeta & \dTTo>{\und{\wh{B}}}
\\
\und\cw^\op & \rTTo^{\und\cw(\_,\1)} & \und\cw
\label{equ-H-0-applied-to-Serre-diagram}
\end{diagram}
Since \(\wh{B_*}\Hom_\cc\cdot\und{\wh{B}}=\Hom_{B_*(\cc)}\) and
\(\wh{B_*}\Hom_{\cc^\op}\cdot\und{\wh{B}}=\Hom_{B_*(\cc)^\op}\), we
obtain a natural transformation
\begin{diagram}
B_*(\cc)^\op\boxt B_*(\cc) & \rTTo^{1\boxt B_*(S)} & B_*(\cc)^\op\boxt
B_*(\cc)
\\
\dTTo<{\Hom^\op_{B_*(\cc)^\op}} &\ldTwoar&\dTTo>{\Hom_{B_*(\cc)}}
\\
\und\cw^\op & \rTTo^{\und\cw(\_,\1)} & \und\cw
\end{diagram}
It is invertible since so are \(\psi\) and \(\zeta\). It follows that
a right Serre \(\cv\)\n-functor \(S:\cc\to\cc\) induces
a right Serre \(\cw\)\n-functor \(B_*(S):B_*(\cc)\to B_*(\cc)\).
\end{proof}

\subsection{From $\ck$-categories to $\gr$-categories.}
Consider the lax symmetric Monoidal base change functor
\((H^\bullet,\kappa^I):\ck\to\gr\),
 \(X\mapsto H^\bullet X=(H^n X)_{n\in\ZZ}\), where for each
\(I\in\Ob\cs\) the morphism
 \(\kappa^I:\tens^{i\in I}H^\bullet X_i\to H^\bullet\tens^{i\in I}X_i\)
is the K\"unneth map. There is a \(\gr\)\n-functor
\(\und{\wh{H^\bullet}}:H^\bullet_*\und\ck\to\und\gr\),
 \(X\mapsto H^\bullet X\), that acts on morphisms via the map
\[
\ck(X[-n],Y)=H^n\und\ck(X,Y)\to\und\gr(H^\bullet X,H^\bullet Y)^n
=\prod_{d\in\ZZ}\und\kMod(H^{d-n}X,H^dY)
\]
which sends the homotopy class of a chain map \(f:X[-n]\to Y\) to
\((H^d(f))_{d\in\ZZ}\). Note that \(H^\bullet\) preserves the unit
object, therefore there is a natural transformation
\begin{diagram}
H^\bullet_*\und\ck^\op &\rTTo^{H^\bullet_*\und\ck(\_,\kk)}
&H^\bullet_*\und\ck
\\
\dTTo<{(\und{\wh{H^\bullet}})^\op} &\ldTwoar^\zeta
&\dTTo>{\und{\wh{H^\bullet}}}
\\
\und\gr^\op &\rTTo^{\und\gr(\_,\kk)} &\und\gr
\end{diagram}
Explicitly, the map
\(\zeta_X=\und{\wh{H^\bullet}}_{X,\kk}:H^\bullet(\und\ck(X,\kk))\to\und\gr(H^\bullet
X,H^\bullet\kk)=\und\gr(H^\bullet X,\kk)\) is given by its components
\[
\ck(X[-n],\kk)=H^n\und\ck(X,\kk)\to\und\gr(H^\bullet
X,\kk)^n=\und\kMod(H^{-n}X,\kk),\quad f\mapsto H^0(f).
\]
In general, \(\zeta\) is not invertible. However, if \(\kk\) is a
field, \(\zeta\) is an isomorphism. In fact, in this case
\(H^\bullet:\ck\to\gr\) is an equivalence. A quasi-inverse is given by
the functor \(F:\gr\to\ck\) which equips a graded \(\kk\)\n-module with
the trivial differential.

\begin{corollary}\label{cor-Serre-Hbullet-Serre}
Suppose \(\kk\) is a field. Let \(S:\cc\to\cc\) be a (right) Serre
\(\ck\)\n-functor. Then
 \(H^\bullet_*(S):H^\bullet_*(\cc)\to H^\bullet_*(\cc)\) is a (right)
Serre \(\gr\)\n-functor. Moreover, \(H^\bullet_*\) reflects (right) Serre
functors: if \(H^\bullet_*(\cc)\) admits a (right) Serre
\(\gr\)\n-functor, then \(\cc\) admits a (right) Serre \(\ck\)\n-functor.
\end{corollary}

\begin{proof}
The first assertion follows from
\propref{prop-K-Serre-exists-iff-k-Serre}. For the proof of the second,
note that the symmetric Monoidal functor \(F:\gr\to\ck\) induces a
symmetric Monoidal \(\Cat\)\n-functor \(F_*:\wh\grCat\to\wh\KCat\). The
corresponding \(\ck\)\n-functor \(\und{\wh{F}}:F_*\und\gr\to\und\ck\)
acts  as the identity on morphisms (the complex \(\und\ck(FX,FY)\)
carries the trivial differential and coincides with \(\und\gr(X,Y)\) as
a graded \(\kk\)\n-module). Furthermore, \(F\) preserves the unit
object, therefore \propref{prop-K-Serre-exists-iff-k-Serre} applies. It
follows that if \(\bar S:H^\bullet_*(\cc)\to H^\bullet_*(\cc)\) is a
right Serre \(\gr\)\n-functor, then
 \(F_*(\bar S):F_*H^\bullet_*(\cc)\to F_*H^\bullet_*(\cc)\) is a right
Serre \(\ck\)\n-functor. Since the \(\ck\)\n-category
 \(F_* H^\bullet_*(\cc)\) is isomorphic to \(\cc\), the right Serre
\(\ck\)\n-functor \(F_*(\bar S)\) translates to a right Serre
\(\ck\)\n-functor on \(\cc\).
\end{proof}

\subsection{From $\gr$-categories to $\kk$-categories.}
Consider a lax symmetric Monoidal base change functor
\((N,\nu^I):\gr\to\kMod\), \(X=(X^n)_{n\in\ZZ}\mapsto X^0\), where for
each \(I\in\Ob\cs\) the map \(\nu^I:\tens^{i\in I}NX_i=\tens^{i\in
I}X^0_i\to N\tens^{i\in I}X_i=\oplus_{\sum_{i\in I}n_i=0}X^{n_i}_i\) is
the natural embedding. The \(\kk\)\n-functor \(\und{\wh
N}:N_*\und\gr\to\und\kMod\), \(X\mapsto NX=X^0\), acts on morphisms via
the projection
\begin{multline*}
N\und\gr(X,Y)=\und\gr(X,Y)^0=\prod_{d\in\ZZ}\und\kMod(X^d,Y^d)
\\
\to\und\kMod(X^0,Y^0)=\und\kMod(NX,NY).
\end{multline*}
The functor \(N\) preserves the unit object, therefore there exists a
natural transformation
\begin{diagram}
N_*\und\gr^\op & \rTTo^{N_*\und\gr(\_,\kk)} & N_*\und\gr
\\
\dTTo<{(\und{\wh N})^\op} &\ldTwoar^\zeta & \dTTo>{\und{\wh N}}
\\
\und\kMod^\op & \rTTo^{\und\kMod(\_,\kk)} & \und\kMod
\end{diagram}
Explicitly, the map \(\zeta_X=\und{\wh N}_{X,\kk}\) is the identity map
\[N\und\gr(X,\kk)=\und\gr(X,\kk)^0\to\und\kMod(X^0,\kk)=\und\kMod(NX,\kk).\]

\begin{corollary}[to \propref{prop-K-Serre-exists-iff-k-Serre}]
 \label{cor-gr-Serre-k-Serre}
Suppose \(S:\cc\to\cc\) is a right Serre \(\gr\)\n-functor. Then
\(N_*(S):N_*(\cc)\to N_*(\cc)\) is a right Serre \(\kk\)\n-functor.
\end{corollary}

If \(N_*(\cc)\) possess a right Serre \(\kk\)\n-functor, it does not
imply, in general, that $\cc$ has a right Serre \(\gr\)\n-functor.
However, this will be the case if $\cc$ is closed under shifts, as
explained in the next section.

\subsection{Categories closed under shifts.}
As in \cite[Chapter~10]{BesLyuMan-book} denote by $\cz$ the following
algebra (strict monoidal category) in the symmetric monoidal category
of $\dg$\n-categories, $\ck$\n-categories or $\gr$\n-categories. As a
graded quiver $\cz$ has $\Ob\cz=\ZZ$ and $\cz(m,n)=\kk[n-m]$. In the
first two cases $\cz$ is supplied with zero differential. Composition
in the category $\cz$ comes from the multiplication in $\kk$:
\begin{align*}
\mu_\cz: \cz(l,m)\tens_\kk\cz(m,n) =\kk[m-l]\tens_\kk\kk[n-m]
&\to \kk[n-l] =\cz(l,n),
\\
1s^{m-l}\tens1s^{n-m} &\mapsto 1s^{n-l}.
\end{align*}
The elements \(1\in\kk=\cz(n,n)\) are identity morphisms of $\cz$.

The object $\cz$ of \((\VCat,\boxt)\) (where $\cv$ is $\dg$, $\ck$ or
$\gr$) is equipped with an algebra (a strict monoidal category)
structure, given by multiplication -- the $\cv$\n-functor
\begin{align*}
\tens_\psi: \cz\boxt\cz \to \cz, \qquad m\times n &\mapsto m+n,
\\
\tens_\psi: (\cz\boxt\cz)(n\times m,k\times l)
= \cz(n,k)\tens\cz(m,l) &\to \cz(n+m,k+l),
\\
1s^{k-n}\tens1s^{l-m} &\mapsto (-1)^{k(m-l)}s^{k+l-n-m}.
\end{align*}

Therefore, for the three mentioned $\cv$ the functor
\(-\boxt\cz:\VCat\to\VCat\), \(\cc\mapsto\cc\boxt\cz\), is a monad. It
takes a $\cv$\n-category $\cc$ to the $\cv$\n-category $\cc\boxt\cz$
with the set of objects \(\Ob\cc\boxt\cz=\Ob\cc\times\ZZ\) and with the
graded modules of morphisms
\((\cc\boxt\cz)((X,n),(Y,m))=\cc(X,Y)\tens\kk[m-n]\). The composition
is given by the following morphism in $\cv$:
\begin{multline}
\mu_{\cc\boxt\cz} =\bigl[
(\cc\boxt\cz)((X,n),(Y,m))\tens(\cc\boxt\cz)((Y,m),(Z,p)) =
\\
\cc(X,Y)\tens\kk[m-n]\tens\cc(Y,Z)\tens\kk[p-m] \rTTo^{1\tens c\tens1}
\cc(X,Y)\tens\cc(Y,Z)\tens\kk[m-n]\tens\kk[p-m]
\\
\rTTo^{\mu_\cc\tens\mu_\cz} \cc(X,Z)\tens\kk[p-n]
=(\cc\boxt\cz)((X,n),(Z,p)) \bigr].
\label{eq-muCZ-1c1-muCmuZ}
\end{multline}
The unit \(u_\sh=\id\boxt\1_\cz:\cc\to\cc\boxt\cz\) of the monad
\(-\boxt\cz:\VCat\to\VCat\) is the natural embedding \(X\mapsto(X,0)\)
bijective on morphisms. Here \(\1_\cz:\kk\to\cz\), \(*\mapsto0\) is the
unit of the algebra $\cz$, whose source is the graded category $\kk$
with one object.

The $\cv$\n-category \(\cc\boxt\cz\) admits an isomorphic form
\(\cc^\sh\) whose set of objects is \(\Ob\cc^\sh=\Ob\cc\times\ZZ\),
likewise \(\cc\boxt\cz\). The graded $\kk$\n-modules of morphisms are
$\cc^\sh((X,n),(Y,m))=\cc(X,Y)[m-n]$. This graded quiver is identified
with \(\cc\boxt\cz\) via the isomorphism
\begin{equation*}
\text\ss=\bigl[ \cc(X,Y)[m-n] \rTTo^{s^{n-m}} \cc(X,Y)
\rTTo^{\lambda^{\msf{I\lar.}}(1\tens s^{m-n})} \cc(X,Y)\tens\cz(n,m)
\bigr]
\end{equation*}
in \cite[Chapter~10]{BesLyuMan-book}. Therefore, in the cases of
\(\cv=\Com\) or \(\cv=\ck\) the graded $\kk$\n-module
\(\cc^\sh((X,n),(Y,m))\) is equipped with the differential
\((-1)^{m-n}s^{n-m}d_\cc s^{m-n}\). Multiplication in \(\cc^\sh\) is
found from \eqref{eq-muCZ-1c1-muCmuZ} as
\begin{multline}
\mu_{\cc^\sh} =\bigl[\cc(X,Y)[m-n]\tens\cc(Y,Z)[p-m]
\rTTo^{\text\ss\tens\text\ss}
\cc(X,Y)\tens\kk[m-n]\tens\cc(Y,Z)\tens\kk[p-m]
\\
\hfill \rTTo^{\mu_{\cc\boxt\cz}} \cc(X,Z)\tens\kk[p-n]
\rTTo^{\text\ss^{-1}} \cc(X,Z)[p-n] \bigr] \quad
\\
\quad =\bigl[\cc(X,Y)[m-n]\tens\cc(Y,Z)[p-m]
 \rTTo^{(s^{m-n}\tens s^{p-m})^{-1}} \cc(X,Y)\tens\cc(Y,Z) \hfill
\\
\rTTo^{\mu_\cc} \cc(X,Z) \rTTo^{s^{p-n}} \cc(X,Z)[p-n] \bigr].
\label{eq-muC[]-ss1-muC-s}
\end{multline}

\begin{definition}\label{def-closed-under-shifts}
We say that a $\cv$\n-category $\cc$ is \emph{closed under shifts} if
every object $(X,n)$ of $\cc^\sh$ is isomorphic in $\cc^\sh$ to some
object $(Y,0)$, $Y=X[n]\in\Ob\cc$.
\end{definition}

Clearly, \(-^\sh:\VCat\to\VCat\) is also a monad, whose unit
\(u_\sh:\cc\to\cc^\sh\) is the natural embedding \(X\mapsto(X,0)\)
identity on morphisms. Immediately one finds that a $\cv$\n-category
$\cc$ is closed under shifts if and only if the functor
\(u_\sh:\cc\to\cc^\sh\) is an equivalence.

The lax symmetric Monoidal base change functor
\((H^\bullet,\kappa^I):\ck\to\gr\) gives, in particular, the K\"unneth
functor
 \(\kappa:H^\bullet\cc\boxt H^\bullet\cz\to H^\bullet(\cc\boxt\cz)\),
identity on objects. It is an isomorphism of $\gr$\n-categories because
$\cz(m,n)=\kk[n-m]$ are flat graded $\kk$\n-modules. Clearly, $\cz$
coincides with \(H^\bullet\cz\) as a graded $\kk$\n-quiver, hence we
have the isomorphism
 \(\kappa:(H^\bullet\cc)\boxt\cz\to H^\bullet(\cc\boxt\cz)\).
Equivalently we may write the isomorphism
 \((H^\bullet\cc)^\sh\simeq H^\bullet(\cc^\sh)\). From the lax
monoidality of \((H^\bullet,\kappa^I)\) we deduce the following
equation:
\[ H^\bullet(u_\sh) = \bigl[ H^\bullet\cc \rTTo^{u_\sh}
(H^\bullet\cc)\boxt\cz \rTTo^\kappa_\sim H^\bullet(\cc\boxt\cz) \bigr].
\]
Therefore, if $\cc$ is a $\ck$\n-category closed under shifts, then
$H^\bullet\cc$ is a $\gr$\n-category closed under shifts.

For a \(\gr\)\n-category \(\cc\), the components of the graded
\(\kk\)\n-module \(\cc(X,Y)\) are denoted by \(\cc(X,Y)^n=\cc^n(X,Y)\),
\(X,Y\in\Ob\cc\), \(n\in\ZZ\). The \(\kk\)\n-category \(N_*(\cc)\) is
denoted by \(\cc^0\).

\begin{proposition}\label{pro-Serre0-gr-Serre}
Let \(\cc\) be a \(\gr\)\n-category closed under shifts. Suppose
\(S^0:\cc^0\to \cc^0\) is a right Serre \(\kk\)\n-functor. Then there
exists a right Serre \(\gr\)\n-functor \(S:\cc\to\cc\) such that
\(N_*(S)=S^0\).
\end{proposition}

\begin{proof}
Let
\(\psi^0=(\psi^0_{X,Y}:\cc^0(X,YS^0)\to\und\kMod(\cc^0(Y,X),\kk))_{X,Y\in\Ob\cc}\)
be a natural isomorphism. Let
\(\phi^0_{X,Y}:\cc^0(Y,X)\tens\cc^0(X,YS)\to\kk\), \(X,Y\in\Ob\cc\),
denote the corresponding pairings from~\eqref{equ-def-pairing-phi}.
Define trace functionals \(\tau^0_X:\cc^0(X,XS)\to\kk\),
\(X\in\Ob\cc\), by formula~\eqref{equ-def-tau}. We are going to apply
\propref{prop-traces-determine-Serre-functor}. For this we need to specify a
map \(\Ob S:\Ob\cc\to\Ob\cc\) and trace functionals
\(\tau_X:\cc(X,XS)\to\kk\), \(X\in\Ob\cc\). Set \(\Ob S=\Ob S^0\).
Let the \(0\)\n-th component of \(\tau_X\) be equal to the map
\(\tau^0_X\), the other components necessarily vanish since \(\kk\) is
concentrated in degree \(0\). Let us prove that the pairings
\(\phi_{X,Y}\) given by \eqref{equ-def-pairing-phi} are perfect. For
\(n\in\ZZ\), the restriction of \(\phi_{X,Y}\) to the summand
\(\cc^n(Y,X)\tens\cc^{-n}(X,YS)\) is given by
\[
\phi_{X,Y}=\bigl[
\cc^n(Y,X)\tens\cc^{-n}(X,YS)\rTTo^{\mu_\cc}\cc^0(Y,YS)\rTTo^{\tau^0_Y}\kk
\bigr].
\]
It can be written as follows:
\begin{multline*}
\phi_{X,Y}=\bigl[ \cc^n(Y,X)\tens\cc^{-n}(X,YS)
=\cc^\sh((Y,0),(X,n))^0\tens\cc^\sh((X,n),(YS,0))^0
\\
\rTTo^{(-)^n\mu_{\cc^\sh}} \cc^\sh((Y,0),(YS,0))^0
=\cc^0(Y,YS)\rTTo^{\tau^0_Y}\kk \bigr].
\end{multline*}
Since \(\cc\) is closed under shifts, there exist an object
\(X[n]\in\Ob\cc\) and an isomorphism \(\alpha:(X,n)\to(X[n],0)\) in
\(\cc^\sh\). Using associativity of \(\mu_{\cc^\sh}\), we obtain:
\begin{multline*}
\phi_{X,Y}=\bigl[ \cc^n(Y,X)\tens\cc^{-n}(X,YS)
=\cc^\sh((Y,0),(X,n))^0\tens\cc^\sh((X,n),(YS,0))^0
\\
\rTTo^{\cc^\sh(1,\alpha)^0\tens\cc^\sh(\alpha^{-1},1)^0}
\cc^\sh((Y,0),(X[n],0))^0\tens\cc^\sh((X[n],0),(YS,0))^0
\\
\hfill\rTTo^{(-)^n\mu_{\cc^\sh}} \cc^\sh((Y,0),(YS,0))^0 =\cc^0(Y,YS)
\rTTo^{\tau^0_Y}\kk \bigr]\quad
\\
\quad=\bigl[ \cc^n(Y,X)\tens\cc^{-n}(X,YS)
=\cc^\sh((Y,0),(X,n))^0\tens\cc^\sh((X,n),(YS,0))^0\hfill
\\
\rTTo^{\cc^\sh(1,\alpha)^0\tens\cc^\sh(\alpha^{-1},1)^0}
\cc^\sh((Y,0),(X[n],0))^0\tens\cc^\sh((X[n],0),(YS,0))^0
\\
\hfill=\cc^0(Y,X[n])\tens\cc^0(X[n],YS) \rTTo^{(-)^n\mu_{\cc^0}}
\cc^0(Y,YS)\rTTo^{\tau^0_Y}\kk \bigr]\quad
\\
\quad=\bigl[ \cc^n(Y,X)\tens\cc^{-n}(X,YS)
=\cc^\sh((Y,0),(X,n))^0\tens\cc^\sh((X,n),(YS,0))^0\hfill
\\
\rTTo^{\cc^\sh(1,\alpha)^0\tens\cc^\sh(\alpha^{-1},1)^0}
\cc^\sh((Y,0),(X[n],0))^0\tens\cc^\sh((X[n],0),(YS,0))^0
\\
=\cc^0(Y,X[n])\tens\cc^0(X[n],YS) \rTTo^{(-)^n\phi^0_{X[n],Y}}\kk
\bigr].
\end{multline*}
Since \(\phi^0_{X[n],Y}\) is a perfect pairing and the maps
\(\cc^\sh(1,\alpha)^0\) and \(\cc^\sh(\alpha^{-1},1)^0\) are
invertible, the pairing \(\phi_{X,Y}\) is perfect as well. Indeed, it
is easy to see that the corresponding maps \(\psi^{-n}_{X,Y}\) and
\(\psi^0_{X[n],Y}\) are related as follows:
\begin{multline*}
\psi^{-n}_{X,Y}=\bigl[ \cc^{-n}(X,YS) \rTTo^{\cc^\sh(\alpha^{-1},1)^0}
\cc^0(X[n],YS) \rTTo^{(-)^n\psi^0_{X[n],Y}}
\\
\Hom_\kk(\cc^0(Y,X[n]),\kk) \rTTo^{\Hom_\kk(\cc^\sh(1,\alpha))^0,1)}
\Hom_\kk(\cc^n(Y,X),\kk) \bigr].
\end{multline*}
\propref{prop-traces-determine-Serre-functor} implies that there is a
right Serre \(\gr\)\n-functor \(S:\cc\to\cc\). Its components are
determined unambiguously by equation~\eqref{dia-S1-c-mutau-mutau}.
Applying the multifunctor \(\wh N\) to it we find that the functor
\(N_*(S):\cc^0\to\cc^0\) satisfies the same equation the functor
\(S^0:\cc^0\to\cc^0\) does. By uniqueness of the solution, \(N_*(S)=S^0\).
\end{proof}

\section{$A_\infty$-categories and $\ck$-categories}
In this section we recall and deepen the relationship between
\ainf-categories and $\ck$\n-categories. It is implemented by a
multifunctor \(\kf:\Ainftyu\to\wh\KCat\) from
\cite[Chapter~13]{BesLyuMan-book}, where $\KCat$ is the symmetric
Monoidal category of $\ck$\n-categories and $\ck$\n-functors, and
\(\wh\KCat\) is the corresponding symmetric multicategory. This
multifunctor extends to non-unital \ainf-categories as a sort of
multifunctor \(\kf:\Ainfty\to\wh{\KCat^{nu}}\), where $\KCat^{nu}$ is
the symmetric Monoidal category of non-unital $\ck$\n-categories and
$\ck$\n-functors, and \(\wh{\KCat^{nu}}\) is the corresponding
symmetric multicategory [loc.~cit.].

\subsection{Opposite $A_\infty$-categories.}
Recall the following definitions from
\cite[Appendix~A]{math.CT/0306018}. Let $\ca$ be a graded
$\kk$\n-quiver. Then its \emph{opposite quiver} $\ca^\op$ is defined as
the quiver with the same class of objects $\Ob\ca^\op=\Ob\ca$, and with
graded $\kk$\n-modules of morphisms $\ca^\op(X,Y)=\ca(Y,X)$.

Let $\gamma:Ts\ca^\op\to Ts\ca$ denote the following anti-isomorphism
of coalgebras and algebras (free categories):
\begin{multline}
\gamma= (-1)^k\omega^0_c: s\ca^\op(X_0,X_1)\tdt s\ca^\op(X_{k-1},X_k)
\\
\to s\ca(X_k,X_{k-1})\tdt s\ca(X_1,X_0),
\label{eq-gamma-anti}
\end{multline}
where $\omega^0=\bigl(
\begin{smallmatrix}
1 & 2 & \dots & k-1 & k \\
k &k-1& \dots &  2  & 1
\end{smallmatrix}
\bigr)\in\SSS_k$.
Clearly,
$\gamma\Delta_0=\Delta_0(\gamma\tens\gamma)c=\Delta_0c(\gamma\tens\gamma)$,
which is the anti-isomorphism property. Notice also that
$(\ca^\op)^\op=\ca$ and $\gamma^2=\id$.

When $\ca$ is an \ainf-category with the codifferential
$b:Ts\ca\to Ts\ca$, then $\gamma b\gamma:Ts\ca^\op\to Ts\ca^\op$ is
also a codifferential. Indeed,
\begin{align*}
\gamma b\gamma\Delta_0 &= \gamma b\Delta_0 c(\gamma\tens\gamma)
=\gamma\Delta_0(1\tens b+b\tens1)c(\gamma\tens\gamma)
\\
&= \Delta_0(\gamma\tens\gamma)c(1\tens b+b\tens1)c(\gamma\tens\gamma)
=\Delta_0(\gamma b\gamma\tens1+1\tens\gamma b\gamma).
\end{align*}
The \emph{opposite \ainf-category} $\ca^\op$ to an \ainf-category $\ca$
is the opposite quiver, equipped with the codifferential $b^\op=\gamma
b\gamma:Ts\ca^\op\to Ts\ca^\op$. The components of $b^\op$ are computed
as follows:
\begin{multline}
b^\op_k = (-)^{k+1}\bigl[s\ca^\op(X_0,X_1)\tdt s\ca^\op(X_{k-1},X_k)
\rTTo^{\omega^0_c}
\\
s\ca(X_k,X_{k-1})\tdt s\ca(X_1,X_0) \rTTo^{b_k}
s\ca(X_k,X_0) =s\ca^\op(X_0,X_k) \bigr].
\label{equ-codiff-opposite}
\end{multline}
The sign $(-1)^k$ in \eqref{eq-gamma-anti} ensures that the above
definition agrees with the definition of the opposite usual category,
meaning that, for an arbitrary \ainf-category \(\ca\),
\(\kf(\ca^\op)=(\kf \ca)^\op\). Indeed, clearly, both categories have
\(\Ob\ca\) as the set of objects. Furthermore, for each pair of objects
\(X,Y\in\Ob\ca\),
\[
m^\op_1=sb^\op_1s^{-1}=sb_1s^{-1}=m_1:\ca^\op(X,Y)=\ca(Y,X)\to\ca(Y,X)=\ca^\op(X,Y),
\]
therefore
\((\kf\ca^\op)(X,Y)=(\ca^\op(X,Y),m^\op_1)=(\ca(Y,X),m_1)=(\kf\ca)^\op(X,Y)\).
Finally, the compositions in both categories coincide:
\[ \mu_{\kf\ca^\op} =m^\op_2 =(s\tens s)b^\op_2s^{-1}
=-(s\tens s)cb_2s^{-1} =c(s\tens s)b_2s^{-1} =cm_2 =\mu_{(\kf\ca)^\op}.
\]
In particular, it follows that \(\ca^\op\) is unital if so is \(\ca\),
with the same unit elements.

For an arbitrary \ainf-functor $f:\boxt^{i\in\mb n}Ts\ca_i\to Ts\cb$
there is another \ainf-functor $f^\op$ defined by the commutative
square
\begin{diagram}[width=4em]
\boxtimes^{i\in\mb n}Ts\ca_i &\rTTo^f &Ts\cb
\\
\dTTo<{\boxtimes^{\mb n}\gamma} &&\dTTo>\gamma
\\
\boxtimes^{i\in\mb n}Ts\ca_i^\op &\rTTo^{f^\op} &Ts\cb^\op
\end{diagram}
Since \(\gamma^2=\id\), the \ainf-functor \(f^\op\) is found as the
composite
\[
f^\op=\bigl[\boxt^{i\in\mb n}Ts\ca^\op \rTTo^{\boxt^{\mb n}\gamma}
\boxt^{i\in\mb n}Ts\ca_i \rTTo^f Ts\cb \rTTo^\gamma Ts\cb^\op\bigr].
\]

A non-unital $\ck$\n-functor
 \(\kf f:\boxt^{i\in\mb n}\kf\ca_i\to\kf\cb\) is associated with $f$ in
\cite[Chapter~13]{BesLyuMan-book}. This $\ck$\n-functor acts on objects
in the same way as $f$. It is determined by the components \(f_{e_j}\),
\(e_j=(0,\dots,0,1,0,\dots,0)\):
\begin{multline*}
\kf f = \bigl[ \tens^{j\in\mb n}\kf\ca_j(X_j,Y_j)
\rTTo^{\tens^{j\in\mb n}(sf_{e_j}s^{-1})}
\\
\tens^{j\in\mb n}\kf\cb\bigl(((Y_i)_{i<j},(X_i)_{i\ge j})f,
 ((Y_i)_{i\le j},(X_i)_{i>j})f\bigr)
\rTTo^{\mu_{\kf\cb}^{\mb n}}
\kf\cb\bigl((X_i)_{i\in\mb n}f,(Y_i)_{i\in\mb n}f\bigr) \bigr],
\end{multline*}
The case of $n=1$ was considered in
\cite[Proposition~8.6]{Lyu-AinfCat}. According to \cite{BesLyuMan-book}
an \ainf-functor $f$ is called \emph{unital}, if the $\ck$\n-functor
\(\kf f\) is unital. The set of unital \ainf-functors
 \(\boxt^{i\in\mb n}Ts\ca_i\to Ts\cb\) is denoted
\(\Ainftyu((\ca_i)_{i\in\mb n};\cb)\). The assignment
\(\ca\mapsto\kf\ca\), \(f\mapsto\kf f\) gives a multifunctor
\(\kf:\Ainftyu\to\wh\KCat\), see \cite[Chapter~13]{BesLyuMan-book}.

\begin{lemma}\label{lem-kf-commutes-with-op-on-functors}
For an arbitrary \ainf-functor \(f:\boxt^{i\in\mb n}Ts\ca_i\to Ts\cb\),
the \(\ck\)\n-functors
\(\kf f^\op,(\kf f)^\op:\boxt^{i\in\mb n}\kf\ca^\op_i\to\kf\cb^\op\) coincide.
\end{lemma}

\begin{proof}
The case \(n=1\) is straightforward. We provide a proof in the case
\(n=2\), which we are going to use later.

Let \(f:Ts\ca\boxt Ts\cb\to Ts\cc\) be an \ainf-functor. The components
of \(f^\op\) are given by
\begin{multline}
f^\op_{kn}=(-)^{k+n-1}\bigl[
s\ca^\op(X_0,X_1)\tdt s\ca^\op(X_{k-1},X_k)\tens
\\
\hfill \tens s\cb^\op(U_0,U_1)\tdt s\cb^\op(U_{n-1},U_n) \quad
\\
=s\ca(X_1,X_0)\tdt s\ca(X_k,X_{k-1})\tens s\cb(U_1,U_0)\tdt
s\cb(U_n,U_{n-1})
\\
\rTTo^{\pi^{kn}_c}s\ca(X_k,X_{k-1})\tdt
s\ca(X_1,X_0)\tens s\cb(U_n,U_{n-1})\tdt
s\cb(U_1,U_0)\\
\rTTo^{f_{kn}}s\cc((X_k,U_n)f,(X_0,U_0)f)=s\cc^\op((X_0,U_0)f,(X_k,U_n)f)
\bigr],
\label{equ-A8-functor-op-components}
\end{multline}
where
\(\pi^{kn}=\bigl(
\begin{smallmatrix}
1 & 2 & \dots & k & k+1 & k+2 & \dots & k+n\\
k & k-1 & \dots & 1 & k+n & k+n-1 & \dots & k+1
\end{smallmatrix}
\bigr)\in\SSS_{k+n}\), and \(\pi^{kn}_c\) is the corresponding signed
permutation.

Clearly, both \(\kf f^\op\) and \((\kf f)^\op\) act as \(\Ob f\) on
objects. Let \(X,Y\in\Ob\ca\), \(U,V\in\Ob\cb\). Then
\begin{multline*}
\kf f^\op=\bigl[
\ca^\op(X,Y)\tens\cb^\op(U,V)\rTTo^{sf^\op_{10}s^{-1}\tens
sf^\op_{01}s^{-1}}\\
\cc^\op((X,U)f,(Y,U)f)\tens\cc^\op((Y,U)f,(Y,V)f)\rTTo^{\mu_{\kf\cc^\op}}
\cc^\op((X,U)f,(Y,V)f)
\bigr].
\end{multline*}
By \eqref{equ-A8-functor-op-components},
\begin{align*}
f^\op_{10}&=f_{10}: s\ca^\op(X,Y) \to
s\cc((Y,U)f,(X,U)f)=s\cc^\op((X,U)f,(Y,U)f),
\\
f^\op_{01}&=f_{01}: s\cb^\op(U,V) \to
s\cc((Y,V)f,(Y,U)f)=s\cc^\op((Y,U)f,(Y,V)f),
\end{align*}
therefore
\begin{multline*}
\kf f^\op=\bigl[
\ca(Y,X)\tens\cb(V,U) \rTTo^{sf_{10}s^{-1}\tens sf_{01}s^{-1}}
\\
\cc((Y,U)f,(X,U)f)\tens\cc((Y,V)f,(Y,U)f) \rTTo^c
\\
\hfill \cc((Y,V)f,(Y,U)f)\tens\cc((Y,U)f,(X,U)f) \rTTo^{\mu_{\kf\cc}}
\cc((Y,V)f,(X,U)f) \bigr]\quad
\\
\quad=\bigl[ \ca(Y,X)\tens\cb(V,U) \rTTo^c \cb(V,U)\tens\ca(Y,X)
\rTTo^{sf_{01}s^{-1}\tens sf_{10}s^{-1}} \hfill
\\
\cc((Y,V)f,(Y,U)f)\tens\cc((Y,U)f,(X,U)f) \rTTo^{\mu_{\kf\cc}}
\cc((Y,V)f,(X,U)f) \bigr].
\end{multline*}
Further,
\begin{multline*}
(\kf f)^\op=\bigl[
\ca(Y,X)\tens\cb(V,U) \rTTo^{sf_{10}s^{-1}\tens sf_{01}s^{-1}}
\\
\cc((Y,V)f,(X,V)f)\tens\cc((X,V)f,(X,U)f) \rTTo^{\mu_{\kf\cc}}
\cc((Y,V)f,(X,U)f) \bigr].
\end{multline*}
We must therefore prove the following equation in \(\ck\):
\begin{multline}
\bigl[
\ca(Y,X)\tens\cb(V,U)\rTTo^c\cb(V,U)\tens\ca(Y,X)
\rTTo^{sf_{01}s^{-1}\tens sf_{10}s^{-1}}\hfill\\
\hfill\cc((Y,V)f,(Y,U)f)\tens\cc((Y,U)f,(X,U)f)\rTTo^{\mu_{\kf\cc}}
\cc((Y,V)f,(X,U)f)
\bigr]\quad\\
\quad=\bigl[
\ca(Y,X)\tens\cb(V,U)\rTTo^{sf_{10}s^{-1}\tens sf_{01}s^{-1}}
\cc((Y,V)f,(X,V)f)\tens\cc((X,V)f,(X,U)f)\hfill\\
\rTTo^{\mu_{\kf\cc}}\cc((Y,V)f,(X,U)f)
\bigr].
\label{equ-kf-fop-kf-f-op}
\end{multline}
By definition of an \ainf-functor the equation \(fb=(b\boxt1+1\boxt
b)f:Ts\ca\boxt Ts\cb\to Ts\cc\) holds.  Restricting it to \(s\ca\boxt
s\cb\) and composing with \(\pr_1:Ts\cc\to s\cc\), we obtain
\begin{multline*}
(f_{10}\tens f_{01})b_2+c(f_{01}\tens f_{10})b_2+f_{11}b_1\\
=(1\tens b_1+b_1\tens1)f_{11}:s\ca(Y,X)\tens s\cb(V,U)\to
s\cc((Y,V)f,(X,U)f).
\end{multline*}
Thus, \((f_{10}\tens f_{01})b_2+c(f_{01}\tens f_{10})b_2\) is a
boundary. Therefore,
\[
(s\tens s)(f_{10}\tens f_{01})b_2=c(s\tens s)(f_{01}\tens f_{10})b_2
\]
in \(\ck\). This implies equation~\eqref{equ-kf-fop-kf-f-op}.
\end{proof}

In particular, \(f^\op\) is a unital \ainf-functor if \(f\) is unital.

\begin{proposition}
The correspondences \(\ca\mapsto\ca^\op\), \(f\mapsto f^\op\) define a
symmetric multifunctor \(-^\op:\Ainfty\to\Ainfty\) which restricts to a
symmetric multifunctor \(-^\op:\Ainftyu\to\Ainftyu\).
\end{proposition}

\begin{proof}
Straightforward.
\end{proof}

As an arbitrary symmetric multifunctor between closed multicategories,
\(-^\op\) possesses a closing transformation
 \(\und\op:\und\Ainfty((\ca_i)_{i\in I};\cb)^\op
 \to\und\Ainfty((\ca^\op_i)_{i\in I};\cb^\op)\)
uniquely determined by the following equation in \(\Ainfty\):
\begin{multline}
\bigl[(\ca^\op_i)_{i\in I},\und\Ainfty((\ca_i)_{i\in I};\cb)^\op
\rTTo^{1,\und\op}
 (\ca^\op_i)_{i\in I},\und\Ainfty((\ca^\op_i)_{i\in I};\cb^\op)
\rTTo^{\ev^{\Ainfty}} \cb^\op\bigr]
\\
=(\ev^{\Ainfty})^\op.
\label{equ-op-closing-transform}
\end{multline}
The \ainf-functor \((\ev^{\Ainfty})^\op\) acts on objects in the same
way as \(\ev^{\Ainfty}\). It follows that \((X_i)_{i\in
I}(f)\und\op=(X_i)_{i\in I}f\) for an arbitrary \ainf-functor
\(f:(\ca_i)_{i\in I}\to\cb\) and a family of objects
\(X_i\in\Ob\ca_i\), \(i\in I\). The components
\begin{multline}
(\ev^{\Ainfty})^\op_{(m_i),m}=-\bigl[
\boxt^{i\in I}T^{m_i}s\ca^\op_i\boxt T^ms\und\Ainfty((\ca_i)_{i\in I};\cb)^\op
\rTTo^{\boxt^I(\gamma)_{I}\boxt\gamma}
\\
\boxt^{i\in I}T^{m_i}s\ca_i\boxt T^ms\und\Ainfty((\ca_i)_{i\in I};\cb)
\rTTo^{\ev^{\Ainfty}_{(m_i),m}}s\cb \bigr]
 \label{equ-ev-op}
\end{multline}
vanish unless \(m=0\) or \(m=1\) since the same holds for
\(\ev^{\Ainfty}_{(m_i),m}\). From
equations~\eqref{equ-op-closing-transform} and \eqref{equ-ev-op} we
infer that
\begin{multline*}
(\tens^{i\in I}1^{\tens m_i}\tens\Ob\und\op)\ev^{\Ainfty}_{(m_i),0}
\\
 \quad =-\bigl[
\tens^{i\in I}\tens^{p_i\in\mb{m_i}}s\ca^\op_i(X^i_{p_i-1},X^i_{p_i})
\tens T^0s\und\Ainfty((\ca_i)_{i\in I};\cb)^\op(f,f) \hfill
\\
\simeq\tens^{i\in I}\tens^{p_i\in\mb{m_i}}s\ca^\op_i(X^i_{p_i-1},X^i_{p_i})
\rTTo^{\tens^{i\in I}(-)^{m_i}\omega^0_c}
\\
\tens^{i\in I}\tens^{p_i\in\mb{m_i}}s\ca_i(X^i_{m_i-p_i},X^i_{m_i-p_i+1})
\rTTo^{f_{(m_i)}}s\cb((X^i_0)_{i\in I}f,(X^i_{m_i})_{i\in I}f)
\bigr],
\end{multline*}
therefore \((f)\und\op=f^\op:(\ca^\op_i)_{i\in I}\to\cb^\op\).
Similarly,
\begin{multline*}
(\tens^{i\in I}1^{\tens m_i}\tens\und\op_1)\ev^{\Ainfty}_{(m_i),1}
\\
 \quad =\bigl[
\tens^{i\in I}\tens^{p_i\in\mb{m_i}}s\ca^\op_i(X^i_{p_i-1},X^i_{p_i})\tens
s\und\Ainfty((\ca_i)_{i\in I};\cb)^\op(f,g) \hfill
\\
 \quad \rTTo^{\tens^{i\in I}(-)^{m_i}\omega^0_c\tens1}
\tens^{i\in I}\tens^{p_i\in\mb{m_i}}s\ca_i(X^i_{m_i-p_i},X^i_{m_i-p_i+1})
 \tens s\und\Ainfty((\ca_i)_{i\in I};\cb)(g,f) \hfill
\\
\quad\rTTo^{\tens^{i\in I}\tens^{p_i\in\mb{m_i}}1\tens\pr}
 \tens^{i\in I}\tens^{p_i\in\mb{m_i}}s\ca_i(X^i_{m_i-p_i},X^i_{m_i-p_i+1})
\hfill
\\
\hfill\tens\uCom(\tens^{i\in I}\tens^{p_i\in\mb{m_i}}s\ca_i(X^i_{m_i-p_i},X^i_{m_i-p_i+1}),
s\cb((X^i_{m_i})_{i\in I}g,(X^i_0)_{i\in I}f))\quad
\\
\rTTo^{\ev^{\Com}}
s\cb((X^i_{m_i})_{i\in I}g,(X^i_0)_{i\in I}f)
=s\cb^\op((X^i_0)_{i\in I}f^\op,(X^i_{m_i})_{i\in I}g^\op)
\bigr].
\end{multline*}
It follows that the map \(\und\op_1:s\und\Ainfty((\ca_i)_{i\in
I};\cb)^\op(f,g)\to s\und\Ainfty((\ca^\op_i)_{i\in
I};\cb^\op)(f^\op,g^\op)\) takes an \ainf-transformation \(r:g\to
f:(\ca_i)_{i\in I}\to\cb\) to the opposite \ainf-transformation
\(r^\op\overset{\text{def}}=(r)\und\op_1:f^\op\to
g^\op:(\ca^\op_i)_{i\in I}\to\cb^\op\) with the components
\begin{multline*}
[(r)\und\op_1]_{(m_i)}=(-)^{m_1+\dots+m_n}\bigl[
\tens^{i\in I}\tens^{p_i\in\mb{m_i}}s\ca^\op_i(X^i_{p_i-1},X^i_{p_i})
\rTTo^{\tens^{i\in I}\omega^0_c}
\\
\tens^{i\in I}\tens^{p_i\in\mb{m_i}}s\ca_i(X^i_{m_i-p_i},X^i_{m_i-p_i+1})
\rTTo^{r_{(m_i)}}
s\cb((X^i_{m_i})_{i\in I}g,(X^i_0)_{i\in I}f)
\bigr].
\end{multline*}
The higher components of \(\und\op\) vanish. Similar computations can
be performed in the multicategory \(\Ainftyu\). They lead to the same
formulas for \(\und\op\), which means that the \ainf-functor
\(\und\op\) restricts to a unital \ainf-functor
\(\und\op:\und\Ainftyu((\ca_i)_{i\in
I};\cb)^\op\to\und\Ainftyu((\ca^\op_i)_{i\in I};\cb^\op)\) if the
\ainf-categories \(\ca_i\), \(i\in I\), \(\cb\) are unital.

As an easy application of the above considerations note that if
\(r:f\to g:(\ca_i)_{i\in I}\to\cb\) is an isomorphism of
\ainf-functors, then \(r^\op:g^\op\to f^\op:(\ca^\op_i)_{i\in
I}\to\cb^\op\) is an isomorphism as well.

\subsection{The $\kCat$-multifunctor $\kf$.}
The multifunctor \(\kf:\Ainftyu\to\wh\KCat\) is defined in
\cite[Chapter~13]{BesLyuMan-book}. Here we construct its extension to
natural \ainf-transformations as follows. Let
 \(f,g:(\ca_i)_{i\in I}\to\cb\) be unital \ainf-functors,
\(r:f\to g:(\ca_i)_{i\in I}\to\cb\) a natural \ainf-transformation. It
gives rise to a natural transformation of \(\ck\)\n-functors
 \(\kf r:\kf f\to\kf g:\boxt^{i\in I}\kf\ca_i\to\kf\cb\). Components of
\(\kf r\) are given by
\[\sS{_{(X_i)_{i\in I}}}\kf r=\sS{_{(X_i)_{i\in I}}}r_0s^{-1}:
\kk\to\cb((X_i)_{i\in I}f,(X_i)_{i\in I}g), \quad X_i\in\Ob\ca_i, \quad
i\in I.
\]
Since \(r_0b_1=0\), \(\kf r\) is a chain map. Naturality is expressed
by the following equation in \(\ck\):
\begin{diagram}
\boxt^{i\in I}\kf\ca_i(X_i,Y_i) & \rTTo^{\kf f} & \kf\cb((X_i)_{i\in I}f,(Y_i)_{i\in
I}f)\\
\dTTo<{\kf g} &=& \dTTo>{(1\tens\sS{_{(Y_i)_{i\in I}}}\kf
r)\mu_{\kf\cb}}\\
\kf\cb((X_i)_{i\in I}g,(Y_i)_{i\in I}g) & \rTTo^{(\sS{_{(X_i)_{i\in I}}}\kf
r\tens1)\mu_{\kf\cb}} & \kf\cb((X_i)_{i\in I}f,(Y_i)_{i\in I}g).
\end{diagram}
Associativity of \(\mu_{\kf\cb}\) allows to write it as follows:
\begin{multline*}
\bigl[ \tens^{i\in I}\kf\ca_i(X_i,Y_i) \rTTo^{\tens^{i\in
I}sf_{e_i}s^{-1}\tens\sS{_{(Y_i)_{i\in I}}}r_0s^{-1}}
\\
\tens^{i\in I}\kf\cb(((Y_j)_{j<i},(X_j)_{j\ge i})f,((Y_j)_{j\le
i},(X_j)_{j>i})f)\tens\kf\cb((Y_i)_{i\in I}f,(Y_i)_{i\in I}g)\\
\hfill\rTTo^{\mu^{I\sqcup\mb1}_{\kf\cb}}
\kf\cb((X_i)_{i\in I}f,(Y_i)_{i\in I}g)
\bigr]\quad\\
\quad=\bigl[ \tens^{i\in I}\kf\ca_i(X_i,Y_i)\rTTo^{\sS{_{(X_i)_{i\in
I}}}r_0s^{-1}\tens\tens^{i\in I}sg_{e_i}s^{-1}}\hfill\\
 \kf\cb((X_i)_{i\in
I}f,(X_i)_{i\in I}g)\tens\tens^{i\in I}\kf\cb(((Y_j)_{j<i},(X_j)_{j\ge
i})g,((Y_j)_{j\le i},(X_j)_{j>i})g)\\
\rTTo^{\mu^{\mb1\sqcup I}_{\kf\cb}}\kf\cb((X_i)_{i\in I}f,(Y_i)_{i\in I}g)
\bigr].
\end{multline*}
This equation is a consequence of the following equation in \(\ck\):
\begin{multline*}
(s(f|^{(Y_j)_{j<i},(X_j)_{j>i}}_i)_1s^{-1}\tens\sS{_{(Y_j)_{j\le
i},(X_j)_{j>i}}}r_0s^{-1})\mu_{\kf\cb}\\
=(\sS{_{(Y_j)_{j<i},(X_j)_{j\ge i}}}r_0s^{-1}\tens
s(g|^{(Y_j)_{j<i},(X_j)_{j>i}}_i)_1s^{-1})\mu_{\kf\cb}:\\
\ca(X_i,Y_i)\to\cb(((Y_j)_{j<i},(X_j)_{j\ge i})f,
((Y_j)_{j\le i},(X_j)_{j>i})g),
\end{multline*}
which in turn follows from the equation \((rB_1)_{e_i}=0\):
\begin{multline*}
(sf_{e_i}s^{-1}\tens r_0s^{-1})m_2-(r_0s^{-1}\tens
sg_{e_i}s^{-1})m_2+sr_{e_i}s^{-1}m_1+m_1sr_{e_i}s^{-1}\\
=s[(f_{e_i}\tens r_0)b_2+(r_0\tens g_{e_i})b_2+r_{e_i}b_1+b_1r_{e_i}]s^{-1}=0:\\
\ca(X_i,Y_i)\to\cb(((Y_j)_{j<i},(X_j)_{j\ge i})f,
((Y_j)_{j\le i},(X_j)_{j>i})g).
\end{multline*}

The 2\n-category \(\KCat\) is naturally a symmetric Monoidal
\(\kCat\)\n-category, therefore \(\wh\KCat\) is a symmetric
\(\kCat\)\n-multicategory. According to the general recipe, for each
map \(\phi:I\to J\), the composition in \(\wh\KCat\) is given by the
\(\kk\)\n-linear functor
\begin{multline*}
\mu^{\wh\KCat}_\phi=\bigl[
\boxt^{J\sqcup\mb1}[(\KCat(\boxt^{i\in\phi^{-1}j}\ca_i,\cb_j))_{j\in
J},\KCat(\boxt^{j\in J}\cb_j,\cc)] \rTTo^{\Lambda^\phi_{\kCat}}_\sim
\\
\boxt^{j\in J}\KCat(\boxt^{i\in\phi^{-1}j}\ca_i,\cb_j)\boxt\KCat(\boxt^{j\in J}\cb_j,\cc)
\rTTo^{\boxt^J\boxt1}\\
\KCat(\boxt^{j\in J}\boxt^{i\in\phi^{-1}j}\ca_i,\boxt^{j\in
J}\cb_j)\boxt\KCat(\boxt^{j\in J}\cb_j,\cc)
\rTTo^{\lambda^\phi\cdot-\cdot-}
\KCat(\boxt^{i\in I}\ca_i,\cc)
\bigr].
\end{multline*}
In particular, the action on natural transformations is given by the
map
\begin{multline*}
\tens^{j\in
J}\KCat(\boxt^{i\in\phi^{-1}j}\ca_i,\cb_j)(f_j,g_j)\tens\KCat(\boxt^{j\in
J}\cb_j,\cc)(h,k)\\
\to\KCat(\boxt^{i\in I}\ca_i,\cc)((f_j)_{j\in J}\cdot h,(g_j)_{j\in J}\cdot
k),\quad
\tens^{j\in J}r^j\tens p\mapsto (r^j)_{j\in J}\cdot p,
\end{multline*}
where for each collection of objects \(X_i\in\ca_i\), \(i\in I\),
\begin{multline*}
\sS{_{(X_i)_{i\in I}}}[(r^j)_{j\in J}\cdot p]=\bigl[ \kk
\rTTo^{\lambda^{\emptyset\to J\sqcup\mb1}}_\sim \tens^{J\sqcup\mb1}\kk
\rTTo^{\tens^{j\in
J}\sS{_{(X_i)_{i\in\phi^{-1}j}}}r^j\tens\sS{_{((X_i)_{i\in\phi^{-1}j}g_j)_{j\in
J}}}p}
\\
\tens^{j\in J}\cb((X_i)_{i\in\phi^{-1}j}f_j,(X_i)_{i\in\phi^{-1}j}g_j)
\tens\cc(((X_i)_{i\in\phi^{-1}j}g_j)_{j\in
J}h,((X_i)_{i\in\phi^{-1}j}g_j)_{j\in J}k)
\\
\rTTo^{h\tens1}\cc((X_i)_{i\in I}(f_j)_{j\in J}h,(X_i)_{i\in I}(g_j)_{j\in
J}h)\tens\cc((X_i)_{i\in I}(g_j)_{j\in
J}h,(X_i)_{i\in I}(g_j)_{j\in
J}k)\\
\rTTo^{\mu_\cc}\cc((X_i)_{i\in I}(f_j)_{j\in J}h,(X_i)_{i\in
I}(g_j)_{j\in J}k) \bigr].
\end{multline*}

The base change functor \(H^0:\Ainftyu\to\kCat\) turns the symmetric
\(\Ainftyu\)\n-multicategory \(\und\Ainftyu\) into a symmetric
\(\kCat\)\n-multicategory, which we denote by \(A^u_\infty\). That is,
the objects of \(A^u_\infty\) are unital \ainf-categories, and for each
collection \((\ca_i)_{i\in I}\), \(\cb\) of unital \ainf-categories,
there is a \(\kk\)\n-linear category \(A^u_\infty((\ca_i)_{i\in
I};\cb)=H^0\und\Ainftyu((\ca_i)_{i\in I};\cb)\), whose objects are
unital \ainf-functors, and whose morphisms are equivalence classes of
natural \ainf-transformations. The composition in \(A^u_\infty\) is
given by the \(\kk\)\n-linear functor
\(\mu^{A^u_\infty}_\phi=H^0(\mu^{\und\Ainftyu}_\phi)=H^0(\kf\mu^{\und\Ainftyu}_\phi)\),
where
\begin{multline*}
\kf\mu^{\und\Ainftyu}_\phi=\bigl[ \tens^{j\in
J}\und\Ainftyu((\ca_i)_{i\in\phi^{-1}j};\cb_j)(f_j,g_j)
\tens\und\Ainftyu((\cb_j)_{j\in J};\cc)(h,k)
\\
\rTTo^{\tens^{j\in J}sM_{e_j0}s^{-1}\tens sM_{0\dots01}s^{-1}}
\\
\quad \tens^{j\in J}\und\Ainftyu((\ca_i)_{i\in I};\cc)
(((g_l)_{l<j},(f_l)_{l\ge j})h,((g_l)_{l\le j},(f_l)_{l>j})h) \hfill
\\
\hfill \tens\und\Ainftyu((\ca_i)_{i\in I};\cc)
 ((g_j)_{j\in J}h,(g_j)_{j\in J}k) \quad
\\
\rTTo^{\mu^{J\sqcup\mb1}_{\kf\und\Ainftyu((\ca_i)_{i\in I};\cc)}}
\und\Ainftyu((\ca_i)_{i\in I};\cc)((f_j)_{j\in J}h,(g_j)_{j\in J}k)
\bigr].
\end{multline*}

\begin{proposition}
There is a symmetric \(\kCat\)\n-multifunctor
\(\kf:A^u_\infty\to\wh\KCat\).
\end{proposition}

\begin{proof}
It remains to prove compatibility of \(\kf\) with
\(\mu^{A^u_\infty}_\phi\) on the level of transformations. Let \(r^j\in
s\und\Ainftyu((\ca_i)_{i\in\phi^{-1}j};\cb_j)(f_j,g_j)\), \(j\in J\),
\(p\in s\und\Ainftyu((\cb_j)_{j\in J};\cc)\) be natural
\ainf-transformations. Then \(((r^js^{-1})_{j\in
J},ps^{-1})\mu^{A^u_\infty}_\phi\) is the equivalence class of the
following \ainf-transformation:
\[
[\tens^{j\in J}((g_l)_{l<j},r^j,(f_l)_{l>j},h)M_{e_j0}s^{-1}\tens((g_j)_{j\in
J},p)M_{0\dots01}s^{-1}]\mu^{J\sqcup\mb1}_{\kf\und\Ainftyu((\ca_i)_{i\in I};\cc)}.
\]
In order to find \(\kf[((r^js^{-1})_{j\in
J},ps^{-1})\mu^{A^u_\infty}_\phi]\) we need the 0-th components of the
above expression. Since \([(t\tens q)B_2]_0=(t_0\tens q_0)b_2\), for
arbitrary composable \ainf-transformations \(t\) and \(q\), it follows that
\begin{multline*}
\sS{_{(X_i)_{i\in I}}}\kf[((r^js^{-1})_{j\in J},ps^{-1})\mu^{A^u_\infty}_\phi]
\\
\quad =(\tens^{j\in J}
\sS{_{(X_i)_{i\in\phi^{-1}j}}}[((g_l)_{l<j},r^j,(f_l)_{l>j},h)M_{e_j0}]_0s^{-1}
\tens \hfill
\\
\hfill \tens \sS{_{(X_i)_{i\in I}}}[((g_j)_{j\in
J},p)M_{0\dots01}]_0s^{-1}) \mu^{J\sqcup\mb1}_{\kf\cc} \quad
\\
\quad =(\tens^{j\in J}
\sS{_{(X_i)_{i\in\phi^{-1}j}}}r^j_0s^{-1}s(h|_j^{((X_i)_{i\in\phi^{-1}l}g_l)_{l<j},
((X_i)_{i\in\phi^{-1}l}f_l)_{l>j}})_1s^{-1} \tens \hfill
\\
\tens\sS{_{((X_i)_{i\in\phi^{-1}j}g_j)_{j\in J}}}p_0s^{-1})
\mu^{J\sqcup\mb1}_{\kf\cc}.
\end{multline*}
By associativity of \(\mu_{\kf\cc}\), this equals
\begin{multline*}
(\tens^{j\in J}\sS{_{(X_i)_{i\in\phi^{-1}j}}}r^j_0s^{-1}\cdot
\tens^{j\in J}s(h|_j^{((X_i)_{i\in\phi^{-1}l}g_l)_{l<j},
((X_i)_{i\in\phi^{-1}l}f_l)_{l>j}})_1s^{-1}\mu^J_{\kf\cc}
\\
\hfill \tens\sS{_{((X_i)_{i\in\phi^{-1}j}g_j)_{j\in J}}}p_0s^{-1})
\mu_{\kf\cc} \quad
\\
=(\tens^{j\in J}\sS{_{(X_i)_{i\in\phi^{-1}j}}}\kf r^j\cdot\kf h\tens
\sS{_{((X_i)_{i\in\phi^{-1}j}g_j)_{j\in J}}}\kf p)\mu_{\kf\cc}
=\sS{_{(X_i)_{i\in I}}}[(\kf r^j)_{j\in J}\cdot\kf p].
\end{multline*}
Therefore, \(\kf[((r^js^{-1})_{j\in J},ps^{-1})\mu^{A^u_\infty}_\phi]
=((\kf r^j)_{j\in J},\kf p)\mu^{\wh\KCat}_\phi\), hence \(\kf\) is a
multifunctor.
\end{proof}

The quotient functor \(Q:\dg=\Com\to\ck\) equipped with the identity
transformation \(\tens^{i\in I}QX_i\to Q\tens^{i\in I}X_i\) is a
symmetric Monoidal functor. It gives rise to a symmetric Monoidal
\(\Cat\)\n-functor \(Q_*:\dgCat\to\KCat\). Let
\(\wh{Q_*}:\wh{\dgCat}\to\wh{\KCat}\) denote the corresponding
symmetric \(\Cat\)\n-multifunctor.

\begin{proposition}\label{pro-multinatural-isomorphism-xi}
There is a multinatural isomorphism
\begin{diagram}[height=2.3em]
\Ainftyu\times\wh\dgCat & \rTTo^{\boxdot} & \Ainftyu
\\
\dTTo<{\kf\times\wh{Q_*}} &\ldTwoar^{\xi} &\dTTo>{\kf}
\\
\wh\KCat\times\wh\KCat & \rTTo^\boxt & \wh\KCat
\end{diagram}
where \(\boxdot:\Ainftyu\boxt\wh\dgCat\to\Ainftyu\) is the action of
differential graded categories on unital \ainf-categories constructed
in \cite[Appendices~C.10--C.13]{BesLyuMan-book}.
\end{proposition}

\begin{proof}
Given a unital \ainf-category \(\ca\) and a differential graded
category \(\cc\), define an isomorphism of \(\ck\)\n-quivers
\(\xi:\kf(\ca\boxdot\cc)\to\kf\ca\boxt Q_*\cc\), identity on objects,
as follows. For \(X,Y\in\Ob\ca\), \(U,V\in\Ob\cc\), we have
\begin{align*}
\kf(\ca\boxdot\cc)((X,U),(Y,V))&=((\ca\boxdot\cc)((X,U),(Y,V)),m^{\ca\boxdot\cc}_1)
\\
&=((s\ca(X,Y)\tens\cc(U,V))[-1],sb^{\ca\boxdot\cc}s^{-1}),
\\
(\kf\ca\boxt Q_*\cc)((X,U),(Y,V))&=\kf\ca(X,Y)\tens Q_*\cc(U,V)\\
&=(\ca(X,Y)\tens\cc(U,V),m_1\tens1+1\tens
d).
\end{align*}
Define \(\xi\) by
\[
\xi=s(s^{-1}\tens1):(s\ca(X,Y)\tens\cc(U,V))[-1]\to\ca(X,Y)\tens\cc(U,V).
\]
The morphism \(\xi\) commutes with the differential since
\begin{multline*}
m^{\ca\boxdot\cc}_1\cdot\xi=sb^{\ca\boxdot\cc}s^{-1}\cdot s(s^{-1}\tens1)
=s(b_1\tens1-1\tens d)(s^{-1}\tens1)
\\
=s(s^{-1}\tens1)(sb_1s^{-1}\tens1+1\tens d)
=\xi\cdot(m_1\tens1+1\tens d),
\end{multline*}
therefore it is an isomorphism of \(\ck\)\n-quivers. We claim that
it also respects the composition. Indeed, suppose \(X,Y,Z\in\Ob\ca\),
\(U,V,W\in\Ob\cc\). From \cite[(C.10.1)]{BesLyuMan-book} we find that
\begin{multline*}
\mu_{\kf(\ca\boxdot\cc)}=m^{\ca\boxdot\cc}_2=
\bigl[
(\ca\boxdot\cc)((X,U),(Y,V))\tens(\ca\boxdot\cc)((Y,V),(Z,W))\rTTo^{s\tens s}
\\
(s\ca(X,Y)\tens\cc(U,V))\tens(s\ca(Y,Z)\tens\cc(V,W))\rTTo^{\sigma_{(12)}}
\\
(s\ca(X,Y)\tens s\ca(Y,Z))\tens(\cc(U,V)\tens\cc(V,W))\rTTo^{b_2\tens\mu_\cc}
\\
\hfill s\ca(X,Z)\tens\cc(U,W)\rTTo^{s^{-1}}(\ca\boxdot\cc)((X,U),(Z,W))
\bigr],\quad
\\
\quad\mu_{\kf\ca\boxt Q_*\cc}=\bigl[
(\ca(X,Y)\tens\cc(U,V))\tens(\ca(Y,Z)\tens\cc(V,W))\rTTo^{\sigma_{(12)}}\hfill
\\
(\ca(X,Y)\tens\ca(Y,Z))\tens(\cc(U,V)\tens\cc(V,W))\rTTo^{m_2\tens\mu_\cc}
\ca(X,Z)\tens\cc(U,W)
\bigr].
\end{multline*}
It follows that
\begin{multline*}
\mu_{\kf(\ca\boxdot\cc)}\cdot\xi=\bigl[
(\ca\boxdot\cc)((X,U),(Y,V))\tens(\ca\boxdot\cc)((Y,V),(Z,W))\rTTo^{s\tens s}
\\
(s\ca(X,Y)\tens\cc(U,V))\tens(s\ca(Y,Z)\tens\cc(V,W))\rTTo^{\sigma_{(12)}}
\\
\hfill(s\ca(X,Y)\tens s\ca(Y,Z))\tens(\cc(U,V)\tens\cc(V,W))
\rTTo^{b_2s^{-1}\tens\mu_\cc} \ca(X,Z)\tens\cc(U,W) \bigr],\quad
\\
\quad(\xi\tens\xi)\cdot\mu_{\kf\ca\boxt Q_*\cc}=\bigl[
(\ca\boxdot\cc)((X,U),(Y,V))\tens(\ca\boxdot\cc)((Y,V),(Z,W)) \hfill
\\
\rTTo^{s(s^{-1}\tens1)\tens s(s^{-1}\tens1)}
(\ca(X,Y)\tens\cc(U,V))\tens(\ca(Y,Z)\tens\cc(V,W))
\rTTo^{\sigma_{(12)}}
\\
\hfill(\ca(X,Y)\tens\ca(Y,Z))\tens(\cc(U,V)\tens\cc(V,W))
\rTTo^{m_2\tens\mu_\cc} \ca(X,Z)\tens\cc(U,W) \bigr]\quad
\\
\quad=\bigl[
(\ca\boxdot\cc)((X,U),(Y,V))\tens(\ca\boxdot\cc)((Y,V),(Z,W))
\rTTo^{s\tens s}\hfill
\\
(s\ca(X,Y)\tens\cc(U,V))\tens(s\ca(Y,Z)\tens\cc(V,W))
\rTTo^{\sigma_{(12)}}
\\
(s\ca(X,Y)\tens s\ca(Y,Z))\tens(\cc(U,V)\tens\cc(V,W))
\rTTo^{-(s^{-1}\tens s^{-1})m_2\tens\mu_\cc} \ca(X,Z)\tens\cc(U,W)
\bigr],
\end{multline*}
therefore
\(\mu_{\kf(\ca\boxdot\cc)}\cdot\xi=(\xi\tens\xi)\cdot\mu_{\kf\ca\boxt
Q_*\cc}\), as \(b_2s^{-1}=-(s^{-1}\tens s^{-1})m_2\). The morphism
\(\xi\) also respects the identity morphisms since
\[
1^{\kf(\ca\boxdot\cc)}_{(X,U)}\xi=(\sS{_X}\uni^\ca_0\tens
1^\cc_U)s^{-1}\cdot s(s^{-1}\tens1)=(\sS{_X}\uni^\ca_0s^{-1}\tens
1^\cc_U)=(1^{\kf\ca}_X\tens 1^\cc_U)=1^{\kf\ca\boxt Q_*\cc}_{(X,U)}.
\]
Thus, \(\xi\) is an isomorphism of \(\ck\)\n-categories.

Multinaturality of \(\xi\) reduces to the following problem. Let
\(f:\boxt^{i\in I}Ts\ca_i\to Ts\cb\) be an \ainf-functor,
\(g:\boxt^{i\in I}\cc_i\to\cd\) a differential graded functor. Then the
diagram
\begin{diagram}
\boxt^{i\in I}\kf(\ca_i\boxdot\cc_i) & \rTTo^{\boxt^I\xi\;}
&\boxt^{i\in I}(\kf\ca_i\boxt Q_*\cc_i)
\\
\dTTo<{\kf(f\boxdot g)} && \dTTo>{\sigma_{(12)}\cdot (\kf f\boxt\wh{Q_*}g)}
\\
\kf(\cb\boxdot\cd) & \rTTo^\xi & \kf\cb\boxt Q_*\cd
\end{diagram}
must commute; let us prove this. Let \(X_i,Y_i\in\Ob\ca_i\),
\(U_i,V_i\in\Ob\cc_i\), \(i\in I\), be families of objects. Then
\begin{multline*}
\kf(f\boxdot g)=\bigl[
\tens^{i\in I}(\ca_i\boxdot\cc_i)((X_i,U_i),(Y_i,V_i))
\rTTo^{\tens^{i\in I}s(f\boxdot g)_{e_i}s^{-1}}
\\
\quad\tens^{i\in I}(\cb\boxdot\cd)\bigl((((Y_j)_{j<i},(X_j)_{j\ge i})f,((V_j)_{j<i},(U_j)_{j\ge i})g),\hfill
\\
\hfill(((Y_j)_{j\le
i},(X_j)_{j>i})f,((V_j)_{j\le
i},(U_j)_{j>i})g)\bigr)\quad
\\
\rTTo^{\mu^I_{\kf(\cb\boxdot\cd)}}(\cb\boxdot\cd)\bigl(((X_i)_{i\in I}f,(U_i)_{i\in I}g),
((Y_i)_{i\in I}f,(V_i)_{i\in I}g)\bigr)
\bigr].
\end{multline*}
    From \cite[(C.5.1)]{BesLyuMan-book} we infer that
\begin{multline*}
(f\boxdot g)_{e_i}=\bigl[
s\ca_i(X_i,Y_i)\tens\cc(U_i,V_i)\rTTo^{f_{e_i}\tens g_{e_i}}
\\
\quad s\cb(((Y_j)_{j<i},(X_j)_{j\ge i})f,((Y_j)_{j\le i},(X_j)_{j>i})f)
\tens \hfill
\\
\quad \tens
\cd(((V_j)_{j<i},(U_j)_{j\ge i})g,((V_j)_{j\le i},(U_j)_{j>i})g)
\bigr],
\end{multline*}
where
\begin{multline*}
g_{e_i}=\bigl[ \cc_i(U_i,V_i)
 \rTTo^{\lambda^{i:\mb1\hookrightarrow I}\cdot
 \tens^{j\in I}[(1_{V_j})_{j<i},\id,(1_{U_j})_{j>i}]}
\\
\tens^{j\in I}
[(\cc_j(V_j,V_j))_{j<i},\cc_i(U_i,V_i),(\cc_j(U_j,U_j))_{j>i}]
\\
\rTTo^g
\cd(((V_j)_{j<i},(U_j)_{j\ge i})g,((V_j)_{j\le i},(U_j)_{j>i})g)
\bigr].
\end{multline*}
Therefore
\begin{multline*}
\kf(f\boxdot g)\cdot\xi=\bigl[
\tens^{i\in I}(\ca_i\boxdot\cc_i)((X_i,U_i),(Y_i,V_i))
\rTTo^{\tens^{i\in I}s(f_{e_i}\tens g_{e_i})s^{-1}}
\\
\quad\tens^{i\in I}(\cb\boxdot\cd)\bigl((((Y_i)_{j<i},(X_j)_{j\ge
i})f,((V_i)_{j<i},(U_j)_{j\ge i})g),\hfill
\\
\hfill(((Y_j)_{j\le
i},(X_j)_{j>i})f,((V_j)_{j\le
i},(U_j)_{j>i})g)\bigr)\quad
\\
\quad\rTTo^{\tens^{i\in I}\xi\;}
\tens^{i\in I}[\cb(((Y_i)_{j<i},(X_j)_{j\ge i})f,((Y_j)_{j\le
i},(X_j)_{j>i})f)\hfill
\\
\hfill\tens\cd(((V_i)_{j<i},(U_j)_{j\ge i})g,((V_j)_{j\le
i},(U_j)_{j>i})g)]\quad
\\
\rTTo^{\mu^I_{\kf\cb\boxt Q_*\cd}}
\cb((X_i)_{i\in I}f,(Y_i)_{i\in I}f)\tens\cd((U_i)_{i\in I}g,(V_i)_{i\in I}g)
\bigr]
\end{multline*}
since \(\xi\) respects the composition. The above expression can be
transformed as follows:
\begin{multline*}
\kf(f\boxdot g)\cdot\xi=\bigl[
\tens^{i\in I}(\ca_i\boxdot\cc_i)((X_i,U_i),(Y_i,V_i))
 \rTTo^{\tens^{i\in I}s(s^{-1}\tens1)}
\\
\tens^{i\in I}(\ca_i(X_i,Y_i)\tens\cc_i(U_i,V_i))
\rTTo^{\tens^{i\in I}(sf_{e_i}s^{-1}\tens g_{e_i})}
\\
\quad \tens^{i\in I}[\cb(((Y_i)_{j<i},(X_j)_{j\ge i})f,((Y_j)_{j\le
i},(X_j)_{j>i})f) \hfill
\\
\hfill\tens\cd(((V_i)_{j<i},(U_j)_{j\ge i})g,((V_j)_{j\le
i},(U_j)_{j>i})g)]\quad
\\
\quad\rTTo^{\sigma_{(12)}}
\bigl(\tens^{i\in I}\cb(((Y_j)_{j<i},(X_j)_{j\ge i})f,((Y_j)_{j\le
i},(X_j)_{j>i})f)\bigr)\hfill
\\
\hfill\tens\bigl(\tens^{i\in I}\cd(((V_j)_{j<i},(U_j)_{j\ge i})g,((V_j)_{j\le
i},(U_j)_{j>i})g)\bigr)\quad
\\
\hfill\rTTo^{\mu^I_{\kf\cb}\tens\mu^I_\cd}
\cb((X_i)_{i\in I}f,(Y_i)_{i\in I}f)\tens\cd((U_i)_{i\in I}g,(V_i)_{i\in I}g)
\bigr]\quad
\\
\quad=\bigl[
\tens^{i\in I}(\ca_i\boxdot\cc_i)((X_i,U_i),(Y_i,V_i))\rTTo^{\tens^{i\in
I}\xi}\hfill
\\
\tens^{i\in I}(\ca_i(X_i,Y_i)\tens\cc_i(U_i,V_i))
\rTTo^{\sigma_{(12)}}
(\tens^{i\in I}\ca_i(X_i,Y_i))\tens(\tens^{i\in I}\cc_i(U_i,V_i))
\\
\hfill\rTTo^{(\tens^{i\in I}sf_{e_i}s^{-1})\mu^I_{\kf\cb}\tens(\tens^{i\in I}g_{e_i})\mu_\cd}
\cb((X_i)_{i\in I}f,(Y_i)_{i\in I}f)\tens\cd((U_i)_{i\in I}g,(V_i)_{i\in I}g)
\bigr]\quad
\\
\quad=\tens^{i\in I}\xi\cdot\sigma_{(12)}\cdot(\kf f\tens g),\hfill
\end{multline*}
due to the definition of \(\kf f\) and the identity
\begin{multline*}
\bigl[
\tens^{i\in I}\cc_i(U_i,V_i)\rTTo^{\tens^{i\in I}g_{e_i}\;}
\tens^{i\in I}\cd(((V_j)_{j<i},(U_j)_{j\ge i})g,((V_j)_{j\le i},(U_j)_{j>i})g)
\\
\rTTo^{\mu^I_\cd}\cd((U_i)_{i\in I}g,(V_i)_{i\in I}g)
\bigr]=g,
\end{multline*}
which is a consequence of \(g\) being a functor and associativity of
\(\mu_\cd\). The proposition is proven.
\end{proof}

\subsection{$A_\infty$-categories closed under shifts.}
 \label{sec-ainf-categories-closed-under-shifts}
Unital \ainf-categories closed under shifts are defined in
\cite[Chapter~15]{BesLyuMan-book} similarly to
\defref{def-closed-under-shifts}. A unital \ainf-category $\cc$ is
closed under shifts if and only if the \ainf-functor
\(u_\sh:\cc\to\cc^\sh\simeq\cc\boxdot\cz=\cc\boxt\cz\) is an
equivalence.

For an arbitrary \ainf-category $\cc$ the operations in $\cc^\sh$ are
described explicitly in \cite[(15.2.2)]{BesLyuMan-book}. In particular,
\begin{multline*}
b^{\cc^\sh}_2 =(-)^{p-n}\bigl[
 s\cc^\sh((X,n),(Y,m))\tens s\cc^\sh((Y,m),(Z,p))
\\
=s\cc(X,Y)[m-n]\tens s\cc(Y,Z)[p-m]
 \rTTo^{(s^{m-n}\tens s^{p-m})^{-1}} s\cc(X,Y)\tens s\cc(Y,Z)
\\
\rTTo^{b^\cc_2} s\cc(X,Z) \rTTo^{s^{p-n}} s\cc(X,Z)[p-n]
=s\cc^\sh((X,n),(Z,p)) \bigr].
\end{multline*}
The above proposition implies that the binary operation
\(m^{\cc^\sh}_2=(s\tens s)b^{\cc^\sh}_2s^{-1}\) in $\cc^\sh$ is
homotopic to multiplication in \(\kf\cc^\sh\) given by formula
\eqref{eq-muC[]-ss1-muC-s}. Actually, \(m^{\cc^\sh}_2\) is given
precisely by chain map~\eqref{eq-muC[]-ss1-muC-s}, as one easily
deduces from the above expression for \(b^{\cc^\sh}_2\).

We have denoted the algebra $\cz$ in \(\dgCat\), ``the same'' algebra
\(Q_*\cz\) in \(\KCat\) and ``the same'' algebra \(H^\bullet(Q_*\cz)\)
in \(\grCat\) all by the same letter $\cz$ by abuse of notation. Since
units of the monads \(-\boxdot\cz\) and \(-\boxt\cz\) reduce
essentially to the unit of the algebra $\cz$,
\propref{pro-multinatural-isomorphism-xi} implies the following
relation between them:
\[ \bigl[ \kf\ca \rTTo^{\kf u_\sh} \kf(\ca\boxdot\cz) \rTTo^{\xi}_\sim
\kf\ca\boxt Q_*\cz \bigr] =\bigl[ \kf\ca \rTTo^{u_\sh} \kf\ca\boxt\cz
=\kf\ca\boxt Q_*\cz \bigr].
\]
Thus, if one of the \(\ck\)\n-functors
 \(\kf u_\sh:\kf\ca\to\kf(\ca\boxdot\cz)\) and
\(u_\sh:\kf\ca\to\kf\ca\boxt\cz\) is an equivalence, then so is the
other. The former is a \(\ck\)\n-equivalence if and only if the
\ainf-functor \(u_\sh:\ca\to\ca\boxdot\cz\) is an equivalence.
Therefore, the \ainf-category $\ca$ is closed under shifts if and only
if the \(\ck\)\n-category \(\kf\ca\) is closed under shifts.

\subsection{Shifts as differential graded functors.}
 \label{sec-Shifts-dg-functors}
Let \(f=(f^i:C^i\to D^{i+\deg f})_{i\in\ZZ}\in\uCom(C,D)\) be a
homogeneous element (a $\kk$\n-linear map $f:C\to D$ of certain degree
$d=\deg f$). Define
 $f^{[n]}=(-)^{fn}s^{-n}fs^n=(-)^{dn}\bigl(C[n]^i=C^{i+n}
 \rTTo^{f^{i+n}} D^{i+n+d}=D[n]^{i+d}\bigr)$,
which is an element of \(\uCom(C[n],D[n])\) of the same degree $\deg f$.

Define the shift differential graded functor \([n]:\uCom\to\uCom\) as
follows. It takes a complex \((C,d)\) to the complex
\((C[n],d^{[n]})\), \(d^{[n]}=(-)^ns^{-n}ds^n\). On morphisms it acts
via \(\uCom(s^{-n},1)\cdot\uCom(1,s^n):\uCom(C,D)\to\uCom(C[n],D[n])\),
\(f\mapsto f^{[n]}\). Clearly, \([n]\cdot[m]=[n+m]\).

\section{$A_\infty$-modules.}\label{sec-A8-modules}
Consider the monoidal category \((\cQuiver/S,\tens)\) of graded
$\kk$\n-quivers. When $S=\mb1$ it reduces to the category of graded
$\kk$\n-modules used by Keller~\cite{math.RA/9910179} in his definition
of $A_\infty$\n-modules over $A_\infty$\n-algebras. Let $C$, $D$ be
coassociative counital coalgebras; let \(\psi:C\to D\) be a
homomorphism; let \(\delta:M\to M\tens C\) and \(\delta:N\to N\tens D\)
be counital comodules; let \(f:M\to N\) be a $\psi$\n-comodule
homomorphism, \(f\delta=\delta(f\tens\psi)\); let $\xi:C\to D$ be a
\((\psi,\psi)\)-coderivation,
\(\xi\Delta_0=\Delta_0(\psi\tens\xi+\xi\tens\psi)\). Define a
\emph{\((\psi,f,\xi)\)-connection} as a morphism \(r:M\to N\) of
certain degree such that
\begin{diagram}
M &\rTTo^\delta &M\tens C
\\
\dTTo<r &= &\dTTo>{f\tens\xi+r\tens\psi}
\\
N &\rTTo^\delta &N\tens D
\end{diagram}
compare with Tradler~\cite{xxx0108027}. Let $(C,b^C)$ be a differential
graded coalgebra. Let a counital comodule $M$ have a
\((1,1,b^C)\)-connection \(b^M:M\to M\) of degree 1, that is,
\(b^M\delta=\delta(1\tens b^C+b^M\tens1)\). Its \emph{curvature}
\((b^M)^2:M\to M\) is always a $C$\n-comodule homomorphism of degree 2.
If it vanishes, $b^M$ is called a flat connection (a differential) on
$M$.

Equivalently, we consider the category \((\dQ/S,\tens)\) of
differential graded quivers, and coalgebras and comodules therein. For
\ainf-applications it suffices to consider coalgebras (resp. comodules)
whose underlying graded coalgebra (resp. comodule) has the form
\(Ts\ca\) (resp. \(s\cm\tens Ts\cc\)).

Let \(\cm\in\Ob\cQuiver/S\) be graded quiver such that
\(\cm(X,Y)=\cm(Y)\) depends only on $Y\in S$. For any quiver
\(\cc\in\Ob\cQuiver/S\) the tensor quiver \(C=(Ts\cc,\Delta_0)\) is a
coalgebra. The comodule
 \(\delta=1\tens\Delta_0:M=s\cm\tens Ts\cc\to s\cm\tens Ts\cc\tens
Ts\cc\) is counital. Let $(\cc,b^\cc)$ be an \ainf-category.
Equivalently, we consider augmented coalgebras in \((\dQ/S,\tens)\) of
the form \((Ts\cc,\Delta_0,b^\cc)\). Let \(b^\cm:s\cm\tens Ts\cc\to
s\cm\tens Ts\cc\) be a \((1,1,b^\cc)\)\n-connection. Define the matrix
coefficients of \(b^\cm\) to be
\[
b^\cm_{mn}=(1\tens\inj_m)\cdot
b^\cm\cdot(1\tens\pr_n):s\cm\tens T^ms\cc\to s\cm\tens T^ns\cc,\quad m,n\ge0.
\]
The coefficients \(b^\cm_{m0}:s\cm\tens T^ms\cc\to s\cm\) are
abbreviated to \(b^\cm_m\) and called components of \(b^\cm\).

A version of the following statement occurs in
\cite[Lemme~2.1.2.1]{Lefevre-Ainfty-these}.

\begin{lemma}
Any \((1,1,b^\cc)\)-connection
 \(b^\cm:s\cm\tens Ts\cc\to s\cm\tens Ts\cc\) is determined in a
unique way by its components \(b^\cm_n:s\cm\tens T^ns\cc\to s\cm\),
$n\ge0$. The matrix coefficients of \(b^\cm\) are expressed via
components of \(b^\cm\) and components of the codifferential \(b^\cc\)
as follows:
\begin{equation*}
b^\cm_{mn}=b^\cm_{m-n}\tens1^{\tens n}
 +\sum_{\substack{p+k+q=m\\ p+1+q=n}}
1^{\tens1+p}\tens b^\cc_k\tens 1^{\tens q}: s\cm\tens T^ms\cc\to
s\cm\tens T^ns\cc
\end{equation*}
for \(m\ge n\). If \(m<n\), the matrix coefficient \(b^\cm_{mn}\)
vanishes.
\end{lemma}

Such comodules are particular cases of bimodules discussed below. That
is why statements about comodules are only formulated. We prove more
general results in the next section.

The morphism \((b^\cm)^2:s\cm\tens Ts\cc\to s\cm\tens Ts\cc\) is a
\((1,1,0)\)\n-connection of degree~\(2\), therefore equation
\((b^\cm)^2=0\) is equivalent to its particular case
\((b^\cm)^2(1\tens\pr_0)=0:s\cm\tens Ts\cc\to s\cm\). Thus \(b^\cm\) is
a flat connection if for each \(m\ge0\) the following equation holds:
\begin{multline}
\sum_{n=0}^m(b^\cm_{m-n}\tens1^{\tens n})b^\cm_{n}+\sum_{p+k+q=m}
(1^{\tens1+p}\tens b^\cc_k\tens 1^{\tens q})b^\cm_{p+1+q}=0:
\\
s\cm\tens T^ms\cc\to s\cm.
 \label{equ-bM2-0}
\end{multline}
Equivalently, such a \(Ts\cc\)\n-comodule with a flat connection is the
\(Ts\cc\)\n-comodule \((s\cm\tens Ts\cc,b^\cm)\) in the category
\((\dQ/S,\tens)\). It consists of the following data: a graded
\(\kk\)\n-module \(\cm(X)\) for each object \(X\) of \(\cc\); a family
of \(\kk\)\n-linear maps of degree~\(1\)
\[
b^\cm_n:s\cm(X_0)\tens s\cc(X_0,X_1)\tdt
s\cc(X_{n-1},X_n)\to s\cm(X_n),\quad n\ge0,
\]
subject to equations~\eqref{equ-bM2-0}. Equation~\eqref{equ-bM2-0} for
\(m=0\) implies \((b^\cm_0)^2=0\), that is, \((s\cm(X),b^\cm_0)\) is a
chain complex, for each object \(X\in\Ob\cc\). We call a
\(Ts\cc\)\n-comodule with a flat connection
 \((s\cm\tens Ts\cc,b^\cm)\), \(\cm(*,Y)=\cm(Y)\), a
\emph{$\cc$\n-module} (an \emph{\ainf-module over} $\cc$).
$\cc$\n-modules form a differential graded category \(\cc\modul\). The
notion of a module over some kind of \ainf-category was
introduced by Lef\`evre-Hasegawa under the name of polydule
\cite{Lefevre-Ainfty-these}.

\begin{proposition}\label{prop-ainf-modules-ainf-functors-A-Com}
An arbitrary \ainf-functor \(\phi:\cc\to\uCom\) determines a
\(Ts\cc\)-comodule \(s\cm\tens Ts\cc\) with a flat connection $b^\cm$
by the formulae: \(\cm(X)=X\phi\), for each object \(X\) of \(\cc\),
\(b^\cm_0=s^{-1}ds:s\cm(X)\to s\cm(X)\), where \(d\) is the
differential in the complex \(X\phi\), and for $n>0$
\begin{multline}
b^\cm_n = \bigl[ s\cm(X_0)\tens s\cc(X_0,X_1)\tdt s\cc(X_{n-1},X_n)
\\
\rTTo^{1\tens\phi_n} s\cm(X_0)\tens s\uCom(\cm(X_0),\cm(X_n))
\rTTo^{(s\tens s)^{-1}} \cm(X_0)\tens\uCom(\cm(X_0),\cm(X_n))
\\
\rTTo^{\ev^\Com} \cm(X_n) \rTTo^s s\cm(X_n) \bigr].
\label{equ-ainf-functor-to-ainf-module}
\end{multline}
This mapping from \ainf-functors to $\cc$\n-modules is bijective.
Moreover, the differential graded categories \(\und\Ainfty(\cc;\uCom)\)
and \(\cc\modul\) are isomorphic.
\end{proposition}

A \emph{$\cc$\n-module} (an \emph{\ainf-module over} $\cc$) is defined
as an \ainf-functor \(\phi:\cc\to\uCom\) by Seidel
\cite[Section~1j]{SeidelP-book-Fukaya}. The above proposition shows
that the both definitions of $\cc$\n-modules are equivalent. In the
differential graded case $\cc$\n-modules are actively used by
Drinfeld~\cite{Drinf:DGquot}.

\begin{definition}
Let $\cc$ be a unital \ainf-category. A $\cc$\n-module $\cm$ determined
by an \ainf-functor \(\phi:\cc\to\uCom\) is called \emph{unital} if
\(\phi\) is unital.
\end{definition}

\begin{proposition}\label{prop-unital-ainf-module}
A $\cc$\n-module $\cm$ is unital if and only if for each
\(X\in\Ob\cc\) the composition
\[ \bigl[s\cm(X) \simeq s\cm(X)\tens\kk \rTTo^{1\tens\sS{_X}\uni^\cc_0}
s\cm(X)\tens s\cc(X,X) \rTTo^{b^\cm_1} s\cm(X) \bigr]
\]
is homotopic to identity map.
\end{proposition}

\begin{proof}
The second statement expands to the property that
\begin{multline*}
\bigl[ s\cm(X) \simeq s\cm(X)\tens\kk
\rTTo^{s^{-1}\tens\sS{_X}\uni^\cc_0} \cm(X)\tens s\cc(X,X)
\\
\rTTo^{1\tens\phi_1s^{-1}} \cm(X)\tens\uCom(\cm(X),\cm(X))
\rTTo^{\ev^\Com} \cm(X) \rTTo^s s\cm(X) \bigr]
\end{multline*}
is homotopic to identity. That is,
\[ \sS{_X}\uni^\cc_0\phi_1s^{-1} =1_{s\cm(X)} +vm^\uCom_1,
\quad\text{or,}\quad \sS{_X}\uni^\cc_0\phi_1 =1_{\cm(X)}s +vsb^\uCom_1.
\]
In other words, \ainf-functor $\phi$ is unital.
\end{proof}

\section{$A_\infty$-bimodules.}\label{sec-A8-bimodules}
Consider monoidal category \((\cQuiver/S,\tens)\) of graded
$\kk$\n-quivers. When $S=\mb1$ it reduces to the category of graded
$\kk$\n-modules used by
Tradler~\cite{xxx0108027,math.QA/0210150} in his definition
of $A_\infty$\n-bimodules over $A_\infty$\n-algebras. We extend his definitions of
$A_\infty$-bimodules improved in \cite{math.QA/0305052} from graded
$\kk$\n-modules to graded $\kk$\n-quivers.
The notion of a bimodule over some kind of \ainf-categories was
introduced by Lef\`evre-Hasegawa under the name of bipolydule
\cite{Lefevre-Ainfty-these}.

\begin{definition}
Let $A$, $C$ be coassociative counital coalgebras in
\((\cQuiver/R,\tens)\) resp. \((\cQuiver/S,\tens)\). A \emph{counital
$(A,C)$\n-bicomodule} $(P,\delta^{P})$ consists of a graded
$\kk$\n-span ($\gr$\n-span) $P$ with \(\Ob_sP=R\), \(\Ob_tP=S\),
 \(\Par P=\Ob_sP\times\Ob_tP\), \(\src=\pr_1\), \(\tgt=\pr_2\) and a
coaction
$\delta^{P}=(\delta',\delta''):P\to(A\tens_RP)\oplus(P\tens_SC)$ of
degree $0$ such that the following diagram commutes
\begin{diagram}[height=2.5em]
P &\rTTo^\delta &(A\tens_RP)\oplus(P\tens_SC)
\\
\dTTo<\delta &&\dTTo>{(\Delta\tens1)\oplus(\delta\tens1)}
\\
(A\tens_RP)\oplus(P\tens_SC)
&\rTTo^{(1\tens\delta)\oplus(1\tens\Delta)}
&
\begin{array}{r}
(A\tens_RA\tens_RP)\oplus(A\tens_RP\tens_SC) \\
\oplus(P\tens_SC\tens_SC)
\end{array}
\end{diagram}
and \(\delta'\cdot(\eps\tens1)=1=\delta''\cdot(1\tens\eps):P\to P\).
\end{definition}

The equation presented on the diagram consists in fact of three
equations claiming that $P$ is a left $A$\n-comodule, a right
$C$\n-comodule and the coactions commute.

Let $A$, $B$, $C$, $D$ be coassociative counital coalgebras; let
\(\phi:A\to B\), \(\psi:C\to D\) be homomorphisms; let $\chi:A\to B$ be
a \((\phi,\phi)\)-coderivation and let $\xi:C\to D$ be a
\((\psi,\psi)\)-coderivation of certain degree, that is,
\(\chi\Delta=\Delta(\phi\tens\chi+\chi\tens\phi)\),
\(\xi\Delta=\Delta(\psi\tens\xi+\xi\tens\psi)\). Let
\(\delta:P\to(A\tens P)\oplus(P\tens C)\) be a counital
$(A,C)$\n-bicomodule and let
 \(\delta:Q\to(B\tens Q)\oplus(Q\tens D)\) be a counital
$(B,D)$\n-bicomodule. A $\kk$\n-span morphism \(f:P\to Q\) of degree 0
with \(\Ob_sf=\Ob\phi\), \(\Ob_tf=\Ob\psi\) is a
\emph{$(\phi,\psi)$\n-bicomodule homomorphism} if
\(f\delta'=\delta'(\phi\tens f):P\to B\tens Q\) and
\(f\delta''=\delta''(f\tens\psi):P\to Q\tens D\). Define a
\emph{\((\phi,\psi,f,\chi,\xi)\)-connection} as a $\kk$\n-span morphism
\(r:P\to Q\) of certain degree with \(\Ob_sr=\Ob\phi\),
\(\Ob_tr=\Ob\psi\) such that
\begin{diagram}
P &\rTTo^\delta &(A\tens P)\oplus(P\tens C)
\\
\dTTo<r
 &&\dTTo>{(\phi\tens r+\chi\tens f)\oplus(f\tens\xi+r\tens\psi)}
\\
Q &\rTTo^\delta &(B\tens Q)\oplus(Q\tens D)
\end{diagram}

Let $(A,b^A)$, $(C,b^C)$ be differential graded coalgebras and let $P$
be an $(A,C)$\n-bicomodule with an
\((\id_A,\id_C,\id_P,b^A,b^C)\)-connection \(b^P:P\to P\) of degree 1,
that is, \(b^P\delta'=\delta'(1\tens b^P+b^A\tens1)\) and
\(b^P\delta''=\delta''(1\tens b^C+b^P\tens1)\). Its \emph{curvature}
\((b^P)^2:P\to P\) is always an $(A,C)$\n-bicomodule homomorphism of
degree 2. If it vanishes, $b^P$ is called a \emph{flat connection} (a
differential) on $P$.

Taking for $(A,b^A)$ the trivial differential graded coalgebra $\kk$
with the trivial coactions we recover the notions introduced in
\secref{sec-A8-modules}. Namely, an $(A,C)$\n-bicomodule $P$ with an
\((\id_A,\id_C,\id_P,b^A,b^C)\)-connection \(b^P:P\to P\) of degree 1
is the same as a $C$\n-comodule with a \((1,1,b^C)\)-connection, both
flatness conditions coincide, etc.

Equivalently, bicomodules with flat connections are bicomodules which
live in the category of differential graded spans. The set of
\(A\)\n-\(C\)-bicomodules becomes the set of objects of a differential
graded category \(A\text-C\bicomod\). For differential graded
bicomodules \(P\), \(Q\), the \(k\)\n-th component of the graded
\(\kk\)\n-module \(A\text-C\bicomod(P,Q)\) consists of
\((\id_A,\id_C)\)\n-bicomodule homomorphisms  \(t:P\to Q\) of degree
\(k\). The differential of \(t\) is the commutator
\(tm_1=tb^Q-(-)^tb^Pt:P\to Q\), which is again a homomorphism of
bicomodules, naturally of degree \(k+1\). Composition of homomorphisms
of bicomodules is the ordinary composition of \(\kk\)\n-span morphisms.

The main example of a bicomodule is the following. Let $\ca$, $\cb$,
$\cc$, $\cd$ be graded $\kk$\n-quivers. Let $\cp$, $\cq$ be
$\gr$\n-spans with \(\Ob_s\cp=\Ob\ca\), \(\Ob_t\cp=\Ob\cc\),
\(\Ob_s\cq=\Ob\cb\), \(\Ob_t\cq=\Ob\cd\),
\(\Par\cp=\Ob_s\cp\times\Ob_t\cp\), \(\Par\cq=\Ob_s\cq\times\Ob_t\cq\),
\(\src=\pr_1\), \(\tgt=\pr_2\). Take coalgebras \(A=Ts\ca\),
\(B=Ts\cb\), \(C=Ts\cc\), \(D=Ts\cd\) and bicomodules
 \(P=Ts\ca\tens s\cp\tens Ts\cc\), \(Q=Ts\cb\tens s\cq\tens Ts\cd\)
equipped with the cut comultiplications (coactions)
\begin{align*}
\Delta_0(a_{1},\dots,a_{n})
&=\sum_{i=0}^{n}(a_{1},\dots,a_{i})\tens(a_{i+1},\dots,a_{n}),
\\
\delta(a_{1},\ldots,a_{k},p,c_{k+1},\ldots,c_{k+l}) &=
\sum_{i=0}^{k}(a_{1},\dots,a_{i})\tens(a_{i+1},\ldots,p,\ldots,c_{k+l})
\\
&\qquad +\sum_{i=k}^{k+l}
(a_{1},\ldots,p,\ldots,c_{i})\tens(c_{i+1},\ldots,c_{k+l}).
\end{align*}

Notice that a graded quiver \(\cm\in\Ob\cQuiver/S\) such that
\(\cm(X,Y)=\cm(Y)\) depends only on $Y\in S$ is nothing else but a
$\gr$\n-span $\cm$ with \(\Ob_s\cm=\{*\}\), \(\Ob_t\cm=S\). Thus,
$Ts\cc$-comodules of the form \(s\cm\tens Ts\cc\) from
\secref{sec-A8-modules} are nothing else but
$Ts\ca$-$Ts\cc$-bicomodules \(Ts\ca\tens s\cm\tens Ts\cc\) for the
graded quiver $\ca=\1_u$ with one object $*$ and with \(\1_u(*,*)=0\).
Furthermore, \ainf-modules $\cm$ over an \ainf-category $\cc$ are the
same as $\1_u$-$\cc$-bimodules, as defined before
\propref{pro-dg-cat-ACbimod-A(ACC)-isomorphic}.

Let \(\phi:Ts\ca\to Ts\cb\), \(\psi:Ts\cc\to Ts\cd\) be augmented
coalgebra morphisms. Let \(g:P\to Q\) be a \(\kk\)\n-span morphism of
certain degree with \(\Ob_s g=\Ob\phi\), \(\Ob_t g=\Ob\psi\). Define
the matrix coefficients of \(g\) to be
\begin{multline*}
g_{kl;mn}=(\inj_k\tens1\tens\inj_l)\cdot
g\cdot(\pr_m\tens1\tens\pr_n):\\
T^ks\ca\tens s\cp\tens T^ls\cc\to
T^ms\cb\tens s\cq\tens T^ns\cd, \quad k,l,m,n\ge0.
\end{multline*}
The coefficients \(g_{kl;00}:T^ks\ca\tens s\cp\tens T^ls\cc\to s\cq\)
are abbreviated to \(g_{kl}\) and called components of \(g\). Denote by
\(\check g\) the composite
 \(g\cdot(\pr_0\tens1\tens\pr_0):Ts\ca\tens s\cp\tens Ts\cc\to s\cq\).
The restriction of \(\check g\) to the summand \(T^ks\ca\tens s\cp\tens
T^ls\cc\) is precisely the component \(g_{kl}\).

Let \(f:P\to Q\) be a \((\phi,\psi)\)\n-bicomodule homomorphism. It is
uniquely recovered from its components similarly to
Tradler~\cite[Lemma~4.2]{xxx0108027}. Let us supply the details. The
coaction \(\delta^P\) has two components,
\begin{align*}
\delta'=\Delta_0\tens1\tens1&:Ts\ca\tens s\cp\tens Ts\cc\to Ts\ca\tens
Ts\ca\tens s\cp\tens Ts\cc,\\
\delta''=1\tens1\tens\Delta_0&:Ts\ca\tens s\cp\tens Ts\cc\to Ts\ca\tens
s\cp\tens Ts\cc\tens Ts\cc,
\end{align*}
and similarly for \(\delta^Q\). As \(f\) is a
\((\phi,\psi)\)\n-bicomodule homomorphism, it satisfies the equations
\begin{align*}
f(\Delta_0\tens1\tens1)=(\Delta_0\tens1\tens1)(\phi\tens f)&:Ts\ca\tens
s\cp\tens Ts\cc\to Ts\cb\tens Ts\cb\tens s\cq\tens Ts\cd,\\
f(1\tens1\tens\Delta_0)=(1\tens1\tens\Delta_0)(f\tens\psi)&:Ts\ca\tens
s\cp\tens Ts\cc\to Ts\cb\tens s\cq\tens Ts\cd\tens Ts\cd.
\end{align*}
It follows that
\begin{multline*}
f(\Delta_0\tens1\tens\Delta_0)=(\Delta_0\tens1\tens\Delta_0)(\phi\tens
f\tens\psi):\\
Ts\ca\tens
s\cp\tens Ts\cc\to Ts\cb\tens Ts\cb\tens s\cq\tens Ts\cd\tens Ts\cd.
\end{multline*}
Composing both sides with the morphism
\begin{equation}
1\tens\pr_0\tens1\tens\pr_0\tens1:Ts\cb\tens Ts\cb\tens s\cq\tens
Ts\cd\tens Ts\cd\to Ts\cb\tens s\cq\tens Ts\cd,
\label{equ-1-pr0-1-pr0-1}
\end{equation}
and taking into account the identities \(\Delta_0(1\tens\pr_0)=1\),
\(\Delta_0(\pr_0\tens1)=1\), we obtain
\begin{equation}
f=(\Delta_0\tens1\tens\Delta_0)(\phi\tens\check f\tens\psi).
\label{equ-f-Delta1Delta-phi-f-psi}
\end{equation}
This equation implies the following formulas for the matrix
coefficients of \(f\):
\begin{multline}
f_{kl;mn}=\sum_{\substack{i_1+\dots+i_m+p=k\\ j_1+\dots+j_n+q=l}}
(\phi_{i_1}\tdt\phi_{i_m}\tens f_{pq}\tens\psi_{j_1}\tdt\psi_{j_n}):
\\
T^ks\ca\tens s\cp\tens T^ls\cc\to T^ms\cb\tens s\cq\tens T^ns\cd,\quad
k,l,m,n\ge0.
\label{eq-fklmn-phiphifpsipsi}
\end{multline}
In particular, if \(k<m\) or \(l<n\), the matrix coefficient
\(f_{kl;mn}\) vanishes.

Let \(r:P\to Q\) be a \((\phi,\psi,f,\chi,\xi)\)\n-connection. It
satisfies the following equations:
\begin{multline*}
r(\Delta_0\tens1\tens1)=(\Delta_0\tens1\tens1)(\phi\tens r+\chi\tens
f):\\ \hfill Ts\ca\tens
s\cp\tens Ts\cc\to Ts\cb\tens Ts\cb\tens s\cq\tens Ts\cd,\quad\\
\quad r(1\tens1\tens\Delta_0)=(1\tens1\tens\Delta_0)(f\tens\xi+r\tens\psi):\hfill\\
Ts\ca\tens
s\cp\tens Ts\cc\to Ts\cb\tens s\cq\tens Ts\cd\tens Ts\cd.
\end{multline*}
They imply that
\begin{multline*}
r(\Delta_0\tens1\tens\Delta_0)=
(\Delta_0\tens1\tens\Delta_0)(\phi\tens f\tens\xi+\phi\tens r\tens\psi+\chi\tens
f\tens\psi):\\
Ts\ca\tens
s\cp\tens Ts\cc\to Ts\cb\tens Ts\cb\tens s\cq\tens Ts\cd\tens Ts\cd.
\end{multline*}
Composing both side with the morphism~\eqref{equ-1-pr0-1-pr0-1} we
obtain
\begin{equation}
r=(\Delta_0\tens1\tens\Delta_0)
(\phi\tens\check f\tens\xi+\phi\tens\check r\tens\psi+\chi\tens\check
f\tens\psi).
\label{equ-connection-via-check-r}
\end{equation}
    From this equation we find the following expression for the matrix
coefficient \(r_{kl;mn}\):
\begin{multline}
\sum^{p+1+q=n}_{\substack{i_1+\dots+i_m+i=k\\ j+j_1+\dots+j_p+t+j_{p+1}+\dots+j_{p+q}=l}}
 \hspace*{-2em} \phi_{i_1}\tdt\phi_{i_m}\tens f_{ij}\tens\psi_{j_1}
\tdt\psi_{j_p}\tens\xi_t\tens\psi_{j_{p+1}}\tdt\psi_{j_{p+q}}
\\
+\sum_{\substack{i_1+\dots+i_m+i=k\\ j+j_1+\dots+j_n=l}}
\phi_{i_1}\tdt\phi_{i_m}\tens r_{ij}\tens\psi_{j_1}\tdt\psi_{j_n}
\\
+\sum^{a+1+c=m}_{\substack{i_1+\dots+i_a+u+i_{a+1}+\dots+i_{a+c}+i=k\\ j+j_1+\dots+j_n=l}}
 \hspace*{-2em} \phi_{i_1}\tdt\phi_{i_a}\tens\chi_u\tens\phi_{i_{a+1}}
\tdt\phi_{i_{a+c}}\tens f_{ij}\tens\psi_{j_1}\tdt\psi_{j_n}:
\\
T^ks\ca\tens s\cp\tens T^ls\cc\to T^ms\cb\tens s\cq\tens T^ns\cd,
 \quad k,l,m,n\ge0.
\label{equ-phi-psi-f-chi-xi-connection-comp}
\end{multline}

Let $\ca$, $\cc$ be \ainf-categories and let $\cp$ be a $\gr$\n-span
with \(\Ob_s\cp=\Ob\ca\), \(\Ob_t\cp=\Ob\cc\),
\(\Par\cp=\Ob_s\cp\times\Ob_t\cp\), \(\src=\pr_1\), \(\tgt=\pr_2\). Let
\(A=Ts\ca\), \(C=Ts\cc\), and consider the bicomodule \(P=Ts\ca\tens
s\cp\tens Ts\cc\). The set of \((1,1,1,b^\cA,b^\cC)\)-connections
\(b^\cP:P\to P\) of degree 1 with \((b^\cp_{00})^2=0\) is in bijection
with the set of augmented coalgebra homomorphisms
\(\phi^\cp:Ts\ca^\op\boxt Ts\cc\to Ts\uCom\). Indeed, collections of
complexes \((\phi^\cp(X,Y),d)^{X\in\Ob\ca}_{Y\in\Ob\cc}\) are
identified with the \(\dg\)\n-spans \((\cp,sb^\cp_{00}s^{-1})\). In
particular, for each pair of objects \(X\in\Ob\ca\), \(Y\in\Ob\cc\)
holds \((\phi^\cp(X,Y))[1]=(s\cp(X,Y),-b^\cp_{00})\). The components
\(b^\cp_{kn}\) and \(\phi^\cp_{kn}\) are related for \((k,n)\ne(0,0)\)
by the formula
\begin{multline*}
b^\cp_{kn} = \bigl[
 s\ca(X_k,X_{k-1})\tdt s\ca(X_1,X_0)\tens s\cp(X_0,Y_0)\tens
 s\cc(Y_0,Y_1)\tdt s\cc(Y_{n-1},Y_n)
\\
\rTTo^{\tilde{\gamma}\tens1^{\tens n}}
\\
s\cp(X_0,Y_0)\tens s\ca^\op(X_0,X_1)\tdt s\ca^\op(X_{k-1},X_k)
 \tens s\cc(Y_0,Y_1)\tdt s\cc(Y_{n-1},Y_n)
\\
\rTTo^{1\tens \phi^\cp_{kn}\;{}}
 s\cp(X_0,Y_0)\tens s\uCom(\cp(X_0,Y_0),\cp(X_k,Y_n))
\\
\rTTo^{1\tens s^{-1}}
s\cp(X_0,Y_0)\tens\uCom(\cp(X_0,Y_0),\cp(X_k,Y_n))
\\
\rTTo^{1\tens[1]} s\cp(X_0,Y_0)\tens\uCom(s\cp(X_0,Y_0),s\cp(X_k,Y_n))
\rTTo^{\ev^\Com} s\cp(X_k,Y_n) \bigr],
\end{multline*}
where $\tilde{\gamma}=(12\dots k+1)\cdot\gamma$, and anti-isomorphism
$\gamma$ is defined by \eqref{eq-gamma-anti}.

The components of \(b^\cp\) can be
written in a more concise form. Given objects \(X,Y\in\Ob\ca\),
\(Z,W\in\Ob\cc\), define
\begin{multline}
\check b^\cp_+=\bigl[Ts\ca(Y,X)\tens s\cp(X,Z)\tens Ts\cc(Z,W)
\\
\rTTo^{c\tens1} s\cp(X,Z)\tens Ts\ca(Y,X)\tens Ts\cc(Z,W)
\\
\rTTo^{1\tens\gamma\tens1}
 s\cp(X,Z)\tens Ts\ca^\op(X,Y)\tens Ts\cc(Z,W)
\\
\rTTo^{1\tens\check\phi^\cp} s\cp(X,Z)\tens s\uCom(\cp(X,Z),\cp(Y,W))
\\
\hfill \rTTo^{(s\tens s)^{-1}} \cp(X,Z)\tens\uCom(\cp(X,Z),\cp(Y,W))
\rTTo^{\ev^{\Com}} \cp(Y,W) \rTTo^s s\cp(Y,W) \bigr] \quad
\\
\quad =\bigl[Ts\ca(Y,X)\tens s\cp(X,Z)\tens Ts\cc(Z,W) \rTTo^{c\tens1}
s\cp(X,Z)\tens Ts\ca(Y,X)\tens Ts\cc(Z,W) \hfill
\\
\rTTo^{1\tens\gamma\tens1}
 s\cp(X,Z)\tens Ts\ca^\op(X,Y)\tens Ts\cc(Z,W)
\\
\rTTo^{1\tens\check\phi^\cp} s\cp(X,Z)\tens s\uCom(\cp(X,Z),\cp(Y,W))
\\
\rTTo^{1\tens s^{-1}[1]} s\cp(X,Z)\tens\uCom(s\cp(X,Z),s\cp(Y,W))
\rTTo^{\ev^{\Com}} s\cp(Y,W) \bigr],
 \label{equ-check-b-P-+}
\end{multline}
where \(\gamma:Ts\ca\to Ts\ca^\op\) is the coalgebra anti-isomorphism
\eqref{eq-gamma-anti}, and
\(\check\phi^\cp=\phi^\cp\pr_1:Ts\ca^\op\boxt Ts\cc\to s\uCom\).
Conversely, components of the \ainf-functor \(\phi^\cp\) can be found
as
\begin{multline}
\check\phi^\cp=\bigl[
Ts\ca^\op(X,Y)\tens Ts\cc(Z,W)\rTTo^{\gamma\tens1}Ts\ca(Y,X)\tens Ts\cc(Z,W)
\\
\rTTo^{\coev^\Com}
\uCom(s\cp(X,Z),s\cp(X,Z)\tens Ts\ca(Y,X)\tens Ts\cc(Z,W))
\\
\rTTo^{\uCom(1,(c\tens1)\check b^\cp_+)}
\uCom(s\cp(X,Z),s\cp(Y,W))
\rTTo^{[-1]s}s\uCom(\cp(X,Z),\cp(Y,W)) \bigr].
\label{eq-check-phi-P-coev-bP}
\end{multline}
Define also
\begin{equation}
\check b^\cp_0=\bigl[Ts\ca(Y,X)\tens s\cp(X,Z)\tens Ts\cc(Z,W)
\rTTo^{\pr_0\tens1\tens\pr_0} s\cp(X,Z)\rTTo^{b^\cp_{00}}
s\cp(X,Z)\bigr].
 \label{equ-check-b-P-0}
\end{equation}
Note that \(\check b^\cp_+\) vanishes on \(T^0s\ca(Y,X)\tens
s\cp(X,Z)\tens T^0s\cc(Z,W)\) since \(\check\phi^\cp\) vanishes on
\(T^0s\ca^\op(X,Y)\tens T^0s\cc(Z,W)\). It follows that \(\check
b^\cp=\check b^\cp_+ +\check b^\cp_0\).

The following statement was proven by Lef\`evre-Hasegawa in assumption
that the ground ring is a field \cite[Lemme~5.3.0.1]{Lefevre-Ainfty-these}.

\begin{proposition}\label{prop-flat-connection-vs-A-inf-functor}
$b^\cp$ is a flat connection, that is, \((Ts\ca\tens s\cp\tens
Ts\cc,b^\cp)\) is a bicomodule in \(\dQ\), if and only if the
corresponding augmented coalgebra homomorphism
\(\phi^\cp:Ts\ca^\op\boxt Ts\cc\to Ts\uCom\)  is an \ainf-functor.
\end{proposition}

    \ifx\chooseClass1
\straightForward.
\proofInArXiv.
    \else
\begin{proof}
According to \eqref{equ-connection-via-check-r},
\begin{align*}
b^\cp&=(\Delta_0\tens1\tens\Delta_0)
(1\tens\pr_0\tens1\tens\pr_0\tens b^\cc+1\tens\check b^\cp\tens1
+b^\ca\tens\pr_0\tens1\tens\pr_0\tens1)
\\
&=1\tens1\tens b^\cc+(\Delta_0\tens1\tens\Delta_0)(1\tens\check
b^\cp\tens1)+b^\ca\tens1\tens1.
\end{align*}
The \(\kk\)\n-span morphism \((b^\cp)^2:P\to P\) is a
\((1,1,1,0,0)\)\n-connection of degree 2, therefore the equation
\((b^\cp)^2=0\) is equivalent to its particular case
\((b^\cp)^2(\pr_0\tens1\tens\pr_0)=0:Ts\ca\tens s\cp\tens Ts\cc\to
s\cp\). In terms of \(\check b^\cp\), the latter reads as follows:
\begin{equation}
(b^\ca\tens1\tens1+1\tens1\tens b^\cc)\check b^\cp+(\Delta_0\tens1\tens\Delta_0)(1\tens\check
b^\cp\tens1)\check b^\cp=0.
\label{equ-bP-0}
\end{equation}
Note that \((b^\ca\tens1\tens1+1\tens1\tens b^\cc)\check
b^\cp=(b^\ca\tens1\tens1+1\tens1\tens b^\cc)\check b^\cp_+\), since
\(b^\ca\pr_0=0\), \(b^\cc\pr_0=0\). The second term in the above
equation splits into the sum of four summands, which we are going to
compute separately. First of all,
\begin{multline*}
(\Delta_0\tens1\tens\Delta_0)(1\tens\check b^\cp_0\tens1)\check
b^\cp_0=(\Delta_0(\pr_0\tens\pr_0)\tens
1\tens\Delta(\pr_0\tens\pr_0))(b^\cp_{00})^2\\
=(\pr_0\tens1\tens\pr_0)s^{-1}d^2s=0:
Ts\ca\tens s\cp\tens Ts\cc\to s\cp.
\end{multline*}
Secondly,
\begin{multline*}
(\Delta_0\tens1\tens\Delta_0)(1\tens\check b^\cp_+\tens1)\check
b^\cp_0+(\Delta_0\tens1\tens\Delta_0)(1\tens\check b^\cp_0\tens1)\check
b^\cp_+\\
\quad=(\Delta_0(\pr_0\tens1)\tens1\tens\Delta_0(1\tens\pr_0))\check
b^\cp_+b^\cp_{00}+(\Delta_0(1\tens\pr_0)\tens b^\cp_{00}\tens\Delta(\pr_0\tens1))\check
b^\cp_+\hfill\\
\quad=\check b^\cp_+ b^\cp_{00} +(1\tens b^\cp_{00}\tens 1)\check b^\cp_+\hfill\\
\quad=\bigl[
Ts\ca(X,Y)\tens s\cp(Y,Z)\tens Ts\cc(Z,W)\rTTo^{c\tens1}
s\cp(Y,Z)\tens Ts\ca(X,Y)\tens Ts\cc(Z,W)\hfill
\\
\rTTo^{1\tens\gamma\tens1}s\cp(Y,Z)\tens Ts\ca^\op(Y,X)\tens
Ts\cc(Z,W)\rTTo^{1\tens\check\phi^\cp} s\cp(Y,Z)\tens
s\uCom(\cp(Y,Z),\cp(X,W))
\\
\rTTo^{1\tens s^{-1}}
s\cp(Y,Z)\tens\uCom(\cp(Y,Z),\cp(X,W))
\rTTo^{1\tens[1]}
s\cp(Y,Z)\tens\uCom(s\cp(Y,Z),s\cp(X,W))\\
\hfill\rTTo^{\ev^{\Com}}
s\cp(X,W)\rTTo^{b^\cp_{00}} s\cp(X,W)
\bigr]\quad\\
\quad+\bigl[
Ts\ca(X,Y)\tens s\cp(Y,Z)\tens Ts\cc(Z,W)\rTTo^{1\tens b^\cp_{00}\tens1}
Ts\ca(X,Y)\tens s\cp(Y,Z)\tens Ts\cc(Z,W)\hfill\\
\rTTo^{c\tens1}s\cp(Y,Z)\tens Ts\ca(X,Y)\tens Ts\cc(Z,W)\rTTo^{1\tens\gamma\tens1}
s\cp(Y,Z)\tens Ts\ca^\op(Y,X)\tens Ts\cc(Z,W)\\
\rTTo^{1\tens\check\phi^\cp}s\cp(Y,Z)\tens s\uCom(\cp(Y,Z),\cp(X,W))
\rTTo^{1\tens s^{-1}}
s\cp(Y,Z)\tens\uCom(\cp(Y,Z),\cp(X,W))\\
\hfill\rTTo^{1\tens[1]}
s\cp(Y,Z)\tens\uCom(s\cp(Y,Z),s\cp(X,W))
\rTTo^{\ev^{\Com}}s\cp(X,W)
\bigr]\quad\\
\quad=\bigl[
Ts\ca(X,Y)\tens s\cp(Y,Z)\tens Ts\cc(Z,W)\rTTo^{c\tens1}
s\cp(Y,Z)\tens Ts\ca(X,Y)\tens Ts\cc(Z,W)\hfill
\\
\rTTo^{1\tens\gamma\tens1}s\cp(Y,Z)\tens Ts\ca^\op(Y,X)\tens
Ts\cc(Z,W)\rTTo^{1\tens\check\phi^\cp} s\cp(Y,Z)\tens
s\uCom(\cp(Y,Z),\cp(X,W))
\\
\rTTo^{1\tens s^{-1}}
s\cp(Y,Z)\tens\uCom(\cp(Y,Z),\cp(X,W))
\rTTo^{1\tens[1]}
s\cp(Y,Z)\tens\uCom(s\cp(Y,Z),s\cp(X,W))\\
\rTTo^{\ev^{\Com}b^\cp_{00}-(b^\cp_{00}\tens1)\ev^{\Com}}
s\cp(X,W)
\bigr].
\end{multline*}
The complexes \(s\cp(Y,Z)\) and \(s\cp(X,W)\) carry the differential
\(-b^\cp_{00}\). Since \(\ev^{\Com}\) is a chain map, it follows that
\(\ev^{\Com}b^\cp_{00}-(b^\cp_{00}\tens1)\ev^{\Com}=-(1\tens m^{\uCom}_1)\ev^{\Com}\),
therefore
\begin{multline*}
(\Delta_0\tens1\tens\Delta_0)(1\tens\check b^\cp_+\tens1)\check
b^\cp_0+(\Delta_0\tens1\tens\Delta_0)(1\tens\check b^\cp_0\tens1)\check
b^\cp_+\\
\quad=-\bigl[
Ts\ca(X,Y)\tens s\cp(Y,Z)\tens Ts\cc(Z,W)\rTTo^{c\tens1}
s\cp(Y,Z)\tens Ts\ca(X,Y)\tens Ts\cc(Z,W)\hfill
\\
\rTTo^{1\tens\gamma\tens1}s\cp(Y,Z)\tens Ts\ca^\op(Y,X)\tens Ts\cc(Z,W)
\rTTo^{1\tens\check\phi^\cp} s\cp(Y,Z)\tens s\uCom(\cp(Y,Z),\cp(X,W))
\\
\rTTo^{1\tens s^{-1}} s\cp(Y,Z)\tens\uCom(\cp(Y,Z),\cp(X,W))
\rTTo^{1\tens[1]} s\cp(Y,Z)\tens\uCom(s\cp(Y,Z),s\cp(X,W))
\\
\rTTo^{(1\tens m^{\uCom}_1)\ev^{\Com}}s\cp(X,W)\bigr].
\end{multline*}
Since \([1]\) is a differential graded functor, it follows that
\([1]m^\uCom_1=m^\uCom_1[1]\). Together with the relation
\(b^\uCom_1s^{-1}=s^{-1}m^\uCom_1\) this implies that
\begin{multline*}
(\Delta_0\tens1\tens\Delta_0)(1\tens\check b^\cp_+\tens1)\check
b^\cp_0+(\Delta_0\tens1\tens\Delta_0)(1\tens\check b^\cp_0\tens1)\check
b^\cp_+\\
\quad=-\bigl[
Ts\ca(X,Y)\tens s\cp(Y,Z)\tens Ts\cc(Z,W)\rTTo^{c\tens1}
s\cp(Y,Z)\tens Ts\ca(X,Y)\tens Ts\cc(Z,W)\hfill
\\
\rTTo^{1\tens\gamma\tens1}s\cp(Y,Z)\tens Ts\ca^\op(Y,X)\tens Ts\cc(Z,W)
\\
\rTTo^{1\tens\check\phi^\cp b^{\uCom}_1}
 s\cp(Y,Z)\tens s\uCom(\cp(Y,Z),\cp(X,W))
\rTTo^{1\tens s^{-1}} s\cp(Y,Z)\tens\uCom(\cp(Y,Z),\cp(X,W))
\\
\rTTo^{1\tens[1]} s\cp(Y,Z)\tens\uCom(s\cp(Y,Z),s\cp(X,W))
\rTTo^{\ev^{\Com}}s\cp(X,W)\bigr].
\end{multline*}
Next, let us compute
\begin{multline*}
(\Delta_0\tens1\tens\Delta_0)(1\tens\check b^\cp_+\tens1)\check b^\cp_+
=\bigl[ Ts\ca(X,Y)\tens s\cp(Y,Z)\tens Ts\cc(Z,W)
\rTTo^{\Delta_0\tens1\tens\Delta_0}
\\
\bigoplus_{U\in\Ob\ca,V\in\Ob\cc} \hspace*{-2em}
 Ts\ca(X,U)\tens Ts\ca(U,Y)\tens s\cp(Y,Z)\tens Ts\cc(Z,V)\tens Ts\cc(V,W)
\rTTo^{\sum 1\tens c\tens1\tens1}
\\
\bigoplus_{U\in\Ob\ca,V\in\Ob\cc} \hspace*{-2em}
 Ts\ca(X,U)\tens s\cp(Y,Z)\tens Ts\ca(U,Y)\tens Ts\cc(Z,V)\tens Ts\cc(V,W)
\rTTo^{\sum 1\tens1\tens\gamma\tens1\tens1}
\\
\bigoplus_{U\in\Ob\ca,V\in\Ob\cc} \hspace*{-2em}
 Ts\ca(X,U)\tens s\cp(Y,Z)\tens Ts\ca^\op(Y,U)\tens Ts\cc(Z,V)\tens Ts\cc(V,W)
\rTTo^{\sum 1\tens1\tens\check\phi^\cp\tens1}
\\
\bigoplus_{U\in\Ob\ca,V\in\Ob\cc} \hspace*{-2em}
 Ts\ca(X,U)\tens s\cp(Y,Z)\tens s\uCom(\cp(Y,Z),\cp(U,V))\tens Ts\cc(V,W)
\rTTo^{\sum 1\tens 1\tens s^{-1}\tens1}
\\
\bigoplus_{U\in\Ob\ca,V\in\Ob\cc} \hspace*{-2em}
 Ts\ca(X,U)\tens s\cp(Y,Z)\tens\uCom(\cp(Y,Z),\cp(U,V))\tens Ts\cc(V,W)
\rTTo^{\sum 1\tens 1\tens [1]\tens1}
\\
\bigoplus_{U\in\Ob\ca,V\in\Ob\cc} \hspace*{-2em}
 Ts\ca(X,U)\tens s\cp(Y,Z)\tens\uCom(s\cp(Y,Z),s\cp(U,V))\tens Ts\cc(V,W)
\rTTo^{\sum 1\tens\ev^{\Com}\tens1}
\\
\bigoplus_{U\in\Ob\ca,V\in\Ob\cc}Ts\ca(X,U)\tens
s\cp(U,V)\tens Ts\cc(V,W)
\rTTo^{\sum c\tens 1}\\
\bigoplus_{U\in\Ob\ca,V\in\Ob\cc}s\cp(U,V)\tens Ts\ca(X,U)\tens Ts\cc(V,W)
\rTTo^{\sum 1\tens\gamma\tens1}\\
\bigoplus_{U\in\Ob\ca,V\in\Ob\cc}s\cp(U,V)\tens Ts\ca^\op(U,X)\tens Ts\cc(V,W)
\rTTo^{\sum 1\tens\check\phi^\cp}\\
\bigoplus_{U\in\Ob\ca,V\in\Ob\cc}s\cp(U,V)\tens s\uCom(\cp(U,V),\cp(X,W))
\rTTo^{\sum 1\tens s^{-1}}\\
\bigoplus_{U\in\Ob\ca,V\in\Ob\cc}s\cp(U,V)\tens \uCom(\cp(U,V),\cp(X,W))
\rTTo^{\sum 1\tens [1]}\\
\bigoplus_{U\in\Ob\ca,V\in\Ob\cc}s\cp(U,V)\tens \uCom(s\cp(U,V),s\cp(X,W))
\rTTo^{\sum \ev^{\Com}}s\cp(X,W)
\bigr].
\end{multline*}
The latter can be written as
\begin{multline*}
\bigl[
Ts\ca(X,Y)\tens s\cp(Y,Z)\tens Ts\cc(Z,W)\rTTo^{c\tens1}
s\cp(Y,Z)\tens Ts\ca(X,Y)\tens Ts\cc(Z,W)
\\
\rTTo^{1\tens\Delta_0 c(\gamma\tens\gamma)\tens\Delta_0} \hspace*{-2.3em}
\bigoplus_{U\in\Ob\ca,V\in\Ob\cc} \hspace*{-2.2em}
 s\cp(Y,Z)\tens Ts\ca^\op(Y,U)\tens Ts\ca^\op(U,X)\tens Ts\cc(Z,V)\tens Ts\cc(V,W)
\\
\rTTo^{\sum 1\tens1\tens c\tens1} \hspace*{-2em}
\bigoplus_{U\in\Ob\ca,V\in\Ob\cc} \hspace*{-2em}
 s\cp(Y,Z)\tens Ts\ca^\op(Y,U)\tens Ts\cc(Z,V)\tens Ts\ca^\op(U,X)\tens Ts\cc(V,W)
\\
\rTTo^{\sum 1\tens\check\phi^\cp\tens\check\phi^\cp} \hspace*{-1em}
\bigoplus_{U\in\Ob\ca,V\in\Ob\cc} \hspace*{-1em}
 s\cp(Y,Z)\tens s\uCom(\cp(Y,Z),\cp(U,V))\tens s\uCom(\cp(U,V),\cp(X,W))
\\
\rTTo^{\sum 1\tens s^{-1}\tens s^{-1}} \hspace*{-1em}
\bigoplus_{U\in\Ob\ca,V\in\Ob\cc} \hspace*{-1em}
 s\cp(Y,Z)\tens\uCom(\cp(Y,Z),\cp(U,V))\tens\uCom(\cp(U,V),\cp(X,W))
\\
\rTTo^{\sum 1\tens[1]\tens[1]} \hspace*{-1em}
\bigoplus_{U\in\Ob\ca,V\in\Ob\cc} \hspace*{-1em}
 s\cp(Y,Z)\tens\uCom(s\cp(Y,Z),s\cp(U,V))\tens\uCom(s\cp(U,V),s\cp(X,W))
\\
\rTTo^{\sum (\ev^{\Com}\tens1)\ev^{\Com}} s\cp(X,W) \bigr].
\end{multline*}
Using the identities \(\Delta_0c(\gamma\tens\gamma)=\gamma\Delta_0\)
(see \cite[Section~A.4]{math.CT/0306018}),
\((\ev^{\Com}\tens1)\ev^{\Com}=(1\tens m^{\uCom}_2)\ev^{\Com}\), and
\((s^{-1}\tens s^{-1})m^{\uCom}_2=-b^{\uCom}_2s^{-1}\), we
transform the above expression as follows:
\begin{multline*}
(\Delta_0\tens1\tens\Delta_0)(1\tens\check b^\cp_+\tens1)\check b^\cp_+
=-\bigl[ Ts\ca(X,Y)\tens s\cp(Y,Z)\tens Ts\cc(Z,W) \rTTo^{c\tens1}
\\
s\cp(Y,Z)\tens Ts\ca(X,Y)\tens Ts\cc(Z,W) \rTTo^{1\tens\gamma\tens1}
s\cp(Y,Z)\tens Ts\ca^\op(Y,X)\tens Ts\cc(Z,W)
\\
\rTTo^{1\tens\Delta_0\tens\Delta_0} \hspace*{-2em}
\bigoplus_{U\in\Ob\ca,V\in\Ob\cc} \hspace*{-2em}
 s\cp(Y,Z)\tens Ts\ca^\op(Y,U)\tens Ts\ca^\op(U,X)\tens Ts\cc(Z,V)\tens Ts\cc(V,W)
\\
\rTTo^{\sum 1\tens1\tens c\tens 1} \hspace*{-2em}
\bigoplus_{U\in\Ob\ca,V\in\Ob\cc} \hspace*{-2em}
 s\cp(Y,Z)\tens Ts\ca^\op(Y,U)\tens Ts\cc(Z,V)\tens Ts\ca^\op(U,X)\tens Ts\cc(V,W)
\\
\rTTo^{\sum 1\tens\check\phi^\cp\tens\check\phi^\cp} \hspace*{-1em}
\bigoplus_{U\in\Ob\ca,V\in\Ob\cc} \hspace*{-1em}
s\cp(Y,Z)\tens s\uCom(\cp(Y,Z),\cp(U,V))\tens s\uCom(\cp(U,V),\cp(X,W))
\\
\rTTo^{\sum 1\tens b^{\uCom}_2}
s\cp(Y,Z)\tens s\uCom(\cp(Y,Z),\cp(X,W))
\rTTo^{1\tens s^{-1}}
s\cp(Y,Z)\tens \uCom(\cp(Y,Z),\cp(X,W))\\
\rTTo^{1\tens [1]}
s\cp(Y,Z)\tens \uCom(s\cp(Y,Z),s\cp(X,W))
\rTTo^{\ev^{\Com}}s\cp(X,W)
\bigr].
\end{multline*}
Finally,
\begin{multline*}
(b^\ca\tens1\tens1+1\tens1\tens b^\cc)\check b^\cp_+
=\bigl[
Ts\ca(X,Y)\tens s\cp(Y,Z)\tens Ts\cc(Z,W)\rTTo^{c\tens1}\\
s\cp(Y,Z)\tens Ts\ca(X,Y)\tens Ts\cc(Z,W)\rTTo^{1\tens\gamma\tens1}
s\cp(Y,Z)\tens Ts\ca^\op(Y,X)\tens Ts\cc(Z,W)\\
\rTTo^{1\tens (b^{\ca^\op}\tens1+1\tens b^\cc)\check\phi^\cp}
s\cp(Y,Z)\tens s\uCom(\cp(Y,Z),\cp(X,W))
\\
\rTTo^{(s\tens s)^{-1}} \cp(Y,Z)\tens\uCom(\cp(Y,Z),\cp(X,W))
\rTTo^{\ev^{\Com}} \cp(X,W) \rTTo^s s\cp(X,W) \bigr],
\end{multline*}
since \(b^{\ca^\op}=\gamma b^\ca\gamma:Ts\ca^\op(Y,X)\to
Ts\ca^\op(Y,X)\). We conclude that the left hand side of
\eqref{equ-bP-0} equals
 \((c\tens1)(1\tens\gamma\tens1)(1\tens R)
 (s\tens s)^{-1}\ev^{\Com}s^{-1}\),
where
\begin{multline*}
R=\bigl[Ts\ca^{\op}(Y,X)\tens Ts\cc(Z,W)\rTTo^{b^{\ca^\op}\tens1+1\tens b^\cc}
Ts\ca^{\op}(Y,X)\tens Ts\cc(Z,W)\\
\hfill\rTTo^{\check\phi^\cp}s\uCom(\phi^\cp(Y,Z),\phi^\cp(X,W))
\bigr]\quad\\
\quad-\check\phi^\cp b^{\uCom}_1
-\bigl[
Ts\ca^\op(Y,X)\tens Ts\cc(Z,W)\rTTo^{\Delta_0\tens\Delta_0}\hfill\\
\bigoplus_{U\in\Ob\ca,V\in\Ob\cc}
Ts\ca^\op(Y,U)\tens Ts\ca^\op(U,X)\tens Ts\cc(Z,V)\tens Ts\cc(V,W)
\rTTo^{\sum 1\tens c\tens 1}\\
\bigoplus_{U\in\Ob\ca,V\in\Ob\cc}
Ts\ca^\op(Y,U)\tens Ts\cc(Z,V)\tens Ts\ca^\op(U,X)\tens Ts\cc(V,W)
\rTTo^{\sum \check\phi^\cp\tens\check\phi^\cp}\\
\bigoplus_{U\in\Ob\ca,V\in\Ob\cc}
s\uCom(\phi^\cp(Y,Z),\phi^\cp(U,V))\tens
s\uCom(\phi^\cp(U,V),\phi^\cp(X,W))\\
\rTTo^{\sum b^{\uCom}_2}
s\uCom(\phi^\cp(Y,Z),\phi^\cp(X,W))
\bigr].
\end{multline*}
By closedness, \(b^\cp\) is a flat connection if and only if \(R=0\),
for all objects \(X,Y\in\Ob\ca\), \(Z,W\in\Ob\cc\), that is, if
\(\phi^\cp\) is an \ainf-functor.
\end{proof}
\fi

Let \(\ca\), \(\cc\) be \ainf-categories. The full subcategory of the
differential graded category \(Ts\ca\text-Ts\cc\bicomod\) consisting of
$\dg$\n-bicomodules whose underlying graded bicomodule has the form
\(Ts\ca\tens s\cp\tens Ts\cc\) is denoted by \(\ca\text-\cc\bimod\).
Its objects are called \emph{\ainf-bimodules}, extending the
terminology of Tradler \cite{xxx0108027}.

\begin{proposition}\label{pro-dg-cat-ACbimod-A(ACC)-isomorphic}
The differential graded categories \(\ca\text-\cc\bimod\) and
\(\und\Ainfty(\ca^\op,\cc;\uCom)\) are isomorphic.
\end{proposition}

    \ifx\chooseClass1
\propref{prop-flat-connection-vs-A-inf-functor} establishes a bijection
between the sets of objects of the differential graded categories
\(\und\Ainfty(\ca^\op,\cc;\uCom)\) and \(\ca\text-\cc\bimod\). It
extends to an isomorphism of differential graded categories.
\proofInArXiv.
    \else
\begin{proof}
\propref{prop-flat-connection-vs-A-inf-functor} establishes a bijection
between the sets of objects of the differential graded categories
\(\und\Ainfty(\ca^\op,\cc;\uCom)\) and \(\ca\text-\cc\bimod\). Let us
extend it to an isomorphism of differential graded categories. Let
\(\phi,\psi:\ca^\op,\cc\to\uCom\) be \ainf-functors, \(\cp\), \(\cq\)
the corresponding \(\ca\text-\cc\)\n-bimodules. Define a
\(\kk\)\n-linear map
\(\Phi:\und\Ainfty(\ca^\op,\cc;\uCom)(\phi,\psi)\to\ca\text-\cc\bimod(\cp,\cq)\)
of degree \(0\) as follows. An element
\(rs^{-1}\in\und\Ainfty(\ca^\op,\cc;\uCom)(\phi,\psi)\) is mapped to an
\((\id_{Ts\ca},\id_{Ts\cc})\)\n-bicomodule homomorphism
 \(t=(rs^{-1})\Phi:
 Ts\ca\tens s\cp\tens Ts\cc\to Ts\ca\tens s\cq\tens Ts\cc\)
given by its components
\begin{multline*}
t_{kn}=(-)^{r+1}\bigl[
\\
s\ca(X_k,X_{k-1})\tdt s\ca(X_1,X_0)\tens s\cp(X_0,Y_0)\tens s\cc(Y_0,Y_1)\tdt s\cc(Y_{n-1},Y_n)
\rTTo^{\tilde\gamma\tens1^{\tens n}}
\\
s\cp(X_0,Y_0)\tens s\ca^\op(X_0,X_1)\tdt s\ca^\op(X_{k-1},X_k)\tens
s\cc(Y_0,Y_1)\tdt s\cc(Y_{n-1},Y_n)
\\
\rTTo^{1\tens r_{kn}}
 s\cp(X_0,Y_0)\tens s\uCom(\cp(X_0,Y_0),\cq(X_k,Y_n))
\rTTo^{1\tens s^{-1}[1]}
\\
s\cp(X_0,Y_0)\tens\uCom(s\cp(X_0,Y_0),s\cq(X_k,Y_n)) \rTTo^{\ev^\Com}
s\cq(X_k,Y_n) \bigr], \quad k,n\ge0,
\end{multline*}
or more concisely,
\begin{multline*}
\check t=(-)^{r+1}\bigl[
Ts\ca(X,Y)\tens s\cp(Y,Z)\tens Ts\cc(Z,W)
\\
\rTTo^{c\tens 1} s\cp(Y,Z)\tens Ts\ca(X,Y)\tens Ts\cc(Z,W)
\rTTo^{1\tens\gamma\tens1}
\\
s\cp(Y,Z)\tens Ts\ca^\op(Y,X)\tens Ts\cc(Z,W)\rTTo^{1\tens\check r}
s\cp(Y,Z)\tens s\uCom(\cp(Y,Z),\cq(X,W))
\\
\rTTo^{1\tens s^{-1}[1]}
s\cp(Y,Z)\tens \uCom(s\cp(Y,Z),s\cq(X,W))\rTTo^{\ev^\Com}s\cq(X,W)
\bigr],
\end{multline*}
where \(\check r=r\cdot\pr_1\in\sspan(Ts\ca^\op\boxt
Ts\cc,s\uCom)(\Ob\phi,\Ob\psi)\). Closedness of \(\gr\) implies that
the map
\(\Phi:\und\Ainfty(\ca^\op,\cc;\uCom)(\phi,\psi)\to\ca\text-\cc\bimod(\cp,\cq)\)
is an isomorphism. Let us prove that it also commutes with the
differential. We must prove that
\(((rs^{-1})\Phi)d=(rs^{-1}m^{\und\Ainfty(\ca^\op,\cc;\uCom)}_1)\Phi=((rB_1)s^{-1})\Phi\)
for each element \(r\in s\und\Ainfty(\ca^\op,\cc;\uCom)(\phi,\psi)\).
Since the both sides of the equation are
\((\id_{Ts\ca},\id_{Ts\cc})\)\n-bimodule homomorphisms of degree \(\deg
r+1\), it suffices to prove the equation
\([((rs^{-1})\Phi)d]^\vee=[((rB_1)s^{-1})\Phi]^\vee\). Using
\eqref{equ-connection-via-check-r}, we obtain:
\begin{align}
[((rs^{-1})\Phi)d]^\vee&=(td)^\vee=(t\cdot b^\cq)^\vee-(-)^t(b^\cp\cdot t)^\vee
=t\cdot\check b^\cq-(-)^tb^\cp\cdot\check t\nonumber
\\
&=t\cdot\check b^\cq_+\label{equ-t-check-b+}
\\
&-(-)^t(\Delta_0\tens1\tens\Delta_0)(1\tens\check b^\cp_+\tens1)\check
t\label{equ-Delta1Delta-1b+1-t}
\\
&+t\cdot\check b^\cq_0-(-)^t(\Delta_0\tens1\tens\Delta_0)(1\tens\check b^\cp_0\tens1)\check t
\label{equ-tb0-Delta1Delta-1b01-t}
\\
&-(-)^t(b^\ca\tens1\tens1+1\tens1\tens b^\cc)\check t.
\label{equ-b11-11b-t}
\end{align}
Let us compute summands \eqref{equ-t-check-b+}--\eqref{equ-b11-11b-t}
separately. According to \eqref{equ-f-Delta1Delta-phi-f-psi},
expression~\eqref{equ-t-check-b+} equals
\begin{multline*}
t\cdot\check
b^\cq_+=(\Delta_0\tens1\tens\Delta_0)(1\tens\check t\tens1)\check b^\cq_+
\\
\qquad =(-)^{r+1}\bigl[ Ts\ca(X,Y)\tens s\cp(Y,Z)\tens
Ts\cc(Z,W)\rTTo^{\Delta_0\tens1\tens\Delta_0} \hfill
\\
\bigoplus_{U\in\Ob\ca,V\in\Ob\cc} \hspace*{-2em}
 Ts\ca(X,U)\tens Ts\ca(U,Y)\tens s\cp(Y,Z)\tens Ts\cc(Z,V)\tens Ts\cc(V,W)
\rTTo^{\sum 1\tens c\tens1\tens1}
\\
\bigoplus_{U\in\Ob\ca,V\in\Ob\cc} \hspace*{-2em}
 Ts\ca(X,U)\tens s\cp(Y,Z)\tens Ts\ca(U,Y)\tens Ts\cc(Z,V)\tens Ts\cc(V,W)
\rTTo^{\sum 1\tens1\tens\gamma\tens1\tens1}
\\
\bigoplus_{U\in\Ob\ca,V\in\Ob\cc} \hspace*{-2em}
 Ts\ca(X,U)\tens s\cp(Y,Z)\tens Ts\ca^\op(Y,U)\tens Ts\cc(Z,V)\tens Ts\cc(V,W)
\rTTo^{\sum 1\tens1\tens\check r\tens1}
\\
\bigoplus_{U\in\Ob\ca,V\in\Ob\cc} \hspace*{-2em}
 Ts\ca(X,U)\tens s\cp(Y,Z)\tens s\uCom(\cp(Y,Z),\cq(U,V))\tens Ts\cc(V,W)
\rTTo^{\sum 1\tens 1\tens s^{-1}[1]\tens1}
\\
\bigoplus_{U\in\Ob\ca,V\in\Ob\cc} \hspace*{-2em}
 Ts\ca(X,U)\tens s\cp(Y,Z)\tens\uCom(s\cp(Y,Z),s\cq(U,V))\tens Ts\cc(V,W)
\rTTo^{\sum 1\tens\ev^{\Com}\tens1}
\\
\bigoplus_{U\in\Ob\ca,V\in\Ob\cc}Ts\ca(X,U)\tens s\cq(U,V)\tens
Ts\cc(V,W) \rTTo^{\sum c\tens 1}
\\
\bigoplus_{U\in\Ob\ca,V\in\Ob\cc}s\cq(U,V)\tens Ts\ca(X,U)\tens Ts\cc(V,W)
\rTTo^{\sum 1\tens\gamma\tens1}
\\
\bigoplus_{U\in\Ob\ca,V\in\Ob\cc}s\cq(U,V)\tens Ts\ca^\op(U,X)\tens Ts\cc(V,W)
\rTTo^{\sum 1\tens\check\psi}
\\
\bigoplus_{U\in\Ob\ca,V\in\Ob\cc}s\cq(U,V)\tens s\uCom(\cq(U,V),\cq(X,W))
\rTTo^{\sum 1\tens s^{-1}[1]}
\\
\bigoplus_{U\in\Ob\ca,V\in\Ob\cc}s\cq(U,V)\tens \uCom(s\cq(U,V),s\cq(X,W))
\rTTo^{\sum \ev^{\Com}}s\cq(X,W)
\bigr].
\end{multline*}
As in the proof of \propref{prop-flat-connection-vs-A-inf-functor}, the
above composite can be transformed as follows:
\begin{multline*}
(-)^{r+1}\bigl[ Ts\ca(X,Y)\tens s\cp(Y,Z)\tens Ts\cc(Z,W)
\\
\rTTo^{c\tens1} s\cp(Y,Z)\tens Ts\ca(X,Y)\tens Ts\cc(Z,W)
\\
\rTTo^{1\tens\Delta_0 c(\gamma\tens\gamma)\tens\Delta_0} \hspace*{-2.3em}
\bigoplus_{U\in\Ob\ca,V\in\Ob\cc} \hspace*{-2.2em}
 s\cp(Y,Z)\tens Ts\ca^\op(Y,U)\tens Ts\ca^\op(U,X)\tens Ts\cc(Z,V)\tens Ts\cc(V,W)
\\
\rTTo^{\sum 1\tens1\tens c\tens1} \hspace*{-2em}
\bigoplus_{U\in\Ob\ca,V\in\Ob\cc} \hspace*{-2em}
 s\cp(Y,Z)\tens Ts\ca^\op(Y,U)\tens Ts\cc(Z,V)\tens Ts\ca^\op(U,X)\tens Ts\cc(V,W)
\\
\rTTo^{\sum 1\tens\check r\tens\check\psi} \hspace*{-1em}
\bigoplus_{U\in\Ob\ca,V\in\Ob\cc} \hspace*{-1em}
 s\cp(Y,Z)\tens s\uCom(\cp(Y,Z),\cq(U,V))\tens s\uCom(\cq(U,V),\cq(X,W))
\\
\rTTo^{\sum 1\tens s^{-1}[1]\tens s^{-1}[1]} \hspace*{-2.2em}
\bigoplus_{U\in\Ob\ca,V\in\Ob\cc} \hspace*{-2em}
 s\cp(Y,Z)\tens \uCom(s\cp(Y,Z),s\cq(U,V))\tens\uCom(s\cq(U,V),s\cq(X,W))
\\
\rTTo^{\sum (\ev^{\Com}\tens1)\ev^{\Com}} s\cq(X,W) \bigr].
\end{multline*}
Applying the already mentioned identities
\(\Delta_0c(\gamma\tens\gamma)=\gamma\Delta_0\),
\((\ev^{\Com}\tens1)\ev^{\Com}=(1\tens m^{\uCom}_2)\ev^{\Com}\), and
\((s^{-1}\tens s^{-1})m^{\uCom}_2=-b^{\uCom}_2s^{-1}\), we find:
\begin{multline*}
t\cdot\check b^\cq_+ =(-)^r\bigl[ Ts\ca(X,Y)\tens s\cp(Y,Z)\tens
Ts\cc(Z,W) \rTTo^{c\tens1}
\\
s\cp(Y,Z)\tens Ts\ca(X,Y)\tens Ts\cc(Z,W) \rTTo^{1\tens\gamma\tens1}
s\cp(Y,Z)\tens Ts\ca^\op(Y,X)\tens Ts\cc(Z,W)
\\
\rTTo^{1\tens\Delta_0\tens\Delta_0} \hspace*{-2em}
\bigoplus_{U\in\Ob\ca,V\in\Ob\cc} \hspace*{-2em}
 s\cp(Y,Z)\tens Ts\ca^\op(Y,U)\tens Ts\ca^\op(U,X)\tens Ts\cc(Z,V)\tens Ts\cc(V,W)
\\
\rTTo^{\sum 1\tens1\tens c\tens 1} \hspace*{-2em}
\bigoplus_{U\in\Ob\ca,V\in\Ob\cc} \hspace*{-2em}
 s\cp(Y,Z)\tens Ts\ca^\op(Y,U)\tens Ts\cc(Z,V)\tens Ts\ca^\op(U,X)\tens Ts\cc(V,W)
\\
\rTTo^{\sum 1\tens\check r\tens\check\psi} \hspace*{-2em}
\bigoplus_{U\in\Ob\ca,V\in\Ob\cc} \hspace*{-2em}
 s\cp(Y,Z)\tens s\uCom(\cp(Y,Z),\cq(U,V))\tens s\uCom(\cq(U,V),\cq(X,W))
\\
\rTTo^{\sum 1\tens b^{\uCom}_2}
 s\cp(Y,Z)\tens s\uCom(\cp(Y,Z),\cq(X,W))
\\
\rTTo^{1\tens s^{-1}[1]} s\cp(Y,Z)\tens \uCom(s\cp(Y,Z),s\cq(X,W))
\rTTo^{\ev^{\Com}}s\cq(X,W) \bigr].
\end{multline*}
Similarly, composite~\eqref{equ-Delta1Delta-1b+1-t} equals
\begin{multline*}
-(-)^t(\Delta_0\tens1\tens\Delta_0)(1\tens\check b^\cp_+\tens1)\check t
=-\bigl[ Ts\ca(X,Y)\tens s\cp(Y,Z)\tens
Ts\cc(Z,W)\rTTo^{\Delta_0\tens1\tens\Delta_0}
\\
\bigoplus_{U\in\Ob\ca,V\in\Ob\cc} \hspace*{-2em}
 Ts\ca(X,U)\tens Ts\ca(U,Y)\tens s\cp(Y,Z)\tens Ts\cc(Z,V)\tens Ts\cc(V,W)
\rTTo^{\sum 1\tens c\tens1\tens1}
\\
\bigoplus_{U\in\Ob\ca,V\in\Ob\cc} \hspace*{-2em}
 Ts\ca(X,U)\tens s\cp(Y,Z)\tens Ts\ca(U,Y)\tens Ts\cc(Z,V)\tens Ts\cc(V,W)
\rTTo^{\sum 1\tens1\tens\gamma\tens1\tens1}
\\
\bigoplus_{U\in\Ob\ca,V\in\Ob\cc} \hspace*{-2em}
 Ts\ca(X,U)\tens s\cp(Y,Z)\tens Ts\ca^\op(Y,U)\tens Ts\cc(Z,V)\tens Ts\cc(V,W)
\rTTo^{\sum 1\tens1\tens\check\phi\tens1}
\\
\bigoplus_{U\in\Ob\ca,V\in\Ob\cc} \hspace*{-2em}
 Ts\ca(X,U)\tens s\cp(Y,Z)\tens s\uCom(\cp(Y,Z),\cp(U,V))\tens Ts\cc(V,W)
\rTTo^{\sum 1\tens 1\tens s^{-1}[1]\tens1}
\\
\bigoplus_{U\in\Ob\ca,V\in\Ob\cc} \hspace*{-2em}
 Ts\ca(X,U)\tens s\cp(Y,Z)\tens\uCom(s\cp(Y,Z),s\cp(U,V))\tens Ts\cc(V,W)
\rTTo^{\sum 1\tens\ev^{\Com}\tens1}
\\
\bigoplus_{U\in\Ob\ca,V\in\Ob\cc}Ts\ca(X,U)\tens
s\cp(U,V)\tens Ts\cc(V,W)
\rTTo^{\sum c\tens 1}
\\
\bigoplus_{U\in\Ob\ca,V\in\Ob\cc}s\cp(U,V)\tens Ts\ca(X,U)\tens Ts\cc(V,W)
\rTTo^{\sum 1\tens\gamma\tens1}
\\
\bigoplus_{U\in\Ob\ca,V\in\Ob\cc}s\cp(U,V)\tens Ts\ca^\op(U,X)\tens Ts\cc(V,W)
\rTTo^{\sum 1\tens\check r}
\\
\bigoplus_{U\in\Ob\ca,V\in\Ob\cc}s\cp(U,V)\tens s\uCom(\cp(U,V),\cq(X,W))
\rTTo^{\sum 1\tens s^{-1}[1]}
\\
\hfill\bigoplus_{U\in\Ob\ca,V\in\Ob\cc}s\cp(U,V)\tens \uCom(s\cp(U,V),s\cq(X,W))
\rTTo^{\sum \ev^{\Com}}s\cq(X,W)
\bigr]\quad
\\
\quad =(-)^{r+1}\bigl[ Ts\ca(X,Y)\tens s\cp(Y,Z)\tens Ts\cc(Z,W) \hfill
\\
\rTTo^{c\tens1} s\cp(Y,Z)\tens Ts\ca(X,Y)\tens Ts\cc(Z,W)
\\
\rTTo^{1\tens\Delta_0 c(\gamma\tens\gamma)\tens\Delta_0} \hspace*{-2.3em}
\bigoplus_{U\in\Ob\ca,V\in\Ob\cc} \hspace*{-2.2em}
 s\cp(Y,Z)\tens Ts\ca^\op(Y,U)\tens Ts\ca^\op(U,X)\tens Ts\cc(Z,V)\tens Ts\cc(V,W)
\\
\rTTo^{\sum 1\tens1\tens c\tens1} \hspace*{-2em}
\bigoplus_{U\in\Ob\ca,V\in\Ob\cc} \hspace*{-2em}
 s\cp(Y,Z)\tens Ts\ca^\op(Y,U)\tens Ts\cc(Z,V)\tens Ts\ca^\op(U,X)\tens Ts\cc(V,W)
\\
\rTTo^{\sum 1\tens\check\phi\tens\check r} \hspace*{-1em}
\bigoplus_{U\in\Ob\ca,V\in\Ob\cc} \hspace*{-1em}
 s\cp(Y,Z)\tens s\uCom(\cp(Y,Z),\cp(U,V))\tens s\uCom(\cp(U,V),\cq(X,W))
\\
\rTTo^{\sum 1\tens s^{-1}[1]\tens s^{-1}[1]} \hspace*{-2em}
\bigoplus_{U\in\Ob\ca,V\in\Ob\cc} \hspace*{-2em}
 s\cp(Y,Z)\tens\uCom(s\cp(Y,Z),s\cp(U,V))\tens\uCom(s\cp(U,V),s\cq(X,W))
\\
\hfill\rTTo^{\sum (\ev^{\Com}\tens1)\ev^{\Com}} s\cq(X,W) \bigr]\quad
\\
=(-)^r\bigl[
Ts\ca(X,Y)\tens s\cp(Y,Z)\tens Ts\cc(Z,W)
\rTTo^{c\tens1}\hfill
\\
s\cp(Y,Z)\tens Ts\ca(X,Y)\tens Ts\cc(Z,W) \rTTo^{1\tens\gamma\tens1}
s\cp(Y,Z)\tens Ts\ca^\op(Y,X)\tens Ts\cc(Z,W)
\\
\rTTo^{1\tens\Delta_0\tens\Delta_0} \hspace*{-2em}
\bigoplus_{U\in\Ob\ca,V\in\Ob\cc} \hspace*{-2em}
 s\cp(Y,Z)\tens Ts\ca^\op(Y,U)\tens Ts\ca^\op(U,X)\tens Ts\cc(Z,V)\tens Ts\cc(V,W)
\\
\rTTo^{\sum 1\tens1\tens c\tens 1} \hspace*{-2em}
\bigoplus_{U\in\Ob\ca,V\in\Ob\cc} \hspace*{-2em}
 s\cp(Y,Z)\tens Ts\ca^\op(Y,U)\tens Ts\cc(Z,V)\tens Ts\ca^\op(U,X)\tens Ts\cc(V,W)
\\
\rTTo^{\sum 1\tens\check\phi\tens\check r} \hspace*{-1em}
\bigoplus_{U\in\Ob\ca,V\in\Ob\cc} \hspace*{-1em}
 s\cp(Y,Z)\tens s\uCom(\cp(Y,Z),\cp(U,V))\tens s\uCom(\cp(U,V),\cq(X,W))
\\
\rTTo^{\sum 1\tens b^{\uCom}_2}
s\cp(Y,Z)\tens s\uCom(\cp(Y,Z),\cq(X,W))
\\
\rTTo^{1\tens s^{-1}[1]} s\cp(Y,Z)\tens \uCom(s\cp(Y,Z),s\cq(X,W))
\rTTo^{\ev^{\Com}}s\cq(X,W) \bigr].
\end{multline*}
By \eqref{equ-f-Delta1Delta-phi-f-psi},
expression~\eqref{equ-tb0-Delta1Delta-1b01-t} can be written as
follows:
\begin{multline*}
t\cdot\check b^\cq_0-(-)^t(\Delta_0\tens1\tens\Delta_0)(1\tens\check
b^\cp_0\tens1)\check t
\\
\quad=(\Delta_0\tens1\tens\Delta_0)(1\tens\check
t\tens1)(\pr_0\tens1\tens\pr_0)b^\cq_{00}\hfill
\\
\hfill-(-)^t(\Delta_0\tens1\tens\Delta_0)(1\tens\pr_0\tens1\tens\pr_0\tens1)(1\tens
b^\cp_{00}\tens1)\check t\quad
\\
=\check t\cdot b^\cq_{00}-(-)^t(1\tens b^\cp_{00}\tens1)\check t,
\end{multline*}
therefore
\begin{multline*}
t\cdot\check b^\cq_0-(-)^t(\Delta_0\tens1\tens\Delta_0)(1\tens\check
b^\cp_0\tens1)\check t
\\
\quad =(-)^{r+1}\bigl[ Ts\ca(X,Y)\tens s\cp(Y,Z)\tens Ts\cc(Z,W) \hfill
\\
\rTTo^{c\tens1} s\cp(Y,Z)\tens Ts\ca(X,Y)\tens Ts\cc(Z,W)
\rTTo^{1\tens\gamma\tens1}
\\
s\cp(Y,Z)\tens Ts\ca^\op(Y,X)\tens Ts\cc(Z,W)\rTTo^{1\tens\check r}
s\cp(Y,Z)\tens s\uCom(\cp(Y,Z),\cq(X,W))
\\
\hfill\rTTo^{1\tens s^{-1}[1]} s\cp(Y,Z)\tens\uCom(s\cp(Y,Z),s\cq(X,W))
\rTTo^{\ev^\Com} s\cq(X,W) \rTTo^{b^\cq_{00}} s\cq(X,W) \bigr]\quad
\\
-\bigl[ Ts\ca(X,Y)\tens s\cp(Y,Z)\tens Ts\cc(Z,W)\rTTo^{1\tens
b^\cp_{00}\tens1} Ts\ca(X,Y)\tens s\cp(Y,Z)\tens
Ts\cc(Z,W)\hfill
\\
\rTTo^{c\tens1}s\cp(Y,Z)\tens Ts\ca(X,Y)\tens Ts\cc(Z,W)
\rTTo^{1\tens\gamma\tens1}
\\
s\cp(Y,Z)\tens Ts\ca^\op(Y,X)\tens Ts\cc(Z,W) \rTTo^{1\tens\check r}
s\cp(Y,Z)\tens s\uCom(\cp(Y,Z),\cq(X,W))
\\
\hfill \rTTo^{1\tens s^{-1}[1]}
s\cp(Y,Z)\tens\uCom(s\cp(Y,Z),s\cq(X,W)) \rTTo^{\ev^\Com}s\cq(X,W)
\bigr]\quad
\\
\quad =(-)^{r+1}\bigl[ Ts\ca(X,Y)\tens s\cp(Y,Z)\tens Ts\cc(Z,W) \hfill
\\
\rTTo^{c\tens1} s\cp(Y,Z)\tens Ts\ca(X,Y)\tens Ts\cc(Z,W)
\rTTo^{1\tens\gamma\tens1}
\\
s\cp(Y,Z)\tens Ts\ca^\op(Y,X)\tens Ts\cc(Z,W) \rTTo^{1\tens\check r}
s\cp(Y,Z)\tens s\uCom(\cp(Y,Z),\cq(X,W))
\\
\rTTo^{1\tens s^{-1}[1]} s\cp(Y,Z)\tens\uCom(s\cp(Y,Z),s\cq(X,W))
\rTTo^{\ev^\Com\cdot b^\cq_{00}-(b^\cp_{00}\tens1)\ev^\Com} s\cq(X,W)
\bigr].
\end{multline*}
The complexes \(s\cp(Y,Z)\) and \(s\cq(X,W)\) carry the differential
\(-b^\cp_{00}\). Since \(\ev^{\Com}\) is a chain map, it follows that
 \(\ev^{\Com}b^\cp_{00}-(b^\cp_{00}\tens1)\ev^{\Com}
 =-(1\tens m^{\uCom}_1)\ev^{\Com}\).
Together with the relation \(b^\uCom_1s^{-1}[1]=s^{-1}[1]m^\uCom_1\)
this implies that
\begin{multline*}
t\cdot\check b^\cq_0-(-)^t(\Delta_0\tens1\tens\Delta_0)(1\tens\check
b^\cp_0\tens1)\check t
\\
\quad =(-)^r\bigl[ Ts\ca(X,Y)\tens s\cp(Y,Z)\tens Ts\cc(Z,W) \hfill
\\
\rTTo^{c\tens1} s\cp(Y,Z)\tens Ts\ca(X,Y)\tens Ts\cc(Z,W)
\rTTo^{1\tens\gamma\tens1}
\\
s\cp(Y,Z)\tens Ts\ca^\op(Y,X)\tens Ts\cc(Z,W) \rTTo^{1\tens\check r b^\uCom_1}
s\cp(Y,Z)\tens s\uCom(\cp(Y,Z),\cq(X,W))
\\
\rTTo^{1\tens s^{-1}[1]}
s\cp(Y,Z)\tens\uCom(s\cp(Y,Z),s\cq(X,W))\rTTo^{\ev^\Com}s\cq(X,W)
\bigr].
\end{multline*}
Finally, notice that
\begin{multline*}
-(-)^t(b^\ca\tens1\tens1+1\tens1\tens b^\cc)\check t
\\
=-\bigl[
Ts\ca(X,Y)\tens s\cp(Y,Z)\tens Ts\cc(Z,W)\rTTo^{c\tens1}
s\cp(Y,Z)\tens Ts\ca(X,Y)\tens Ts\cc(Z,W)
\\
\rTTo^{1\tens\gamma\tens1}
 s\cp(Y,Z)\tens Ts\ca^\op(Y,X)\tens Ts\cc(Z,W)
\\
\rTTo^{1\tens(b^{\ca^\op}\tens1+1\tens b^\cc)\check r}
 s\cp(Y,Z)\tens s\uCom(\cp(Y,Z),\cq(X,W))
\\
\rTTo^{1\tens s^{-1}[1]}
s\cp(Y,Z)\tens\uCom(s\cp(Y,Z),s\cq(X,W))\rTTo^{\ev^\Com}s\cq(X,W)
\bigr].
\end{multline*}
Summing up, we conclude that
\((td)^\vee=(-)^r(c\tens1)(1\tens\gamma\tens1)(1\tens R)(1\tens
s^{-1}[1])\ev^\Com\), where
\begin{multline*}
R=\bigl[
Ts\ca^\op(Y,X)\tens Ts\cc(Z,W)\rTTo^{\Delta_0\tens\Delta_0}
\\
\bigoplus_{U\in\Ob\ca,V\in\Ob\cc}Ts\ca^\op(Y,U)\tens
Ts\ca^\op(U,X)\tens Ts\cc(Z,V)\tens Ts\cc(V,W)
\rTTo^{\sum 1\tens c\tens1}
\\
\bigoplus_{U\in\Ob\ca,V\in\Ob\cc} \hspace*{-1em}
 Ts\ca^\op(Y,U)\tens Ts\cc(Z,V)\tens Ts\ca^\op(U,X)\tens Ts\cc(V,W)
\rTTo^{\sum(\check\phi\tens\check r+\check r\tens\check\psi)b^\uCom_2}
\\
\hfill s\uCom(\cp(Y,Z),\cq(X,W))
\bigr]
+\check rb^\uCom_1
-(-)^r(b^{\ca^\op}\tens1+1\tens b^\cc)\check r\quad
\\
\quad=[rb^\uCom-(-)^r(b^{\ca^\op}\tens1+1\tens
b^\cc)r]^\vee=[rB_1]^\vee.\hfill
\end{multline*}
The claim follows.

Let us prove that the constructed chain maps are compatible with the composition.
Let \(\phi,\psi,\chi:\ca^\op,\cc\to\uCom\) be \ainf-functors, \(\cp\),
\(\cq\), \(\ct\) the corresponding \(\ca\text-\cc\)\n-bimodules. Pick
arbitrary \(r\in s\und\Ainfty(\ca^\op,\cc;\uCom)(\phi,\psi)\) and \(q\in
s\und\Ainfty(\ca^\op,\cc;\uCom)(\psi,\chi)\), and denote by
\(t=(rs^{-1})\Phi\) and \(u=(qs^{-1})\Phi\) the corresponding
bicomodule homomorphisms. We must show that
\[
t\cdot u=((rs^{-1}\tens
qs^{-1})m^{\und\Ainfty(\ca^\op,\cc;\uCom)}_2)\Phi=(-)^{q+1}((r\tens
q)B_2s^{-1})\Phi.
\]
Again, it suffices to prove the equation
 \((t\cdot u)^\vee=(-)^{q+1}[((r\tens q)B_2s^{-1})\Phi]^\vee\).
We have:
\begin{multline*}
(t\cdot u)^\vee=t\cdot\check u
=(\Delta_0\tens1\tens\Delta_0)(1\tens\check t\tens1)\check u
\\
\quad=(-)^{r+q}\bigl[ Ts\ca(X,Y)\tens s\cp(Y,Z)\tens Ts\cc(Z,W) \hfill
\\
\rTTo^{\Delta_0\tens1\tens\Delta_0} \hspace*{-2em}
\bigoplus_{U\in\Ob\ca,V\in\Ob\cc} \hspace*{-2em}
 Ts\ca(X,U)\tens Ts\ca(U,Y)\tens s\cp(Y,Z)\tens Ts\cc(Z,V)\tens Ts\cc(V,W)
\\
\rTTo^{\sum 1\tens c\tens1\tens1} \hspace*{-2em}
\bigoplus_{U\in\Ob\ca,V\in\Ob\cc} \hspace*{-2em}
 Ts\ca(X,U)\tens s\cp(Y,Z)\tens Ts\ca(U,Y)\tens Ts\cc(Z,V)\tens Ts\cc(V,W)
\\
\rTTo^{\sum 1\tens1\tens\gamma\tens1\tens1} \hspace*{-2em}
\bigoplus_{U\in\Ob\ca,V\in\Ob\cc} \hspace*{-2em}
 Ts\ca(X,U)\tens s\cp(Y,Z)\tens Ts\ca^\op(Y,U)\tens Ts\cc(Z,V)\tens Ts\cc(V,W)
\\
\rTTo^{\sum 1\tens1\tens\check r\tens1} \hspace*{-2em}
\bigoplus_{U\in\Ob\ca,V\in\Ob\cc} \hspace*{-2em}
 Ts\ca(X,U)\tens s\cp(Y,Z)\tens s\uCom(\cp(Y,Z),\cq(U,V))\tens Ts\cc(V,W)
\\
\rTTo^{\sum 1\tens 1\tens s^{-1}[1]\tens1} \hspace*{-2em}
\bigoplus_{U\in\Ob\ca,V\in\Ob\cc} \hspace*{-2em}
 Ts\ca(X,U)\tens s\cp(Y,Z)\tens\uCom(s\cp(Y,Z),s\cq(U,V))\tens Ts\cc(V,W)
\\
\rTTo^{\sum 1\tens\ev^{\Com}\tens1}
\bigoplus_{U\in\Ob\ca,V\in\Ob\cc}Ts\ca(X,U)\tens
s\cq(U,V)\tens Ts\cc(V,W)
\rTTo^{\sum c\tens 1}
\\
\bigoplus_{U\in\Ob\ca,V\in\Ob\cc}s\cq(U,V)\tens Ts\ca(X,U)\tens Ts\cc(V,W)
\rTTo^{\sum 1\tens\gamma\tens1}
\\
\bigoplus_{U\in\Ob\ca,V\in\Ob\cc}s\cq(U,V)\tens Ts\ca^\op(U,X)\tens Ts\cc(V,W)
\rTTo^{\sum 1\tens\check q}
\\
\bigoplus_{U\in\Ob\ca,V\in\Ob\cc}s\cq(U,V)\tens s\uCom(\cq(U,V),\ct(X,W))
\rTTo^{\sum 1\tens s^{-1}[1]}
\\
\hfill\bigoplus_{U\in\Ob\ca,V\in\Ob\cc}s\cq(U,V)\tens \uCom(s\cq(U,V),s\ct(X,W))
\rTTo^{\sum \ev^{\Com}}s\ct(X,W)
\bigr]\quad
\\
=(-)^r\bigl[
Ts\ca(X,Y)\tens s\cp(Y,Z)\tens Ts\cc(Z,W)\rTTo^{c\tens1}
s\cp(Y,Z)\tens Ts\ca(X,Y)\tens Ts\cc(Z,W)
\\
\rTTo^{1\tens\Delta_0 c(\gamma\tens\gamma)\tens\Delta_0} \hspace*{-2.3em}
\bigoplus_{U\in\Ob\ca,V\in\Ob\cc} \hspace*{-2.2em}
 s\cp(Y,Z)\tens Ts\ca^\op(Y,U)\tens Ts\ca^\op(U,X)\tens Ts\cc(Z,V)\tens Ts\cc(V,W)
\\
\rTTo^{\sum 1\tens1\tens c\tens1} \hspace*{-2em}
\bigoplus_{U\in\Ob\ca,V\in\Ob\cc} \hspace*{-2em}
 s\cp(Y,Z)\tens Ts\ca^\op(Y,U)\tens Ts\cc(Z,V)\tens Ts\ca^\op(U,X)\tens Ts\cc(V,W)
\\
\rTTo^{\sum 1\tens\check r\tens\check q} \hspace*{-2em}
\bigoplus_{U\in\Ob\ca,V\in\Ob\cc} \hspace*{-2em}
 s\cp(Y,Z)\tens s\uCom(\cp(Y,Z),\cq(U,V))\tens s\uCom(\cq(U,V),\ct(X,W))
\\
\rTTo^{\sum 1\tens s^{-1}[1]\tens s^{-1}[1]} \hspace*{-2em}
\bigoplus_{U\in\Ob\ca,V\in\Ob\cc} \hspace*{-2em}
 s\cp(Y,Z)\tens \uCom(s\cp(Y,Z),s\cq(U,V))\tens\uCom(s\cq(U,V),s\ct(X,W))
\\
\hfill \rTTo^{\sum(\ev^{\Com}\tens1)\ev^{\Com}} s\ct(X,W) \bigr]\quad
\\
=(-)^{r+1}\bigl[ Ts\ca(X,Y)\tens s\cp(Y,Z)\tens Ts\cc(Z,W)
\rTTo^{c\tens1} \hfill
\\
s\cp(Y,Z)\tens Ts\ca(X,Y)\tens Ts\cc(Z,W)
\rTTo^{1\tens\gamma\tens1}
s\cp(Y,Z)\tens Ts\ca^\op(Y,X)\tens Ts\cc(Z,W)
\\
\rTTo^{1\tens\Delta_0\tens\Delta_0} \hspace*{-2em}
\bigoplus_{U\in\Ob\ca,V\in\Ob\cc} \hspace*{-2em}
 s\cp(Y,Z)\tens Ts\ca^\op(Y,U)\tens Ts\ca^\op(U,X)\tens Ts\cc(Z,V)\tens Ts\cc(V,W)
\\
\rTTo^{\sum 1\tens1\tens c\tens 1} \hspace*{-2em}
\bigoplus_{U\in\Ob\ca,V\in\Ob\cc} \hspace*{-2em}
 s\cp(Y,Z)\tens Ts\ca^\op(Y,U)\tens Ts\cc(Z,V)\tens Ts\ca^\op(U,X)\tens Ts\cc(V,W)
\\
\rTTo^{\sum 1\tens\check r\tens\check q} \hspace*{-2em}
\bigoplus_{U\in\Ob\ca,V\in\Ob\cc} \hspace*{-2em}
 s\cp(Y,Z)\tens s\uCom(\cp(Y,Z),\cq(U,V))\tens s\uCom(\cq(U,V),\ct(X,W))
\\
\rTTo^{\sum 1\tens b^{\uCom}_2}
 s\cp(Y,Z)\tens s\uCom(\cp(Y,Z),\ct(X,W))
\\
\rTTo^{1\tens s^{-1}[1]} s\cp(Y,Z)\tens \uCom(s\cp(Y,Z),s\ct(X,W))
\rTTo^{\ev^{\Com}}s\ct(X,W) \bigr].
\end{multline*}
It remains to note that
\begin{multline*}
[(r\tens q)B_2]^\vee=\bigl[
Ts\ca^\op(Y,X)\tens Ts\cc(Z,W)\rTTo^{\Delta_0\tens\Delta_0}
\\
\bigoplus_{U\in\Ob\ca,V\in\Ob\cc}Ts\ca^\op(Y,U)\tens
Ts\ca^\op(U,X)\tens Ts\cc(Z,V)\tens Ts\cc(V,W)
\rTTo^{\sum 1\tens c\tens1}
\\
\bigoplus_{U\in\Ob\ca,V\in\Ob\cc}Ts\ca^\op(Y,U)\tens Ts\cc(Z,V)\tens
Ts\ca^\op(U,X)\tens Ts\cc(V,W)
\rTTo^{\sum\check r\tens\check q}
\\
\bigoplus_{U\in\Ob\ca,V\in\Ob\cc}s\uCom(\cp(Y,Z),\cq(U,V))\tens s\uCom(\cq(U,V),\ct(X,W))
\\
\rTTo^{\sum b^\uCom_2} s\uCom(\cp(Y,Z),\ct(X,W)) \bigr],
\end{multline*}
and
 \((-)^{r+1}=(-)^{q+1}(-)^{(r+q+1)+1}
 =(-)^{q+1}(-)^{\deg[(r\tens q)B_2]+1}\).
The claim follows from the definition of \(\Phi\).

Both $\dg$\n-categories \(\ca\text-\cc\bimod\) and
\(\und\Ainfty(\ca^\op,\cc;\uCom)\) are unital. The units are the
identity morphisms in the ordinary categories
\(Z^0(\ca\text-\cc\bimod)\) and
\(Z^0(\und\Ainfty(\ca^\op,\cc;\uCom))\). The $\dg$\n-functor $\Phi$
induces an isomorphism $Z^0\Phi$ of these categories. Hence, $Z^0\Phi$
is unital. In other words, $\Phi$ is unital. The proposition is proven.
\end{proof}
\fi

Let us write explicitly the inverse map
\(\Phi^{-1}:\ca\text-\cc\bimod(\cp,\cq)\to\und\Ainfty(\ca^\op,\cc;\uCom)(\phi,\psi)\).
It takes a bicomodule homomorphism \(t:Ts\ca\tens s\cp\tens Ts\cc\to Ts\ca\tens s\cq\tens
Ts\cc\) to an \ainf-transformation
\(rs^{-1}\in\und\Ainfty(\ca^\op,\cc;\uCom)(\phi,\psi)\) given by its
components
\begin{multline}
\check r=(-)^t\bigl[Ts\ca^\op(Y,X)\tens Ts\cc(Z,W)\rTTo^{\gamma\tens1}
Ts\ca(X,Y)\tens Ts\cc(Z,W)\rTTo^{\coev^\Com}
\\
\uCom(s\cp(Y,Z),s\cp(Y,Z)\tens Ts\ca(X,Y)\tens Ts\cc(Z,W))
\rTTo^{\uCom(1,(c\tens1)\check t)}
\\
\uCom(s\cp(Y,Z),s\cq(X,W))\rTTo^{[-1]s}
s\uCom(\cp(Y,Z),\cq(X,W))
\bigr].
\label{equ-a-inf-trans-via-bicomod-homo}
\end{multline}

\subsection{Regular $A_\infty$-bimodule.}
 \label{sec-regular-A-inf-bimodule}
Let \(\ca\) be an \ainf-category. Extending the notion of regular
\ainf-bimodule given by Tradler \cite[Lemma~5.1(a)]{xxx0108027} from
the case of \ainf-algebras to \ainf-categories, define the
\emph{regular \(\ca\text-\ca\)\n-bimodule} \(\cR=\cR_\ca\) as follows.
Its underlying quiver coincides with \(\ca\). Components of the
codifferential \(b^\cR\) are given by
\[
\check b^\cR=\bigl[ Ts\ca\tens s\ca\tens Ts\ca \rTTo^{\mu_{Ts\ca}}
Ts\ca \rTTo^{\check b^\ca} s\ca \bigr],
\]
where \(\mu_{Ts\ca}\) is the multiplication in the tensor quiver
\(Ts\ca\). Equivalently, \(b^\cR_{kn}=b^\ca_{k+1+n}\), \(k,n\ge0\).
Flatness of \(b^\cR\) in
    \ifx\chooseClass1
the form
\begin{equation}
(b^\ca\tens1\tens1+1\tens1\tens b^\cc)\check b^\cp+(\Delta_0\tens1\tens\Delta_0)(1\tens\check
b^\cp\tens1)\check b^\cp=0.
\label{equ-bP-0}
\end{equation}
    \else
form~\eqref{equ-bP-0}
    \fi
is equivalent to the
\ainf-identity \(b^\ca\cdot\check b^\ca=0\). Indeed, the three summands
of the left hand side of \eqref{equ-bP-0} correspond to three kinds of
subintervals of the interval \([1,k+1+n]\cap\ZZ\). Subintervals of the
first two types miss the point \(k+1\) and those of the third type
contain it.

\begin{definition}
Define an \ainf-functor \(\Hom_\ca:\ca^\op,\ca\to\uCom\) as the
\ainf-functor \(\phi^\cR\) that corresponds to the regular
\(\ca\text-\ca\)\n-bimodule $\cR=\cR_\ca$.
\end{definition}

The \ainf-functor \(\Hom_\ca\) takes a pair of objects \(X,Z\in\Ob\ca\)
to the chain complex \((\ca(X,Z),m_1)\). The components of \(\Hom_\ca\)
are found from equation~\eqref{eq-check-phi-P-coev-bP}:
\begin{multline}
(\Hom_\ca)_{kn} =\bigl[ T^ks\ca^\op(X,Y)\tens T^ns\ca(Z,W)
\rTTo^{\gamma\tens1} T^ks\ca(Y,X)\tens T^ns\ca(Z,W)
\\
\rTTo^{\coev^\Com}
\uCom(s\ca(X,Z),s\ca(X,Z)\tens T^ks\ca(Y,X)\tens T^ns\ca(Z,W))
\\
\rTTo^{\uCom(1,(c\tens1)b^\ca_{k+1+n})} \uCom(s\ca(X,Z),s\ca(Y,W))
\rTTo^{[-1]s} s\uCom(\ca(X,Z),\ca(Y,W)) \bigr].
\label{equ-ainf-Hom-components}
\end{multline}

Closedness of the multicategory $\Ainfty$
\cite[Theorem~12.19]{BesLyuMan-book} implies that there exists a unique
\ainf-functor \(\Yo:\ca\to\und\Ainfty(\ca^\op;\uCom)\) (called the
\emph{Yoneda \ainf-functor}) such that
\[ \Hom_\ca=\bigl[\ca^\op,\ca \rTTo^{1,\Yo}
\ca^\op,\und\Ainfty(\ca^\op;\uCom) \rTTo^{\ev^{\Ainfty}} \uCom\bigr].
\]
Explicit formula \cite[(12.25.4)]{BesLyuMan-book} for evaluation
component \(\ev^{\Ainfty}_{k0}\) shows that the value of $\Yo$ on an
object $Z$ of $\ca$ is given by the restriction \ainf-functor
\[ Z\Yo =H^Z =H_\ca^Z =\Hom_\ca\big|^Z: \ca^\op \to \uCom,
\quad X \mapsto (\ca(X,Z),m_1) =\Hom_\ca(X,Z)
\]
with the components
\begin{multline}
H^Z_k =(\Hom_\ca)_{k0} =(-1)^k\bigl[ T^ks\ca^\op(X,Y)
\rTTo^{\coev^\Com}
\\
\uCom(s\ca(X,Z),s\ca(X,Z)\tens T^ks\ca^\op(X,Y))
\rTTo^{\uCom(1,\omega^0_cb^\ca_{k+1})}
\\
\uCom(s\ca(X,Z),s\ca(Y,Z))
\rTTo^{[-1]s} s\uCom(\ca(X,Z),\ca(Y,Z)) \bigr],
\label{equ-components-of-H-X}
\end{multline}
where
 \(\omega^0=\left(
\begin{smallmatrix}
0 & 1 & \dots & k-1 & k\\
k & k-1 & \dots & 1 & 0
\end{smallmatrix}
 \right)\in\mathfrak S_{k+1}\),
and \(\omega^0_c\) is the corresponding signed permutation.
Restrictions of \ainf-functors in general are defined in
\cite[Section~12.18]{BesLyuMan-book}, in particular, the $k$\n-th
component of \(\Hom_\ca\big|^Z_1\) described by [\textit{loc.~cit.},
(12.18.2)] equals \((1,\Ob\Yo)\ev^{\Ainfty}_{k0}\). Equivalently,
components of the \ainf-functor $H^Z:\ca^\op\to\uCom$ are determined by
the equation
\begin{multline*}
s^{\tens k}H^Z_ks^{-1} =(-1)^{k(k+1)/2+1}\bigl[ T^k\ca^\op(X,Y)
\rTTo^{\coev^\Com}
\\
\uCom(\ca(X,Z),\ca(X,Z)\tens T^k\ca^\op(X,Y))
\rTTo^{\uCom(1,\omega^0_cm^\ca_{k+1})} \uCom(\ca(X,Z),\ca(Y,Z)) \bigr].
\end{multline*}

Notice that $\ev^{\Ainfty}_{km}$ vanishes unless \(m\le1\).
Formula~\cite[(12.25.4)]{BesLyuMan-book} for the component
\(\ev^{\Ainfty}_{k1}\) implies that the component \((\Hom_\ca)_{kn}\)
is determined for \(n\ge1\), \(k\ge0\) by \(\Yo_{nk}\) which is the
composition of $\Yo_n$ with
\[ \pr_k: s\und\Ainfty(\ca^\op;\uCom)(H^Z,H^W)
\to\uCom(T^ks\ca^\op(X,Y),s\uCom(XH^Z,YH^W))
\]
as follows:
\begin{multline*}
(\Hom_\ca)_{kn} =\bigl[ T^ks\ca^\op(X,Y)\tens T^ns\ca(Z,W)
\\
\rTTo^{1\tens\Yo_{nk}}
T^ks\ca^\op(X,Y)\tens\uCom\bigl(T^ks\ca^\op(X,Y),s\uCom(\ca(X,Z),\ca(Y,W))\bigr)
\\
\rTTo^{\ev^\Com} s\uCom(\ca(X,Z),\ca(Y,W)) \bigr].
\end{multline*}
Conversely, the component \(\Yo_n\) is determined by the components
\((\Hom_\ca)_{kn}\) for all $k\ge0$ via the formula
\begin{multline*}
\Yo_{nk} =\bigl[ T^ns\ca(Z,W) \rTTo^{\coev^\Com}
\uCom(T^ks\ca^\op(X,Y),T^ks\ca^\op(X,Y)\tens T^ns\ca(Z,W))
\\
\rTTo^{\uCom(1,(\Hom_\ca)_{kn})}
\uCom\bigl(T^ks\ca^\op(X,Y),s\uCom(\ca(X,Z),\ca(Y,W))\bigr) \bigr].
\end{multline*}
Plugging in expression~\eqref{equ-ainf-Hom-components} we get
\begin{multline*}
\Yo_{nk} =\bigl[ T^ns\ca(Z,W) \rTTo^{\coev^\Com}
\uCom(T^ks\ca^\op(X,Y),T^ks\ca^\op(X,Y)\tens T^ns\ca(Z,W))
\\
\rTTo^{\uCom(1,\coev^\Com)}
\\
\uCom\bigl(T^ks\ca^\op(X,Y),
 \uCom(s\ca(X,Z),s\ca(X,Z)\tens T^ks\ca^\op(X,Y)\tens T^ns\ca(Z,W))\bigr)
\\
\rTTo^{\uCom(1,\uCom(1,(1\tens\gamma\tens1)(c\tens1)b^\ca_{k+1+n}))}
\uCom\bigl(T^ks\ca^\op(X,Y),\uCom(s\ca(X,Z),s\ca(Y,W))\bigr)
\\
\rTTo^{\uCom(1,[-1]s)}
\uCom\bigl(T^ks\ca^\op(X,Y),s\uCom(\ca(X,Z),\ca(Y,W))\bigr) \bigr].
\end{multline*}

Another kind of the Yoneda \ainf-functor
\(Y:\ca\to\und\Ainfty(\ca^\op;\uCom)\) was introduced in
\cite[Appendix~A]{math.CT/0306018}. Actually, it was defined there as
an \ainf-functor from \(\ca^\op\) to \(\und\Ainfty(\ca;\uCom)\). It
turns out that $Y$ which we shall call the \emph{shifted Yoneda
\ainf-functor} differs from $\Yo$ by a shift:
\begin{equation}
Y=\Yo\cdot\und\Ainfty(1;[1]):\ca\to\und\Ainfty(\ca^\op;\uCom).
\label{eq-Y-Yo-[1]}
\end{equation}
Indeed, an object $Z$ of $\ca$ is taken by $Y$ to the \ainf-functor
\(ZY=h^Z:\ca^\op\to\uCom\),
\(X\mapsto(s\ca(X,Z),-b_1)=(\ca(X,Z),m_1)[1]=(XH^Z)[1]\). The
components of \(H^Z\cdot[1]\)
\begin{multline*}
H^Z_ks^{-1}[1]s =(-1)^k\bigl[ T^ks\ca^\op(X,Y) \rTTo^{\coev^\Com}
\\
\uCom(s\ca(X,Z),s\ca(X,Z)\tens T^ks\ca^\op(X,Y))
\rTTo^{\uCom(1,\omega^0_cb^\ca_{k+1})}
\\
\uCom(s\ca(X,Z),s\ca(Y,Z)) \rTTo^s s\uCom(s\ca(X,Z),s\ca(Y,Z)) \bigr]
\end{multline*}
coincide with the components \(h^Z_k\) by
\cite[Appendix~A]{math.CT/0306018}. Therefore, \(h^Z=H^Z\cdot[1]\).
Furthermore, the components $Y_n$ are determined by
\(Y_{nk}=Y_n\cdot\pr_k\), which turn out [\textit{loc. cit.}] to
coincide with
\begin{multline*}
\Yo_{nk}\cdot\uCom(1,s^{-1}[1]s) =\bigl[ T^ns\ca(Z,W)
\rTTo^{\coev^\Com}
\\
\uCom(T^ks\ca^\op(X,Y),T^ks\ca^\op(X,Y)\tens T^ns\ca(Z,W))
\rTTo^{\uCom(1,\coev^\Com)}
\\
\uCom\bigl(T^ks\ca^\op(X,Y),
 \uCom(s\ca(X,Z),s\ca(X,Z)\tens T^ks\ca^\op(X,Y)\tens T^ns\ca(Z,W))\bigr)
\\
\rTTo^{\uCom(1,\uCom(1,(1\tens\gamma\tens1)(c\tens1)b^\ca_{k+1+n}))}
\uCom\bigl(T^ks\ca^\op(X,Y),\uCom(s\ca(X,Z),s\ca(Y,W))\bigr)
\\
\hfill \rTTo^{\uCom(1,s)}
\uCom\bigl(T^ks\ca^\op(X,Y),s\uCom(s\ca(X,Z),s\ca(Y,W))\bigr) \bigr]
\quad
\\
\quad =\bigl[ T^ns\ca(Z,W) \rTTo^{\coev^\Com} \hfill
\\
\uCom(s\ca(X,Z)\tens T^ks\ca^\op(X,Y),
 s\ca(X,Z)\tens T^ks\ca^\op(X,Y)\tens T^ns\ca(Z,W))
\\
\rTTo^{\uCom(1,(1\tens\gamma\tens1)(c\tens1)b^\ca_{k+1+n})}
\uCom(s\ca(X,Z)\tens T^ks\ca^\op(X,Y),s\ca(Y,W))
\\
\rTTo^{\und\varphi^{-1}}_\sim
\uCom\bigl(T^ks\ca^\op(X,Y),\uCom(s\ca(X,Z),s\ca(Y,W))\bigr)
\\
\rTTo^{\uCom(1,s)}
\uCom\bigl(T^ks\ca^\op(X,Y),s\uCom(s\ca(X,Z),s\ca(Y,W))\bigr) \bigr].
\end{multline*}
The natural isomorphism
\(\und\varphi^{\Com}:\uCom(A,\uCom(B,C))\to\uCom(B\tens A,C)\) is found
from the equation
\begin{multline*}
\bigl[ B\tens A\tens \uCom(A,\uCom(B,C)) \rTTo^{1\tens\ev^{\Com}}
B\tens\uCom(B,C) \rTTo^{\ev^{\Com}} C \bigr]
\\
=\bigl[ B\tens A\tens\uCom(A,\uCom(B,C)) \rTTo^{\und\varphi^{\Com}}
B\tens A\tens\uCom(B\tens A,C) \rTTo^{\ev^{\Com}} C \bigr].
\end{multline*}
Its solvability is implied by closedness of \(\Com\). Summing up,
\eqref{eq-Y-Yo-[1]} holds and the two Yoneda \ainf-functors agree.

\subsection{Restriction of scalars.}\label{sec-Restriction-scalars}
Let \(f:\ca\to\cb\), \(g:\cc\to\cd\) be \ainf-functors. Let \(\cp\) be
a \(\cb\text-\cd\)\n-bimodule, \(\phi:\cb^\op,\cd\to\uCom\) the
corresponding \ainf-functor. Define an \(\ca\text-\cc\)\n-bimodule
\(\sS{_f}\cp_g\) as the bimodule corresponding to the composite
\[
\ca^\op,\cc\rTTo^{f^\op,g}\cb^\op,\cd\rTTo^\phi\uCom.
\]
Its underlying \(\gr\)\n-span is given by
\(\sS{_f}\cp_g(X,Y)=\cp(Xf,Yg)\), \(X\in\Ob\ca\), \(Y\in\Ob\cc\).
Components of the codifferential \(b^{\sS{_f}\cp_g}\) are found using
formulas~\eqref{equ-check-b-P-+} and \eqref{equ-check-b-P-0}:
\begin{multline*}
\check b^{\sS{_f}\cp_g}_+=\bigl[
Ts\ca(X,Y)\tens s\cp(Yf,Zg)\tens Ts\cc(Z,W)
\\
\rTTo^{c\tens1}
s\cp(Yf,Zg)\tens Ts\ca(X,Y)\tens Ts\cc(Z,W)
\\
\rTTo^{1\tens\gamma\tens1}
s\cp(Yf,Zg)\tens Ts\ca^\op(Y,X)\tens Ts\cc(Z,W)
\\
\rTTo^{1\tens[(f^\op,g)\phi]^\vee}
s\cp(Yf,Zg)\tens s\uCom(\cp(Yf,Zg),\cp(Xf,Wg))
\\
\hfill\rTTo^{1\tens s^{-1}[1]}
s\cp(Yf,Zg)\tens \uCom(s\cp(Yf,Zg),s\cp(Xf,Wg))\rTTo^{\ev^\Com}s\cp(Xf,Wg)
\bigr]\quad
\\
\quad=\bigl[
Ts\ca(X,Y)\tens s\cp(Yf,Zg)\tens Ts\cc(Z,W)\rTTo^{f\tens1\tens g}\hfill
\\
Ts\cb(Xf,Yf)\tens s\cp(Yf,Zg)\tens Ts\cd(Zg,Wg)\rTTo^{c\tens1}
\\
s\cp(Yf,Zg)\tens Ts\cb(Xf,Yf)\tens Ts\cd(Zg,Wg)\rTTo^{1\tens\gamma\tens1}
\\
s\cp(Yf,Zg)\tens Ts\cb^\op(Yf,Xf)\tens Ts\cd(Zg,Wg)\rTTo^{1\tens\check\phi}
\\
s\cp(Yf,Zg)\tens s\uCom(\cp(Yf,Zg),\cp(Xf,Wg))\rTTo^{1\tens s^{-1}[1]}
\\
\hfill s\cp(Yf,Zg)\tens\uCom(s\cp(Yf,Zg),s\cp(Xf,Wg))\rTTo^{\ev^\Com}s\cp(Xf,Wg)
\bigr],\quad
\\
\quad\check b^{\sS{_f}\cp_g}_0=\bigl[
 Ts\ca(X,Y)\tens s\cp(Yf,Zg)\tens Ts\cc(Z,W)
\rTTo^{\pr_0\tens1\tens\pr_0} s\cp(Yf,Zg) \hfill
\\
\rTTo^{b^\cp_{00}}s\cp(Yf,Zg) \bigr].
\end{multline*}
These equations can be combined into a single formula
\begin{multline}
\check b^{\sS{_f}\cp_g}=\bigl[
Ts\ca(X,Y)\tens s\cp(Yf,Zg)\tens Ts\cc(Z,W)\rTTo^{f\tens1\tens g}
\\
Ts\cb(Xf,Yf)\tens s\cp(Yf,Zg)\tens Ts\cd(Zg,Wg)\rTTo^{\check b^\cp}
s\cp(Xf,Wg)
\bigr].
\label{eq-bfPg-f1g-bP}
\end{multline}

Let \(f:\ca\to\cb\) be an \ainf-functor. Define an
\((\id_{Ts\ca},\id_{Ts\ca})\)\n-bicomodule homomorphism
\(t^f:\cR_\ca=\ca\to\sS{_f}\cb_f=\sS{_f}(\cR_\cb)_f\) of degree \(0\)
by its components
\begin{equation*}
\check t^f=\bigl[
Ts\ca(X,Y)\tens s\ca(Y,Z)\tens Ts\ca(Z,W)\rTTo^{\mu_{Ts\ca}}
Ts\ca(X,W)\rTTo^{\check f} s\cb(Xf,Wf)
\bigr],
\end{equation*}
or in extended form,
\begin{multline}
t^f_{kn}=\bigl[ s\ca(X_k,X_{k-1})\tdt s\ca(X_1,X_0)\tens s\ca(X_0,Z_0)
\tens
\\
\tens s\ca(Z_0,Z_1)\tdt s\ca(Z_{n-1},Z_n)
 \rTTo^{f_{k+1+n}} s\cb(X_kf,Z_nf) \bigr].
\label{equ-A-to-fBf-components}
\end{multline}
We claim that \(t^fd=0\). As usual, it suffices to show that
\((t^fd)^\vee=0\). From the identity
\begin{multline*}
(t^fd)^\vee=t^f\cdot\check b^{\sS{_f}(\cR_\cb)_f}-b^{\cR_\ca}\cdot\check t^f
=(\Delta_0\tens1\tens\Delta_0)(1\tens\check t^f\tens1)\check b^{\sS{_f}(\cR_\cb)_f}
\\
-(b^\ca\tens1\tens1+1\tens1\tens
b^\ca)\check t^f-(\Delta_0\tens1\tens\Delta_0)(1\tens\check
b^{\cR_\ca}\tens1)\check t^f
\end{multline*}
it follows that
\begin{multline*}
(t^fd)^\vee=\bigl[
Ts\ca(X,Y)\tens s\ca(Y,Z)\tens
Ts\ca(Z,W)\rTTo^{\Delta_0\tens1\tens\Delta_0}\hfill
\\
\bigoplus_{U,V\in\Ob\ca} \hspace*{-1em}
 Ts\ca(X,U)\tens Ts\ca(U,Y)\tens s\ca(Y,Z)\tens Ts\ca(Z,V)\tens Ts\ca(V,W)
\rTTo^{\sum 1\tens\mu_{Ts\ca}\tens1}
\\
\bigoplus_{U,V\in\Ob\ca}Ts\ca(X,U)\tens Ts\ca(U,V)\tens Ts\ca(V,W)
\rTTo^{\sum 1\tens\check f\tens1}
\\
\bigoplus_{U,V\in\Ob\ca}Ts\ca(X,U)\tens s\cb(Uf,Vf)\tens Ts\ca(V,W)
\rTTo^{\sum f\tens1\tens f}
\\
\bigoplus_{U,V\in\Ob\ca}Ts\cb(Xf,Uf)\tens s\cb(Uf,Vf)\tens Ts\cb(Vf,Wf)
\rTTo^{\mu_{Ts\cb}}
\\
\hfill Ts\cb(Xf,Wf)\rTTo^{\check b^\cb}s\cb(Xf,Wf)
\bigr]\quad
\\
\quad-\bigl[
Ts\ca(X,Y)\tens s\ca(Y,Z)\tens Ts\ca(Z,W)\rTTo^{b^\ca\tens1\tens1+1\tens1\tens b^\ca}
\hfill
\\
\hfill Ts\ca(X,Y)\tens s\ca(Y,Z)\tens Ts\ca(Z,W)\rTTo^{\mu_{Ts\ca}}
Ts\ca(X,W)\rTTo^{\check f}s\cb(Xf,Wf)
\bigr]\quad
\\
\quad-\bigl[
Ts\ca(X,Y)\tens s\ca(Y,Z)\tens
Ts\ca(Z,W)\rTTo^{\Delta_0\tens1\tens\Delta_0}\hfill
\\
\bigoplus_{U,V\in\Ob\ca}
 Ts\ca(X,U)\tens Ts\ca(U,Y)\tens s\ca(Y,Z)\tens Ts\ca(Z,V)\tens Ts\ca(V,W)
\\
\rTTo^{\sum 1\tens\mu_{Ts\ca}\tens1}
\bigoplus_{U,V\in\Ob\ca}Ts\ca(X,U)\tens Ts\ca(U,V)\tens Ts\ca(V,W)
\rTTo^{\sum 1\tens\check b^\ca\tens1}
\\
\bigoplus_{U,V\in\Ob\ca}Ts\ca(X,U)\tens s\ca(U,V)\tens Ts\ca(V,W)
\rTTo^{\sum\mu_{Ts\ca}} Ts\ca(X,W)
\\
\rTTo^{\check f} s\cb(Xf,Wf) \bigr].
\end{multline*}
Likewise \secref{sec-regular-A-inf-bimodule} we see that the equation
\((t^fd)^\vee=0\) is equivalent to \(f\cdot\check b^\cb=b^\ca\cdot\check f\).

\begin{corollary}\label{cor-Hom-A-fop-f-Hom-B}
Let \(f:\ca\to\cb\) be an \ainf-functor. There is a natural
\ainf-transformation
\(r^f:\Hom_\ca\to(f^\op,f)\cdot\Hom_\cb:\ca^\op,\ca\to\uCom\) depicted
as follows:
\begin{diagram}
\ca^\op,\ca && \rTTo^{\Hom_\ca} && \uCom
\\
&\rdTTo_{f^\op,f} &\dTwoar>{r^f} &\ruTTo_{\Hom_\cb}
\\
&&\cb^\op,\cb &&
\end{diagram}
It is invertible if \(f\) is homotopy full and faithful.
\end{corollary}

\begin{proof}
Define
 \(r^f=(t^f)\Phi^{-1}s\in
 s\und\Ainfty(\ca^\op,\ca;\uCom)(\Hom_\ca,(f^\op,f)\Hom_\cb)\),
where \(t^f:\ca\to\sS{_f}\cb_f\) is the closed bicomodule homomorphism
defined above. Since \(\Phi\) is an invertible chain map, it follows
that \(r^f\) is a natural \ainf-transformation. Suppose \(f\) is homotopy
full and faithful. That is, its first component \(f_1\) is homotopy
invertible. This implies that the \((0,0)\)\n-component
\begin{multline*}
\sS{_{X,Z}}r^f_{00}=\bigl[
\kk\rTTo^\coev\uCom(s\ca(X,Z),s\ca(X,Z))\rTTo^{\uCom(1,f_1)}
\\
\uCom(s\ca(X,Z),s\cb(Xf,Zf))\rTTo^{[-1]s}s\uCom(\ca(X,Z),\cb(Xf,Zf))
\bigr],
\end{multline*}
found from \eqref{equ-a-inf-trans-via-bicomod-homo} and
\eqref{equ-A-to-fBf-components},
is invertible modulo boundaries in \(s\uCom(\ca(X,Z),\cb(Xf,Zf))\),
thus \(r^f\) is invertible by \cite[Lemma~13.9]{BesLyuMan-book}. The
corollary is proven.
\end{proof}

\subsection{Opposite bimodule.}
Let \(\cp\) be an \(\ca\text-\cc\)\n-bimodule,
\(\phi:\ca^\op,\cc\to\uCom\) the corresponding \ainf-functor. Define
an \emph{opposite bimodule \(\cp^\op\)} as
the \(\cc^\op\text-\ca^\op\)\n-bimodule corresponding to the
\ainf-functor
\[
(\id_{\cc},\id_{\ca^\op},\phi)\mu^{\Ainfty}_{\msf{X}:\mb2\to\mb2}=\bigl[Ts\cc\boxt
Ts\ca^\op\rTTo^c Ts\ca^\op\boxt Ts\cc\rTTo^\phi Ts\uCom\bigr].
\]
Its underlying \(\gr\)\n-span is given by \(\cp^\op(Y,X)=\cp(X,Y)\),
\(X\in\Ob\ca\), \(Y\in\Ob\cc\). Components of the differential
\(b^{\cp^\op}\) are found from equations~\eqref{equ-check-b-P-+} and
\eqref{equ-check-b-P-0}:
\begin{multline*}
\check b^{\cp^\op}_+=\bigl[
Ts\cc^\op(W,Z)\tens s\cp^\op(Z,Y)\tens Ts\ca^\op(Y,X)\rTTo^{c\tens 1}
\\
s\cp(Y,Z)\tens Ts\cc^\op(W,Z)\tens Ts\ca^\op(Y,X)\rTTo^{1\tens\gamma\tens1}
s\cp(Y,Z)\tens Ts\cc(Z,W)\tens Ts\ca^\op(Y,X)
\\
\rTTo^{1\tens c}
s\cp(Y,Z)\tens Ts\ca^\op(Y,X)\tens Ts\cc(Z,W)\rTTo^{1\tens\check\phi}
s\cp(Y,Z)\tens s\uCom(\cp(Y,Z),\cp(X,W))
\\
\hfill\rTTo^{1\tens s^{-1}[1]}
s\cp(Y,Z)\tens\uCom(s\cp(Y,Z),s\cp(X,W))\rTTo^{\ev^\Com}s\cp(X,W)=s\cp^\op(W,X)
\bigr]\quad
\\
\quad =\bigl[ Ts\cc^\op(W,Z)\tens s\cp^\op(Z,Y)\tens Ts\ca^\op(Y,X)
\rTTo^{(13)^\sim} \hfill
\\
Ts\ca^\op(Y,X)\tens s\cp(Y,Z)\tens Ts\cc^\op(W,Z)
\rTTo^{\gamma\tens1\tens\gamma}
\\
Ts\ca(X,Y)\tens s\cp(Y,Z)\tens Ts\cc(Z,W) \rTTo^{\check b^\cp_+}
s\cp(X,W)=s\cp^\op(W,X) \bigr],\quad
\\
\quad\check b^{\cp^\op}_0=\bigl[
 Ts\cc^\op(W,X)\tens s\cp^\op(Z,Y)\tens Ts\ca^\op(Y,X)
\rTTo^{\pr_0\tens1\tens\pr_0} s\cp(Y,Z) \hfill
\\
\rTTo^{b^\cp_{00}}s\cp(Y,Z)
\bigr].
\end{multline*}
These equations are particular cases of a single formula
\begin{multline}
\check b^{\cp^\op}=\bigl[
Ts\cc^\op(W,Z)\tens s\cp^\op(Z,Y)\tens Ts\ca^\op(Y,X)\rTTo^{(13)^\sim}
\\
Ts\ca^\op(Y,X)\tens s\cp(Y,Z)\tens Ts\cc^\op(W,Z)
\rTTo^{\gamma\tens1\tens\gamma}
\\
\hfill Ts\ca(X,Y)\tens s\cp(Y,Z)\tens Ts\cc(Z,W) \rTTo^{\check b^\cp}
s\cp(X,W)=s\cp^\op(W,X) \bigr]\quad
\\
\quad=-\bigl[ Ts\cc^\op(W,Z)\tens s\cp^\op(Z,Y)\tens Ts\ca^\op(Y,X)
\rTTo^{(13)^\sim} \hfill
\\
Ts\ca^\op(Y,X)\tens s\cp^\op(Z,Y)\tens Ts\cc^\op(W,Z)
\rTTo^{\gamma\tens\gamma\tens\gamma}
\\
Ts\ca(X,Y)\tens s\cp(Y,Z)\tens Ts\cc(Z,W) \rTTo^{\check b^\cp}
s\cp(X,W)=s\cp^\op(W,X) \bigr].
 \label{equ-op-bimod}
\end{multline}

\begin{proposition}
Let \(\ca\) be an \ainf-category. Then \(\cR^\op_\ca=\cR_{\ca^\op}\) as
\(\ca^\op\text-\ca^\op\)\n-bimodules.
\end{proposition}

\begin{proof}
Clearly, the underlying \(\gr\)\n-spans of the both bimodules coincide.
Computing \(\check b^{\cR^\op_\ca}\) by formula~\eqref{equ-op-bimod} yields
\begin{multline*}
\check b^{\cR^\op_\ca}=-\bigl[
Ts\ca^\op(W,Z)\tens s\ca^\op(Z,Y)\tens Ts\ca^\op(Y,X)\rTTo^{(13)^\sim}
\\
\qquad Ts\ca^\op(Y,X)\tens s\ca^\op(Z,Y)\tens
Ts\ca^\op(W,Z)\rTTo^{\gamma\tens\gamma\tens\gamma}\hfill
\\
\hfill Ts\ca(X,Y)\tens s\ca(Y,Z)\tens
Ts\ca(Z,W)\rTTo^{\mu_{Ts\ca}}Ts\ca(X,W)\rTTo^{\check b^\ca}s\ca(X,W)
\bigr]\quad
\\
\quad=-\bigl[
Ts\ca^\op(W,Z)\tens s\ca^\op(Z,Y)\tens Ts\ca^\op(Y,X)\rTTo^{\mu_{Ts\ca^\op}}
Ts\ca^\op(W,X)\hfill
\\
\rTTo^\gamma Ts\ca(X,W)\rTTo^{\check b^\ca}s\ca(X,W)
\bigr]
\end{multline*}
since \(\gamma:Ts\ca^\op\to Ts\ca\) is a category anti-isomorphism.
Since \(b^{\ca^\op}=\gamma b^\ca\gamma\), it follows that
 \(\check b^{\ca^\op}=-\gamma\check b^\ca:Ts\ca^\op(W,X)\to s\ca(W,X)\),
therefore
\begin{multline*}
\check b^{\cR^\op_\ca}=\bigl[
Ts\ca^\op(W,Z)\tens s\ca^\op(Z,Y)\tens Ts\ca^\op(Y,X)\rTTo^{\mu_{Ts\ca^\op}}
\\
Ts\ca^\op(W,X)\rTTo^{\check b^{\ca^\op}}s\ca^\op(W,X)
\bigr]=\check b^{\cR_{\ca^\op}}.
\end{multline*}
The proposition is proven.
\end{proof}

\begin{corollary}\label{cor-Hom-opposite-cat}
Let \(\ca\) be an \ainf-category. Then
\[
\Hom_{\ca^\op}=\bigl[
Ts\ca\boxt Ts\ca^\op\rTTo^c Ts\ca^\op\boxt Ts\ca\rTTo^{\Hom_\ca}Ts\uCom
\bigr].
\]
\end{corollary}

\begin{proposition}\label{prop-kf-commutes-Hom}
For an arbitrary \ainf-category \(\ca\),
\(\kf\Hom_\ca=\Hom_{\kf\ca}:\kf\ca^\op\boxt\kf\ca\to\und\ck\).
\end{proposition}

\begin{proof}
Let \(X,Y,U,V\in\Ob\ca\). Then
\begin{multline*}
\kf\Hom_\ca=\bigl[ \ca^\op(X,Y)\tens\ca(U,V)
\rTTo^{s(\Hom_\ca)_{10}s^{-1}\tens s(\Hom_\ca)_{01}s^{-1}}
\\
\uCom(\ca(X,U),\ca(Y,U))\tens\uCom(\ca(Y,U),\ca(Y,V)) \rTTo^{m^\uCom_2}
\uCom(\ca(X,U),\ca(Y,V)) \bigr].
\end{multline*}
According to \eqref{equ-components-of-H-X},
\begin{multline}
s(\Hom_\ca)_{10}s^{-1}=-\bigl[ \ca(Y,X) \rTTo^s s\ca(Y,X)
\\
\rTTo^{\coev^\Com} \uCom(s\ca(X,U),s\ca(X,U)\tens s\ca(Y,X))
\\
\hfill\rTTo^{\uCom(1,cb_2)} \uCom(s\ca(X,U),s\ca(Y,U))
\rTTo^{[-1]}~\|_{\uCom(s,1)\cdot\uCom(1,s^{-1})}
\uCom(\ca(X,U),\ca(Y,U)) \bigr]\quad
\\
\quad=-\bigl[ \ca(Y,X) \rTTo^{\coev^\Com}
\uCom(\ca(X,U),\ca(X,U)\tens\ca(Y,X)) \rTTo^{\uCom(1,s\tens s)} \hfill
\\
\hfill\uCom(\ca(X,U),s\ca(X,U)\tens s\ca(Y,X))
\rTTo^{\uCom(1,cb_2s^{-1})} \uCom(\ca(X,U),\ca(Y,U)) \bigr]\quad
\\
\quad=\bigl[ \ca(Y,X) \rTTo^{\coev^\Com}
\uCom(\ca(X,U),\ca(X,U)\tens\ca(Y,X)) \rTTo^{\uCom(1,c)} \hfill
\\
\uCom(\ca(X,U),\ca(Y,X)\tens\ca(X,U)) \rTTo^{\uCom(1,m_2)}
\uCom(\ca(X,U),\ca(Y,U)) \bigr].
 \label{equ-Hom-10-comp}
\end{multline}
Similarly we obtain from equation~\eqref{equ-ainf-Hom-components}
\begin{multline*}
s(\Hom_\ca)_{01}s^{-1}=\bigl[ \ca(U,V) \rTTo^s s\ca(U,V)
\\
\rTTo^{\coev^\Com} \uCom(s\ca(Y,U),s\ca(Y,U)\tens s\ca(U,V))
\\
\hfill \rTTo^{\uCom(1,b_2)} \uCom(s\ca(Y,U),s\ca(Y,V)) \rTTo^{[-1]}
\uCom(\ca(Y,U),\ca(Y,V)) \bigr]\quad
\\
=\bigl[ \ca(U,V) \rTTo^{\coev^\Com}
\uCom(\ca(Y,U),\ca(Y,U)\tens\ca(U,V)) \rTTo^{\uCom(1,m_2)}
\uCom(\ca(Y,U),\ca(Y,V)) \bigr].
\end{multline*}
It follows that
\begin{multline*}
\kf\Hom_\ca=\bigl[ \ca(Y,X)\tens\ca(U,V)
\rTTo^{\coev^\Com\tens\coev^\Com}
\\
\uCom(\ca(X,U),\ca(X,U)\tens\ca(Y,X))\tens\uCom(\ca(Y,U),\ca(Y,U)\tens\ca(U,V))
\rTTo^{\uCom(1,cm_2)\tens\uCom(1,m_2)}
\\
\uCom(\ca(X,U),\ca(Y,U))\tens\uCom(\ca(Y,U),\ca(Y,V)) \rTTo^{m^\uCom_2}
\uCom(\ca(X,U),\ca(Y,V)) \bigr].
\end{multline*}
Equation~(A.1.2) of \cite{math.CT/0306018} allows to write the above
expression as follows:
\begin{multline*}
\kf\Hom_\ca=\bigl[ \ca(Y,X)\tens\ca(U,V) \rTTo^{\coev^\Com}
\uCom(\ca(X,U),\ca(X,U)\tens\ca(Y,X)\tens\ca(U,V))
\\
\hfill\rTTo^{\uCom(1,cm_2\tens1)} \uCom(\ca(X,U),\ca(Y,U)\tens\ca(U,V))
\rTTo^{\uCom(1,m_2)} \uCom(\ca(X,U),\ca(Y,V)) \bigr]\quad
\\
=\bigl[ \ca(Y,X)\tens\ca(U,V) \rTTo^{\coev^\Com}
\uCom(\ca(X,U),\ca(X,U)\tens\ca(Y,X)\tens\ca(U,V))
\rTTo^{\uCom(1,c\tens1)}
\\
\uCom(\ca(X,U),\ca(Y,X)\tens\ca(X,U)\tens\ca(U,V))
\rTTo^{\uCom(1,(m_2\tens1)m_2)} \uCom(\ca(X,U),\ca(Y,V)) \bigr]
\\
=\Hom_{\kf\ca}.
\end{multline*}
The proposition is proven.
\end{proof}

\subsection{Duality $A_\infty$-functor.}\label{sec-duality-A8-functor}
The regular module \(\kk\), viewed as a complex concentrated in degree
0, determines the duality \ainf-functor
\(D=H^\kk=h^\kk\cdot[-1]:\uCom^\op\to\uCom\). It maps a complex \(M\)
to its dual \((\uCom(M,\kk),m_1)=(\uCom(M,\kk),-\uCom(d,1))\). Since
\(\uCom\) is a differential graded category, the components \(D_k\)
vanish if \(k>1\), due to \eqref{equ-components-of-H-X}. The component
\(D_1\) is given by
\begin{multline*}
D_1=-\bigl[ s\uCom^\op(M,N)=s\uCom(N,M) \rTTo^{\coev^\Com}
\\
\uCom(s\uCom(M,\kk),s\uCom(M,\kk)\tens s\uCom(N,M))
\rTTo^{\uCom(1,cb^\uCom_2)} \uCom(s\uCom(M,\kk),s\uCom(N,\kk))
\\
\rTTo^{[-1]} \uCom(\uCom(M,\kk),\uCom(N,\kk)) \rTTo^s
s\uCom(\uCom(M,\kk),\uCom(N,\kk)) \bigr].
\end{multline*}
It follows that
\begin{multline*}
sD_1s^{-1}=\bigl[ \uCom^\op(M,N)=\uCom(N,M) \rTTo^{\coev^\Com}
\\
\uCom(\uCom(M,\kk),\uCom(M,\kk)\tens\uCom(N,M))
\rTTo^{\uCom(1,cm^\uCom_2)} \uCom(\uCom(M,\kk),\uCom(N,\kk)) \bigr],
\end{multline*}
\emph{cf.}~\eqref{equ-Hom-10-comp}. It follows from
\eqref{eq-V(-1)=coev-V(1cmu)} that
 \(\kf D=\und\ck(\_,\kk):\kf\uCom^\op=\und\ck^\op\to\kf\uCom=\und\ck\).

\subsection{Dual $A_\infty$-bimodule.}
Let $\ca$, $\cc$ be \ainf-categories, and let $\cp$ be an
$\ca$-$\cc$-bimodule with a flat \((1,1,1,b^\cA,b^\cC)\)-connection
\(b^\cP:Ts\ca\tens s\cp\tens Ts\cc\to Ts\ca\tens s\cp\tens Ts\cc\), and
let \(\phi^\cp:Ts\ca^\op\boxt Ts\cc\to Ts\uCom\) be the corresponding
\ainf-functor. Define the dual $\cc$-$\ca$-bimodule $\cp^*$ as the
bimodule that corresponds to the  following \ainf-functor:
\[ \phi^{\cp^*}=
((\phi^\cp)^\op,D)\mu^{\Ainfty}_{\mb2\to\mb1}\bull\Ainfty(\msf{X};\uCom):
\cc^\op,\ca\to\uCom,
\]
where \(\msf X:\mb2\to\mb2\), \(1\mapsto2\), \(2\mapsto1\), and the map
\(\Ainfty(\msf{X};\uCom)\) is given by the composite
\begin{multline*}
\Ainfty(\ca,\cc^\op;\uCom) \rTTo^{\id_{\cc^\op}\times\id_{\ca}\times1}
\Ainfty(\cc^\op;\cc^\op)\times\Ainfty(\ca;\ca)\times\Ainfty(\ca,\cc^\op;\uCom)
\\
\rTTo^{\mu^{\Ainfty}_{\msf{X}:\mb2\to\mb2}} \Ainfty(\cc^\op,\ca;\uCom).
\end{multline*}
Equivalently,
\begin{multline}
\phi^{\cp^*}=\bigl(Ts\cc^\op\boxt Ts\ca\rTTo^c Ts\ca\boxt Ts\cc^\op
\rTTo^{\gamma\boxt\gamma}
\\
Ts\ca^\op\boxt Ts\cc \rTTo^{\phi^\cp} Ts\uCom \rTTo^\gamma
Ts\uCom^\op\rTTo^D Ts\uCom\bigr).
 \label{equ-A8-functor-dual-bimod}
\end{multline}
The underlying \(\gr\)\n-span of \(\cp^*\) is given by
\(\Ob_s\cp^*=\Ob_t\cp=\Ob\cc\), \(\Ob_t\cp^*=\Ob_s\cp=\Ob\ca\),
\(\Par\cp^*=\Ob_s\cp^*\times\Ob_t\cp^*\), \(\src=\pr_1\),
\(\tgt=\pr_2\), \(\cp^*(X,Y)=\uCom(\cp(Y,X),\kk)\), \(X\in\Ob\cc\),
\(Y\in\Ob\ca\). Moreover,
\begin{multline*}
\check\phi^{\cp^*}=\phi^{\cp^*}\pr_1=\bigl[
Ts\cc^\op(Y,Z)\tens Ts\ca(X,W) \rTTo^{c(\gamma\tens\gamma)}
\\
Ts\ca^\op(W,X)\tens Ts\cc(Z,Y) \rTTo^{\check\phi^\cp}
s\uCom(\cp(W,Z),\cp(X,Y)) \rTTo^{\coev^\Com}
\\
\uCom(s\uCom(\cp(X,Y),\kk),s\uCom(\cp(X,Y),\kk)\tens
s\uCom(\cp(W,Z),\cp(X,Y)))
\rTTo^{\uCom(1,cb^\uCom_2)}
\\
\uCom(s\uCom(\cp(X,Y),\kk),s\uCom(\cp(W,Z),\kk)) \rTTo^{[-1]}
\uCom(\uCom(\cp(X,Y),\kk),\uCom(\cp(W,Z),\kk))
\\
\rTTo^s s\uCom(\uCom(\cp(X,Y),\kk),\uCom(\cp(W,Z),\kk)) \bigr]
\end{multline*}
(the minus sign present in \(\gamma:s\uCom^\op\to s\uCom\) cancels that
present in \(D_1\)). According to \eqref{equ-check-b-P-+},
\begin{multline*}
\check b^{\cp^*}_+=\bigl[ Ts\cc(X,Y)\tens s\cp^*(Y,Z)\tens Ts\ca(Z,W)
\rTTo^{c\tens1}
\\
s\cp^*(Y,Z)\tens Ts\cc(X,Y)\tens Ts\ca(Z,W) \rTTo^{1\tens\gamma\tens1}
s\cp^*(Y,Z)\tens Ts\cc^\op(Y,X)\tens Ts\ca(Z,W)
\\
\rTTo^{1\tens c} s\cp^*(Y,Z)\tens Ts\ca(Z,W)\tens Ts\cc^\op(Y,X)
\rTTo^{1\tens\gamma\tens\gamma}
\\
s\cp^*(Y,Z)\tens Ts\ca^\op(W,Z)\tens Ts\cc(X,Y)
\rTTo^{1\tens\check\phi^\cp} s\cp^*(Y,Z)\tens s\uCom(\cp(W,X),\cp(Z,Y))
\\
\rTTo^{1\tens\coev^\Com}
s\cp^*(Y,Z)\tens\uCom(s\uCom(\cp(Z,Y),\kk),s\uCom(\cp(Z,Y),\kk)\tens
s\uCom(\cp(W,X),\cp(Z,Y)))
\\
\hfill\rTTo^{1\tens\uCom(1,cb^\uCom_2)}
s\cp^*(Y,Z)\tens\uCom(s\cp^*(Y,Z),s\cp^*(X,W)) \rTTo^{\ev^\Com}
s\cp^*(X,W) \bigr]\quad
\\
=\bigl[ Ts\cc(X,Y)\tens s\cp^*(Y,Z)\tens Ts\ca(Z,W) \rTTo^{(123)}
Ts\ca(Z,W)\tens Ts\cc(X,Y)\tens s\cp^*(Y,Z)
\\
 \rTTo^{\gamma\tens1\tens1} Ts\ca^\op(W,Z)\tens Ts\cc(X,Y)\tens s\cp^*(Y,Z)
\rTTo^{\check\phi^\cp\tens1}
\\
s\uCom(\cp(W,X),\cp(Z,Y))\tens s\uCom(\cp(Z,Y),\kk)
\rTTo^{b^\uCom_2}s\uCom(\cp(W,X),\kk)=s\cp^*(X,W) \bigr],
\end{multline*}
by properties of closed monoidal categories. Similarly, by
\eqref{equ-check-b-P-0}
\begin{multline*}
\check b^{\cp^*}_0=\bigl[ Ts\cc(X,Y)\tens s\cp^*(Y,Z)\tens Ts\ca(Z,W)
\rTTo^{\pr_0\tens1\tens\pr_0} s\cp^*(Y,Z)
\\
\rTTo^{-s^{-1}\uCom(d,1)s} s\cp^*(Y,Z) \bigr],
\end{multline*}
where \(d\) is the differential in the complex
\(\phi^\cp(Z,Y)=\cp(Z,Y)\).

\begin{proposition}
The map \(\check b^{\cp^*}:Ts\cc\tens s\cp^*\tens Ts\ca\to s\cp^*\)
satisfies the following equation:
\begin{multline*}
\bigl[ s\cp(W,X)\tens Ts\cc(X,Y)\tens s\cp^*(Y,Z)\tens Ts\ca(Z,W)
\rTTo^{1\tens\check b^{\cp^*}} s\cp(W,X)\tens s\cp^*(X,W)
\\
\hfill\rTTo^{s^{-1}\tens s^{-1}} \cp(W,X)\tens\uCom(\cp(W,X),\kk)
\rTTo^{\ev^{\Com}} \kk \bigr]\quad
\\
\quad=-\bigl[
s\cp(W,X)\tens Ts\cc(X,Y)\tens s\cp^*(Y,Z)\tens Ts\ca(Z,W)
\rTTo^{(1234)}\hfill
\\
Ts\ca(Z,W)\tens s\cp(W,X)\tens Ts\cc(X,Y)\tens s\cp^*(Y,Z)
\rTTo^{\check b^\cp\tens1} s\cp(Z,Y)\tens s\cp^*(Y,Z)
\\
\rTTo^{s^{-1}\tens s^{-1}} \cp(Z,Y)\tens\uCom(\cp(Z,Y),\kk)
\rTTo^{\ev^{\Com}} \kk \bigr].
\end{multline*}
\end{proposition}

\begin{proof}
It suffices to prove that the pairs of morphisms \(\check b^\cp_+\),
\(\check b^{\cp^*}_+\) and \(\check b^\cp_0\), \(\check b^{\cp^*}_0\)
are related by similar equations. The corresponding equation for
\(\check b^\cp_+\), \(\check b^{\cp^*}_+\) is given below:
\begin{multline}
\bigl[ s\cp(W,X)\tens Ts\cc(X,Y)\tens s\cp^*(Y,Z)\tens Ts\ca(Z,W)
\rTTo^{1\tens\check b^{\cp^*}_+} s\cp(W,X)\tens s\cp^*(X,W)
\\
\hfill\rTTo^{s^{-1}\tens s^{-1}} \cp(W,X)\tens\uCom(\cp(W,X),\kk)
\rTTo^{\ev^{\Com}} \kk \bigr]\quad
\\
\quad=-\bigl[
s\cp(W,X)\tens Ts\cc(X,Y)\tens s\cp^*(Y,Z)\tens Ts\ca(Z,W)
\rTTo^{(1234)}\hfill
\\
Ts\ca(Z,W)\tens s\cp(W,X)\tens Ts\cc(X,Y)\tens s\cp^*(Y,Z)
\rTTo^{\check b^\cp_+\tens1} s\cp(Z,Y)\tens s\cp^*(Y,Z)
\\
\rTTo^{s^{-1}\tens s^{-1}} \cp(Z,Y)\tens\uCom(\cp(Z,Y),\kk)
\rTTo^{\ev^{\Com}} \kk \bigr].
 \label{equ-check-b-P-+-check-b-P-*-+}
\end{multline}
Its left hand side equals
\begin{multline*}
\bigl[ s\cp(W,X)\tens Ts\cc(X,Y)\tens s\cp^*(Y,Z)\tens Ts\ca(Z,W)
\rTTo^{1\tens(123)}
\\
s\cp(W,X)\tens Ts\ca(Z,W)\tens Ts\cc(X,Y)\tens s\cp^*(Y,Z)
\rTTo^{1\tens\gamma\tens1\tens1}
\\
s\cp(W,X)\tens Ts\ca^\op(W,Z)\tens Ts\cc(X,Y)\tens s\cp^*(Y,Z)
\rTTo^{1\tens\check\phi^\cp\tens1}
\\
s\cp(W,X)\tens s\uCom(\cp(W,X),\cp(Z,Y))\tens s\uCom(\cp(Z,Y),\kk)
\rTTo^{1\tens b^\uCom_2}
\\
s\cp(W,X)\tens s\uCom(\cp(W,X),\kk) \rTTo^{s^{-1}\tens s^{-1}}
\cp(W,X)\tens\uCom(\cp(W,X),\kk) \rTTo^{\ev^\Com} \kk \bigr].
\end{multline*}
Note that \((1234)(c\tens1\tens1)=1\tens(123)\). Using
\eqref{equ-check-b-P-+}, the right hand side can be written as follows:
\begin{multline*}
-\bigl[ s\cp(W,X)\tens Ts\cc(X,Y)\tens s\cp^*(Y,Z)\tens Ts\ca(Z,W)
\rTTo^{1\tens(123)}
\\
s\cp(W,X)\tens Ts\ca(Z,W)\tens Ts\cc(X,Y)\tens s\cp^*(Y,Z)
\rTTo^{1\tens\gamma\tens1\tens1}
\\
s\cp(W,X)\tens Ts\ca^\op(W,Z)\tens Ts\cc(X,Y)\tens s\cp^*(Y,Z)
\rTTo^{1\tens\check\phi^\cp\tens1}
\\
s\cp(W,X)\tens s\uCom(\cp(W,X),\cp(Z,Y))\tens s\uCom(\cp(Z,Y),\kk)
\rTTo^{1\tens s^{-1}\tens 1}
\\
s\cp(W,X)\tens\uCom(\cp(W,X),\cp(Z,Y))\tens s\uCom(\cp(Z,Y),\kk)
\rTTo^{1\tens[1]\tens1}
\\
s\cp(W,X)\tens\uCom(s\cp(W,X),s\cp(Z,Y))\tens s\uCom(\cp(Z,Y),\kk)
\rTTo^{\ev^\Com\tens1}
\\
s\cp(Z,Y)\tens s\uCom(\cp(Z,Y),\kk) \rTTo^{s^{-1}\tens s^{-1}}
\cp(Z,Y)\tens\uCom(\cp(Z,Y),\kk) \rTTo^{\ev^\Com} \kk \bigr].
\end{multline*}
Equation~\eqref{equ-check-b-P-+-check-b-P-*-+} follows from the
following equation, which we are going to prove:
\begin{multline*}
\bigl[
s\cp(W,X)\tens s\uCom(\cp(W,X),\cp(Z,Y))\tens s\uCom(\cp(Z,Y),\kk)
\rTTo^{1\tens b^\uCom_2}
\\
\hfill s\cp(W,X)\tens s\uCom(\cp(W,X),\kk) \rTTo^{s^{-1}\tens s^{-1}}
\cp(W,X)\tens\uCom(\cp(W,X),\kk) \rTTo^{\ev^\Com} \kk \bigr]\quad
\\
\quad=-\bigl[
s\cp(W,X)\tens s\uCom(\cp(W,X),\cp(Z,Y))\tens s\uCom(\cp(Z,Y),\kk)
\rTTo^{1\tens s^{-1}\tens 1}\hfill
\\
s\cp(W,X)\tens\uCom(\cp(W,X),\cp(Z,Y))\tens s\uCom(\cp(Z,Y),\kk)
\rTTo^{1\tens[1]\tens1}
\\
s\cp(W,X)\tens\uCom(s\cp(W,X),s\cp(Z,Y))\tens s\uCom(\cp(Z,Y),\kk)
\rTTo^{\ev^\Com\tens1}
\\
s\cp(Z,Y)\tens s\uCom(\cp(Z,Y),\kk) \rTTo^{s^{-1}\tens s^{-1}}
\cp(Z,Y)\tens\uCom(\cp(Z,Y),\kk) \rTTo^{\ev^\Com} \kk \bigr].
\end{multline*}
Composing both sides with the morphism \(-s\tens s\tens s\) we obtain
an equivalent equation:
\begin{multline*}
\bigl[ \cp(W,X)\tens\uCom(\cp(W,X),\cp(Z,Y))\tens\uCom(\cp(Z,Y),\kk)
\rTTo^{1\tens m^\uCom_2}
\\
\hfill\cp(W,X)\tens\uCom(\cp(W,X),\kk) \rTTo^{\ev^\Com} \kk \bigr]\quad
\\
\quad=-\bigl[
\cp(W,X)\tens\uCom(\cp(W,X),\cp(Z,Y))\tens\uCom(\cp(Z,Y),\kk)
\rTTo^{s\tens1\tens s}\hfill
\\
s\cp(W,X)\tens\uCom(\cp(W,X),\cp(Z,Y))\tens s\uCom(\cp(Z,Y),\kk)
\rTTo^{1\tens[1]\tens1}
\\
s\cp(W,X)\tens\uCom(s\cp(W,X),s\cp(Z,Y))\tens s\uCom(\cp(Z,Y),\kk)
\rTTo^{\ev^\Com\tens1}
\\
s\cp(Z,Y)\tens s\uCom(\cp(Z,Y),\kk) \rTTo^{s^{-1}\tens s^{-1}}
\cp(Z,Y)\tens\uCom(\cp(Z,Y),\kk) \rTTo^{\ev^\Com} \kk \bigr].
\end{multline*}
It follows from the definition of $\dg$\n-functor \([1]\) that
\[
(s\tens1)(1\tens[1])\ev^\Com s^{-1}
=(s\tens1)(s^{-1}\tens1)\ev^{\Com}ss^{-1} =\ev^{\Com},
\]
therefore the right hand side of the above equation is equal to
\((\ev^\Com\tens1)\ev^\Com\), and it equals the left hand side by the
definition of \(m^\uCom_2\).

The corresponding equation for \(\check b^\cp_0\) and
\(\check b^{\cp^*}_0\) reads as follows:
\begin{multline*}
\bigl[ s\cp(W,X)\tens T^0s\cc(X,Y)\tens s\cp^*(Y,Z)\tens T^0s\ca(Z,W)
\rTTo^{1\tens\check b^{\cp^*}_0}
\\
\hfill s\cp(W,X)\tens s\cp^*(X,W) \rTTo^{s^{-1}\tens s^{-1}}
\cp(W,X)\tens\uCom(\cp(W,X),\kk) \rTTo^{\ev^{\Com}} \kk \bigr]\quad
\\
\quad=-\bigl[
s\cp(W,X)\tens T^0s\cc(X,Y)\tens s\cp^*(Y,Z)\tens T^0s\ca(Z,W)
\rTTo^{(1234)}\hfill
\\
T^0s\ca(Z,W)\tens s\cp(W,X)\tens T^0s\cc(X,Y)\tens s\cp^*(Y,Z)
\rTTo^{\check b^\cp_0\tens1}
\\
s\cp(Z,Y)\tens s\cp^*(Y,Z) \rTTo^{s^{-1}\tens s^{-1}}
\cp(Z,Y)\tens\uCom(\cp(Z,Y),\kk) \rTTo^{\ev^{\Com}} \kk \bigr].
\end{multline*}
It is non-trivial only if \(X=Y\) and \(Z=W\). In this case, up to
obvious isomorphism the left hand side equals
\begin{multline*}
\bigl[ s\cp(Z,X)\tens s\cp^*(X,Z) \rTTo^{s^{-1}\tens s^{-1}}
\cp(Z,X)\tens\uCom(\cp(Z,X),\kk)
\\
\hfill\rTTo^{1\tens\uCom(d,1)} \cp(Z,X)\tens\uCom(\cp(Z,X),\kk)
\rTTo^{\ev^\Com} \kk\bigr]\quad
\\
\quad=\bigl[ s\cp(Z,X)\tens s\cp^*(X,Z) \rTTo^{s^{-1}\tens s^{-1}}
\cp(Z,X)\tens\uCom(\cp(Z,X),\kk) \hfill
\\
\hfill \rTTo^{d\tens1} \cp(Z,X)\tens\uCom(\cp(Z,X),\kk)
\rTTo^{\ev^\Com} \kk \bigr] \quad
\\
\quad=-\bigl[ s\cp(Z,X)\tens s\cp^*(X,Z) \rTTo^{s^{-1}ds\tens 1}
s\cp(Z,X)\tens s\uCom(\cp(Z,X),\kk) \hfill
\\
\rTTo^{s^{-1}\tens s^{-1}} \cp(Z,X)\tens\uCom(\cp(Z,X),\kk)
\rTTo^{\ev^\Com} \kk \bigr],
\end{multline*}
which is the right hand side (up to obvious isomorphism).
\end{proof}

\begin{proposition}
Let \(\ca\) be an \ainf-category. Denote by $\cR$ the regular
$\ca$\n-bimodule. Then
\[ \phi^{\cR^*}
=(\Hom^\op_{\ca^\op},D)\mu^{\Ainfty}_{\mb2\to\mb1}
=\Hom^\op_{\ca^\op}\cdot D: \ca^\op,\ca\to\uCom.
\]
\end{proposition}

\begin{proof}
Formula~\eqref{equ-A8-functor-dual-bimod} implies that
\begin{multline*}
\phi^{\cR^*}=\bigl[ Ts\ca^\op\boxt Ts\ca \rTTo^c Ts\ca\boxt Ts\ca^\op
\rTTo^{\gamma\boxt\gamma} Ts\ca^\op\boxt Ts\ca
\\
\hfill\rTTo^{\Hom_\ca} Ts\uCom \rTTo^\gamma Ts\uCom^\op \rTTo^D Ts\uCom
\bigr]\quad
\\
\quad=\bigl[
Ts\ca^\op\boxt Ts\ca \rTTo^{\gamma\boxt\gamma} Ts\ca\boxt Ts\ca^\op
\rTTo^c Ts\ca^\op\boxt Ts\ca\hfill
\\
\rTTo^{\Hom_\ca} Ts\uCom \rTTo^\gamma Ts\uCom^\op \rTTo^D Ts\uCom
\bigr].
\end{multline*}
By \corref{cor-Hom-opposite-cat},
 \(\Hom_{\ca^\op}=[Ts\ca\boxt Ts\ca^\op\rTTo^c Ts\ca^\op\boxt
 Ts\ca\rTTo^{\Hom_\ca}Ts\uCom]\),
therefore
\begin{multline*}
\phi^{\cR^*}=\bigl[ Ts\ca^\op\boxt Ts\ca \rTTo^{\gamma\boxt\gamma}
Ts\ca\boxt Ts\ca^\op \rTTo^{\Hom_{\ca^\op}} Ts\uCom \rTTo^\gamma
Ts\uCom^\op \rTTo^D Ts\uCom \bigr]
\\
=\bigl[ Ts\ca^\op\boxt Ts\ca \rTTo^{\Hom^\op_{\ca^\op}}
Ts\uCom^\op\rTTo^D Ts\uCom \bigr].
\end{multline*}
The proposition is proven.
\end{proof}

\subsection{Unital $A_\infty$-bimodules.}
\begin{definition}
An $\ca$-$\cc$-bimodule $\cp$ corresponding to an \ainf-functor
\(\phi:\ca^\op,\cc\to\uCom\) is called \emph{unital} if the
\ainf-functor \(\phi\) is unital.
\end{definition}

\begin{proposition}
An $\ca$-$\cc$-bimodule $\cp$ is unital if and only if for all
\(X\in\Ob\ca\), \(Y\in\Ob\cc\) the compositions
\begin{gather*}
\bigl[ s\cp(X,Y) =s\cp(X,Y)\tens\kk \rTTo^{1\tens\sS{_Y}\uni^\cc_0}
s\cp(X,Y)\tens s\cc(Y,Y) \rTTo^{b^\cp_{01}} s\cp(X,Y) \bigr],
\\
-\bigl[ s\cp(X,Y) =\kk\tens s\cp(X,Y) \rTTo^{\sS{_X}\uni^\ca_0\tens1}
s\ca(X,X)\tens s\cp(X,Y) \rTTo^{b^\cp_{10}} s\cp(X,Y) \bigr]
\end{gather*}
are homotopic to the identity map.
\end{proposition}

\begin{proof}
The second statement expands to the property that
\begin{gather*}
\bigl[ s\cp(X,Y) \rTTo^{s^{-1}\tens\sS{_Y}\uni^\cc_0\phi_{01}s^{-1}}
\cp(X,Y)\tens\uCom(\cp(X,Y),\cp(X,Y)) \rTTo^{\ev^\Com s} s\cp(X,Y)
\bigr],
\\
\bigl[ s\cp(X,Y) \rTTo^{s^{-1}\tens\sS{_X}\uni^\ca_0\phi_{10}s^{-1}}
\cp(X,Y)\tens\uCom(\cp(X,Y),\cp(X,Y)) \rTTo^{\ev^\Com s} s\cp(X,Y)
\bigr]
\end{gather*}
are homotopic to identity. That is,
\[
\sS{_Y}\uni^\cc_0\bigl(\phi_{01}|^X_\cc\bigr)-1_{s\cp(X,Y)}s\in\im b_1,
\qquad \sS{_X}\uni^\ca_0\bigl(\phi_{10}|^Y_{\ca^\op}\bigr)
-1_{s\cp(X,Y)}s \in \im b_1
\]
for all \(X\in\Ob\ca\), \(Y\in\Ob\cc\). By
\cite[Proposition~13.6]{BesLyuMan-book} the \ainf-functor $\phi$ is
unital.
\end{proof}

\begin{remark}\label{rem-unital-Hom-rep-Yo}
Suppose \(\ca\) is a unital \ainf-category. By the above criterion, the
regular $\ca$\n-bimodule $\cR_\ca$ is unital, and therefore the
\ainf-functor \(\Hom_\ca:\ca^\op,\ca\to\uCom\) is unital. In
particular, for each object \(Z\) of \(\ca\), the representable
\ainf-functor
\[
H^Z=\Hom_\ca|^Z_1:\ca^\op\to\uCom
\]
is unital, by \cite[Proposition~13.6]{BesLyuMan-book}. Thus, the Yoneda
\ainf-functor \(\Yo:\ca\to\und\Ainfty(\ca^\op;\uCom)\) takes values in
the full \ainf-subcategory \(\und\Ainftyu(\ca^\op;\uCom)\) of
\(\und\Ainfty(\ca^\op;\uCom)\). Furthermore, by the closedness of the
multicategory \(\Ainftyu\), the \ainf-functor \(\Yo\) is unital.
\end{remark}

Let \(\ca\), \(\cb\) be unital \ainf-categories. A unital \ainf-functor
\(f:\ca\to\cb\) is called \emph{homotopy fully faithful} if the
corresponding \(\ck\)\n-functor \(\kf f:\kf\ca\to\kf\cb\) is fully
faithful. That is, \(f\) is homotopy fully faithful if and only if its
first component is homotopy invertible. Equivalently, \(f\) is homotopy
fully faithful if and only if it admits a factorization
\begin{equation}
\ca\rTTo^g \ci\rMono^{e} \cb,
\label{equ-f-g-e}
\end{equation}
where \(\ci\) is a full \ainf-subcategory of \(\cb\), \(e:\ci\to\cb\)
is the embedding, and \(g:\ca\to\ci\) is an \ainf-equivalence.

\begin{lemma}\label{lem-homotopy-fully-faithful}
Suppose \(f:\ca\to\cb\) is a homotopy fully faithful \ainf-functor.
Then for an arbitrary \ainf-category \(\cc\) the \ainf-functor
\(\und\Ainfty(1;f):\und\Ainfty(\cc;\ca)\to\und\Ainfty(\cc;\cb)\) is
homotopy fully faithful.
\end{lemma}

\begin{proof}
Factorize \(f\) as in \eqref{equ-f-g-e}. Then the \ainf-functor
\(\und\Ainfty(1;f)\) factorizes as
\[
\und\Ainfty(\cc;\ca) \rTTo^{\und\Ainfty(1;g)} \und\Ainfty(\cc;\ci)
\rTTo^{\und\Ainfty(1;e)} \und\Ainfty(\cc;\cb).
\]
The \ainf-functor \(\und\Ainfty(1;g)\) is an \ainf-equivalence since
\(\und\Ainfty(\cc;-)\) is a \ainfu-2\n-functor, so it suffices to show
that \(\und\Ainfty(1;e)\) is a full embedding. Since \(e\) is a strict
\ainf-functor, so is \(\und\Ainfty(1;e)\). Its first component is given
by
\begin{multline*}
s\und\Ainfty(\cc;\ci)(\phi,\psi) =\prod^{X,Y\in\Ob\cc}_{n\ge0}
\uCom(T^ns\cc(X,Y),s\ci(X\phi,Y\psi)) \rTTo^{\uCom(1,e_1)}
\\
s\und\Ainfty(\cc;\cb)(\phi e,\psi e)
=\prod^{X,Y\in\Ob\cc}_{n\ge0} \uCom(T^ns\cc(X,Y),s\cb(X\phi,Y\psi)),
\end{multline*}
that is, \(r=(r_n)\mapsto re=(r_ne_1)\). Since
\(s\ci(X\phi,Y\psi)=s\cb(X\phi,Y\psi)\) and
 \(e_1:s\ci(X\phi,Y\psi)\to s\cb(X\phi,Y\psi)\) is the identity
morphism, the above map is the identity morphism, and the proof is
complete.
\end{proof}

\begin{example}
Let \(g:\cc\to\ca\) be an \ainf-functor. Then an $\ca$-$\cc$-bimodule
\(\ca^g\) is associated to it via the \ainf-functor
\begin{align}
\ca^g &= \bigl[ \ca^\op,\cc \rTTo^{1,g} \ca^\op,\ca \rTTo^{1,\Yo}
\ca^\op,\und\Ainfty(\ca^\op,\uCom) \rTTo^{\ev^{\Ainfty}} \uCom \bigr]
\notag
\\
&= \bigl[ \ca^\op,\cc \rTTo^{1,g} \ca^\op,\ca \rTTo^{\Hom_\ca\;}
\uCom \bigr].
\label{eq-Ag-AopC-C}
\end{align}
\end{example}

As we already noticed, $\ca$-$\cc$-bimodules are objects of the
differential graded category
\(\und\Ainfty(\ca^\op,\cc;\uCom)\simeq\und\Ainfty(\cc;\und\Ainfty(\ca^\op;\uCom))\).
Thus it makes sense to speak about their isomorphisms in the homotopy
category \(H^0(\und\Ainfty(\ca^\op,\cc;\uCom))\).

\begin{proposition}\label{pro-g-Ag-unital}
There is an \ainf-functor
\(\und\Ainfty(\cc;\ca)\to\und\Ainfty(\ca^\op,\cc;\uCom)\),
\(g\mapsto\ca^g\), homotopy fully faithful if $\ca$ is a unital
\ainf-category. In that case, \ainf-functors \(g,h:\cc\to\ca\) are
isomorphic if and only if the bimodules \(\ca^g\) and \(\ca^h\) are
isomorphic. If both $\cc$ and $\ca$ are unital, then the above
\ainf-functor restricts to
\(\und\Ainftyu(\cc;\ca)\to\und\Ainftyu(\ca^\op,\cc;\uCom)\). Moreover,
$g$ is unital if and only if $\ca^g$ is unital.
\end{proposition}

\begin{proof}
The functor in question is the composite
\[
\und\Ainfty(\cc;\ca) \rTTo^{\und\Ainfty(1;\Yo)}
\und\Ainfty(\cc;\und\Ainfty(\ca^\op;\uCom))
\rTTo^{\und\varphi^{\Ainfty}}_\sim \und\Ainfty(\ca^\op,\cc;\uCom),
\]
where \(\Yo\) is the Yoneda \ainf-functor. When $\ca$ is unital,
\(\Yo:\ca\to\und\Ainftyu(\ca^\op;\uCom)\) is homotopy fully faithful by
\corref{cor-Yoneda-equivalence}, see also
\cite[Theorem~A.11]{math.CT/0306018}. Thus, the first claim follows
from \lemref{lem-homotopy-fully-faithful}.

Assume that \ainf-categories $\cc$ and $\ca$ are unital, and
\ainf-functor~\eqref{eq-Ag-AopC-C} is unital. Let us prove that
\(g:\cc\to\ca\) is unital. Denote
 \(f=g\cdot\Yo:\cc\to\und\Ainfty(\ca^\op;\uCom)\). The \ainf-functor
\[ f' = \bigl[ \ca^\op,\cc \rTTo^{1,f}
\ca^\op,\und\Ainftyu(\ca^\op,\uCom) \rTTo^\ev \uCom \bigr]
\]
is unital by assumption. The bijection
\[ \varphi^{\Ainfty}: \Ainfty(\cc;\und\Ainfty(\ca^\op;\uCom))
\to \Ainfty(\ca^\op,\cc;\uCom)
\]
shows that given $f'$ can be obtained from a unique $f$. The bijection
\[ \varphi^{\Ainftyu}: \Ainftyu(\cc;\und\Ainftyu(\ca^\op;\uCom))
\to \Ainftyu(\ca^\op,\cc;\uCom)
\]
shows that such \ainf-functor $f$ is unital.

Thus, the composition of \(g:\cc\to\ca\) with the unital homotopy fully
faithful \ainf-functor \(\Yo:\ca\to\und\Ainftyu(\ca^\op;\uCom)\) is
unital. Denote by \(\Rep(\ca^\op,\uCom)\) the essential image of $\Yo$,
the full differential graded subcategory of
\(\und\Ainftyu(\ca^\op;\uCom)\) whose objects are representable
\ainf-functors \((X)\Yo=H^X\). The composition of $g$ with the
\ainf-equivalence \(\Yo:\ca\to\Rep(\ca^\op,\uCom)\) is unital. Denote
by \(\Psi:\Rep(\ca^\op,\uCom)\to\ca\) a quasi-inverse to $\Yo$. We find
that $g$ is isomorphic to a unital \ainf-functor
\(g\cdot\Yo\cdot\Psi:\cc\to\ca\). Thus, $g$ is unital itself by
\cite[(8.2.4)]{Lyu-AinfCat}.
\end{proof}

\begin{proposition}\label{prop-bimodul-isom-A-g}
Let \(\ca,\cc\) be \ainf-categories, and suppose \(\ca\) is unital. Let
\(\cp\) be an $\ca$-$\cc$-bimodule, \(\phi^\cp:\ca^\op,\cc\to\uCom\)
the corresponding \ainf-functor. The $\ca$-$\cc$-bimodule \(\cp\) is
isomorphic to \(\ca^g\) for some \ainf-functor \(g:\cc\to\ca\) if and
only if for each object \(Y\in\Ob\cc\) the \ainf-functor
\(\phi^\cp|^Y_1:\ca^\op\to\uCom\) is representable.
\end{proposition}

\begin{proof}
The ``only if'' part is obvious. For the proof of ``if'', consider the
\ainf-functor
\(f=(\varphi^{\Ainfty})^{-1}(\phi^\cp):\cc\to\und\Ainfty(\ca^\op;\uCom)\).
It acts on objects by \(Y\mapsto\phi^\cp|^Y_1\), \(Y\in\Ob\cc\),
therefore it takes values in the \ainf-subcategory
\(\Rep(\ca^\op,\uCom)\) of representable \ainf-functors. Denote by
\(\Psi:\Rep(\ca^\op,\uCom)\to\ca\) a quasi-inverse to \(\Yo\). Let
\(g\) denote the \ainf-functor \(f\cdot\Psi:\cc\to\ca\). Then the
composite
\(g\cdot\Yo=f\cdot\Psi\cdot\Yo:\cc\to\und\Ainfty(\ca^\op;\uCom)\) is
isomorphic to \(f\), therefore the \ainf-functor
\[
\varphi^{\Ainfty}(g\cdot\Yo)=\bigl[\ca^\op,\cc \rTTo^{1,g\cdot \Yo}
\ca^\op,\und\Ainfty(\ca^\op;\uCom) \rTTo^{\ev^{\Ainfty}} \uCom\bigr],
\]
corresponding to the bimodule \(\ca^g\), is isomorphic to
\(\varphi^{\Ainfty}(f)=\phi^\cp\). Thus, \(\ca^g\) is isomorphic to
\(\cp\).
\end{proof}

\begin{lemma}\label{lem-dual-bimodule-unital}
If $\ca$-$\cc$-bimodule $\cp$ is unital, then the dual
$\cc$-$\ca$-bimodule $\cp^*$ is unital as well.
\end{lemma}

\begin{proof}
The \ainf-functor \(\phi^{\cp^*}\) is the composite of two
\ainf-functors, \((\phi^\cp)^\op:\ca,\cc^\op\to\uCom^\op\) and
\(D=H^\kk:\uCom^\op\to\uCom\). The latter is unital by
\remref{rem-unital-Hom-rep-Yo}. The former is unital if and only if
\(\phi^\cp\) is unital.
\end{proof}

\section{Serre $A_\infty$-functors}\label{sec-Serre-A8-functors}
Here we present Serre \ainf-functors as an application of
\ainf-bimodules.

\begin{definition}[cf. Soibelman~\cite{MR2095670},
 Kontsevich and Soibelman, sequel to \cite{math.RA/0606241}]
A \emph{right Serre \ainf-functor} \(S:\ca\to\ca\) in a unital
\ainf-category $\ca$ is an \ainf-functor for which the
$\ca$\n-bimodules
 \(\ca^S=\bigl[\ca^\op,\ca\rTTo^{1,S} \ca^\op,\ca\rTTo^{\Hom_\ca\;}
 \uCom\bigr]\)
and \(\ca^*\) are isomorphic. If, moreover, $S$ is an
\ainf-equivalence, it is called a \emph{Serre \ainf-functor}.
\end{definition}

By \lemref{lem-dual-bimodule-unital} and \propref{pro-g-Ag-unital}, if
a right Serre \ainf-functor exists, then it is unital and unique up to
an isomorphism.

\begin{proposition}\label{prop-Serre-ainf-implies-Serre-K}
If \(S:\ca\to\ca\) is a (right) Serre \ainf-functor, then
 \(\kf S:\kf\ca\to\kf\ca\) is a (right) Serre \(\ck\)\n-functor.
\end{proposition}

\begin{proof}
Let \(p:\ca^S\to\ca^*\) be an isomorphism. More precisely, \(p\) is an
isomorphism
\[
(1,S,\Hom_\ca)\mu^{\Ainftyu}_{\id:\mb2\to\mb2}
\to(\Hom^\op_{\ca^\op},D)\mu^{\Ainftyu}_{\mb2\to\mb1}:\ca^\op,\ca\to\uCom.
\]
We visualize this by the following diagram:
\begin{diagram}
\ca^\op,\ca & \rTTo^{1,S} & \ca^\op,\ca\\
\dTTo<{\Hom^\op_{\ca^\op}} & \ldTwoar^p & \dTTo>{\Hom_\ca}\\
\uCom^\op & \rTTo^D & \uCom
\end{diagram}
Applying the \(\kCat\)\n-multifunctor \(\kf\), and using
\lemref{lem-kf-commutes-with-op-on-functors},
\propref{prop-kf-commutes-Hom}, and results of
\secref{sec-duality-A8-functor}, we get a similar diagram in \(\KCat\):
\begin{diagram}
\kf\ca^\op\boxt\kf\ca & \rTTo^{1\boxt\kf S} & \kf\ca^\op\boxt\kf\ca\\
\dTTo<{\Hom^\op_{\kf\ca^\op}} &\ldTwoar^{\kf p} &\dTTo>{\Hom_{\kf\ca}}\\
\und\ck^\op & \rTTo_{\und\ck(\_,\kk)} & \und\ck
\end{diagram}
Since \(\kf p\) is an isomorphism, it follows that \(\kf S\) is a right
Serre \(\ck\)\n-functor.

The \ainf-functor $S$ is an equivalence if and only if \(\kf S\) is a
\(\ck\)\n-equivalence.
\end{proof}

When $\ca$ is an \ainf-algebra and $S$ is its identity endomorphism,
the natural transformation \(p:\ca\to\ca^*\) identifies with an
$\infty$\n-inner-product on $\ca$, as defined by
Tradler~\cite[Definition~5.3]{xxx0108027}.

\begin{corollary}\label{cor-ainf-Serre-exists-iff-K-Serre}
Let \(\ca\) be a unital \ainf-category. Then \(\ca\) admits a (right)
Serre \ainf-functor if and only if \(\kf\ca\) admits a (right) Serre
\(\ck\)\n-functor.
\end{corollary}

\begin{proof}
The ``only if'' part is proven above. Suppose \(\kf\ca\) admits a Serre
\(\ck\)\n-functor. By \propref{prop-existence-Serre-K-functor} this
implies representability of the \(\ck\)\n-functor
\[
\Hom_{\kf\ca}(Y,\_)^\op\cdot\und\ck(\_,\kk)=\kf[\Hom_\ca(Y,\_)^\op\cdot
D]:\kf\ca^\op\to\und\ck=\kf\uCom,
\]
for each object \(Y\in\Ob\ca\). By \corref{cor-ainf-equiv-repr-K} the
\ainf-functor
\[\Hom_\ca(Y,\_)^\op\cdot D
=(\Hom^\op_{\ca^\op},D)\mu^{\Ainfty}_{\mb2\to\mb1}|^Y_1:\ca^\op\to\uCom
\]
is representable for each \(Y\in\Ob\ca\). By
\propref{prop-bimodul-isom-A-g} the bimodule \(\ca^*\) corresponding to
the \ainf-functor \((\Hom^\op_{\ca^\op},D)\mu^{\Ainfty}_{\mb2\to\mb1}\)
is isomorphic to \(\ca^S\) for some \ainf-functor \(S:\ca\to\ca\).
\end{proof}

\begin{corollary}
Suppose \(\ca\) is a Hom\n-reflexive \ainf-category, \emph{i.e.}, the
complex \(\ca(X,Y)\) is reflexive in $\ck$ for each pair of objects
\(X,Y\in\Ob\ca\). If \(S:\ca\to\ca\) is a right Serre \ainf-functor,
then \(S\) is homotopy fully faithful.
\end{corollary}

\begin{proof}
The $\ck$\n-functor \(\kf S\) is fully faithful by
\propref{prop-e-K-psi-S-psi}.
\end{proof}

In particular, the above corollary applies if \(\kk\) is a field and
all homology spaces \(H^n(\ca(X,Y))\) are finite dimensional. If
\(\ca\) is closed under shifts, the latter condition is equivalent to
requiring that \(H^0(\ca(X,Y))\) be finite dimensional for each pair
\(X,Y\in\Ob\ca\). Indeed,
\(H^n(\ca(X,Y))=H^n(\kf\ca(X,Y))=H^0(\kf\ca(X,Y)[n])=H^0((\kf\ca)^\sh((X,0),(Y,n)))\).
The \(\ck\)\n-category \(\kf\ca\) is closed under shifts by results of
\secref{sec-ainf-categories-closed-under-shifts}, therefore there
exists an isomorphism \(\alpha:(Y,n)\to(Z,0)\) in \((\kf\ca)^\sh\), for
some \(Z\in\Ob\ca\). It induces an isomorphism
\((\kf\ca)^\sh(1,\alpha):(\kf\ca)^\sh((X,0),(Y,n))\to(\kf\ca)^\sh((X,0),(Z,0))=\kf\ca(X,Z)\)
in \(\ck\), thus an isomorphism in homology
\[ H^n(\ca(X,Y))= H^0((\kf\ca)^\sh((X,0),(Y,n)))
\simeq H^0(\kf\ca(X,Z)) =H^0(\ca(X,Z)).
\]
The latter space is finite dimensional by assumption.

\begin{theorem}\label{thm-Serre-A8-k-functor}
Suppose \(\kk\) is a field, \(\ca\) is a unital \ainf-category closed
under shifts. Then the following conditions are equivalent:
\begin{itemize}
\item[(a)] \(\ca\) admits a (right) Serre \ainf-functor;
\item[(b)] \(\kf\ca\) admits a (right) Serre \(\ck\)\n-functor;
\item[(c)] \(H^\bullet\ca\overset{\text{def}}=H^\bullet(\kf\ca)\)
    admits a (right) Serre \(\gr\)\n-functor;
\item[(d)] \(H^0(\ca)\) admits (right) Serre \(\kk\)\n-linear functor.
\end{itemize}
\end{theorem}

\begin{proof}
Equivalence of (a) and (b) is proven in
\corref{cor-ainf-Serre-exists-iff-K-Serre}. Conditions (b) and (c) are
equivalent due to \corref{cor-Serre-Hbullet-Serre}, because
\(H^\bullet:\ck\to\gr\) is an equivalence of symmetric monoidal
categories. Condition (c) implies (d) for arbitrary $\gr$\n-category by
\corref{cor-gr-Serre-k-Serre}, in particular, for \(H^\bullet\ca\).
Note that \(H^\bullet\ca\) is closed under shifts by
\secref{sec-ainf-categories-closed-under-shifts} and the discussion
preceding \propref{pro-Serre0-gr-Serre}. Therefore, (d) implies (c) due
to \propref{pro-Serre0-gr-Serre}.
\end{proof}

An application of this theorem is the following. Let $\kk$ be a field.
Drinfeld's construction of quotients of pretriangulated
$\dg$\n-categories \cite{Drinf:DGquot} allows to find a pretriangulated
$\dg$\n-category $\ca$ such that \(H^0(\ca)\) is some familiar derived
category (e.g. category \(D^b_{\text{coh}}(X)\) of complexes of
coherent sheaves on a projective variety $X$). If a right Serre functor
exists for \(H^0(\ca)\), then $\ca$ admits a right Serre \ainf-functor
$S$ by the above theorem. That is the case of
 \(H^0(\ca)\simeq D^b_{\text{coh}}(X)\), where $X$ is a smooth
projective variety \cite[Example~3.2]{MR1039961}. Notice that
\(S:\ca\to\ca\) does not have to be a $\dg$\n-functor.

\begin{proposition}
Let \(S:\ca\to\ca\), \(S':\cb\to\cb\) be right Serre \ainf-functors.
Let \(g:\cb\to\ca\) be an \ainf-equivalence. Then the \ainf-functors
\(S'g:\cb\to\ca\) and \(gS:\cb\to\ca\) are isomorphic.
\end{proposition}

\begin{proof}
Consider the following diagram:
\begin{diagram}
\cb^\op,\cb &&& \rTTo^{1,S'} &&& \cb^\op,\cb
\\
&\rdTTo^{g^\op,g}&&&&\ldTTo^{g^\op,g} &
\\
\dTTo<{\Hom^\op_{\cb^\op}} &\lTwoar^{(r^g)^\op}_\sim & \ca^\op,\ca
&\rTTo^{1,S} & \ca^\op,\ca &\lTwoar^{r^g}_\sim &\dTTo>{\Hom_\cb}
\\
&\ldTTo_{\Hom^\op_{\ca^\op}} &&&& \rdTTo_{\Hom_\ca} &
\\
\uCom^\op &&&\rTTo^D &&& \uCom
\end{diagram}
Here the natural \ainf-isomorphism \(r^g\) is that constructed in
\corref{cor-Hom-A-fop-f-Hom-B}. The exterior and the lower trapezoid
commute up to natural \ainf-isomorphisms by definition of right Serre
functor. It follows that the \ainf-functors
\((g^\op,S'g)\Hom_\ca:\cb^\op,\cb\to\uCom\) and
\((g^\op,gS)\Hom_\ca:\cb^\op,\cb\to\uCom\) are isomorphic. Consider the
\ainf-functors
\[
\cb \rTTo^{S'g} \ca \rTTo^\Yo \und\Ainftyu(\ca^\op;\uCom)
\rTTo^{\und\Ainftyu(g^\op;1)} \und\Ainftyu(\cb^\op;\uCom)
\]
and
\[
\cb \rTTo^{gS} \ca \rTTo^\Yo \und\Ainftyu(\ca^\op;\uCom)
\rTTo^{\und\Ainftyu(g^\op;1)} \und\Ainftyu(\cb^\op;\uCom)
\]
that correspond to \((g^\op,S'g)\Hom_\ca:\cb^\op,\cb\to\uCom\) and
\((g^\op,gS)\Hom_\ca:\cb^\op,\cb\to\uCom\) by closedness. More
precisely, the upper line is equal to
\((\und{\varphi}^{\Ainftyu})^{-1}((g^\op,S'g)\Hom_\ca)\) and the bottom
line is equal to
\((\und{\varphi}^{\Ainftyu})^{-1}((g^\op,gS)\Hom_\ca)\), where
\[
\und{\varphi}^{\Ainftyu}:\und\Ainftyu(\cb;\und\Ainftyu(\cb^\op;\uCom))
\to\und\Ainftyu(\cb^\op,\cb;\uCom)
\]
is the natural isomorphism of \ainf-categories coming from the closed
structure. It follows that the \ainf-functors
 \(S'\cdot g\cdot\Yo\cdot\und\Ainftyu(g^\op;1)\) and
\(g\cdot S\cdot\Yo\cdot\und\Ainftyu(g^\op;1)\) are isomorphic.
Obviously, the \ainf-functor \(g^\op\) is an equivalence, therefore so
is the \ainf-functor \(\und\Ainftyu(g^\op;1)\) since
\(\und\Ainftyu(\_;\uCom)\) is an \ainfu-2-functor. Therefore, the
\ainf-functors \(S'\cdot g\cdot\Yo\) and \(g\cdot S\cdot\Yo\) are
isomorphic. However, this implies that the \ainf-functors
\((1,S'g)\Hom_\ca=\und\varphi^{\Ainftyu}(S'\cdot g\cdot\Yo)\) and
\((1,gS)\Hom_\ca=\und\varphi^{\Ainftyu}(g\cdot S\cdot\Yo)\) are
isomorphic as well. These \ainf-functors correspond to
\((\ca,\cb)\)\n-\ainf-bimodules \(\ca^{S'g}\) and \(\ca^{gS}\)
respectively. \propref{pro-g-Ag-unital} implies an isomorphism between
the \ainf-functors \(S'g\) and \(gS\).
\end{proof}

\begin{remark}
Let $\ca$ be an \ainf-category, let \(S:\ca\to\ca\) be an
\ainf-functor. The \((0,0)\)-component of a cycle
\(p\in\und\Ainfty(\ca^\op,\ca;\uCom)(\ca^S,\ca^*)[1]^{-1}\) determines
for all objects $X$, $Y$ of $\ca$ a degree 0 map
\begin{multline*}
\kk \simeq T^0s\ca^\op(X,X)\tens T^0s\ca(Y,Y) \rTTo^{p_{00}}
s\uCom(\ca(X,YS),\ca^*(X,Y))
\\
\rTTo^{s^{-1}} \uCom(\ca(X,YS),\uCom(\ca(Y,X),\kk)).
\end{multline*}
The obtained mapping \(\ca(X,YS)\to\uCom(\ca(Y,X),\kk)\) is a chain
map, since \(p_{00}s^{-1}m_1=0\). Its homotopy class gives
\(\psi_{X,Y}\) from \eqref{equ-Serre-natural-transform} when the pair
\((S,p)\) is projected to \((\kf S,\psi=\kf p)\) via the multifunctor
$\kf$.
\end{remark}

\subsection{A strict case of a Serre $A_\infty$-functor.}
Let us consider a particularly simple case of an \ainf-category $\ca$
with a right Serre functor \(S:\ca\to\ca\) which is a strict
\ainf-functor (only the first component does not vanish) and with the
invertible natural \ainf-transformation
\(p:\ca^S\to\ca^*:\ca^\op,\ca\to\uCom\) whose only non-vanishing
component is
 \(p_{00}:T^0s\ca^\op(X,X)\tens T^0s\ca(Y,Y)\to
 s\uCom(\ca(X,YS),\uCom(\ca(Y,X),\kk))\).
Invertibility of $p$, equivalent to invertibility of \(p_{00}\), means
that the induced chain maps \(r_{00}:\ca(X,YS)\to\uCom(\ca(Y,X),\kk)\)
are homotopy invertible for all objects $X$, $Y$ of $\ca$. General
formulae for \(pB_1\) give the components \((pB_1)_{00}=p_{00}b_1\) and
\begin{equation}
(pB_1)_{kn} =([(1,S)\Hom_\ca]_{kn}\tens p_{00})b^{\uCom}_2
+(p_{00}\tens[\Hom_{\ca^\op}^\op\cdot D]_{kn})b^{\uCom}_2
\label{eq-(pB1)kn}
\end{equation}
for \(k+n>0\). Since $p$ is natural, \(pB_1=0\), thus the right hand
side of \eqref{eq-(pB1)kn} has to vanish. Expanding the first summand
we get
\begin{multline*}
(-)^k\bigl[ T^ks\ca^\op(X_0,X_k)\tens T^ns\ca(Y_0,Y_n) \rTTo^\coev
\\
 \uCom\bigl(s\ca(X_0,Y_0S),
s\ca(X_0,Y_0S)\tens T^ks\ca^\op(X_0,X_k)\tens T^ns\ca(Y_0,Y_n)\bigr)
\\
 \rTTo^{\uCom(1,\rho_c(1^{\tens k+1}\tens S_1^{\tens n})b_{k+1+n}r_{00})}
\uCom\bigl(s\ca(X_0,Y_0S),s\uCom(\ca(Y_n,X_k),\kk)\bigr)
\\
\rTTo^{[-1]s} s\uCom\bigl(\ca(X_0,Y_0S),\uCom(\ca(Y_n,X_k),\kk)\bigr)
\bigr].
\end{multline*}
Expanding the second summand we obtain
\begin{multline*}
-(-)^k\bigl[ T^ks\ca^\op(X_0,X_k)\tens T^ns\ca(Y_0,Y_n)
\\
\rTTo^\coev
 \uCom\bigl(s\ca(Y_n,X_k),
s\ca(Y_n,X_k)\tens T^ks\ca^\op(X_0,X_k)\tens T^ns\ca(Y_0,Y_n)\bigr)
\\
\rTTo^{\uCom(1,(123)_c(1\tens1\tens\omega^0_c))}
 \uCom\bigl(s\ca(Y_n,X_k),
T^ns\ca(Y_0,Y_n)\tens s\ca(Y_n,X_k)\tens T^ks\ca(X_k,X_0)\bigr)
\\
 \rTTo^{\uCom(1,b_{n+1+k})}
\uCom\bigl(s\ca(Y_n,X_k),s\ca(Y_0,X_0)\bigr) \rTTo^{[-1]s}
s\uCom\bigl(\ca(Y_n,X_k),\ca(Y_0,X_0)\bigr)
\\
\rTTo^\coev
 \uCom\bigl(s\ca(X_0,Y_0S),
s\ca(X_0,Y_0S)\tens s\uCom(\ca(Y_n,X_k),\ca(Y_0,X_0))\bigr)
\rTTo^{\uCom(1,(r_{00}\tens1)cb_2)}
\\
\uCom\bigl(s\ca(X_0,Y_0S),s\uCom(\ca(Y_n,X_k),\kk)\bigr) \rTTo^{[-1]s}
s\uCom\bigl(\ca(X_0,Y_0S),\uCom(\ca(Y_n,X_k),\kk)\bigr) \bigr].
\end{multline*}
Sum of the two above expressions must vanish. The obtained equation can
be simplified further by closedness of $\Com$. The homotopy isomorphism
\(r_{00}\) induces the pairing
\[ q_{00} =\bigl[ \ca(Y,X)\tens\ca(X,YS) \rTTo^{1\tens r_{00}}
\ca(Y,X)\tens\uCom(\ca(Y,X),\kk) \rTTo^\ev \kk \bigr].
\]
Using it we write down the naturality condition for $p$ as follows: for
all \(k\ge0\), \(n\ge0\)
\begin{multline}
\bigl[
 \ca(Y_n,X_k)\tens T^k\ca(X_k,X_0)\tens\ca(X_0,Y_0S)\tens T^n\ca(Y_0,Y_n)
\\
\hfill\rTTo^{(1^{\tens3}\tens(sS_1s^{-1})^{\tens n})(1\tens m_{k+1+n})}
\ca(Y_n,X_k)\tens\ca(X_k,Y_nS) \rTTo^{q_{00}} \kk \bigr] \quad
\\
\quad =(-)^{(k+1)(n+1)} \bigl[
 \ca(Y_n,X_k)\tens T^k\ca(X_k,X_0)\tens\ca(X_0,Y_0S)\tens T^n\ca(Y_0,Y_n)
\hfill
\\
 \rTTo^{(1234)_c}
T^n\ca(Y_0,Y_n)\tens\ca(Y_n,X_k)\tens T^k\ca(X_k,X_0)\tens\ca(X_0,Y_0S)
\\
\rTTo^{m_{n+1+k}\tens1} \ca(Y_0,X_0)\tens\ca(X_0,Y_0S) \rTTo^{q_{00}}
\kk \bigr].
\label{eq-Snmk1nq00-1234mn1kq00}
\end{multline}
Let us give a sufficient condition for this equation to hold true.

\begin{proposition}
Let $\ca$ be an \ainf-category, and let \(S:\ca\to\ca\) be a strict
\ainf-functor. Suppose given a pairing
\(q_{00}:\ca(Y,X)\tens\ca(X,YS)\to\kk\) for all objects $X$, $Y$ of
$\ca$. Assume that for all \(X,Y\in\Ob\ca\)
\begin{enumerate}
\renewcommand{\labelenumi}{(\alph{enumi})}
\item $q_{00}$ is a chain map;

\item the induced chain map
\[ r_{00} = \bigl[ \ca(X,YS) \rTTo^\coev
\uCom(\ca(Y,X),\ca(Y,X)\tens\ca(X,YS)) \rTTo^{\uCom(1,q_{00})}
\uCom(\ca(Y,X),\kk) \bigr]
\]
is homotopy invertible;

\item the pairing $q_{00}$ is symmetric in a sense similar to
diagram~\eqref{dia-1Sphi-cphi}, namely, the following diagram of chain
maps commutes:
\begin{diagram}[LaTeXeqno]
\ca(X,YS)\tens\ca(Y,X) &\rTTo^{1\tens sS_1s^{-1}}
&\ca(X,YS)\tens\ca(YS,XS)
\\
\dTTo<c &= &\dTTo>{q_{00}}
\\
\ca(Y,X)\tens\ca(X,YS) &\rTTo^{q_{00}} &\kk
\label{dia-1Sq00-cq00}
\end{diagram}

\item the following equation holds for all $k\ge0$ and all objects
$X_0$, \dots, $X_k$, $Y$
\begin{multline}
\bigl[ \ca(Y,X_k)\tens T^k\ca(X_k,X_0)\tens\ca(X_0,YS)
\rTTo^{1\tens m_{k+1}}
\ca(Y,X_k)\tens\ca(X_k,YS) \rTTo^{q_{00}} \kk \bigr]
\\
\quad =(-)^{k+1} \bigl[ \ca(Y,X_k)\tens T^k\ca(X_k,X_0)\tens\ca(X_0,YS)
\hfill
\\
\rTTo^{m_{1+k}\tens1} \ca(Y,X_0)\tens\ca(X_0,YS) \rTTo^{q_{00}} \kk
\bigr].
 \label{eq-1mk1q00-m1k1q00}
\end{multline}
\end{enumerate}
Then the natural \ainf-transformation
\(p:\ca^S\to\ca^*:\ca^\op,\ca\to\uCom\) with the only non-vanishing
component \(p_{00}:1\mapsto r_{00}\) is invertible and \(S:\ca\to\ca\)
is a Serre \ainf-functor.
\end{proposition}

Notice that \eqref{eq-1mk1q00-m1k1q00} is precisely the case of
\eqref{eq-Snmk1nq00-1234mn1kq00} with $n=0$. On the other hand,
diagram~\eqref{dia-1Sphi-cphi} written for $\ck$\n-category
\(\cc=\kf\ca\) and the pairing \(\phi=[q_{00}]\) says that
\eqref{dia-1Sq00-cq00} has to commute only up to homotopy. Thus,
condition (c) is sufficient but not necessary.

\begin{proof}
We have to prove equation~\eqref{eq-Snmk1nq00-1234mn1kq00} for all
$k\ge0$, $n\ge0$. The case of $n=0$ holds by condition~(d). Let us
proceed by induction on $n$. Assume that
\eqref{eq-Snmk1nq00-1234mn1kq00} holds true for all $k\ge0$,
 \(0\le n<N\). Let us prove equation~\eqref{eq-Snmk1nq00-1234mn1kq00}
for $k\ge0$, $n=N$. We have
\begin{align*}
&(-)^{(k+1)(n+1)} (13524)_c \cdot (m_{n+1+k}\tens1) \cdot q_{00}
\\
&\overset{(d)}= (-)^{(k+1)(n+1)+k+n+1} (13524)_c \cdot
 (1\tens m_{n+k+1}) \cdot q_{00}
\\
&=(-)^{kn} (12345)_c \cdot (m_{n+k+1}\tens1) \cdot c \cdot q_{00}
\\
&\overset{(c)}= (-)^{kn} (12345)_c \cdot (m_{n+k+1}\tens sS_1s^{-1})
\cdot q_{00}
\\
&=(-)^{(k+2)n} (1^{\tens3}\tens sS_1s^{-1}\tens1) \cdot (12345)_c \cdot
(m_{n+k+1}\tens1) \cdot q_{00}
\\
&\underset{\text{for }k+1,n-1}{\overset{\text{by \eqref{eq-Snmk1nq00-1234mn1kq00}}}=}
(1^{\tens3}\tens sS_1s^{-1}\tens1) \cdot
(1^{\tens4}\tens T^{n-1}(sS_1s^{-1})) \cdot (1\tens m_{k+1+n}) \cdot
q_{00}:
\\
&\ca(Y_n,X_k)\tens T^k\ca(X_k,X_0)\tens\ca(X_0,Y_0S)\tens\ca(Y_0,Y_1)
\tens T^{n-1}\ca(Y_1,Y_n) \to \kk.
\end{align*}
This is just equation~\eqref{eq-Snmk1nq00-1234mn1kq00} for $k$, $n$.
\end{proof}

Some authors like to consider a special case of the above in which
\(S=[d]\) is the shift functor (when it makes sense), the paring
\(q_{00}\) is symmetric and cyclically symmetric with respect to
$n$\n-ary compositions, cf. \cite[Section~6.2]{math.QA/0412149}. Then
$\ca$ is called a Calabi--Yau \ainf-category. General Serre
\ainf-functors cover wider scope, although they require more data to
work with.

\appendix
\section{The Yoneda Lemma.}\label{sec-The-Yoneda-Lemma}
A version of the classical Yoneda Lemma is presented in Mac~Lane's
book~\cite[Section~III.2]{MacLane} as the following statement. For any
category $\cc$ there is an isomorphism of functors
\[ \ev^{\Cat} \simeq \bigl[ \cc\times\und\Cat(\cc,\Set)
\rTTo^{\Yo^\op\times1} \und\Cat(\cc,\Set)^\op\times\und\Cat(\cc,\Set)
\rTTo^{\Hom_{\und\Cat(\cc,\Set)}} \Set \bigr],
\]
where \(\Yo:\cc^\op\to\und\Cat(\cc,\Set)\), \(X\mapsto\cc(X,\_)\), is
the Yoneda embedding. Here we generalize this to \ainf-setting.

\begin{theorem}[The Yoneda Lemma]\label{thm-Yoneda-Lemma}
For any \ainf-category $\ca$ there is a natural \ainf-transformation
\[ \Omega: \ev^{\Ainfty} \to \bigl[ \ca,\und\Ainfty(\ca;\uCom)
\rTTo^{\Yo^\op,1}
\und\Ainfty(\ca;\uCom)^\op,\und\Ainfty(\ca;\uCom)
\rTTo^{\Hom_{\und\Ainfty(\ca;\uCom)}} \uCom \bigr].
\]
If the \ainf-category $\ca$ is unital, $\Omega$ restricts to an
\emph{invertible} natural \ainf-transformation
\begin{diagram}[width=6em]
\ca,\und\Ainftyu(\ca;\uCom) &&\rTTo^{\ev^{\Ainftyu}} &&\uCom
\\
&\rdTTo<{\Yo^\op,1} &\dTwoar<\Omega
&\ruTTo>{\Hom_{\und\Ainftyu(\ca;\uCom)}} &
\\
&&\und\Ainftyu(\ca;\uCom)^\op,\und\Ainftyu(\ca;\uCom) &&
\end{diagram}
\end{theorem}

Previously published \ainf-versions of Yoneda Lemma contented with the
statement that for unital \ainf-category $\ca$, the Yoneda
\ainf-functor \(\Yo:\ca^\op\to\und\Ainftyu(\ca;\uCom)\) is homotopy
full and faithful \cite[Theorem~9.1]{Fukaya:FloerMirror-II},
\cite[Theorem~A.11]{math.CT/0306018}. A more general form of the Yoneda
Lemma is considered by Seidel \cite[Lemma~2.12]{SeidelP-book-Fukaya}
over a ground field $\kk$. We shall see that these are corollaries of
the above theorem.

\begin{proof}
First of all we describe the \ainf-transformation $\Omega$ for an
arbitrary \ainf-category $\ca$. The discussion of
\secref{sec-Restriction-scalars} applied to the \ainf-functor
\[ \psi =\bigl[ \ca,\und\Ainfty(\ca;\uCom) \rTTo^{\Yo^\op,1}
\und\Ainfty(\ca;\uCom)^\op,\und\Ainfty(\ca;\uCom)
\rTTo^{\Hom_{\und\Ainfty(\ca;\uCom)}} \uCom \bigr]
\]
presents the corresponding
$\ca^\op$-\(\und\Ainfty(\ca;\uCom)\)-bimodule
\(\cq=\sS{_\Yo}{\und\Ainfty(\ca;\uCom)}_1\) via the regular
\ainf-bimodule. Thus,
\[ (\cq(X,f),sb^\cq_{00}s^{-1})
=(\und\Ainfty(\ca;\uCom)(H^X,f),sB_1s^{-1}).
\]
According to \eqref{equ-components-of-H-X}
\(H^X=\ca^\op(\_,X)=\ca(X,\_)\) has the components
\begin{multline}
H^X_k =(\Hom_{\ca^\op})_{k0} =\bigl[ T^ks\ca(Y,Z) \rTTo^{\coev^\Com}
\uCom(s\ca(X,Y),s\ca(X,Y)\tens T^ks\ca(Y,Z))
\\
\rTTo^{\uCom(1,b^\ca_{1+k})} \uCom(s\ca(X,Y),s\ca(X,Z)) \rTTo^{[-1]s}
s\uCom(\ca(X,Y),\ca(X,Z)) \bigr].
\label{eq-HXk-HomAopk0}
\end{multline}
We have \(b^\cq_{00}=B_1\) and, moreover, by \eqref{eq-bfPg-f1g-bP}
\begin{multline*}
\check{b}^\cq =\bigl[
 Ts\ca^\op(Y,X)\tens s\cq(X,f)\tens Ts\und\Ainfty(\ca;\uCom)(f,g)
\\
\rTTo^{\Yo\tens1\tens1}
 Ts\und\Ainfty(\ca;\uCom)(H^Y,H^X)\tens
 s\und\Ainfty(\ca;\uCom)(H^X,f)\tens Ts\und\Ainfty(\ca;\uCom)(f,g)
\\
\rTTo^{\check{B}} s\und\Ainfty(\ca;\uCom)(H^Y,g) =s\cq(Y,g) \bigr].
\end{multline*}
Since \(\und\Ainfty(\ca;\uCom)\) is a differential graded category,
\(B_p=0\) for $p>2$. Therefore, \(b^\cq_{kn}=0\) if $n>1$, and
\(b^\cq_{k1}=0\) if $k>0$. The non-trivial components are (for $k>0$)
\begin{multline}
b^\cq_{k0} =\bigl[
 T^ks\ca^\op(Y,X)\tens s\cq(X,f)\tens
 T^0s\und\Ainfty(\ca;\uCom)(f,f)\rTTo^{\Yo_k\tens1\tens1}
\\
s\und\Ainfty(\ca;\uCom)(H^Y,H^X)\tens s\und\Ainfty(\ca;\uCom)(H^X,f)
\rTTo^{B_2} s\und\Ainfty(\ca;\uCom)(H^Y,f) =s\cq(Y,f) \bigr],
\\
\hskip\multlinegap b^\cq_{01} =\bigl[
 T^0s\ca^\op(X,X)\tens s\cq(X,f)\tens
 s\und\Ainfty(\ca;\uCom)(f,g)\hfill
\\
\rTTo^{B_2} s\und\Ainfty(\ca;\uCom)(H^X,g) =s\cq(X,g) \bigr].
\label{equ-components-bQ}
\end{multline}

Denote by $\ce$ the $\ca^\op$-\(\und\Ainfty(\ca;\uCom)\)-bimodule
corresponding to the \ainf-functor
\(\ev^{\Ainfty}:\ca,\und\Ainfty(\ca;\uCom)\to\uCom\). For any object
$X$ of $\ca$ and any \ainf-functor \(f:\ca\to\uCom\) the complex
\((\ce(X,f),sb^\ce_{00}s^{-1})\) is \((Xf,d)\). According to
\eqref{equ-check-b-P-+}
\begin{multline*}
\check b^\ce_+ =\bigl[
 Ts\ca^\op(Y,X)\tens s\ce(X,f)\tens Ts\und\Ainfty(\ca;\uCom)(f,g)
\\
\rTTo^{c\tens1}
 s\ce(X,f)\tens Ts\ca^\op(Y,X)\tens Ts\und\Ainfty(\ca;\uCom)(f,g)
\rTTo^{1\tens\gamma\tens1}
\\
s\ce(X,f)\tens Ts\ca(X,Y)\tens Ts\und\Ainfty(\ca;\uCom)(f,g)
\rTTo^{1\tens\check\ev^{\Ainfty}} s\ce(X,f)\tens s\uCom(Xf,Yg)
\\
\rTTo^{1\tens s^{-1}[1]} Xf[1]\tens\uCom(Xf[1],Yg[1])
\rTTo^{\ev^{\Com}} Yg[1] =s\ce(Y,g) \bigr].
\end{multline*}
Explicit formula~\eqref{eq-evA80-evA81} for \(\ev^{\Ainfty}\)
shows that \(b^\ce_{kn}=0\) if $n>1$. The remaining components are
described as
\begin{multline*}
b^\ce_{k0} =\bigl[
 T^ks\ca^\op(Y,X)\tens s\ce(X,f)\tens T^0s\und\Ainfty(\ca;\uCom)(f,f)
\\
\rTTo^{c\tens1} s\ce(X,f)\tens T^ks\ca^\op(Y,X) \rTTo^{1\tens\gamma}
s\ce(X,f)\tens T^ks\ca(X,Y) \rTTo^{1\tens f_k}
\\
Xf[1]\tens s\uCom(Xf,Yf) \rTTo^{1\tens s^{-1}[1]}
Xf[1]\tens\uCom(Xf[1],Yf[1]) \rTTo^{\ev^{\Com}} Yf[1] =s\ce(Y,f)
\bigr]
\end{multline*}
for $k>0$, and if $k\ge0$ there is
\begin{multline*}
b^\ce_{k1} =\bigl[
 T^ks\ca^\op(Y,X)\tens s\ce(X,f)\tens s\und\Ainfty(\ca;\uCom)(f,g)
\\
\rTTo^{c\tens1}
 s\ce(X,f)\tens T^ks\ca^\op(Y,X)\tens s\und\Ainfty(\ca;\uCom)(f,g)
\\
\rTTo^{1\tens\gamma\tens1}
 s\ce(X,f)\tens T^ks\ca(X,Y)\tens s\und\Ainfty(\ca;\uCom)(f,g)
\\
\rTTo^{1\tens1\tens\pr_k}
 s\ce(X,f)\tens T^ks\ca(X,Y)\tens\uCom(T^ks\ca(X,Y),s\uCom(Xf,Yg))
\rTTo^{1\tens\ev^{\uCom}}
\\
Xf[1]\tens s\uCom(Xf,Yg) \rTTo^{1\tens s^{-1}[1]}
Xf[1]\tens\uCom(Xf[1],Yg[1]) \rTTo^{\ev^{\Com}} Yg[1] =s\ce(Y,g)
\bigr].
\end{multline*}

The \ainf-transformation $\Omega$ in question is constructed via a
homomorphism
\[ \mho =(\Omega s^{-1})\Phi:
Ts\ca^\op\tens s\ce\tens Ts\und\Ainfty(\ca;\uCom)
\to Ts\ca^\op\tens s\cq\tens Ts\und\Ainfty(\ca;\uCom)
\]
of $Ts\ca^\op$-\(Ts\und\Ainfty(\ca;\uCom)\)-bicomodules thanks to
\propref{pro-dg-cat-ACbimod-A(ACC)-isomorphic}. Its matrix coefficients
are recovered from the components via
formula~\eqref{eq-fklmn-phiphifpsipsi} as
\begin{multline*}
\mho_{kl;mn}=\sum_{\substack{m+p=k\\ q+n=l}}
(1^{\tens m}\tens\mho_{pq}\tens1^{\tens n}):
\\
T^ks\ca^\op\tens s\ce\tens T^ls\und\Ainfty(\ca;\uCom) \to
T^ms\ca^\op\tens s\cq\tens T^ns\und\Ainfty(\ca;\uCom).
\end{multline*}

The composition of the morphism
\[ \mho_{pq}:
T^ps\ca^\op(X_0,X_p)\tens X_pf_0[1]\tens T^qs\und\Ainfty(\ca;\uCom)(f_0,f_q)
\to s\und\Ainfty(\ca;\uCom)(H^{X_0},f_q)
\]
with the projection
\begin{equation}
\pr_n:s\und\Ainfty(\ca;\uCom)(H^{X_0},f_q)\to
\uCom(T^ns\ca(Z_0,Z_n),s\uCom(\ca(X_0,Z_0),Z_nf_q))
\label{eq-prn-sAinfty}
\end{equation}
is given by the composite
\begin{multline}
\mho_{pq;n}\overset{\text{def}}=\mho_{pq}\cdot\pr_n=(-)^{p+1}\bigl[
T^ps\ca^\op(X_0,X_p)\tens X_pf_0[1]\tens T^qs\und\Ainfty(\ca;\uCom)(f_0,f_q)
\\
\quad\rTTo^{\coev^\Com}
\uCom(s\ca(X_0,Z_0)\tens T^ns\ca(Z_0,Z_n),\hfill
\\
s\ca(X_0,Z_0)\tens T^ns\ca(Z_0,Z_n)\tens
T^ps\ca^\op(X_0,X_p)\tens X_pf_0[1]\tens
T^qs\und\Ainfty(\ca;\uCom)(f_0,f_q))
\\
\quad\rTTo^{\uCom(1,\text{perm})}
\uCom(s\ca(X_0,Z_0)\tens T^ns\ca(Z_0,Z_n),\hfill
\\
X_pf_0[1]\tens T^ps\ca(X_p,X_0)\tens
s\ca(X_0,Z_0)\tens T^ns\ca(Z_0,Z_n)\tens
T^qs\und\Ainfty(\ca;\uCom)(f_0,f_q))
\\
\quad\rTTo^{\uCom(1,1\tens\ev^{\Ainfty}_{p+1+n,q})}
\uCom(s\ca(X_0,Z_0)\tens T^ns\ca(Z_0,Z_n),X_pf_0[1]\tens
s\uCom(X_pf_0,Z_nf_q))\hfill
\\
\quad\rTTo^{\uCom(1,1\tens s^{-1}[1])}
\uCom(s\ca(X_0,Z_0)\tens
T^ns\ca(Z_0,Z_n),X_pf_0[1]\tens\uCom(X_pf_0[1],Z_nf_q[1]))\hfill
\\
\quad\rTTo^{\uCom(1,\ev^\Com)}
\uCom(s\ca(X_0,Z_0)\tens T^ns\ca(Z_0,Z_n),Z_nf_q[1])\hfill
\\
\quad\rTTo^{(\und{\varphi}^\Com)^{-1}}
\uCom(T^ns\ca(Z_0,Z_n),\uCom(s\ca(X_0,Z_0),Z_nf_q[1]))\hfill
\\
\rTTo^{\uCom(1,[-1]s)}
\uCom(T^ns\ca(Z_0,Z_n),s\uCom(\ca(X_0,Z_0),Z_nf_q))
\bigr].
\label{eq-agemo-pqn}
\end{multline}
Thus, an element
 \(x_1\tdt x_p\tens y\tens r_1\tdt r_q\in
 T^ps\ca^\op(X_0,X_p)\tens X_pf_0[1]\tens
 T^qs\und\Ainfty(\ca;\uCom)(f_0,f_q)\)
is mapped to an \ainf-transformation
 \((x_1\tdt x_p\tens y\tens r_1\tdt r_q)\mho_{pq}\in
 s\und\Ainfty(\ca;\uCom)(H^{X_0},f_q)\)
with components
\begin{gather*}
[(x_1\tdt x_p\tens y\tens r_1\tdt r_q)\mho_{pq}]_n:T^ns\ca(Z_0,Z_n)\to
s\uCom(\ca(X_0,Z_0),Z_nf_q),\\
z_1\tdt z_n\mapsto(z_1\tdt z_n\tens x_1\tdt x_p\tens y\tens r_1\tdt
r_q)\mho'_{pq;n},
\end{gather*}
where
 \(\mho'_{pq;n}\overset{\text{def}}=(1^{\tens n}\tens\mho_{pq;n})\ev^\Com
 =(1^{\tens n}\tens\mho_{pq}\cdot\pr_n)\ev^\Com=(1^{\tens n}\tens\mho_{pq})\ev^{\Ainfty}_{n1}\)
is given by
\begin{multline*}
\mho'_{pq;n}=(-)^{p+1}\bigl[
T^ns\ca(Z_0,Z_n)\tens T^ps\ca^\op(X_0,X_p)\tens X_pf_0[1]\tens
T^qs\und\Ainfty(\ca;\uCom)(f_0,f_q)
\\
\quad\rTTo^{\coev^\Com}
\uCom(s\ca(X_0,Z_0),s\ca(X_0,Z_0)\tens T^ns\ca(Z_0,Z_n)\hfill
\\
\hfill\tens T^ps\ca^\op(X_0,X_p)\tens X_pf_0[1]
\tens T^qs\und\Ainfty(\ca;\uCom)(f_0,f_q))\quad
\\
\quad\rTTo^{\uCom(1,\text{perm})}
\uCom(s\ca(X_0,Z_0),X_pf_0[1]\tens T^ps\ca(X_p,X_0)\hfill
\\
\hfill\tens s\ca(X_0,Z_0)\tens T^ns\ca(Z_0,Z_n)
\tens T^qs\und\Ainfty(\ca;\uCom)(f_0,f_q))\quad
\\
\quad\rTTo^{\uCom(1,1\tens\ev^{\Ainfty}_{p+1+n,q})}
\uCom(s\ca(X_0,Z_0),X_pf_0[1]\tens s\uCom(X_pf_0,Z_nf_q))\hfill
\\
\quad\rTTo^{\uCom(1,1\tens s^{-1}[1])}
\uCom(s\ca(X_0,Z_0),X_pf_0[1]\tens\uCom(X_pf_0[1],Z_nf_q[1]))\hfill
\\
\rTTo^{\uCom(1,\ev^\Com)}
\uCom(s\ca(X_0,Z_0),Z_nf_q[1])
\rTTo^{[-1]s}
s\uCom(\ca(X_0,Z_0),Z_nf_q)
\bigr].
\end{multline*}
It follows that
 \(\mho_{pq}:T^ps\ca^\op\tens s\ce\tens T^qs\und\Ainfty(\ca;\uCom)
 \to s\cq\)
vanishes if \(q>1\). The other components are given by
\begin{multline}
\mho'_{p0;n}=(-)^{p+1}\bigl[
T^ns\ca(Z_0,Z_n)\tens T^ps\ca^\op(X_0,X_p)\tens X_pf[1]
\\
\quad\rTTo^{\coev^\Com}
\uCom(s\ca(X_0,Z_0),s\ca(X_0,Z_0)\tens T^ns\ca(Z_0,Z_n)\tens T^ps\ca^\op(X_0,X_p)\tens
X_pf[1])\hfill
\\
\quad\rTTo^{\uCom(1,\text{perm})}
\uCom(s\ca(X_0,Z_0),X_pf[1]\tens T^ps\ca(X_p,X_0)\tens s\ca(X_0,Z_0)\tens
T^ns\ca(Z_0,Z_n))\hfill
\\
\quad\rTTo^{\uCom(1,1\tens f_{p+1+n})}
\uCom(s\ca(X_0,Z_0),X_pf[1]\tens s\uCom(X_pf,Z_nf))\hfill
\\
\quad\rTTo^{\uCom(1,1\tens s^{-1}[1])}
\uCom(s\ca(X_0,Z_0),X_pf[1]\tens\uCom(X_pf[1],Z_nf[1]))\hfill
\\
\rTTo^{\uCom(1,\ev^\Com)}
\uCom(s\ca(X_0,Z_0),Z_nf[1])
\rTTo^{[-1]s}
s\uCom(\ca(X_0,Z_0),Z_nf)
\bigr]
\label{eq-O'p0n}
\end{multline}
and
\begin{multline*}
\mho'_{p1;n}=(-)^{p+1}\bigl[
T^ns\ca(Z_0,Z_n)\tens T^ps\ca^\op(X_0,X_p)\tens X_pf[1]\tens
s\und\Ainfty(\ca;\uCom)(f,g)
\\
\quad\rTTo^{\coev^\Com}
\uCom(s\ca(X_0,Z_0),s\ca(X_0,Z_0)\tens T^ns\ca(Z_0,Z_n)\hfill
\\
\hfill\tens T^ps\ca^\op(X_0,X_p)\tens X_pf[1]
\tens s\und\Ainfty(\ca;\uCom)(f,g))\quad
\\
\quad\rTTo^{\uCom(1,\text{perm})}
\uCom(s\ca(X_0,Z_0),X_pf[1]\tens T^ps\ca(X_p,X_0)\hfill
\\
\hfill\tens s\ca(X_0,Z_0)\tens T^ns\ca(Z_0,Z_n)
\tens s\und\Ainfty(\ca;\uCom)(f,g))\quad
\\
\quad\rTTo^{\uCom(1,1\tens1^{\tens p+1+n}\tens\pr_{p+1+n})}
\uCom(s\ca(X_0,Z_0),X_pf[1]\tens T^ps\ca(X_p,X_0)\tens s\ca(X_0,Z_0)\hfill
\\
\tens T^ns\ca(Z_0,Z_n)\tens\uCom(T^ps\ca(X_p,X_0)\tens s\ca(X_0,Z_0)\tens
T^ns\ca(Z_0,Z_n),s\uCom(X_pf,Z_ng)))
\\
\quad\rTTo^{\uCom(1,1\tens\ev^\Com)}
\uCom(s\ca(X_0,Z_0),X_pf[1]\tens s\uCom(X_pf,Z_ng))\hfill
\\
\quad\rTTo^{\uCom(1,1\tens s^{-1}[1])}
\uCom(s\ca(X_0,Z_0),X_pf[1]\tens\uCom(X_pf[1],Z_ng[1]))\hfill
\\
\rTTo^{\uCom(1,\ev^\Com)}
\uCom(s\ca(X_0,Z_0),Z_ng[1])
\rTTo^{[-1]s}
s\uCom(\ca(X_0,Z_0),Z_ng)
\bigr].
\end{multline*}

Naturality of the \ainf-transformation $\Omega$ is implied by the
following

\begin{lemma}\label{lem-bicomodule-homomorphism-mho-chain-map}
The bicomodule homomorphism $\mho$ is a chain map.
\end{lemma}

    \ifx\chooseClass1
\straightForward.
\proofInArXiv.
    \else
\begin{proof}
Equivalently, we have to prove the equation
\(\mho\check{b}^\cq=b^\ce\check{\mho}\). In components, the expressions
\begin{align}
(\mho\check{b}^\cq)_{kl} &=\sum_{\substack{m+p=k\\ q+n=l}}
(1_{s\ca^\op}^{\tens m}\tens\mho_{pq}
\tens1_{s\und\Ainfty(\ca;\uCom)}^{\tens n})b^\cq_{mn},
\label{equ-1-mho-1-bQ}
\\
(b^\ce\check{\mho})_{kl} &=\sum_{p+t+q=l}
 (1_{s\ca^\op}^{\tens k}\tens1_{s\ce}
 \tens1_{s\und\Ainfty(\ca;\uCom)}^{\tens p}
 \tens b^{\und\Ainfty(\ca;\uCom)}_t
 \tens1_{s\und\Ainfty(\ca;\uCom)}^{\tens q})\mho_{k,p+1+q}
\label{equ-111-B-1-mho}
\\
&+\sum_{\substack{m+i=k\\ j+n=l}}
 (1_{s\ca^\op}^{\tens m}\tens b^\ce_{ij}
 \tens1_{s\und\Ainfty(\ca;\uCom)}^{\tens n})\mho_{mn}
\label{equ-1-bE-1-mho}
\\
&+\sum_{a+u+c=k}
 (1_{s\ca^\op}^{\tens a}\tens b^{\ca^\op}_u\tens1_{s\ca^\op}^{\tens c}
\tens1_{s\ce}\tens1_{s\und\Ainfty(\ca;\uCom)}^{\tens l})\mho_{a+1+c,l}
\label{equ-1-bop-111-mho}
\end{align}
have to coincide for all \(k,l\ge0\). Let us analyze the details of
this equation. Since \(b^\cq_{mn}=0\) unless \(n=0\) or
\((m,n)=(0,1)\), it follows that the right hand side of
\eqref{equ-1-mho-1-bQ} reduces to
\[
\sum_{m=0}^k(1^{\tens m}_{s\ca^\op}\tens\mho_{k-m,l})b^\cq_{m0}
+(\mho_{k,l-1}\tens1_{s\und\Ainfty(\ca;\uCom)})b^\cq_{01}.
\]
Since \(\und\Ainfty(\ca;\uCom)\) is a differential graded category,
sum~\eqref{equ-111-B-1-mho} reduces to
\begin{multline*}
\sum_{p=1}^{l}(1^{\tens k}_{s\ca^\op}\tens1_{s\ce}
\tens1^{\tens(p-1)}_{s\und\Ainfty(\ca;\uCom)}\tens B_1
\tens1^{\tens(l-p)}_{s\und\Ainfty(\ca;\uCom)})\mho_{kl}
\\
+\sum_{p=1}^{l-1}(1^{\tens k}_{s\ca^\op}\tens1_{s\ce}
\tens1^{\tens(p-1)}_{s\und\Ainfty(\ca;\uCom)}\tens B_2
\tens1^{\tens(l-p-1)}_{s\und\Ainfty(\ca;\uCom)})\mho_{k,l-1}.
\end{multline*}
Since \(b^\ce_{ij}=0\) if \(j>1\), sum~\eqref{equ-1-bE-1-mho} equals
\begin{multline*}
\sum_{m=0}^{k-1}(1^{\tens m}_{s\ca^\op}\tens b^\ce_{k-m,0}\tens1^{\tens
l}_{s\und\Ainfty(\ca;\uCom)})\mho_{ml} +(1^{\tens k}_{s\ca^\op}\tens
b^\ce_{00}\tens1^{\tens l}_{s\und\Ainfty(\ca;\uCom)})\mho_{kl}
\\
+\sum_{m=0}^k(1^{\tens m}_{s\ca^\op}\tens b^\ce_{k-m,1}\tens1^{\tens
(l-1)}_{s\und\Ainfty(\ca;\uCom)})\mho_{m,l-1}.
\end{multline*}
Sum~\eqref{equ-1-bop-111-mho} does not allow further simplification.
Therefore, the equation to prove is
\begin{multline*}
\sum_{m=0}^k(1^{\tens m}_{s\ca^\op}\tens\mho_{k-m,l})b^\cq_{m0}
+(\mho_{k,l-1}\tens1_{s\und\Ainfty(\ca;\uCom)})b^\cq_{01}
\\
\quad=\sum_{p=1}^{l}(1^{\tens k}_{s\ca^\op}\tens1_{s\ce}
\tens1^{\tens(p-1)}_{s\und\Ainfty(\ca;\uCom)}\tens B_1
\tens1^{\tens(l-p)}_{s\und\Ainfty(\ca;\uCom)})\mho_{kl}\hfill
\\
\quad+\sum_{p=1}^{l-1}(1^{\tens k}_{s\ca^\op}\tens1_{s\ce}
\tens1^{\tens(p-1)}_{s\und\Ainfty(\ca;\uCom)}\tens B_2
\tens1^{\tens(l-p-1)}_{s\und\Ainfty(\ca;\uCom)})\mho_{k,l-1}\hfill
\\
\quad+\sum_{m=0}^k(1^{\tens m}_{s\ca^\op}\tens b^\ce_{k-m,0}
\tens1^{\tens l}_{s\und\Ainfty(\ca;\uCom)})\mho_{ml}
+\sum_{m=0}^k(1^{\tens m}_{s\ca^\op}\tens b^\ce_{k-m,1}
\tens1^{\tens(l-1)}_{s\und\Ainfty(\ca;\uCom)})\mho_{m,l-1}\hfill
\\
+\sum_{a+u+c=k}
 (1_{s\ca^\op}^{\tens a}\tens b^{\ca^\op}_u\tens1_{s\ca^\op}^{\tens c}
\tens1_{s\ce}\tens1_{s\und\Ainfty(\ca;\uCom)}^{\tens l})\mho_{a+1+c,l}.
\end{multline*}
Write it in more detailed form using explicit
formulas~\eqref{equ-components-bQ} for components of \(b^\cq\):
\begin{multline*}
S\overset{\text{def}}=\sum_{m=1}^k\bigl[
\bT^ks\ca^\op(X_0,X_k)\tens
X_kf_0[1]\tens\bT^ls\und\Ainfty(\ca;\uCom)(f_0,f_l)
\rTTo^{\Yo_m\tens\mho_{k-m,l}}
\\
\hfill s\und\Ainfty(\ca;\uCom)(H^{X_0},H^{X_m})\tens s\und\Ainfty(\ca;\uCom)(H^{X_m},f_l)
\rTTo^{B_2}s\und\Ainfty(\ca;\uCom)(H^{X_0},f_l)
\bigr]\hskip\multlinegap
\\
\hskip\multlinegap\hphantom{S}+\bigl[ \bT^ks\ca^\op(X_0,X_k)\tens
X_kf_0[1]\tens\bT^ls\und\Ainfty(\ca;\uCom)(f_0,f_l)
\rTTo^{\mho_{kl}}\hfill
\\
\hfill s\und\Ainfty(\ca;\uCom)(H^{X_0},f_l)
\rTTo^{B_1}
s\und\Ainfty(\ca;\uCom)(H^{X_0},f_l)
\bigr]\hskip\multlinegap
\\
\hskip\multlinegap\hphantom{S}+\bigl[ \bT^ks\ca^\op(X_0,X_k)\tens
X_kf_0[1]\tens\bT^ls\und\Ainfty(\ca;\uCom)(f_0,f_l)
\rTTo^{\mho_{k,l-1}\tens1}\hfill
\\
\hfill s\und\Ainfty(\ca;\uCom)(H^{X_0},f_{l-1})\tens s\und\Ainfty(\ca;\uCom)(f_{l-1},f_l)
\rTTo^{B_2}
s\und\Ainfty(\ca;\uCom)(H^{X_0},f_l)
\bigr]\hskip\multlinegap
\\
\hskip\multlinegap -\sum_{p=1}^l\bigl[
\bT^ks\ca^\op(X_0,X_k)\tens
X_kf_0[1]\tens\bT^ls\und\Ainfty(\ca;\uCom)(f_0,f_l)
\rTTo^{1^{\tens k}\tens1\tens1^{\tens p-1}\tens B_1\tens1^{\tens l-p}}
\hfill
\\
\hfill
\bT^ks\ca^\op(X_0,X_k)\tens
X_kf_0[1]\tens\bT^ls\und\Ainfty(\ca;\uCom)(f_0,f_l)
\rTTo^{\mho_{kl}}
s\und\Ainfty(\ca;\uCom)(H^{X_0},f_l)
\bigr]\hskip\multlinegap
\\
-\sum_{p=1}^{l-1}\bigl[ \bT^ks\ca^\op(X_0,X_k)\tens
X_kf_0[1]\tens\bT^ls\und\Ainfty(\ca;\uCom)(f_0,f_l)
 \rTTo^{1^{\tens k}\tens1\tens1^{\tens p-1}\tens B_2\tens1^{\tens l-p-1}}
\\
\bT^ks\ca^\op(X_0,X_k)\tens X_kf_0[1]\tens
\bT^{l-1}s\und\Ainfty(\ca;\uCom)(f_0,\dots,f_{p-1},f_{p+1},\dots,f_l)
\\
\hfill\rTTo^{\mho_{k,l-1}}
s\und\Ainfty(\ca;\uCom)(H^{X_0},f_l)
\bigr]\hskip\multlinegap
\\
\hskip\multlinegap\hphantom{S}-\sum_{m=0}^k\bigl[
\bT^ks\ca^\op(X_0,X_k)\tens
X_kf_0[1]\tens\bT^ls\und\Ainfty(\ca;\uCom)(f_0,f_l)
\rTTo^{1^{\tens m}\tens b^\ce_{k-m,0}\tens1^{\tens l}}
\hfill
\\
\hfill
\bT^ms\ca^\op(X_0,X_m)\tens
X_mf_0[1]\tens\bT^ls\und\Ainfty(\ca;\uCom)(f_0,f_l)
\rTTo^{\mho_{ml}}
s\und\Ainfty(\ca;\uCom)(H^{X_0},f_l)
\bigr]
\\
\hskip\multlinegap\hphantom{S}-\sum_{m=0}^k\bigl[
\bT^ks\ca^\op(X_0,X_k)\tens
X_kf_0[1]\tens\bT^ls\und\Ainfty(\ca;\uCom)(f_0,f_l)
\rTTo^{1^{\tens m}\tens b^\ce_{k-m,1}\tens1^{\tens l-1}}
\hfill
\\
\bT^ms\ca^\op(X_0,X_m)\tens
X_mf_1[1]\tens\bT^{l-1}s\und\Ainfty(\ca;\uCom)(f_1,f_l)
\rTTo^{\mho_{m,l-1}}
s\und\Ainfty(\ca;\uCom)(H^{X_0},f_l)
\bigr]
\\
-\sum_{a+u+c=k}\bigl[ \bT^ks\ca^\op(X_0,X_k)\tens
X_kf_0[1]\tens\bT^ls\und\Ainfty(\ca;\uCom)(f_0,f_l)
 \rTTo^{1^{\tens a}\tens b^{\ca^\op}_u\tens1^{\tens c}\tens1\tens1^{\tens l}}
\\
\bT^{a+1+c}s\ca^\op(X_0,\dots,X_a,X_{a+u},\dots,X_k)
\tens X_kf_0[1]\tens\bT^ls\und\Ainfty(\ca;\uCom)(f_0,f_l)
\\
\rTTo^{\mho_{a+1+c,l}}
s\und\Ainfty(\ca;\uCom)(H^{X_0},f_l)
\bigr]=0.
\end{multline*}
The above equation is equivalent to the system of equations
\begin{multline*}
S\cdot\pr_n=0:\bT^ks\ca^\op(X_0,X_k)\tens
X_kf_0[1]\tens\bT^ls\und\Ainfty(\ca;\uCom)(f_0,f_l)
\\
\to s\uCom(\bT^ns\ca(Z_0,Z_n),s\uCom(\ca(X_0,Z_0),Z_nf_l)),
\end{multline*}
where \(n\ge0\), and \(Z_0,\dots,Z_n\in\Ob\ca\). By closedness, each of
these equations is equivalent to
\begin{multline*}
(1^{\tens n}\tens S\cdot\pr_n)\ev^\Com=(1^{\tens n}\tens
S)\ev^{\Ainfty}_{n1}=0:
\\
\bT^ns\ca(Z_0,Z_n)\tens\bT^ks\ca^\op(X_0,X_k)\tens
X_kf_0[1]\tens\bT^ls\und\Ainfty(\ca;\uCom)(f_0,f_l)
\\
\to s\uCom(\ca(X_0,Z_0),Z_nf_l).
\end{multline*}
The fact that \(\ev^{\Ainfty}\) is an \ainf-functor combined with
explicit formula~\eqref{eq-evA80-evA81} for components of
\(\ev^{\Ainfty}\) allows to derive certain identities. Specifically,
restricting the identity
 \([\ev^{\Ainfty}b^{\uCom}-(b^\ca\boxt1+1\boxt B)\ev^{\Ainfty}]\pr_1=0:
 Ts\ca\boxt Ts\und\Ainfty(\ca;\uCom)\to s\uCom\)
to the summand
 \(\bT^ns\ca(Z_0,Z_n)\tens s\und\Ainfty(\ca;\uCom)(\phi,\psi)
 \tens s\und\Ainfty(\ca;\uCom)(\psi,\chi)\)
yields
\begin{multline}
\bigl[
\bT^ns\ca(Z_0,Z_n)\tens s\und\Ainfty(\ca;\uCom)(\phi,\psi)
\tens s\und\Ainfty(\ca;\uCom)(\psi,\chi)
\\
\qquad\rTTo^{1^{\tens n}\tens B_2}
\bT^ns\ca(Z_0,Z_n)\tens s\und\Ainfty(\ca;\uCom)(\phi,\chi)
\rTTo^{1^{\tens n}\tens\pr_n}\hfill
\\
\hfill\bT^ns\ca(Z_0,Z_n)\tens\uCom(\bT^ns\ca(Z_0,Z_n),s\uCom(Z_0\phi,Z_n\chi))
\rTTo^{\ev^\Com}
s\uCom(Z_0\phi,Z_n\chi)
\bigr]\quad
\\
\quad=\bigl[
\bT^ns\ca(Z_0,Z_n)\tens s\und\Ainfty(\ca;\uCom)(\phi,\psi)
\tens s\und\Ainfty(\ca;\uCom)(\psi,\chi)\hfill
\\
\hfill\rTTo^{1^{\tens n}\tens B_2}
\bT^ns\ca(Z_0,Z_n)\tens s\und\Ainfty(\ca;\uCom)(\phi,\chi)
\rTTo^{\ev^{\Ainfty}_{n1}}
s\uCom(Z_0\phi,Z_n\chi)
\bigr]\quad
\\
\quad=\sum_{p+q=n}\bigl[
\bT^ns\ca(Z_0,Z_n)\tens s\und\Ainfty(\ca;\uCom)(\phi,\psi)
\tens s\und\Ainfty(\ca;\uCom)(\psi,\chi)\hfill
\\
\quad\rTTo^{\perm}
\bT^ps\ca(Z_0,Z_p)\tens s\und\Ainfty(\ca;\uCom)(\phi,\psi)
\tens\bT^qs\ca(Z_p,Z_n)\tens
s\und\Ainfty(\ca;\uCom)(\psi,\chi)\hfill
\\
\rTTo^{\ev^{\Ainfty}_{p1}\tens\ev^{\Ainfty}_{q1}}
s\uCom(Z_0\phi,Z_p\psi)\tens s\uCom(Z_p\psi,Z_n\chi)
\rTTo^{b^{\uCom}_2}
s\uCom(Z_0\phi,Z_n\chi)
\bigr].
\label{equ-1-B2-ev}
\end{multline}
Restricting the same identity to the summand
 \(\bT^ns\ca(Z_0,Z_n)\tens s\und\Ainfty(\ca;\uCom)(\phi,\psi)\) yields
\begin{multline}
\bigl[
\bT^ns\ca(Z_0,Z_n)\tens s\und\Ainfty(\ca;\uCom)(\phi,\psi)
\rTTo^{1^{\tens n}\tens B_1}
\bT^ns\ca(Z_0,Z_n)\tens s\und\Ainfty(\ca;\uCom)(\phi,\psi)
\\
\rTTo^{1^{\tens n}\tens\pr_n}
\bT^ns\ca(Z_0,Z_n)\tens\uCom(\bT^ns\ca(Z_0,Z_n),s\uCom(Z_0\phi,Z_n\psi))
\rTTo^{\ev^\Com}
s\uCom(Z_0\phi,Z_n\psi)
\bigr]
\\
\quad=\bigl[
\bT^ns\ca(Z_0,Z_n)\tens s\und\Ainfty(\ca;\uCom)(\phi,\psi)
\rTTo^{1^{\tens n}\tens B_1}\hfill
\\
\hfill\bT^ns\ca(Z_0,Z_n)\tens s\und\Ainfty(\ca;\uCom)(\phi,\psi)
\rTTo^{\ev^{\Ainfty}_{n1}}
s\uCom(Z_0\phi,Z_n\psi)
\bigr]\quad
\\
=\bigl[
\bT^ns\ca(Z_0,Z_n)\tens s\und\Ainfty(\ca;\uCom)(\phi,\psi)
\rTTo^{\ev^{\Ainfty}_{n1}}
s\uCom(Z_0\phi,Z_n\psi)
\rTTo^{b^{\uCom}_1}
s\uCom(Z_0\phi,Z_n\psi)
\bigr]
\\
\quad+\sum_{p+q=n}^{q>0}\bigl[
\bT^ns\ca(Z_0,Z_n)\tens s\und\Ainfty(\ca;\uCom)(\phi,\psi)
\rTTo^{\perm}\hfill
\\
\bT^ps\ca(Z_0,Z_p)\tens s\und\Ainfty(\ca;\uCom)(\phi,\psi)
\tens\bT^qs\ca(Z_p,Z_n)\tens
T^0s\und\Ainfty(\ca;\uCom)(\psi,\psi)
\\
\hfill\rTTo^{\ev^{\Ainfty}_{p1}\tens\ev^{\Ainfty}_{q0}}
s\uCom(Z_0\phi,Z_p\psi)\tens s\uCom(Z_p\psi,Z_n\psi)
\rTTo^{b^{\uCom}_2}
s\uCom(Z_0\phi,Z_n\psi)
\bigr]\quad
\\
\quad+\sum_{p+q=n}^{p>0}\bigl[
\bT^ns\ca(Z_0,Z_n)\tens s\und\Ainfty(\ca;\uCom)(\phi,\psi)
\rTTo^{\sim}\hfill
\\
\bT^ps\ca(Z_0,Z_p)\tens T^0s\und\Ainfty(\ca;\uCom)(\phi,\phi)
\tens\bT^qs\ca(Z_p,Z_n)\tens s\und\Ainfty(\ca;\uCom)(\phi,\psi)
\\
\hfill\rTTo^{\ev^{\Ainfty}_{p0}\tens\ev^{\Ainfty}_{q1}}
s\uCom(Z_0\phi,Z_p\phi)\tens s\uCom(Z_p\phi,Z_n\psi)
\rTTo^{b^{\uCom}_2}
s\uCom(Z_0\phi,Z_n\psi)
\bigr]\quad
\\
\quad-\sum_{\alpha+t+\beta=n}\bigl[
\bT^ns\ca(Z_0,Z_n)\tens s\und\Ainfty(\ca;\uCom)(\phi,\psi)
\rTTo^{1^{\tens\alpha}\tens b_t\tens1^{\tens\beta}\tens1}\hfill
\\
\bT^{\alpha+1+\beta}\ca(Z_0,\dots,Z_\alpha,Z_{\alpha+t},\dots,Z_n)\tens
 s\und\Ainfty(\ca;\uCom)(\phi,\psi)
\rTTo^{\ev^{\Ainfty}_{\alpha+1+\beta,1}}
s\uCom(Z_0\phi,Z_n\psi)
\bigr]
\\
\quad=\bigl[
\bT^ns\ca(Z_0,Z_n)\tens s\und\Ainfty(\ca;\uCom)(\phi,\psi)
\rTTo^{\ev^{\Ainfty}_{n1}}
s\uCom(Z_0\phi,Z_n\psi)
\rTTo^{b^{\uCom}_1}
s\uCom(Z_0\phi,Z_n\psi)
\bigr]\hfill
\\
\quad+\sum_{p+q=n}^{q>0}\bigl[
\bT^ns\ca(Z_0,Z_n)\tens s\und\Ainfty(\ca;\uCom)(\phi,\psi)
\rTTo^{\perm}\hfill
\\
\bT^ps\ca(Z_0,Z_p)\tens s\und\Ainfty(\ca;\uCom)(\phi,\psi)
\tens\bT^qs\ca(Z_p,Z_n)\tens
T^0s\und\Ainfty(\ca;\uCom)(\psi,\psi)
\\
\hfill\rTTo^{\ev^{\Ainfty}_{p1}\tens\ev^{\Ainfty}_{q0}}
s\uCom(Z_0\phi,Z_p\psi)\tens s\uCom(Z_p\psi,Z_n\psi)
\rTTo^{b^{\uCom}_2}
s\uCom(Z_0\phi,Z_n\psi)
\bigr]\quad
\\
\quad+\sum_{p+q=n}^{p>0}\bigl[
\bT^ns\ca(Z_0,Z_n)\tens s\und\Ainfty(\ca;\uCom)(\phi,\psi)\hfill
\\
\hfill\rTTo^{\phi_p\tens\ev^{\Ainfty}_{q1}}
s\uCom(Z_0\phi,Z_p\phi)\tens s\uCom(Z_p\phi,Z_n\psi)
\rTTo^{b^{\uCom}_2}
s\uCom(Z_0\phi,Z_n\psi)
\bigr]\quad
\\
\quad-\sum_{\alpha+t+\beta=n}\bigl[
\bT^ns\ca(Z_0,Z_n)\tens s\und\Ainfty(\ca;\uCom)(\phi,\psi)
\rTTo^{1^{\tens\alpha}\tens b_t\tens1^{\tens\beta}\tens1}\hfill
\\
\bT^{\alpha+1+\beta}s\ca(Z_0,\dots,Z_\alpha,Z_{\alpha+t},\dots,Z_n)\tens
s\und\Ainfty(\ca;\uCom)(\phi,\psi)
\\
\rTTo^{\ev^{\Ainfty}_{\alpha+1+\beta,1}}
s\uCom(Z_0\phi,Z_n\psi)
\bigr].
\label{equ-1-B1-ev}
\end{multline}
With identities~\eqref{equ-1-B2-ev} and \eqref{equ-1-B1-ev} in hand, it
is the matter of straightforward verification to check that
 \((1^{\tens n}\tens S\cdot\pr_n)\ev^\Com\) admits the following
presentation:
\begin{align}
&\sum_{m=1}^k\sum_{p+q=n}\bigl[
\bT^ns\ca(Z_0,Z_n)\tens\bT^ks\ca^\op(X_0,X_k)\tens
X_kf_0[1]\tens\bT^ls\und\Ainfty(\ca;\uCom)(f_0,f_l)\nonumber
\\
&\hspace{2em}\rTTo^{\perm}
\bT^ps\ca(Z_0,Z_p)\tens\bT^ms\ca^\op(X_0,X_m)\tens\bT^qs\ca(Z_p,Z_n)\nonumber
\\
&\hspace{7em}\tens\bT^{k-m}s\ca^\op(X_m,X_k)\tens X_kf_0[1]\tens
\bT^ls\und\Ainfty(\ca;\uCom)(f_0,f_l)\nonumber
\\
&\hspace{1em}\rTTo^{(\Hom_{\ca^\op})_{pm}\tens\mho'_{k-m,l;q}}
s\uCom(\ca(X_0,Z_0),\ca(X_m,Z_p))\tens
s\uCom(\ca(X_m,Z_p),Z_nf_l)\nonumber
\\
&\hspace{18em}\rTTo^{b^{\uCom}_2} s\uCom(\ca(X_0,Z_0),Z_nf_l) \bigr]
\label{equ-Hom-Mho-b2}
\\
&+\bigl[ \bT^ns\ca(Z_0,Z_n)\tens\bT^ks\ca^\op(X_0,X_k)\tens
X_kf_0[1]\tens \bT^ls\und\Ainfty(\ca;\uCom)(f_0,f_l)\nonumber
\\
&\hspace{6em}\rTTo^{\mho'_{kl;n}} s\uCom(\ca(X_0,Z_0),Z_nf_l)
\rTTo^{b^{\uCom}_1} s\uCom(\ca(X_0,Z_0),Z_nf_l) \bigr]
\label{equ-Mho-b1}
\\
&+\sum_{p+q=n}^{q>0}\bigl[
\bT^ns\ca(Z_0,Z_n)\tens\bT^ks\ca^\op(X_0,X_k)\tens
X_kf_0[1]\tens \bT^ls\und\Ainfty(\ca;\uCom)(f_0,f_l)\nonumber
\\
&\hspace{2em}\rTTo^{\perm}
\bT^ps\ca(Z_0,Z_p)\tens\bT^ks\ca^\op(X_0,X_k)\tens X_kf_0[1]\nonumber
\\
&\hspace{5em}\tens
\bT^ls\und\Ainfty(\ca;\uCom)(f_0,f_l)\tens\bT^qs\ca(Z_p,Z_n)\tens T^0s\und\Ainfty(\ca;\uCom)(f_l,f_l)\nonumber
\\
&\rTTo^{\mho'_{kl;p}\tens\ev^{\Ainfty}_{q0}}
s\uCom(\ca(X_0,Z_0),Z_pf_l)\tens s\uCom(Z_pf_l,Z_nf_l)
\rTTo^{b^{\uCom}_2} s\uCom(\ca(X_0,Z_0),Z_nf_l) \bigl]
\label{equ-Mho-evA8-b2}
\\
&+\sum_{p+q=n}^{p>0}\bigl[
\bT^ns\ca(Z_0,Z_n)\tens\bT^ks\ca^\op(X_0,X_k)\tens
X_kf_0[1]\tens \bT^ls\und\Ainfty(\ca;\uCom)(f_0,f_l)
\nonumber
\\
&\hspace{3em} \rTTo^{H^{X_0}_p\tens\mho'_{kl;q}}
s\uCom(\ca(X_0,Z_0),\ca(X_0,Z_p))\tens s\uCom(\ca(X_0,Z_p),Z_nf_l)
\nonumber
\\
&\hspace{17em} \rTTo^{b^{\uCom}_2} s\uCom(\ca(X_0,Z_0),Z_nf_l) \bigr]
\label{equ-H-Mho-b2}
\\
&-\sum_{\alpha+t+\beta=n}\bigl[
\bT^ns\ca(Z_0,Z_n)\tens\bT^ks\ca^\op(X_0,X_k)
\tens X_kf_0[1]\tens \bT^ls\und\Ainfty(\ca;\uCom)(f_0,f_l)\nonumber
\\
&\hspace{2em}
\rTTo^{1^{\tens\alpha}\tens b_t\tens1^{\tens\beta}\tens1^{\tens k}\tens1\tens1^{\tens l}}
\bT^{\alpha+1+\beta}s\ca(Z_0,\dots,Z_\alpha,Z_{\alpha+t},\dots,Z_n)
\nonumber
\\
&\hspace{6em}
\tens\bT^ks\ca^\op(X_0,X_k)\tens X_kf_0[1]\tens
\bT^ls\und\Ainfty(\ca;\uCom)(f_0,f_l)\nonumber
\\
&\hspace{16em}\rTTo^{\mho'_{kl;\alpha+1+\beta}}
s\uCom(\ca(X_0,Z_0),Z_nf_l) \bigl]
 \label{equ-1-b-1-1-1-1-Mho}
\\
&+\sum_{p+q=n}\bigl[
\bT^ns\ca(Z_0,Z_n)\tens\bT^ks\ca^\op(X_0,X_k)\tens
X_kf_0[1]\tens \bT^ls\und\Ainfty(\ca;\uCom)(f_0,f_l)\nonumber
\\
&\hspace{2em}\rTTo^{\perm}
\bT^ps\ca(Z_0,Z_p)\tens\bT^ks\ca^\op(X_0,X_k)
\tens X_kf_0[1]\nonumber
\\
&\hspace{4em}\tens\bT^{l-1}s\und\Ainfty(\ca;\uCom)(f_0,f_{l-1})
\tens\bT^qs\ca(Z_p,Z_n)\tens
s\und\Ainfty(\ca;\uCom)(f_{l-1},f_l)\nonumber
\\
&\hspace{2em}\rTTo^{\mho'_{k,l-1;p}\tens\ev^{\Ainfty}_{q1}}
s\uCom(\ca(X_0,Z_0),Z_pf_{l-1})\tens s\uCom(Z_pf_{l-1},Z_nf_l)\nonumber
\\
&\hspace{18em}\rTTo^{b^{\uCom}_2} s\uCom(\ca(X_0,Z_0),Z_nf_l) \bigl]
\label{equ-Mho-evA8-b2-15}
\\
&-\sum_{p=1}^{l}\bigl[
\bT^ns\ca(Z_0,Z_n)\tens\bT^ks\ca^\op(X_0,X_k)\tens
X_kf_0[1]\tens \bT^ls\und\Ainfty(\ca;\uCom)(f_0,f_l)\nonumber
\\
&\hspace{12em}
\rTTo^{1^{\tens n}\tens1^{\tens k}\tens1\tens1^{\tens p-1}\tens B_1\tens1^{\tens l-p}}
\nonumber
\\
&\hspace{4em} \bT^ns\ca(Z_0,Z_n)\tens\bT^ks\ca^\op(X_0,X_k)\tens
X_kf_0[1]\tens\bT^ls\und\Ainfty(\ca;\uCom)(f_0,f_l)\nonumber
\\
&\hspace{18em}\rTTo^{\mho'_{kl;n}} s\uCom(\ca(X_0,Z_0),Z_nf_l) \bigl]
\label{equ-1111-B1-1-Mho}
\\
&-\sum_{p=1}^{l-1}\bigl[
\bT^ns\ca(Z_0,Z_n)\tens\bT^ks\ca^\op(X_0,X_k)\tens
X_kf_0[1]\tens \bT^ls\und\Ainfty(\ca;\uCom)(f_0,f_l)\nonumber
\\
&\hspace{2em}
\rTTo^{1^{\tens n}\tens1^{\tens k}\tens1\tens1^{\tens p-1}\tens B_2\tens1^{\tens l-p-1}}
\bT^ns\ca(Z_0,Z_n)\tens\bT^ks\ca^\op(X_0,X_k)
\nonumber
\\
&\hspace{5em} \tens X_kf_0[1]\tens
\bT^{l-1}s\und\Ainfty(\ca;\uCom)(f_0,\dots,f_{p-1},f_{p+1},\dots,f_l)
\nonumber
\\
&\hspace{17em}\rTTo^{\mho'_{k,l-1;n}} s\uCom(\ca(X_0,Z_0),Z_nf_l)
\bigl]
 \label{equ-1111-B2-1-Mho}
\\
&-\sum_{m=0}^{k}\bigl[
\bT^ns\ca(Z_0,Z_n)\tens\bT^ks\ca^\op(X_0,X_k)\tens
X_kf_0[1]\tens \bT^ls\und\Ainfty(\ca;\uCom)(f_0,f_l)\nonumber
\\
&\hspace{2em}
\rTTo^{1^{\tens n}\tens1^{\tens m}\tens b^\ce_{k-m,0}\tens1^{\tens l}}
\bT^ns\ca(Z_0,Z_n)\tens\bT^ms\ca^\op(X_0,X_m)
\nonumber
\\
&\hspace{3em} \tens
X_mf_0[1]\tens\bT^ls\und\Ainfty(\ca;\uCom)(f_0,f_l)
\rTTo^{\mho'_{ml;n}} s\uCom(\ca(X_0,Z_0),Z_nf_l) \bigr]
\label{equ-11-bE0-1-Mho}
\\
&-\sum_{m=0}^{k}\bigl[
\bT^ns\ca(Z_0,Z_n)\tens\bT^ks\ca^\op(X_0,X_k)\tens
X_kf_0[1]\tens \bT^ls\und\Ainfty(\ca;\uCom)(f_0,f_l)\nonumber
\\
&\hspace{12em}
\rTTo^{1^{\tens n}\tens1^{\tens m}\tens b^\ce_{k-m,1}\tens1^{\tens l-1}}
\nonumber
\\
&\hspace{2em}
\bT^ns\ca(Z_0,Z_n)\tens\bT^ms\ca^\op(X_0,X_m)\tens
X_mf_1[1]\tens\bT^{l-1}s\und\Ainfty(\ca;\uCom)(f_1,f_l)\nonumber
\\
&\hspace{16em}\rTTo^{\mho'_{m,l-1;n}} s\uCom(\ca(X_0,Z_0),Z_nf_l)
\bigr]
 \label{equ-11-bE1-1-Mho}
\\
&-\sum_{a+u+c=k}\bigl[
\bT^ns\ca(Z_0,Z_n)\tens\bT^ks\ca^\op(X_0,X_k)\tens
X_kf_0[1]\tens\bT^ls\und\Ainfty(\ca;\uCom)(f_0,f_l)\nonumber
\\
&\hspace{1em}
\rTTo^{1^{\tens n}\tens1^{\tens a}\tens b^{\ca^\op}_u\tens1^{\tens c}\tens1\tens1^{\tens l}}
\bT^ns\ca(Z_0,Z_n)
\nonumber
\\
&\hspace{1em} \tens\bT^{a+1+c}s\ca^\op(X_0,\dots,X_a,X_{a+u},\dots,X_k)
\tens X_kf_0[1]\tens\bT^ls\und\Ainfty(\ca;\uCom)(f_0,f_l)\nonumber
\\
&\hspace{15em}\rTTo^{\mho'_{a+1+c,l;n}} s\uCom(\ca(X_0,Z_0),Z_nf_l)
\bigl].
 \label{equ-11-bop-111-Mho}
\end{align}
Appearance of the component \((\Hom_{\ca^\op})_{pm}\) in
term~\eqref{equ-Hom-Mho-b2} is explained by the identity
 \((\Hom_{\ca^\op})_{pm}=((1,\Yo)\ev^{\Ainfty})_{pm}
 =(1^{\tens p}\tens\Yo_m)\ev^{\Ainfty}_{p1}\),
which holds true by definition of the Yoneda \ainf-functor
\(\Yo:\ca^\op\to\und\Ainfty(\ca;\uCom)\). Expanding
\((\Hom_{\ca^\op})_{pm}\) according to
formula~\eqref{equ-ainf-Hom-components}, term~\eqref{equ-Hom-Mho-b2}
can be written as follows:
\begin{multline*}
-(-)^{k+1}\bigl[
\bT^ns\ca(Z_0,Z_n)\tens\bT^ks\ca^\op(X_0,X_k)\tens X_kf_0[1]\tens
\bT^ls\und\Ainfty(\ca;\uCom)(f_0,f_l)
\\
\quad\rTTo^{\text{perm}}
\bT^ps\ca(Z_0,Z_p)\tens\bT^ms\ca^\op(X_0,X_m)\tens\bT^qs\ca(Z_p,Z_n)\hfill
\\
\hfill
\tens\bT^{k-m}s\ca^\op(X_m,X_k)\tens X_kf_0[1]\tens
\bT^ls\und\Ainfty(\ca;\uCom)(f_0,f_l)\quad
\\
\rTTo^{\coev^\Com\tens\coev^\Com}
\uCom(s\ca(X_0,Z_0),s\ca(X_0,Z_0)\tens\bT^ps\ca(Z_0,Z_p)\tens
\bT^ms\ca^\op(X_0,X_m))
\\
\tens\uCom(s\ca(X_m,Z_p),s\ca(X_m,Z_p)\tens\bT^qs\ca(Z_p,Z_n)
\\
\hfill\tens\bT^{k-m}s\ca^\op(X_m,X_k)
\tens X_kf_0[1]\tens\bT^ls\und\Ainfty(\ca;\uCom)(f_0,f_l))\quad
\\
\rTTo^{\uCom(1,\text{perm})\tens\uCom(1,\text{perm})}
\uCom(s\ca(X_0,Z_0),\bT^ms\ca(X_m,X_0)\tens s\ca(X_0,Z_0)\tens
\bT^ps\ca(Z_0,Z_p))
\\
\tens\uCom(s\ca(X_m,Z_p),X_kf_0[1]\tens\bT^{k-m}s\ca(X_k,X_m)
\\
\hfill\tens\ca(X_m,Z_p)\tens\bT^qs\ca(Z_p,Z_n)
\tens\bT^ls\und\Ainfty(\ca;\uCom)(f_0,f_l))\quad
\\
\quad\rTTo^{\uCom(1,b_{m+1+p})\tens\uCom(1,1\tens\ev^{\Ainfty}_{k-m+1+q,l})}
\uCom(s\ca(X_0,Z_0),s\ca(X_m,Z_p))\hfill
\\
\hfill\tens\uCom(s\ca(X_m,Z_p),X_kf_0[1]\tens s\uCom(X_kf_0,Z_nf_l))\quad
\\
\quad\rTTo^{1\tens\uCom(1,1\tens s^{-1}[1])}
\uCom(s\ca(X_0,Z_0),s\ca(X_m,Z_p))\hfill
\\
\hfill\tens\uCom(s\ca(X_m,Z_p),X_kf_0[1]\tens\uCom(X_kf_0[1],Z_nf_l[1]))\quad
\\
\rTTo^{1\tens\uCom(1,\ev^\Com)}
\uCom(s\ca(X_0,Z_0),s\ca(X_m,Z_p))\tens\uCom(s\ca(X_m,Z_p),Z_nf_l[1])
\rTTo^{[-1]s\tens[-1]s}
\\
s\uCom(\ca(X_0,Z_0),\ca(X_m,Z_p))\tens s\uCom(\ca(X_m,Z_p),Z_nf_l)
\rTTo^{b^{\uCom}_2}
s\uCom(\ca(X_0,Z_0),Z_nf_l)
\bigr].
\end{multline*}
Replacing the last two arrows by the composite \(m^{\uCom}_2[-1]s\) and
applying identity~\eqref{eq-identity-m2} yields:
\begin{multline*}
-(-)^{k+1}\bigl[
\bT^ns\ca(Z_0,Z_n)\tens \bT^ks\ca^\op(X_0,X_k)\tens X_kf_0[1]\tens
\bT^ls\und\Ainfty(\ca;\uCom)(f_0,f_l)
\\
\quad\rTTo^{\coev^\Com}
\uCom(s\ca(X_0,Z_0),s\ca(X_0,Z_0)\tens \bT^ns\ca(Z_0,Z_n)\hfill
\\
\hfill\tens \bT^ks\ca^\op(X_0,X_k)\tens X_kf_0[1]\tens
\bT^ls\und\Ainfty(\ca;\uCom)(f_0,f_l))\quad
\\
\quad\rTTo^{\uCom(1,\text{perm})}
\uCom(s\ca(X_0,Z_0),X_kf_0[1]\tens \bT^ks\ca(X_k,X_0)\tens
s\ca(X_0,Z_0)\hfill
\\
\hfill\tens \bT^ns\ca(Z_0,Z_n)\tens
\bT^ls\und\Ainfty(\ca;\uCom)(f_0,f_l))\quad
\\
\quad
\rTTo^{\uCom(1,1\tens(1^{\tens(k-m)}\tens b_{m+1+p}\tens
 1^{\tens q}\tens1^{\tens l})\ev^{\Ainfty}_{k-m+1+q,l})}
\hfill
\\
\uCom(s\ca(X_0,Z_0),X_kf_0[1]\tens s\uCom(X_kf_0,Z_nf_l))
\\
\quad\rTTo^{\uCom(1,1\tens s^{-1}[1])}
\uCom(s\ca(X_0,Z_0),X_kf_0[1]\tens\uCom(X_kf_0[1],Z_nf_l[1]))\hfill
\\
\rTTo^{\uCom(1,\ev^\Com)}
\uCom(s\ca(X_0,Z_0),Z_nf_l[1])
\rTTo^{[-1]s}
s\uCom(\ca(X_0,Z_0),Z_nf_l)
\bigr].
\end{multline*}
Denote
\(d'=s^{-1}ds=s^{-1}d^{X_kf_0}s=b^\ce_{00}:X_kf_0[1]=s\ce(X_k,f_0)\to
X_kf_0[1]=s\ce(X_k,f_0)\). Thus, the shifted complex \(X_kf_0[1]\)
carries the differential \(-d'\). Since
 \(m^{\uCom}_1=-\uCom(1,d')+\uCom(b_1,1):\uCom(s\ca(X_0,Z_0),Z_nf_l[1])
 \to\uCom(s\ca(X_0,Z_0),Z_nf_l[1])\),
it follows that term~\eqref{equ-Mho-b1} equals
\begin{multline*}
(-)^{k+1}\bigl[
\bT^ns\ca(Z_0,Z_n)\tens\bT^ks\ca^\op(X_0,X_k)\tens X_kf_0[1]\tens
\bT^ls\und\Ainfty(\ca;\uCom)(f_0,f_l)
\\
\quad\rTTo^{\coev^\Com}
\uCom(s\ca(X_0,Z_0),s\ca(X_0,Z_0)\tens\bT^ns\ca(Z_0,Z_n)\hfill
\\
\hfill\tens\bT^ks\ca^\op(X_0,X_k)\tens X_kf_0[1]\tens
\bT^ls\und\Ainfty(\ca;\uCom)(f_0,f_l))\quad
\\
\quad\rTTo^{\uCom(1,\text{perm})}
\uCom(s\ca(X_0,Z_0),X_kf_0[1]\tens\bT^ks\ca(X_k,X_0)\tens
s\ca(X_0,Z_0)\hfill
\\
\hfill\tens\bT^ns\ca(Z_0,Z_n)\tens
\bT^ls\und\Ainfty(\ca;\uCom)(f_0,f_l))\quad
\\
\quad\rTTo^{\uCom(1,1\tens\ev^{\Ainfty}_{k+1+n,l})}
\uCom(s\ca(X_0,Z_0),X_kf_0[1]\tens s\uCom(X_kf_0,Z_nf_l))\hfill
\\
\quad\rTTo^{\uCom(1,1\tens s^{-1}[1])}
\uCom(s\ca(X_0,Z_0),X_kf_0[1]\tens\uCom(X_kf_0[1],Z_nf_l[1]))
\rTTo^{\uCom(1,\ev^\Com)}\hfill
\\
\uCom(s\ca(X_0,Z_0),Z_nf_l[1])
\rTTo^{m^{\uCom}_1}
\uCom(s\ca(X_0,Z_0),Z_nf_l[1])
\rTTo^{[-1]s}
s\uCom(\ca(X_0,Z_0),Z_nf_l)
\bigr]
\\
\hskip\multlinegap=(-)^{k+1}\bigl[
\bT^ns\ca(Z_0,Z_n)\tens\bT^ks\ca^\op(X_0,X_k)\tens X_kf_0[1]\tens
\bT^ls\und\Ainfty(\ca;\uCom)(f_0,f_l)\hfill
\\
\quad\rTTo^{\coev^\Com}
\uCom(s\ca(X_0,Z_0),s\ca(X_0,Z_0)\tens\bT^ns\ca(Z_0,Z_n)\hfill
\\
\hfill\tens\bT^ks\ca^\op(X_0,X_k)\tens X_kf_0[1]\tens
\bT^ls\und\Ainfty(\ca;\uCom)(f_0,f_l))\quad
\\
\quad\rTTo^{\uCom(1,\text{perm})}
\uCom(s\ca(X_0,Z_0),X_kf_0[1]\tens\bT^ks\ca(X_k,X_0)\tens
s\ca(X_0,Z_0)\hfill
\\
\hfill\tens\bT^ns\ca(Z_0,Z_n)\tens
\bT^ls\und\Ainfty(\ca;\uCom)(f_0,f_l))\quad
\\
\quad\rTTo^{\uCom(1,1\tens\ev^{\Ainfty}_{k+1+n,l})}
\uCom(s\ca(X_0,Z_0),X_kf_0[1]\tens s\uCom(X_kf_0,Z_nf_l))\hfill
\\
\quad\rTTo^{\uCom(1,1\tens s^{-1}[1])}
\uCom(s\ca(X_0,Z_0),X_kf_0[1]\tens\uCom(X_kf_0[1],Z_nf_l[1]))\hfill
\\
\hfill\rTTo^{\uCom(1,-\ev^\Com\cdot d')}
\uCom(s\ca(X_0,Z_0),Z_nf_l[1])
\rTTo^{[-1]s}
s\uCom(\ca(X_0,Z_0),Z_nf_l)
\bigr]\quad
\\
\hskip\multlinegap-(-)^{k+1}\bigl[
\bT^ns\ca(Z_0,Z_n)\tens\bT^ks\ca^\op(X_0,X_k)\tens X_kf_0[1]\tens
\bT^ls\und\Ainfty(\ca;\uCom)(f_0,f_l)\hfill
\\
\quad\rTTo^{\coev^\Com}
\uCom(s\ca(X_0,Z_0),s\ca(X_0,Z_0)\tens\bT^ns\ca(Z_0,Z_n)\hfill
\\
\hfill\tens\bT^ks\ca^\op(X_0,X_k)\tens X_kf_0[1]\tens
\bT^ls\und\Ainfty(\ca;\uCom)(f_0,f_l))\quad
\\
\quad\rTTo^{\uCom(1,\text{perm})}
\uCom(s\ca(X_0,Z_0),X_kf_0[1]\tens\bT^ks\ca(X_k,X_0)\tens
s\ca(X_0,Z_0)\hfill
\\
\hfill\tens\bT^ns\ca(Z_0,Z_n)\tens
\bT^ls\und\Ainfty(\ca;\uCom)(f_0,f_l))\quad
\\
\quad\rTTo^{\uCom(1,1\tens(1^{\tens k}\tens b_1\tens 1^{\tens n}\tens1^{\tens l})\ev^{\Ainfty}_{k+1+n,l})}
\uCom(s\ca(X_0,Z_0),X_kf_0[1]\tens s\uCom(X_kf_0,Z_nf_l))\hfill
\\
\quad\rTTo^{\uCom(1,1\tens s^{-1}[1])}
\uCom(s\ca(X_0,Z_0),X_kf_0[1]\tens\uCom(X_kf_0[1],Z_nf_l[1]))\hfill
\\
\rTTo^{\uCom(1,\ev^\Com)}
\uCom(s\ca(X_0,Z_0),Z_nf_l[1])
\rTTo^{[-1]s}
s\uCom(\ca(X_0,Z_0),Z_nf_l)
\bigr].
\end{multline*}
Since \(\ev^\Com\) is a chain map, it follows that
\[
-\ev^\Com\cdot d'=-(d'\tens1)\ev^\Com+(1\tens
m^{\uCom}_1)\ev^\Com:X_kf_0[1]\tens\uCom(X_kf_0[1],Z_nf_l[1])\to Z_nf_l[1],
\]
therefore term~\eqref{equ-Mho-b1} equals
\begin{multline}
(-)^{k+1}\bigl[
\bT^ns\ca(Z_0,Z_n)\tens\bT^ks\ca^\op(X_0,X_k)\tens X_kf_0[1]\tens
\bT^ls\und\Ainfty(\ca;\uCom)(f_0,f_l)
\\
\quad\rTTo^{\coev^\Com}
\uCom(s\ca(X_0,Z_0),s\ca(X_0,Z_0)\tens\bT^ns\ca(Z_0,Z_n)\hfill
\\
\hfill\tens\bT^ks\ca^\op(X_0,X_k)\tens X_kf_0[1]\tens
\bT^ls\und\Ainfty(\ca;\uCom)(f_0,f_l))\quad
\\
\quad\rTTo^{\uCom(1,\text{perm})}
\uCom(s\ca(X_0,Z_0),X_kf_0[1]\tens\bT^ks\ca(X_k,X_0)\tens
s\ca(X_0,Z_0)\hfill
\\
\hfill\tens\bT^ns\ca(Z_0,Z_n)\tens
\bT^ls\und\Ainfty(\ca;\uCom)(f_0,f_l))\quad
\\
\quad\rTTo^{\uCom(1,1\tens\ev^{\Ainfty}_{k+1+n,l}b^{\uCom}_1)}
\uCom(s\ca(X_0,Z_0),X_kf_0[1]\tens s\uCom(X_kf_0,Z_nf_l))\hfill
\\
\quad\rTTo^{\uCom(1,1\tens s^{-1}[1])}
\uCom(s\ca(X_0,Z_0),X_kf_0[1]\tens\uCom(X_kf_0[1],Z_nf_l[1]))\hfill
\\
\hfill\rTTo^{\uCom(1,\ev^\Com)}
\uCom(s\ca(X_0,Z_0),Z_nf_l[1])
\rTTo^{[-1]s}
s\uCom(\ca(X_0,Z_0),Z_nf_l)
\bigr]\quad
\\
\hskip\multlinegap+(-)^{k+1}\bigl[
\bT^ns\ca(Z_0,Z_n)\tens\bT^ks\ca^\op(X_0,X_k)\tens X_kf_0[1]\tens
\bT^ls\und\Ainfty(\ca;\uCom)(f_0,f_l)\hfill
\\
\quad\rTTo^{\coev^\Com}
\uCom(s\ca(X_0,Z_0),s\ca(X_0,Z_0)\tens\bT^ns\ca(Z_0,Z_n)\hfill
\\
\hfill\tens\bT^ks\ca^\op(X_0,X_k)\tens X_kf_0[1]\tens
\bT^ls\und\Ainfty(\ca;\uCom)(f_0,f_l))\quad
\\
\quad\rTTo^{\uCom(1,\text{perm})}
\uCom(s\ca(X_0,Z_0),X_kf_0[1]\tens\bT^ks\ca(X_k,X_0)\tens
s\ca(X_0,Z_0)\hfill
\\
\hfill\tens\bT^ns\ca(Z_0,Z_n)\tens
\bT^ls\und\Ainfty(\ca;\uCom)(f_0,f_l))\quad
\\
\quad\rTTo^{\uCom(1,d'\tens\ev^{\Ainfty}_{k+1+n,l})}
\uCom(s\ca(X_0,Z_0),X_kf_0[1]\tens s\uCom(X_kf_0,Z_nf_l))\hfill
\\
\quad\rTTo^{\uCom(1,1\tens s^{-1}[1])}
\uCom(s\ca(X_0,Z_0),X_kf_0[1]\tens\uCom(X_kf_0[1],Z_nf_l[1]))\hfill
\\
\hfill\rTTo^{\uCom(1,\ev^\Com)}
\uCom(s\ca(X_0,Z_0),Z_nf_l[1])
\rTTo^{[-1]s}
s\uCom(\ca(X_0,Z_0),Z_nf_l)
\bigr]\quad
\\
\hskip\multlinegap-(-)^{k+1}\bigl[
\bT^ns\ca(Z_0,Z_n)\tens\bT^ks\ca^\op(X_0,X_k)\tens X_kf_0[1]\tens
\bT^ls\und\Ainfty(\ca;\uCom)(f_0,f_l)\hfill
\\
\quad\rTTo^{\coev^\Com}
\uCom(s\ca(X_0,Z_0),s\ca(X_0,Z_0)\tens\bT^ns\ca(Z_0,Z_n)\hfill
\\
\hfill\tens\bT^ks\ca^\op(X_0,X_k)\tens X_kf_0[1]\tens
\bT^ls\und\Ainfty(\ca;\uCom)(f_0,f_l))\quad
\\
\quad\rTTo^{\uCom(1,\text{perm})}
\uCom(s\ca(X_0,Z_0),X_kf_0[1]\tens\bT^ks\ca(X_k,X_0)\tens
s\ca(X_0,Z_0)\hfill
\\
\hfill\tens\bT^ns\ca(Z_0,Z_n)\tens
\bT^ls\und\Ainfty(\ca;\uCom)(f_0,f_l))\quad
\\
\quad\rTTo^{\uCom(1,1\tens(1^{\tens k}\tens b_1\tens 1^{\tens n}\tens1^{\tens l})\ev^{\Ainfty}_{k+1+n,l})}
\uCom(s\ca(X_0,Z_0),X_kf_0[1]\tens s\uCom(X_kf_0,Z_nf_l))\hfill
\\
\quad\rTTo^{\uCom(1,1\tens s^{-1}[1])}
\uCom(s\ca(X_0,Z_0),X_kf_0[1]\tens\uCom(X_kf_0[1],Z_nf_l[1]))\hfill
\\
\rTTo^{\uCom(1,\ev^\Com)}
\uCom(s\ca(X_0,Z_0),Z_nf_l[1])
\rTTo^{[-1]s}
s\uCom(\ca(X_0,Z_0),Z_nf_l)
\bigr].
\label{equ*-ev-b1}
\end{multline}
Using the identity
\[
f=\bigl[X\rTTo^{\coev^\uCom}\uCom(Y,Y\tens X)\rTTo^{\uCom(1,1\tens
f)}\uCom(Y,Y\tens\uCom(Y,Z))\rTTo^{\uCom(1,\ev^\Com)}\uCom(Y,Z)
\bigr]
\]
valid for an  arbitrary \(f\in\uCom(X,\uCom(Y,Z))\) by general
properties of closed monoidal categories, term~\eqref{equ-Mho-evA8-b2}
can be written as follows:
\begin{multline}
-(-)^{k+1}\bigl[
\bT^ns\ca(Z_0,Z_n)\tens\bT^ks\ca^\op(X_0,X_k)\tens X_kf_0[1]\tens
\bT^ls\und\Ainfty(\ca;\uCom)(f_0,f_l)
\\
\quad\rTTo^{\text{perm}}
\bT^ps\ca(Z_0,Z_p)\tens\bT^ks\ca^\op(X_0,X_k)\tens X_kf_0[1]
\tens\bT^ls\und\Ainfty(\ca;\uCom)(f_0,f_l)\hfill
\\
\hfill \tens\bT^qs\ca(Z_p,Z_n)\tens T^0s\und\Ainfty(\ca;\uCom)(f_l,f_l)
\rTTo^{\coev^\Com\tens\coev^\Com} \quad
\\
\quad
\uCom(s\ca(X_0,Z_0),s\ca(X_0,Z_0)\tens\bT^ps\ca(Z_0,Z_p)
\tens\bT^ks\ca^\op(X_0,X_k)\tens X_kf_0[1]
\hfill
\\
\tens\bT^ls\und\Ainfty(\ca;\uCom)(f_0,f_l))
\tens\uCom(Z_pf_l[1],Z_pf_l[1]\tens\bT^qs\ca(Z_p,Z_n)\tens
T^0s\und\Ainfty(\ca;\uCom)(f_l,f_l))
\\
\rTTo^{\uCom(1,\text{perm})\tens1}
\uCom(s\ca(X_0,Z_0),X_kf_0[1]\tens\bT^ks\ca(X_k,X_0)\tens
s\ca(X_0,Z_0)\tens\bT^ps\ca(Z_0,Z_p)
\\
\tens\bT^ls\und\Ainfty(\ca;\uCom)(f_0,f_l))
\tens\uCom(Z_pf_l[1],Z_pf_l[1]\tens\bT^qs\ca(Z_p,Z_n)\tens
T^0s\und\Ainfty(\ca;\uCom)(f_l,f_l))
\\
\quad\rTTo^{\uCom(1,1\tens\ev^{\Ainfty}_{k+1+p,l})\tens\uCom(1,1\tens\ev^{\Ainfty}_{q0})}
\uCom(s\ca(X_0,Z_0),X_kf_0[1]\tens s\uCom(X_kf_0,Z_pf_l))\hfill
\\
\hfill\tens\uCom(Z_pf_l[1],Z_pf_l[1]\tens s\uCom(Z_pf_l,Z_nf_l))\quad
\\
\quad\rTTo^{\uCom(1,1\tens s^{-1}[1])\tens\uCom(1,1\tens s^{-1}[1])}
\uCom(s\ca(X_0,Z_0),X_kf_0[1]\tens\uCom(X_kf_0[1],Z_pf_l[1]))\hfill
\\
\hfill\tens\uCom(Z_pf_l[1],Z_pf_l[1]\tens\uCom(Z_pf_l[1],Z_nf_l[1]))\quad
\\
\quad\rTTo^{\uCom(1,\ev^\Com)\tens\uCom(1,\ev^\Com)}
\uCom(s\ca(X_0,Z_0),Z_pf_l[1])\tens\uCom(Z_pf_l[1],Z_nf_l[1])\rTTo^{[-1]s\tens[-1]s}\hfill
\\
s\uCom(\ca(X_0,Z_0),Z_pf_l)\tens s\uCom(Z_pf_l,Z_nf_l)
\rTTo^{b^{\uCom}_2}
s\uCom(\ca(X_0,Z_0),Z_nf_l)
\bigr].
\label{equ*-ev-ev-b2}
\end{multline}
Replacing the last two arrows by the composite \(m^{\uCom}_2[-1]s\) and
applying identity~\eqref{eq-identity-m2} leads to
\begin{multline}
-(-)^{k+1}\bigl[
\bT^ns\ca(Z_0,Z_n)\tens\bT^ks\ca^\op(X_0,X_k)\tens X_kf_0[1]\tens
\bT^ls\und\Ainfty(\ca;\uCom)(f_0,f_l)
\\
\quad\rTTo^{\coev^\Com}
\uCom(s\ca(X_0,Z_0),s\ca(X_0,Z_0)\tens\bT^ns\ca(Z_0,Z_n)\hfill
\\
\hfill\tens\bT^ks\ca^\op(X_0,X_k)\tens X_kf_0[1]\tens
\bT^ls\und\Ainfty(\ca;\uCom)(f_0,f_l))\quad
\\
\quad\rTTo^{\uCom(1,\text{perm})}
\uCom(s\ca(X_0,Z_0),X_kf_0[1]\tens\bT^ks\ca(X_k,X_0)\tens
s\ca(X_0,Z_0)\hfill
\\
\hfill\tens\bT^ns\ca(Z_0,Z_n)\tens
\bT^ls\und\Ainfty(\ca;\uCom)(f_0,f_l))\quad
\\
\rTTo^{\uCom(1,1\tens\text{perm})}
\uCom(s\ca(X_0,Z_0),X_kf_0[1]\tens\bT^ks\ca(X_k,X_0)\tens
s\ca(X_0,Z_0)\tens\bT^ps\ca(Z_0,Z_p)
\\
\hfill\tens\bT^ls\und\Ainfty(\ca;\uCom)(f_0,f_l)\tens\bT^qs\ca(Z_p,Z_n)\tens
T^0s\und\Ainfty(\ca;\uCom)(f_l,f_l))\quad
\\
\rTTo^{\uCom(1,1\tens\ev^{\Ainfty}_{k+1+p,l}\tens\ev^{\Ainfty}_{q0})}
\uCom(s\ca(X_0,Z_0),X_kf_0[1]\tens s\uCom(X_kf_0,Z_pf_l)\tens
s\uCom(Z_pf_l,Z_nf_l))
\\
\quad \rTTo^{\uCom(1,1\tens s^{-1}[1]\tens s^{-1}[1])} \hfill
\\
\uCom(s\ca(X_0,Z_0),X_kf_0[1]\tens\uCom(X_kf_0[1],Z_pf_l[1])\tens\uCom(Z_pf_l[1],Z_nf_l[1]))
\\
\rTTo^{\uCom(1,(\ev^\Com\tens1)\ev^\Com)}
\uCom(s\ca(X_0,Z_0),Z_nf_l[1])
\rTTo^{[-1]s}
s\uCom(\ca(X_0,Z_0),Z_nf_l)
\bigr].
\label{equ*-ev-ev-ev-1-ev}
\end{multline}
Since
\begin{multline*}
(\ev^\Com\tens1)\ev^\Com=(1\tens m^{\uCom}_2)\ev^\Com:
\\
X_kf_0[1]\tens\uCom(X_kf_0[1],Z_pf_l[1])\tens
\Com(Z_pf_l[1],Z_nf_l[1])\to Z_nf_l[1],
\end{multline*}
and \((s^{-1}[1]\tens s^{-1}[1])m^{\uCom}_2=-b^{\uCom}_2s^{-1}[1]\), we
infer that term~\eqref{equ-Mho-evA8-b2} equals
\begin{multline*}
(-)^{k+1}\bigl[
\bT^ns\ca(Z_0,Z_n)\tens\bT^ks\ca^\op(X_0,X_k)\tens X_kf_0[1]\tens
\bT^ls\und\Ainfty(\ca;\uCom)(f_0,f_l)
\\
\quad\rTTo^{\coev^\Com}
\uCom(s\ca(X_0,Z_0),s\ca(X_0,Z_0)\tens\bT^ns\ca(Z_0,Z_n)\hfill
\\
\hfill\tens\bT^ks\ca^\op(X_0,X_k)\tens X_kf_0[1]\tens
\bT^ls\und\Ainfty(\ca;\uCom)(f_0,f_l))\quad
\\
\quad\rTTo^{\uCom(1,\text{perm})}
\uCom(s\ca(X_0,Z_0),X_kf_0[1]\tens\bT^ks\ca(X_k,X_0)\tens
s\ca(X_0,Z_0)\hfill
\\
\hfill\tens\bT^ns\ca(Z_0,Z_n)\tens
\bT^ls\und\Ainfty(\ca;\uCom)(f_0,f_l))\quad
\\
\quad\rTTo^{\uCom(1,1\tens\text{perm}\cdot(\ev^{\Ainfty}_{k+1+p,l}\tens\ev^{\Ainfty}_{q0})b^{\uCom}_2)}
\uCom(s\ca(X_0,Z_0),X_kf_0[1]\tens s\uCom(X_kf_0,Z_nf_l))\hfill
\\
\quad\rTTo^{\uCom(1,1\tens s^{-1}[1])}
\uCom(s\ca(X_0,Z_0),X_kf_0[1]\tens\uCom(X_kf_0[1],Z_nf_l[1]))\hfill
\\
\rTTo^{\uCom(1,\ev^\Com)}
\uCom(s\ca(X_0,Z_0),Z_nf_l[1])
\rTTo^{[-1]s}
s\uCom(\ca(X_0,Z_0),Z_nf_l)
\bigr].
\end{multline*}
Similarly, term~\eqref{equ-H-Mho-b2} equals
\begin{multline*}
-(-)^{k+1}\bigl[
\bT^ns\ca(Z_0,Z_n)\tens\bT^ks\ca^\op(X_0,X_k)\tens X_kf_0[1]\tens
\bT^ls\und\Ainfty(\ca;\uCom)(f_0,f_l)
\\
\quad\rTTo^{\coev^\Com\tens\coev^\Com}
\uCom(s\ca(X_0,Z_0),s\ca(X_0,Z_0)\tens\bT^ps\ca(Z_0,Z_p))
\tens\uCom(s\ca(X_0,Z_p),\hfill
\\
s\ca(X_0,Z_p)\tens\bT^qs\ca(Z_p,Z_n)\tens\bT^ks\ca^\op(X_0,X_k)
\tens X_kf_0[1]\tens\bT^ls\und\Ainfty(\ca;\uCom)(f_0,f_l))
\\
\quad\rTTo^{1\tens\uCom(1,\text{perm})}
\uCom(s\ca(X_0,Z_0),s\ca(X_0,Z_0)\tens\bT^ps\ca(Z_0,Z_p))\hfill
\\
\tens\uCom(s\ca(X_0,Z_p),X_kf_0[1]\tens\bT^ks\ca(X_k,X_0)
\\
\hfill\tens s\ca(X_0,Z_p)
\tens\bT^qs\ca(Z_p,Z_n)\tens\bT^ls\und\Ainfty(\ca;\uCom)(f_0,f_l))\quad
\\
\quad\rTTo^{\uCom(1,b_{p+1})\tens\uCom(1,1\tens\ev^{\Ainfty}_{k+1+q,l})}
\uCom(s\ca(X_0,Z_0),s\ca(X_0,Z_p))\hfill
\\
\hfill\tens\uCom(s\ca(X_0,Z_p),X_kf_0[1]\tens s\uCom(X_kf_0,Z_nf_l))\quad
\\
\quad \rTTo^{1\tens\uCom(1,1\tens s^{-1}[1])} \hfill
\\
\uCom(s\ca(X_0,Z_0),s\ca(X_0,Z_p))
\tens\uCom(s\ca(X_0,Z_p),X_kf_0[1]\tens\uCom(X_kf_0[1],Z_nf_l[1]))
\\
\quad\rTTo^{1\tens\uCom(1,\ev^\Com)}
\uCom(s\ca(X_0,Z_0),s\ca(X_0,Z_p))
\tens\uCom(s\ca(X_0,Z_p),Z_nf_l[1])\hfill
\\
\quad\rTTo^{[-1]s\tens[-1]s}
s\uCom(\ca(X_0,Z_0),\ca(X_0,Z_p))\tens s\uCom(\ca(X_0,Z_p),Z_nf_l)\hfill
\\
\hfill\rTTo^{b^{\uCom}_2}
s\uCom(\ca(X_0,Z_0),Z_nf_l)
\bigr]\quad
\\
\hskip\multlinegap=-(-)^{k+1}\bigl[
\bT^ns\ca(Z_0,Z_n)\tens\bT^ks\ca^\op(X_0,X_k)\tens X_kf_0[1]\tens
\bT^ls\und\Ainfty(\ca;\uCom)(f_0,f_l)\hfill
\\
\quad\rTTo^{\coev^\Com}\uCom(s\ca(X_0,Z_0),s\ca(X_0,Z_0)\tens
\bT^ns\ca(Z_0,Z_n)\hfill
\\
\hfill\tens\bT^ks\ca^\op(X_0,X_k)\tens X_kf_0[1]\tens
\bT^ls\und\Ainfty(\ca;\uCom)(f_0,f_l))\quad
\\
\quad\rTTo^{\uCom(1,\text{perm})}
\uCom(s\ca(X_0,Z_0),X_kf_0[1]\tens\bT^ks\ca(X_k,X_0)\tens
s\ca(X_0,Z_0)\hfill
\\
\hfill\tens\bT^ns\ca(Z_0,Z_n)\tens\bT^ls\und\Ainfty(\ca;\uCom)(f_0,f_l))\quad
\\
\quad\rTTo^{\uCom(1,1\tens(1^{\tens k}\tens b_{p+1}\tens1^{\tens q}\tens1^{\tens l})\ev^{\Ainfty}_{k+1+q,l})}
\uCom(s\ca(X_0,Z_0),X_kf_0[1]\tens s\uCom(X_kf_0,Z_nf_l))\hfill
\\
\quad\rTTo^{\uCom(1,1\tens s^{-1}[1])}
\uCom(s\ca(X_0,Z_0),X_kf_0[1]\tens\uCom(X_kf_0[1],Z_nf_l[1]))\hfill
\\
\rTTo^{\uCom(1,\ev^\Com)}
\uCom(s\ca(X_0,Z_0),Z_nf_l[1])
\rTTo^{[-1]s}
s\uCom(\ca(X_0,Z_0),Z_nf_l)
\bigr]
\end{multline*}
due to identity~\eqref{eq-identity-m2}. It follows immediately from
naturality of \(\coev^\Com\) that term \eqref{equ-1-b-1-1-1-1-Mho}
equals
\begin{multline*}
(-)^{k+1}\bigl[
\bT^ns\ca(Z_0,Z_n)\tens\bT^ks\ca^\op(X_0,X_k)\tens X_kf_0[1]\tens
\bT^ls\und\Ainfty(\ca;\uCom)(f_0,f_l)
\\
\quad\rTTo^{\coev^\Com}\uCom(s\ca(X_0,Z_0),s\ca(X_0,Z_0)\tens
\bT^ns\ca(Z_0,Z_n)\hfill
\\
\hfill\tens\bT^ks\ca^\op(X_0,X_k)\tens X_kf_0[1]\tens
\bT^ls\und\Ainfty(\ca;\uCom)(f_0,f_l))\quad
\\
\quad\rTTo^{\uCom(1,\text{perm})}
\uCom(s\ca(X_0,Z_0),X_kf_0[1]\tens\bT^ks\ca(X_k,X_0)\tens
s\ca(X_0,Z_0)\hfill
\\
\hfill\tens\bT^ns\ca(Z_0,Z_n)\tens\bT^ls\und\Ainfty(\ca;\uCom)(f_0,f_l))\quad
\\
\quad\rTTo^{\uCom(1,1\tens(1^{\tens k}\tens1\tens1^{\tens\alpha}\tens b_t\tens1^{\tens\beta}\tens1^{\tens l})
\ev^{\Ainfty}_{k+\alpha+2+\beta,l})}\hfill
\\
\uCom(s\ca(X_0,Z_0),X_kf_0[1]\tens s\uCom(X_kf_0,Z_nf_l))
\\
\quad\rTTo^{\uCom(1,1\tens s^{-1}[1])}
\uCom(s\ca(X_0,Z_0),X_kf_0[1]\tens\uCom(X_kf_0[1],Z_nf_l[1]))\hfill
\\
\rTTo^{\uCom(1,\ev^\Com)}
\uCom(s\ca(X_0,Z_0),Z_nf_l[1])
\rTTo^{[-1]s}
s\uCom(\ca(X_0,Z_0),Z_nf_l)
\bigr].
\end{multline*}
Term~\eqref{equ-Mho-evA8-b2-15} can be written as follows:
\begin{multline*}
-(-)^{k+1}\bigl[
\bT^ns\ca(Z_0,Z_n)\tens\bT^ks\ca^\op(X_0,X_k)\tens X_kf_0[1]\tens
\bT^ls\und\Ainfty(\ca;\uCom)(f_0,f_l)
\\
\quad\rTTo^{\text{perm}}
\bT^ps\ca(Z_0,Z_p)\tens\bT^ks\ca^\op(X_0,X_k)\tens X_kf_0[1]
\tens\bT^{l-1}s\und\Ainfty(\ca;\uCom)(f_0,f_{l-1})\hfill
\\
\hfill\tens\bT^qs\ca(Z_p,Z_n)\tens
s\und\Ainfty(\ca;\uCom)(f_{l-1},f_l)\quad
\\
\quad\rTTo^{\coev^\Com\tens\coev^\Com}
\uCom(s\ca(X_0,Z_0),s\ca(X_0,Z_0)\tens\bT^ps\ca(Z_0,Z_p)\tens\bT^ks\ca^\op(X_0,X_k)\hfill
\\
\hfill\tens
X_kf_0[1]\tens\bT^{l-1}s\und\Ainfty(\ca;\uCom)(f_0,f_{l-1}))\quad
\\
\hfill\tens\uCom(Z_pf_{l-1}[1],Z_pf_{l-1}[1]\tens\bT^qs\ca(Z_p,Z_n)\tens
s\und\Ainfty(\ca;\uCom)(f_{l-1},f_l))\quad
\\
\quad\rTTo^{\uCom(1,\text{perm})\tens1}
\uCom(s\ca(X_0,Z_0),X_kf_0[1]\tens\bT^ks\ca(X_k,X_0)\tens
s\ca(X_0,Z_0)\hfill
\\
\hfill\tens\bT^ps\ca(Z_0,Z_p)\tens
\bT^{l-1}s\und\Ainfty(\ca;\uCom)(f_0,f_{l-1}))\quad
\\
\hfill\tens\uCom(Z_pf_{l-1}[1],Z_pf_{l-1}[1]\tens\bT^qs\ca(Z_p,Z_n)\tens
s\und\Ainfty(\ca;\uCom)(f_{l-1},f_l))\quad
\\
\quad\rTTo^{\uCom(1,1\tens\ev^{\Ainfty}_{k+1+p,l-1})\tens\uCom(1,1\tens\ev^{\Ainfty}_{q1})}
\uCom(s\ca(X_0,Z_0),X_kf_0[1]\tens s\uCom(X_kf_0,Z_pf_{l-1}))\hfill
\\
\hfill\tens\uCom(Z_pf_{l-1}[1],Z_pf_{l-1}[1]\tens s\uCom(Z_pf_{l-1},Z_nf_l))\quad
\\
\quad\rTTo^{\uCom(1,1\tens s^{-1}[1])\tens\uCom(1,1\tens s^{-1}[1])}
\uCom(s\ca(X_0,Z_0),X_kf_0[1]\tens\uCom(X_kf_0[1],Z_pf_{l-1}[1]))\hfill
\\
\hfill\tens\uCom(Z_pf_{l-1}[1],Z_pf_{l-1}[1]\tens\uCom(Z_pf_{l-1}[1],Z_nf_l[1]))\quad
\\
\rTTo^{\uCom(1,\ev^\Com)\tens\uCom(1,\ev^\Com)}
\uCom(s\ca(X_0,Z_0),Z_pf_{l-1}[1])
\tens\uCom(Z_pf_{l-1}[1],Z_nf_l[1]) \rTTo^{[-1]s\tens[-1]s}
\\
s\uCom(\ca(X_0,Z_0),Z_pf_{l-1})\tens s\uCom(Z_pf_{l-1},Z_nf_l)
\rTTo^{b^{\uCom}_2}
s\uCom(\ca(X_0,Z_0),Z_nf_l)
\bigr].
\end{multline*}
The obtained expression is of the same type as \eqref{equ*-ev-ev-b2}.
Transform it in the same manner to conclude that
term~\eqref{equ-Mho-evA8-b2-15} equals
\begin{multline*}
(-)^{k+1}\bigl[
\bT^ns\ca(Z_0,Z_n)\tens\bT^ks\ca^\op(X_0,X_k)\tens X_kf_0[1]\tens
\bT^ls\und\Ainfty(\ca;\uCom)(f_0,f_l)
\\
\quad\rTTo^{\coev^\Com}\uCom(s\ca(X_0,Z_0),s\ca(X_0,Z_0)\tens
\bT^ns\ca(Z_0,Z_n)\hfill
\\
\hfill\tens\bT^ks\ca^\op(X_0,X_k)\tens X_kf_0[1]\tens
\bT^ls\und\Ainfty(\ca;\uCom)(f_0,f_l))\quad
\\
\quad\rTTo^{\uCom(1,\text{perm})}
\uCom(s\ca(X_0,Z_0),X_kf_0[1]\tens\bT^ks\ca(X_k,X_0)\tens
s\ca(X_0,Z_0)\hfill
\\
\hfill\tens\bT^ns\ca(Z_0,Z_n)\tens\bT^ls\und\Ainfty(\ca;\uCom)(f_0,f_l))\quad
\\
\quad\rTTo^{\uCom(1,1\tens\perm\cdot(\ev^{\Ainfty}_{k+1+p,l-1}\tens\ev^{\Ainfty}_{q1})b^{\uCom}_2)}
\uCom(s\ca(X_0,Z_0),X_kf_0[1]\tens s\uCom(X_kf_0,Z_nf_l))\hfill
\\
\quad\rTTo^{\uCom(1,1\tens s^{-1}[1])}
\uCom(s\ca(X_0,Z_0),X_kf_0[1]\tens\uCom(X_kf_0[1],Z_nf_l[1]))\hfill
\\
\rTTo^{\uCom(1,\ev^\Com)}
\uCom(s\ca(X_0,Z_0),Z_nf_l[1])
\rTTo^{[-1]s}
s\uCom(\ca(X_0,Z_0),Z_nf_l)
\bigr].
\end{multline*}
Obviously, term~\eqref{equ-1111-B1-1-Mho} can be written as follows:
\begin{multline*}
(-)^{k+1}\bigl[
\bT^ns\ca(Z_0,Z_n)\tens\bT^ks\ca^\op(X_0,X_k)\tens X_kf_0[1]\tens
\bT^ls\und\Ainfty(\ca;\uCom)(f_0,f_l)
\\
\quad\rTTo^{\coev^\Com}\uCom(s\ca(X_0,Z_0),s\ca(X_0,Z_0)\tens
\bT^ns\ca(Z_0,Z_n)\hfill
\\
\hfill\tens\bT^ks\ca^\op(X_0,X_k)\tens X_kf_0[1]\tens
\bT^ls\und\Ainfty(\ca;\uCom)(f_0,f_l))\quad
\\
\quad\rTTo^{\uCom(1,\text{perm})}
\uCom(s\ca(X_0,Z_0),X_kf_0[1]\tens\bT^ks\ca(X_k,X_0)\tens
s\ca(X_0,Z_0)\hfill
\\
\hfill\tens\bT^ns\ca(Z_0,Z_n)\tens\bT^ls\und\Ainfty(\ca;\uCom)(f_0,f_l))\quad
\\
\quad\rTTo^{\uCom(1,1\tens(1^{\tens k}\tens1\tens1^{\tens
n}\tens1^{\tens p-1}\tens B_1\tens1^{\tens
l-p})\ev^{\Ainfty}_{k+1+n,l})} \hfill
\\
\uCom(s\ca(X_0,Z_0),X_kf_0[1]\tens s\uCom(X_kf_0,Z_nf_l))
\\
\quad\rTTo^{\uCom(1,1\tens s^{-1}[1])}
\uCom(s\ca(X_0,Z_0),X_kf_0[1]\tens\uCom(X_kf_0[1],Z_nf_l[1]))\hfill
\\
\rTTo^{\uCom(1,\ev^\Com)}
\uCom(s\ca(X_0,Z_0),Z_nf_l[1])
\rTTo^{[-1]s}
s\uCom(\ca(X_0,Z_0),Z_nf_l)
\bigr].
\end{multline*}
Similar presentation holds for term~\eqref{equ-1111-B2-1-Mho}:
\begin{multline*}
(-)^{k+1}\bigl[
\bT^ns\ca(Z_0,Z_n)\tens\bT^ks\ca^\op(X_0,X_k)\tens X_kf_0[1]\tens
\bT^ls\und\Ainfty(\ca;\uCom)(f_0,f_l)
\\
\quad\rTTo^{\coev^\Com}\uCom(s\ca(X_0,Z_0),s\ca(X_0,Z_0)\tens
\bT^ns\ca(Z_0,Z_n)\hfill
\\
\hfill\tens\bT^ks\ca^\op(X_0,X_k)\tens X_kf_0[1]\tens
\bT^ls\und\Ainfty(\ca;\uCom)(f_0,f_l))\quad
\\
\quad\rTTo^{\uCom(1,\text{perm})}
\uCom(s\ca(X_0,Z_0),X_kf_0[1]\tens\bT^ks\ca(X_k,X_0)\tens
s\ca(X_0,Z_0)\hfill
\\
\hfill\tens\bT^ns\ca(Z_0,Z_n)\tens\bT^ls\und\Ainfty(\ca;\uCom)(f_0,f_l))\quad
\\
\quad\rTTo^{\uCom(1,1\tens(1^{\tens k}\tens1\tens1^{\tens
n}\tens1^{\tens p-1}\tens B_2\tens1^{\tens
l-p-1})\ev^{\Ainfty}_{k+1+n,l-1})} \hfill
\\
\uCom(s\ca(X_0,Z_0),X_kf_0[1]\tens s\uCom(X_kf_0,Z_nf_l))
\\
\quad\rTTo^{\uCom(1,1\tens s^{-1}[1])}
\uCom(s\ca(X_0,Z_0),X_kf_0[1]\tens\uCom(X_kf_0[1],Z_nf_l[1]))\hfill
\\
\rTTo^{\uCom(1,\ev^\Com)}
\uCom(s\ca(X_0,Z_0),Z_nf_l[1])
\rTTo^{[-1]s}
s\uCom(\ca(X_0,Z_0),Z_nf_l)
\bigr].
\end{multline*}
If \(m=k\), term~\eqref{equ-11-bE0-1-Mho} equals
\begin{multline*}
(-)^{k+1}\bigl[
\bT^ns\ca(Z_0,Z_n)\tens\bT^ks\ca^\op(X_0,X_k)\tens X_kf_0[1]\tens
\bT^ls\und\Ainfty(\ca;\uCom)(f_0,f_l)
\\
\quad\rTTo^{\coev^\Com}
\uCom(s\ca(X_0,Z_0),s\ca(X_0,Z_0)\tens\bT^ns\ca(Z_0,Z_n)\hfill
\\
\hfill\tens\bT^ks\ca^\op(X_0,X_k)\tens X_kf_0[1]\tens
\bT^ls\und\Ainfty(\ca;\uCom)(f_0,f_l))\quad
\\
\quad\rTTo^{\uCom(1,\text{perm})}
\uCom(s\ca(X_0,Z_0),X_kf_0[1]\tens\bT^ks\ca(X_k,X_0)\tens
s\ca(X_0,Z_0)\hfill
\\
\hfill\tens\bT^ns\ca(Z_0,Z_n)\tens
\bT^ls\und\Ainfty(\ca;\uCom)(f_0,f_l))\quad
\\
\quad\rTTo^{\uCom(1,d'\tens\ev^{\Ainfty}_{k+1+n,l})}
\uCom(s\ca(X_0,Z_0),X_kf_0[1]\tens s\uCom(X_kf_0,Z_nf_l))\hfill
\\
\quad\rTTo^{\uCom(1,1\tens s^{-1}[1])}
\uCom(s\ca(X_0,Z_0),X_kf_0[1]\tens\uCom(X_kf_0[1],Z_nf_l[1]))\hfill
\\
\hfill\rTTo^{\uCom(1,\ev^\Com)}
\uCom(s\ca(X_0,Z_0),Z_nf_l[1])
\rTTo^{[-1]s}
s\uCom(\ca(X_0,Z_0),Z_nf_l)
\bigr],
\end{multline*}
it cancels one of the summands present in \eqref{equ*-ev-b1}. Suppose
that \(0\le m\le k-1\). Expressing \(b^\ce_{k-m,0}\) via
\(\ev^{\Ainfty}_{k-m,0}\), we write term~\eqref{equ-11-bE0-1-Mho} as
follows:
\begin{multline*}
(-)^{k+1}\bigl[
\bT^ns\ca(Z_0,Z_n)\tens\bT^ks\ca^\op(X_0,X_k)\tens X_kf_0[1]\tens
\bT^ls\und\Ainfty(\ca;\uCom)(f_0,f_l)
\\
\quad\rTTo^{\text{perm}}
\bT^ns\ca(Z_0,Z_n)\tens\bT^ms\ca^\op(X_0,X_m)\tens X_kf_0[1]\tens
\bT^{k-m}s\ca(X_k,X_m)\hfill
\\
\hfill\tens T^0s\und\Ainfty(\ca;\uCom)(f_0,f_0)
\tens\bT^ls\und\Ainfty(\ca;\uCom)(f_0,f_l)\quad
\\
\quad\rTTo^{1^{\tens n}\tens1^{\tens
m}\tens1\tens\ev^{\Ainfty}_{k-m,0}\tens1^{\tens l}}
\bT^ns\ca(Z_0,Z_n)\tens\bT^ms\ca^\op(X_0,X_m)\tens X_kf_0[1]\hfill
\\
\hfill\tens s\uCom(X_kf_0,X_mf_0)\tens\bT^ls\und\Ainfty(\ca;\uCom)(f_0,f_l)\quad
\\
\quad\rTTo^{1^{\tens n}\tens1^{\tens m}\tens1\tens s^{-1}[1]\tens1^{\tens l}}
\bT^ns\ca(Z_0,Z_n)\tens\bT^ms\ca^\op(X_0,X_m)\tens X_kf_0[1]\hfill
\\
\hfill\tens\uCom(X_kf_0[1],X_mf_0[1])\tens
\bT^ls\und\Ainfty(\ca;\uCom)(f_0,f_l)\quad
\\
\quad\rTTo^{1^{\tens n}\tens1^{\tens m}\tens\ev^\Com\tens1^{\tens l}}
\hfill
\\
\bT^ns\ca(Z_0,Z_n)\tens\bT^ms\ca^\op(X_0,X_m)\tens X_mf_0[1]\tens
\bT^ls\und\Ainfty(\ca;\uCom)(f_0,f_l)
\\
\quad\rTTo^{\coev^\Com}
\uCom(s\ca(X_0,Z_0),s\ca(X_0,Z_0)\tens\bT^ns\ca(Z_0,Z_n)\tens\bT^ms\ca^\op(X_0,X_m)\hfill
\\
\hfill\tens X_mf_0[1]\tens
\bT^ls\und\Ainfty(\ca;\uCom)(f_0,f_l))\quad
\\
\quad\rTTo^{\uCom(1,\text{perm})}
\uCom(s\ca(X_0,Z_0),X_mf_0[1]\tens\bT^ms\ca(X_m,X_0)\tens
s\ca(X_0,Z_0)\hfill
\\
\hfill\tens\bT^ns\ca(Z_0,Z_n)\tens\bT^ls\und\Ainfty(\ca;\uCom)(f_0,f_l))\quad
\\
\quad\rTTo^{\uCom(1,1\tens\ev^{\Ainfty}_{m+1+n,l})}
\uCom(s\ca(X_0,Z_0),X_mf_0[1]\tens s\uCom(X_mf_0,Z_nf_l))\hfill
\\
\quad\rTTo^{\uCom(1,1\tens s^{-1}[1])}
\uCom(s\ca(X_0,Z_0),X_mf_0[1]\tens\uCom(X_mf_0[1],Z_nf_l[1]))\hfill
\\
\hfill\rTTo^{\uCom(1,\ev^\Com)}
\uCom(s\ca(X_0,Z_0),Z_nf_l[1])
\rTTo^{[-1]s}
s\uCom(\ca(X_0,Z_0),Z_nf_l)
\bigr]\quad
\\
\hskip\multlinegap=(-)^{k+1}\bigl[
\bT^ns\ca(Z_0,Z_n)\tens\bT^ks\ca^\op(X_0,X_k)\tens X_kf_0[1]
\tens\bT^ls\und\Ainfty(\ca;\uCom)(f_0,f_l)\hfill
\\
\quad\rTTo^{\coev^\Com}
\uCom(s\ca(X_0,Z_0),s\ca(X_0,Z_0)\tens\bT^ns\ca(Z_0,Z_n)\tens
\bT^ks\ca^\op(X_0,X_k)\hfill
\\
\hfill\tens X_kf_0[1]
\tens\bT^ls\und\Ainfty(\ca;\uCom)(f_0,f_l))\quad
\\
\quad\rTTo^{\uCom(1,\text{perm})}
\uCom(s\ca(X_0,Z_0),X_kf_0[1]\tens\bT^{k-m}s\ca(X_k,X_m)\tens
T^0s\und\Ainfty(\ca;\uCom)(f_0,f_0)\hfill
\\
\hfill\tens\bT^ms\ca(X_m,X_0)\tens
s\ca(X_0,Z_0)\tens\bT^ns\ca(Z_0,Z_n) \tens
\bT^ls\und\Ainfty(\ca;\uCom)(f_0,f_l))\quad
\\
\quad
\rTTo^{\uCom(1,1\tens\ev^{\Ainfty}_{k-m,0}\tens\ev^{\Ainfty}_{m+1+n,l})}
\hfill
\\
\uCom(s\ca(X_0,Z_0),X_kf_0[1]\tens s\uCom(X_kf_0,X_mf_0)\tens
s\uCom(X_mf_0,Z_nf_l))
\\
\quad \rTTo^{\uCom(1,1\tens s^{-1}[1]\tens s^{-1}[1])} \hfill
\\
\uCom(s\ca(X_0,Z_0),X_kf_0[1]\tens\uCom(X_kf_0[1],X_mf_0[1])\tens\uCom(X_mf_0[1],Z_nf_l[1]))
\\
\rTTo^{\uCom(1,(\ev^\Com\tens1)\ev^\Com)}
\uCom(s\ca(X_0,Z_0),Z_nf_l[1])
\rTTo^{[-1]s}
s\uCom(\ca(X_0,Z_0),Z_nf_l)
\bigr].
\end{multline*}
The further transformations are parallel to \eqref{equ*-ev-ev-ev-1-ev}. We conclude that
term~\eqref{equ-11-bE0-1-Mho} equals
\begin{multline*}
-(-)^{k+1}\bigl[
\bT^ns\ca(Z_0,Z_n)\tens\bT^ks\ca^\op(X_0,X_k)\tens X_kf_0[1]\tens
\bT^ls\und\Ainfty(\ca;\uCom)(f_0,f_l)
\\
\quad\rTTo^{\coev^\Com}\uCom(s\ca(X_0,Z_0),s\ca(X_0,Z_0)\tens
\bT^ns\ca(Z_0,Z_n)\hfill
\\
\hfill\tens\bT^ks\ca^\op(X_0,X_k)\tens X_kf_0[1]\tens
\bT^ls\und\Ainfty(\ca;\uCom)(f_0,f_l))\quad
\\
\quad\rTTo^{\uCom(1,\text{perm})}
\uCom(s\ca(X_0,Z_0),X_kf_0[1]\tens\bT^ks\ca(X_k,X_0)\tens
s\ca(X_0,Z_0)\hfill
\\
\hfill\tens\bT^ns\ca(Z_0,Z_n)\tens\bT^ls\und\Ainfty(\ca;\uCom)(f_0,f_l))\quad
\\
\quad\rTTo^{\uCom(1,1\tens(\ev^{\Ainfty}_{k-m,0}\tens\ev^{\Ainfty}_{m+1+n,l})b^{\uCom}_2)}
\uCom(s\ca(X_0,Z_0),X_kf_0[1]\tens s\uCom(X_kf_0,Z_nf_l))\hfill
\\
\quad\rTTo^{\uCom(1,1\tens s^{-1}[1])}
\uCom(s\ca(X_0,Z_0),X_kf_0[1]\tens\uCom(X_kf_0[1],Z_nf_l[1]))\hfill
\\
\rTTo^{\uCom(1,\ev^\Com)}
\uCom(s\ca(X_0,Z_0),Z_nf_l[1])
\rTTo^{[-1]s}
s\uCom(\ca(X_0,Z_0),Z_nf_l)
\bigr].
\end{multline*}
The case of term~\eqref{equ-11-bE1-1-Mho} is quite similar, we only
give the result:
\begin{multline*}
-(-)^{k+1}\bigl[
\bT^ns\ca(Z_0,Z_n)\tens\bT^ks\ca^\op(X_0,X_k)\tens X_kf_0[1]\tens
\bT^ls\und\Ainfty(\ca;\uCom)(f_0,f_l)
\\
\quad\rTTo^{\coev^\Com}\uCom(s\ca(X_0,Z_0),s\ca(X_0,Z_0)\tens
\bT^ns\ca(Z_0,Z_n)\hfill
\\
\hfill\tens\bT^ks\ca^\op(X_0,X_k)\tens X_kf_0[1]\tens
\bT^ls\und\Ainfty(\ca;\uCom)(f_0,f_l))\quad
\\
\quad\rTTo^{\uCom(1,\text{perm})}
\uCom(s\ca(X_0,Z_0),X_kf_0[1]\tens\bT^ks\ca(X_k,X_0)\tens
s\ca(X_0,Z_0)\hfill
\\
\hfill\tens\bT^ns\ca(Z_0,Z_n)\tens\bT^ls\und\Ainfty(\ca;\uCom)(f_0,f_l))\quad
\\
\quad\rTTo^{\uCom(1,1\tens\text{perm}\cdot(\ev^{\Ainfty}_{k-m,1}\tens\ev^{\Ainfty}_{m+1+n,l-1})b^{\uCom}_2)}
\uCom(s\ca(X_0,Z_0),X_kf_0[1]\tens s\uCom(X_kf_0,Z_nf_l))\hfill
\\
\quad\rTTo^{\uCom(1,1\tens s^{-1}[1])}
\uCom(s\ca(X_0,Z_0),X_kf_0[1]\tens\uCom(X_kf_0[1],Z_nf_l[1]))\hfill
\\
\rTTo^{\uCom(1,\ev^\Com)}
\uCom(s\ca(X_0,Z_0),Z_nf_l[1])
\rTTo^{[-1]s}
s\uCom(\ca(X_0,Z_0),Z_nf_l)
\bigr].
\end{multline*}
Finally, using formula~\eqref{equ-codiff-opposite} for
\(b^{\ca^\op}_u\), we find that term~\eqref{equ-11-bop-111-Mho} equals
\begin{multline*}
(-)^{a+2+c}\bigl[
\bT^ns\ca(Z_0,Z_n)\tens\bT^ks\ca^\op(X_0,X_k)\tens X_kf_0[1]\tens
\bT^ls\und\Ainfty(\ca;\uCom)(f_0,f_l)
\\
\quad\rTTo^{\coev^\Com}
\uCom(s\ca(X_0,Z_0),s\ca(X_0,Z_0)\tens\bT^ns\ca(Z_0,Z_n)\tens\bT^ks\ca^\op(X_0,X_k)\hfill
\\
\hfill\tens X_kf_0[1]\tens
\bT^ls\und\Ainfty(\ca;\uCom)(f_0,f_l))\quad
\\
\quad\rTTo^{\uCom(1,1\tens1^{\tens n}\tens1^{\tens a}\tens b^{\ca^\op}_u\tens1^{\tens c}\tens1\tens1^{\tens l})}
\uCom(s\ca(X_0,Z_0),s\ca(X_0,Z_0)\tens\bT^ns\ca(Z_0,Z_n)\hfill
\\
\hfill\tens \bT^{a+1+c}s\ca^\op(X_0,\dots,X_a,X_{a+u},\dots,X_k)\tens X_kf_0[1]\tens
\bT^ls\und\Ainfty(\ca;\uCom)(f_0,f_l))\quad
\\
\quad\rTTo^{\uCom(1,\text{perm})}
\uCom(s\ca(X_0,Z_0),X_kf_0[1]\tens\bT^{c+1+a}s\ca(X_k,\dots,X_{a+u},X_a,\dots,X_0)\hfill
\\
\hfill\tens s\ca(X_0,Z_0)\tens\bT^ns\ca(Z_0,Z_n) \tens
\bT^ls\und\Ainfty(\ca;\uCom)(f_0,f_l))\quad
\\
\quad\rTTo^{\uCom(1,1\tens\ev^{\Ainfty}_{c+2+a+n,l})}
\uCom(s\ca(X_0,Z_0),X_kf_0[1]\tens s\uCom(X_kf_0,Z_nf_l))\hfill
\\
\quad\rTTo^{\uCom(1,\tens s^{-1}[1])}
\uCom(s\ca(X_0,Z_0),X_kf_0[1]\tens\uCom(X_kf_0[1],Z_nf_l[1]))\hfill
\\
\hfill\rTTo^{\uCom(1,\tens\ev^\Com)}
\uCom(s\ca(X_0,Z_0),Z_nf_l[1])
\rTTo^{[-1]s}
s\uCom(\ca(X_0,Z_0),Z_nf_l)
\bigr]\quad
\\
=(-)^{a+u+c+1}\bigl[
\bT^ns\ca(Z_0,Z_n)\tens\bT^ks\ca^\op(X_0,X_k)\tens X_kf_0[1]\tens
\bT^ls\und\Ainfty(\ca;\uCom)(f_0,f_l)
\\
\quad\rTTo^{\coev^\Com}\uCom(s\ca(X_0,Z_0),s\ca(X_0,Z_0)\tens
\bT^ns\ca(Z_0,Z_n)\hfill
\\
\hfill\tens\bT^ks\ca^\op(X_0,X_k)\tens X_kf_0[1]\tens
\bT^ls\und\Ainfty(\ca;\uCom)(f_0,f_l))\quad
\\
\quad\rTTo^{\uCom(1,\text{perm})}
\uCom(s\ca(X_0,Z_0),X_kf_0[1]\tens\bT^ks\ca(X_k,X_0)\tens
s\ca(X_0,Z_0)\hfill
\\
\hfill\tens\bT^ns\ca(Z_0,Z_n)\tens\bT^ls\und\Ainfty(\ca;\uCom)(f_0,f_l))\quad
\\
\quad\rTTo^{\uCom(1,1\tens(1^{\tens c}\tens b_u\tens1^{\tens a}\tens1
 \tens1^{\tens n}\tens1^{\tens l})\ev^{\Ainfty}_{c+2+a+n,l})}
\hfill
\\
\uCom(s\ca(X_0,Z_0),X_kf_0[1]\tens s\uCom(X_kf_0,Z_nf_l))
\\
\quad\rTTo^{\uCom(1,1\tens s^{-1}[1])}
\uCom(s\ca(X_0,Z_0),X_kf_0[1]\tens\uCom(X_kf_0[1],Z_nf_l[1]))\hfill
\\
\rTTo^{\uCom(1,\ev^\Com)}
\uCom(s\ca(X_0,Z_0),Z_nf_l[1])
\rTTo^{[-1]s}
s\uCom(\ca(X_0,Z_0),Z_nf_l)
\bigr].
\end{multline*}
The overall sign is \((-)^{a+u+c+1}=(-)^{k+1}\).

Summing up, we conclude that
\begin{multline*}
(1^{\tens n}\tens S\cdot\pr_n)\ev^\Com=(-)^{k+1}\bigl[
\bT^ns\ca(Z_0,Z_n)\tens\bT^ks\ca^\op(X_0,X_k)\tens
X_kf_0[1]\hfill
\\
\hfill\tens
\bT^ls\und\Ainfty(\ca;\uCom)(f_0,f_l)\quad
\\
\quad\rTTo^{\coev^\Com}\uCom(s\ca(X_0,Z_0),s\ca(X_0,Z_0)\tens
\bT^ns\ca(Z_0,Z_n)\hfill
\\
\hfill\tens\bT^ks\ca^\op(X_0,X_k)\tens X_kf_0[1]\tens
\bT^ls\und\Ainfty(\ca;\uCom)(f_0,f_l))\quad
\\
\quad\rTTo^{\uCom(1,\text{perm})}
\uCom(s\ca(X_0,Z_0),X_kf_0[1]\tens\bT^ks\ca(X_k,X_0)\tens
s\ca(X_0,Z_0)\hfill
\\
\hfill\tens\bT^ns\ca(Z_0,Z_n)\tens\bT^ls\und\Ainfty(\ca;\uCom)(f_0,f_l))\quad
\\
\quad\rTTo^{\uCom(1,1\tens R)}
\uCom(s\ca(X_0,Z_0),X_kf_0[1]\tens s\uCom(X_kf_0,Z_nf_l))\hfill
\\
\quad\rTTo^{\uCom(1,1\tens s^{-1}[1])}
\uCom(s\ca(X_0,Z_0),X_kf_0[1]\tens\uCom(X_kf_0[1],Z_nf_l[1]))\hfill
\\
\rTTo^{\uCom(1,\ev^\Com)}
\uCom(s\ca(X_0,Z_0),Z_nf_l[1])
\rTTo^{[-1]s}
s\uCom(\ca(X_0,Z_0),Z_nf_l)
\bigr],
\end{multline*}
where
\begin{multline*}
R=-\sum_{m=1}^k\sum_{p+q=n}(1^{\tens k-m}\tens b_{m+1+p}\tens1^{\tens q}\tens1^{\tens
l})\ev^{\Ainfty}_{k-m+1+q,l}
+\ev^{\Ainfty}_{k+1+n,l}b^{\uCom}
\\
\hskip\multlinegap\hphantom{R}-(1^{\tens k}\tens b_1\tens1^{\tens
n}\tens1^{\tens l})\ev^{\Ainfty}_{k+1+n,l}
+\sum_{p+q=n}^{q>0}\perm\cdot(\ev^{\Ainfty}_{k+1+p,l}\tens\ev^{\Ainfty}_{q0})b^{\uCom}_2\hfill
\\
\hskip\multlinegap\hphantom{R}-\sum_{p+q=n}^{p>0}(1^{\tens k}\tens
b_{p+1}\tens1^{\tens q}\tens1^{\tens l})\ev^{\Ainfty}_{k+1+q,l}\hfill
\\
\hskip\multlinegap\hphantom{R}-\sum_{\alpha+t+\beta=n}(1^{\tens
k}\tens1\tens1^{\tens\alpha}\tens b_t\tens1^{\tens\beta}\tens1^{\tens
l})\ev^{\Ainfty}_{k+\alpha+2+\beta,l}\hfill
\\
\hskip\multlinegap\hphantom{R}
+\sum_{p+q=n}\perm\cdot(\ev^{\Ainfty}_{k+1+p,l-1}\tens\ev^{\Ainfty}_{q1})b^{\uCom}_2
\hfill
\\
\hskip\multlinegap\hphantom{R}
-\sum_{p=1}^l(1^{\tens k+1+n}\tens1^{\tens p-1}\tens B_1\tens1^{\tens
l-p})\ev^{\Ainfty}_{k+1+n,l}\hfill
\\
\hskip\multlinegap\hphantom{R}-\sum_{p=1}^{l-1}(1^{\tens
k+1+n}\tens1^{\tens p-1}\tens B_2\tens1^{\tens
l-p-1})\ev^{\Ainfty}_{k+1+n,l-1} \hfill
\\
\hskip\multlinegap\hphantom{R}
+\sum_{m=0}^{k-1}(\ev^{\Ainfty}_{k-m,0}\tens\ev^{\Ainfty}_{m+1+n,l})b^{\uCom}_2
+\sum_{m=0}^{k}\perm\cdot(\ev^{\Ainfty}_{k-m,1}\tens\ev^{\Ainfty}_{m+1+n,l-1})b^{\uCom}_2
\hfill
\\
\hskip\multlinegap\hphantom{R}-\sum_{a+u+c=k}(1^{\tens c}\tens
b_u\tens1^{\tens a}\tens1\tens1^{\tens n}\tens1^{\tens
l})\ev^{\Ainfty}_{c+2+a+n,l}:\hfill
\\
\hskip\multlinegap\hphantom{R}
\bT^ks\ca(X_k,X_0)\tens
s\ca(X_0,Z_0)\tens\bT^ns\ca(Z_0,Z_n)\tens\bT^ls\und\Ainfty(\ca;\uCom)(f_0,f_l)
\hfill
\\
\to s\uCom(X_kf_0,Z_nf_l),
\end{multline*}
which is easily seen to be the restriction of
\((\ev^{\Ainfty}b^{\uCom}-(b^\ca\boxt1+1\boxt
B)\ev^{\Ainfty})\pr_1:Ts\ca\boxt Ts\und\Ainfty(\ca;\uCom)\to s\uCom\)
to the summand \(\bT^ks\ca(X_k,X_0)\tens
s\ca(X_0,Z_0)\tens\bT^ns\ca(Z_0,Z_n)\tens\bT^ls\und\Ainfty(\ca;\uCom)(f_0,f_l)\)
of the source. Since \(\ev^{\Ainfty}\) is an \ainf-functor, it follows
that \(R=0\), and the equation is proven. That finishes the proof of
\lemref{lem-bicomodule-homomorphism-mho-chain-map}.
\end{proof}
\fi

Let \(\ca\) be an \ainf-category, and let \(f:\ca\to\uCom\) be an
\ainf-functor. Denote by \(\cm\) the $\ca$-module determined by $f$ in
\propref{prop-ainf-modules-ainf-functors-A-Com}. Denote
\begin{equation}
\Upsilon=\mho_{00}:
 s\cm(X) =s\ce(X,f) =Xf[1] \to s\und\Ainfty(\ca;\uCom)(H^X,f)
\label{eq-Upsilon-mho00}
\end{equation}
for the sake of brevity. The composition of \(\Upsilon\) with the
projection $\pr_n$ from \eqref{eq-prn-sAinfty} is given by the
particular case \(p=q=0\) of \eqref{eq-agemo-pqn}:
\begin{multline*}
\Upsilon_n=-\bigl[ s\cm(X) \rTTo^{\coev^\Com}
\\
\uCom(s\ca(X,Z_0)\tens T^ns\ca(Z_0,Z_n),
 s\ca(X,Z_0)\tens T^ns\ca(Z_0,Z_n)\tens s\cm(X))
\\
\rTTo^{\uCom(1,\tau_c)}
 \uCom(s\ca(X,Z_0)\tens T^ns\ca(Z_0,Z_n),
 s\cm(X)\tens s\ca(X,Z_0)\tens T^ns\ca(Z_0,Z_n))
\\
\rTTo^{\uCom(1,b^\cm_{n+1})}
\uCom(s\ca(X,Z_0)\tens T^ns\ca(Z_0,Z_n),s\cm(Z_n))
\\
\rTTo^{(\und{\varphi}^\Com)^{-1}}
\uCom(T^ns\ca(Z_0,Z_n),\uCom(s\ca(X,Z_0),s\cm(Z_n)))
\\
\rTTo^{\uCom(1,[-1])}
\uCom(T^ns\ca(Z_0,Z_n),\uCom(\ca(X,Z_0),\cm(Z_n)))
\\
\rTTo^{\uCom(1,s)} \uCom(T^ns\ca(Z_0,Z_n),s\uCom(\ca(X,Z_0),\cm(Z_n)))
\bigr],
\end{multline*}
where \(n\ge0\),
 \(\tau=\left(
\begin{smallmatrix}
0 & 1 & \dots & n & n+1\\
1 & 2 & \dots & n+1 & 0
\end{smallmatrix}
\right)\in\SSS_{n+2}
 \).
An element \(r\in s\cm(X)\) is mapped to an \ainf-transformation
\((r)\Upsilon\) with the components
\begin{align*}
(r)\Upsilon_n: T^ns\ca(Z_0,Z_n) &\to s\uCom(\ca(X,Z_0),\cm(Z_n)), \quad
n\ge0,
\\
z_1\tdt z_n &\mapsto(z_1\tdt z_n\tens r)\Upsilon'_n,
\end{align*}
where
\begin{multline*}
\Upsilon'_n=-\bigl[ T^ns\ca(Z_0,Z_n)\tens s\cm(X) \rTTo^{\coev^\Com}
\\
\uCom(s\ca(X,Z_0),s\ca(X,Z_0)\tens T^ns\ca(Z_0,Z_n)\tens s\cm(X))
\rTTo^{\uCom(1,\tau_c)}
\\
\uCom(s\ca(X,Z_0),s\cm(X)\tens s\ca(X,Z_0)\tens T^ns\ca(Z_0,Z_n))
\rTTo^{\uCom(1,b^\cm_{n+1})}
\\
\uCom(s\ca(X,Z_0),s\cm(Z_n)) \rTTo^{[-1]} \uCom(\ca(X,Z_0),\cm(Z_n))
\rTTo^s s\uCom(\ca(X,Z_0),\cm(Z_n)) \bigr].
\end{multline*}

Since $\mho$ is a chain map by
\lemref{lem-bicomodule-homomorphism-mho-chain-map},
    \ifx\chooseClass1
the map
    \else
vanishing of \eqref{equ-tb0-Delta1Delta-1b01-t} on \(s\cm(X)\) implies that the map
    \fi
\[ \mho_{00}=\Upsilon:(s\cm(X),b^\ce_{00}=s^{-1}d^{Xf}s=b^\cm_0)
\to(s\und\Ainfty(\ca;\uCom)(H^X,f),b^\cq_{00}=B_1)
\]
is a chain map as well. The following result generalizes previously
known \ainf-version of the Yoneda Lemma
\cite[Theorem~9.1]{Fukaya:FloerMirror-II},
\cite[Proposition~A.9]{math.CT/0306018}, and gives the latter if
\(f=H^W\) for \(W\in\Ob\ca\). It can also be found in Seidel's book
\cite[Lemma~2.12]{SeidelP-book-Fukaya}, where it is proven assuming
that the ground ring \(\kk\) is a field, and the proof is based on a
spectral sequence argument. The proof presented
    \ifx\chooseClass1
in archive version \cite{LyuMan-AmodSerre} of this article
    \else
here
    \fi
is considerably
longer than that of Seidel, however it works in the case of an
arbitrary commutative ground ring.

\begin{proposition}\label{prop-Upsilon-homotopy-invertible}
Let \(\ca\) be a unital \ainf-category, let \(X\) be an object of
\(\ca\), and let \(f:\ca\to\uCom\) be a unital \ainf-functor. Then the
map \(\Upsilon\) is homotopy invertible.
\end{proposition}

\begin{proof}
The \ainf-module \(\cm\) corresponding to $f$ is unital by
\propref{prop-unital-ainf-module}. The components of \(f\) are
expressed via the components of \(b^\cm\) as follows ($k\ge1$):
\begin{multline}
f_k=\bigl[ T^ks\ca(Z_0,Z_k) \rTTo^{\coev^\Com}
\uCom(s\cm(Z_0),s\cm(Z_0)\tens T^ks\ca(Z_0,Z_k))
\\
\rTTo^{\uCom(1,b^\cm_k)} \uCom(s\cm(Z_0),s\cm(Z_k)) \rTTo^{[-1]}
\uCom(\cm(Z_0),\cm(Z_k))
\\
\rTTo^s s\uCom(\cm(Z_0),\cm(Z_k)) \bigr].
\label{equ-ainf-functor-via-ainf-module}
\end{multline}

Define a map \(\alpha:s\und\Ainfty(\ca;\uCom)(H^X,f)\to s\cm(X)\) as
follows:
\begin{multline}
\alpha=\bigl[ s\und\Ainfty(\ca;\uCom)(H^X,f) \rTTo^{\pr_0}
s\uCom(\ca(X,X),\cm(X)) \rTTo^{s^{-1}} \uCom(\ca(X,X),\cm(X))
\\
\rTTo^{[1]} \uCom(s\ca(X,X),s\cm(X)) \rTTo^{\uCom(\sS{_X}\uni^\ca_0,1)}
\uCom(\kk,s\cm(X))=s\cm(X) \bigr].
\label{eq-alpha-pr0s1[1]C(i1)}
\end{multline}
The map \(\alpha\) is a chain map. Indeed, \(\pr_0\) is a chain map,
and
\begin{align*}
s^{-1}[1] &\uCom(\sS{_X}\uni^\ca_0,1)b^\cm_0
=s^{-1}[1]\uCom(\sS{_X}\uni^\ca_0,1)\uCom(1,b^\cm_0)
\\
&=s^{-1}[1](-\uCom(1,b^\cm_0)+\uCom(b_1,1))\uCom(\sS{_X}\uni^\ca_0,1)
=s^{-1}[1]m^\uCom_1\uCom(\sS{_X}\uni^\ca_0,1)
\\
&=s^{-1}m^\uCom_1[1]\uCom(\sS{_X}\uni^\ca_0,1)=b^\uCom_1
s^{-1}[1]\uCom(\sS{_X}\uni^\ca_0,1),
\end{align*}
since \(\sS{_X}\uni^\ca_0\) is a chain map, and \([1]\) is a
differential graded functor. Let us compute \(\Upsilon\alpha\):
\begin{multline}
\Upsilon\alpha=\Upsilon_0s^{-1}[1]\uCom(\sS{_X}\uni^\ca_0,1) =-\bigl[
s\cm(X) \rTTo^{\coev^\Com} \uCom(s\ca(X,X),s\ca(X,X)\tens s\cm(X))
\\
\rTTo^{\uCom(1,c)} \uCom(s\ca(X,X),s\cm(X)\tens
s\ca(X,X))\rTTo^{\uCom(1,b^\cm_1)} \uCom(s\ca(X,X),s\cm(X))
\\
\hfill\rTTo^{\uCom(\sS{_X}\uni^\ca_0,1)}
\uCom(\kk,s\cm(X))=s\cm(X)
\bigr]\quad\\
\quad=\bigl[
s\cm(X)\rTTo^{\coev^\Com}\uCom(s\ca(X,X),s\ca(X,X)\tens s\cm(X))
\rTTo^{\uCom(\sS{_X}\uni^\ca_0,1)}\hfill\\
\hfill\uCom(\kk,s\ca(X,X)\tens s\cm(X))
\rTTo^{\uCom(1,cb^\cm_1)}\uCom(\kk,s\cm(X))=s\cm(X)
\bigr]\quad\\
\quad=\bigl[
s\cm(X)\rTTo^{\coev^\Com}\uCom(\kk,\kk\tens s\cm(X))\rTTo^{\uCom(1,\sS{_X}\uni^\ca_0\tens1)}
\uCom(\kk,s\ca(X,X)\tens
s\cm(X))\hfill\\
\hfill\rTTo^{\uCom(1,cb^\cm_1)}\uCom(\kk,s\cm(X))=s\cm(X)
\bigr]\quad\\
\quad=\bigl[
s\cm(X)\rTTo^{1\tens\sS{_X}\uni^\ca_0}s\cm(X)\tens
s\ca(X,X)\rTTo^{b^\cm_1}s\cm(X)
\bigr].\hfill
\label{eq-Upsilon-alpha}
\end{multline}
Since \(\cm\) is a unital \ainf-module by
\propref{prop-unital-ainf-module}, it follows that \(\Upsilon\alpha\) is
homotopic to identity. Let us prove that \(\alpha\Upsilon\) is
homotopy invertible.

The graded $\kk$-module $s\und\Ainfty(\ca;\uCom)(H^X,f)$ is
$V=\prod_{n=0}^{\infty}V_n$, where
\[ V_n=\prod_{Z_0,\,\dots,\,Z_n\in\Ob\ca}
\uCom(\bar{T}^ns\ca(Z_0,Z_n),s\uCom(\ca(X,Z_0),\cm(Z_n)))
\]
and all products are taken in the category of graded $\kk$\n-modules.
In other terms, for $d\in\ZZ$, \(V^d = \prod_{n=0}^{\infty}V_n^d\),
where
\[
V_n^d = \prod_{Z_0,\,\dots,\,Z_n\in\Ob\ca}
\uCom(\bar{T}^ns\ca(Z_0,Z_n),s\uCom(\ca(X,Z_0),\cm(Z_n)))^d.
\]
We consider $V_n^d$ as Abelian groups with discrete topology. The
Abelian group $V^d$ is equipped with the topology of the product. Thus,
its basis of neighborhoods of 0 is given by \(\kk\)\n-submodules
$\Phi_m^d=0^{m-1}\times\prod_{n=m}^{\infty}V_n^d$. They form a
filtration $V^d=\Phi_0^d\supset\Phi_1^d\supset\Phi_2^d\supset\dots$. We
call a \(\kk\)\n-linear map $a:V\to V$ of degree $p$ continuous if the
induced maps \(a^{d,d+p}=a\big|_{V^d}:V^d\to V^{d+p}\) are continuous
for all $d\in\ZZ$. This holds if and only if for any $d\in\ZZ$ and
$m\in\NN\overset{\text{def}}=\ZZ_{\ge0}$ there exists an integer
$\kappa=\kappa(d,m)\in\NN$ such that
$(\Phi_\kappa^d)a\subset\Phi_m^{d+p}$. We may assume that
\begin{equation}
m'<m'' \quad \text{ implies } \quad \kappa(d,m') \le \kappa(d,m'').
\label{eq-mm-kk-inequa}
\end{equation}
Indeed, a given function $m\mapsto\kappa(d,m)$ can be replaced with the
function $m\mapsto\kappa'(d,m)=\min_{n\ge m}\kappa(d,n)$ and
\(\kappa'\) satisfies condition~\eqref{eq-mm-kk-inequa}. Continuous
linear maps $a:V\to V$ of degree $p$ are in bijection with families of
$\NN\times\NN$-matrices $(A^{d,d+p})_{d\in\ZZ}$ of linear maps
$A_{nm}^{d,d+p}:V_n^d\to V_m^{d+p}$ with finite number of non-vanishing
elements in each column of $A^{d,d+p}$. Indeed, to each continuous map
\(a^{d,d+p}:V^d\to V^{d+p}\) corresponds the inductive limit over $m$
of $\kappa(d,m)\times m$-matrices of maps
$V^d/\Phi_{\kappa(d,m)}^d\to V^{d+p}/\Phi_m^{d+p}$. On the other hand,
to each family $(A^{d,d+p})_{d\in\ZZ}$ of $\NN\times\NN$-matrices with
finite number of non-vanishing elements in each column correspond
obvious maps $a^{d,d+p}:V^d\to V^{d+p}$, and they are continuous. Thus,
\(a=(a^{d,d+p})_{d\in\ZZ}\) is continuous. A continuous map
\(a:V\to V\) can be completely recovered from one $\NN\times\NN$-matrix
\((a_{nm})_{n,m\in\NN}\) of maps
$a_{nm}=(A_{nm}^{d,d+p})_{d\in\ZZ}:V_n\to V_m$ of degree $p$.
Naturally, not any such matrix determines a continuous map, however, if
the number of non-vanishing elements in each column of \((a_{nm})\) is
finite, then this matrix does determine a continuous map.

The differential $D\overset{\text{def}}=B_1:V\to V$,
$r\mapsto(r)B_1=rb-(-)^rbr$ is continuous and the function \(\kappa\)
for it is simply \(\kappa(d,m)=m\). Its matrix is given by
\begin{align*}
D = B_1 =
\begin{bmatrix}
    D_{0,0} & D_{0,1} & D_{0,2} & \dots\\
     0 & D_{1,1} & D_{1,2} & \dots\\
     0 & 0 & D_{2,2} & \dots\\
     \vdots & \vdots & \vdots &\ddots
\end{bmatrix}
,
\end{align*}
where
\begin{gather*}
D_{k,k} = \uCom(1,b^\uCom_1) -
\uCom\bigl(\sum_{p+1+q=k}
1^{\tens p}\tens b_1\tens1^{\tens q},1\bigr): V_k \to V_k, \\
r_kD_{k,k} = r_kb_1^{\uCom} -(-)^r\sum_{p+1+q=k}
(1^{\tens p}\tens b_1\tens 1^{\tens q})r_k,
\end{gather*}
(one easily recognizes the differential in the complex $V_k$),
\[ r_kD_{k,k+1} = (r_k\tens f_1)b^{\uCom}_2
+(H_1^X\tens r_k)b^{\uCom}_2 -(-)^r\sum_{p+q=k-1}
(1^{\tens p}\tens b_2\tens1^{\tens q})r_k.
\]
Further we shall see that we do not need to compute other components.

Composition of \(\alpha\Upsilon\) with \(\pr_n\) equals
\begin{multline*}
\alpha\Upsilon_n=-\bigl[
s\und\Ainfty(\ca;\uCom)(H^X,f)\rTTo^{\pr_0}s\uCom(\ca(X,X),\cm(X))
\\
\rTTo^{s^{-1}[1]}\uCom(s\ca(X,X),s\cm(X))
\rTTo^{\uCom(\sS{_X}\uni^\ca_0,1)}\uCom(\kk,s\cm(X))=s\cm(X)
\\
\rTTo^{\coev^\Com}
 \uCom(s\ca(X,Z_0)\tens T^ns\ca(Z_0,Z_n),
s\ca(X,Z_0)\tens T^ns\ca(Z_0,Z_n)\tens s\cm(X))
\\
\rTTo^{\uCom(1,\tau_c)}
 \uCom(s\ca(X,Z_0)\tens T^ns\ca(Z_0,Z_n),
s\cm(X)\tens s\ca(X,Z_0)\tens T^ns\ca(Z_0,Z_n))
\\
\rTTo^{\uCom(1,b^\cm_{n+1})}
\uCom(s\ca(X,Z_0)\tens T^ns\ca(Z_0,Z_n),s\cm(Z_n))
\\
\rTTo^{(\und{\varphi}^\uCom)^{-1}}
\uCom(T^ns\ca(Z_0,Z_n),\uCom(s\ca(X,Z_0),s\cm(Z_n)))
\\
\rTTo^{\uCom(1,[-1]s)}
\uCom(T^ns\ca(Z_0,Z_n),s\uCom(\ca(X,Z_0),\cm(Z_n))) \bigr].
\end{multline*}
Clearly, $\alpha\Upsilon$ is continuous (take $\kappa(d,m)=1$). Its
$\NN\times\NN$-matrix has the form
\[ \alpha\Upsilon =
\begin{bmatrix}
    * & * & * &\dots \\
    0 & 0& 0& \dots\\
    0 & 0 & 0 &\dots\\
    \vdots & \vdots &\vdots & \ddots
\end{bmatrix}
 .
\]

\begin{lemma}\label{lem-alphaY1-homotopy-upper-triangular}
The map $\alpha\Upsilon:V\to V$ is homotopic to a continuous map $V\to V$,
whose $\NN\times\NN$-matrix is upper-triangular with the identity maps
$\id:V_k\to V_k$ on the diagonal.
\end{lemma}

    \ifx\chooseClass1
The proof is analogous to proof of \cite[Lemma~A.10]{math.CT/0306018}.
\proofInArXiv.
    \else
\begin{proof}
Define a continuous $\kk$\n-linear map
$K:s\und\Ainfty(\ca;\uCom)(H^X,f)\to s\und\Ainfty(\ca;\uCom)(H^X,f)$ of
degree $-1$ by its matrix
\[ K =
\begin{bmatrix}
    0 & 0 & 0 &\dots \\
    K_{1,0} & 0& 0& \dots\\
    0 & K_{2,1} & 0 &\dots\\
    \vdots & \vdots &\vdots & \ddots
\end{bmatrix}
,
\]
so $\kappa(d,m)=m+1$, where $K_{k+1,k}$ maps the factor indexed by
$(X,Z_0,\dots,Z_k)$ to the factor indexed by $(Z_0,\dots,Z_k)$ as
follows:
\begin{align*}
K_{k+1,k} = \bigl[
&\uCom(s\ca(X,Z_0)\tens\bT s\ca(Z_0,Z_k),s\uCom(\ca(X,X),\cm(Z_k)))
\\
\rTTo^{\uCom(1,s^{-1}[1])} &
\uCom(s\ca(X,Z_0)\tens\bT s\ca(Z_0,Z_k),\uCom(s\ca(X,X),s\cm(Z_k)))
\\
\rTTo^{\uCom(1,\uCom(\unix,1))} &
\uCom(s\ca(X,Z_0)\tens\bT s\ca(Z_0,Z_k),\uCom(\kk,s\cm(Z_k)))
\\
\rTTo^{(\und{\varphi}^\Com)^{-1}}& \uCom(\bT s\ca(Z_0,Z_k),\uCom(s\ca(X,Z_0),s\cm(Z_k)))
\\
\rTTo^{\uCom(1,[-1]s)}&
\uCom(\bT s\ca(Z_0,Z_k),s\uCom(\ca(X,Z_0),\cm(Z_k))) \bigr].
\end{align*}
Other factors are ignored.

The composition of continuous maps \(V\to V\) is continuous as well. In particular,
one finds the matrices of $B_1K$ and $KB_1$:
\begin{gather*}
B_1K =
\begin{bmatrix}
D_{0,1}K_{1,0} & D_{0,2}K_{2,1} & D_{0,3}K_{3,2} &\dots\\
D_{1,1}K_{1,0} & D_{1,2}K_{2,1} & D_{1,3}K_{3,2} &\dots\\
0 & D_{2,2}K_{2,1} & D_{2,3}K_{3,2} &\dots\\
0 & 0 &D_{3,3}K_{3,2}& \dots\\
\vdots & \vdots & \vdots &\ddots
\end{bmatrix}
,\\
KB_1 =
\begin{bmatrix}
0 & 0 & 0 & \dots \\
K_{1,0}D_{0,0} & K_{1,0}D_{0,1} & K_{1,0}D_{0,2} &\dots\\
0 & K_{2,1}D_{1,1} & K_{2,1}D_{1,2} &\dots\\
0 & 0 & K_{3,2}D_{2,2} &\dots\\
\vdots & \vdots & \vdots &\ddots
\end{bmatrix}
 .
\end{gather*}
We have $D_{k+1,k+1}K_{k+1,k}+K_{k+1,k}D_{k,k}=0$ for all $k\ge0$.
Indeed, conjugating the expanded left hand side with $\uCom(1,[-1]s)$ we
come to the following identity:
\begin{multline*}
\bigl[ \uCom(s\ca(X,Z_0)\tens\bT s\ca(Z_0,Z_k),
\uCom(s\ca(X,X),s\cm(Z_k)))
\\
\rTTo^{\uCom(1,m^\uCom_1)
 +\uCom(\sum_{p+q=k}1^{\tens p}\tens b_1\tens1^{\tens q},1)}
\\
\uCom(s\ca(X,Z_0)\tens\bT s\ca(Z_0,Z_k),
\uCom(s\ca(X,X),s\cm(Z_k)))
\\
\rTTo^{\uCom(1,\uCom(\unix,1))}
\uCom(s\ca(X,Z_0)\tens\bT s\ca(Z_0,Z_k),
\uCom(\kk,s\cm(Z_k)))
\\
\hfill\rTTo^{(\und{\varphi}^\uCom)^{-1}}
\uCom(\bT s\ca(Z_0,Z_k),\uCom(s\ca(X,Z_0),s\cm(Z_k)))
\bigr]\quad
\\
\quad+\bigl[
\uCom(s\ca(X,Z_0)\tens\bT s\ca(Z_0,Z_k),\uCom(s\ca(X,X),s\cm(Z_k)))\hfill
\\
\rTTo^{\uCom(1,\uCom(\unix,1))}
\uCom(s\ca(X,Z_0)\tens\bT s\ca(Z_0,Z_k),\uCom(\kk,s\cm(Z_k)))
\\
\rTTo^{(\und{\varphi}^\uCom)^{-1}}
\uCom(\bT s\ca(Z_0,Z_k),\uCom(s\ca(X,Z_0),s\cm(Z_k)))
\\
\rTTo^{\uCom(1,m^\uCom_1)
 +\uCom(\sum_{p+q=k-1}1^{\tens p}\tens b_1\tens1^{\tens q},1)}
\uCom(\bT s\ca(Z_0,Z_k),\uCom(s\ca(X,Z_0),s\cm(Z_k))) \bigr]=0.
\end{multline*}
After reducing all terms to the common form, beginning with
$\uCom(1,\uCom(\unix,1))$, all terms cancel each other, so the identity
is proven.

Therefore, the chain map $a=\alpha\Upsilon+B_1K+KB_1$ is also
represented by an upper-triangular matrix. Its diagonal elements are
chain maps $a_{kk}:V_k\to V_k$. We are going to show that they are
homotopic to identity maps.

Let us compute the matrix element
 $a_{00}:V_0\to V_0=\prod_{Z\in\Ob\ca}s\uCom(\ca(X,Z),\cm(Z))$. We have
\[ a_{00}\pr_Z = (\alpha\Upsilon_0+B_1K_{1,0})\pr_Z:
V_0\to V_0 \rTTo^{\pr_Z} s\uCom(\ca(X,Z),\cm(Z)).
\]
In the expanded form these terms are as follows:
\begin{multline*}
\alpha\Upsilon_0\pr_Z=-\bigl[
V_0\rTTo^{\pr_X}s\uCom(\ca(X,X),\cm(X))\rTTo^{s^{-1}[1]}\uCom(s\ca(X,X),s\cm(X))\\
\rTTo^{\uCom(\unix,1)}\uCom(\kk,s\cm(X))=s\cm(X)\rTTo^{\coev^\Com}
\uCom(s\ca(X,Z),s\ca(X,Z)\tens s\cm(X))\\
\rTTo^{\uCom(1,cb^\cm_1)}
\uCom(s\ca(X,Z),s\cm(Z))\rTTo^{[-1]s}s\uCom(\ca(X,Z),\cm(Z))
\bigr],
\end{multline*}
\[
B_1K_{1,0}\pr_Z=[(1\tens f_1)b^\uCom_2+(H^X_1\tens1)b^\uCom_2]K_{1,0}\pr_Z,
\]
\begin{multline*}
(1\tens f_1)b^\uCom_2=\bigl[ V_0\rTTo^{\pr_X} s\uCom(\ca(X,X),\cm(X))
\\
\quad\rTTo^{\coev^\Com} \uCom(s\ca(X,Z),s\ca(X,Z)\tens s\uCom(\ca(X,X),\cm(X)))\hfill
\\
\rTTo^{\uCom(1,c)} \uCom(s\ca(X,Z),s\uCom(\ca(X,X),\cm(X))\tens
s\ca(X,Z))
\\
\hfill\rTTo^{\uCom(1,1\tens f_1)}
\uCom(s\ca(X,Z),s\uCom(\ca(X,X),\cm(X))\tens s\uCom(\cm(X),\cm(Z)))
\quad
\\
\rTTo^{\uCom(1,b^\uCom_2)} \uCom(s\ca(X,Z),s\uCom(\ca(X,X),\cm(Z)))
\bigr],
\end{multline*}
\begin{multline*}
(H^X_1\tens1)b^\uCom_2=\bigl[
V_0\rTTo^{\pr_Z}s\uCom(\ca(X,Z),\cm(Z))\\
\quad\rTTo^{\coev^\Com}\uCom(s\ca(X,Z),s\ca(X,Z)\tens s\uCom(\ca(X,Z),\cm(Z)))\hfill\\
\rTTo^{\uCom(1,H^X_1\tens1)}
\uCom(s\ca(X,Z),s\uCom(\ca(X,X),\ca(X,Z))\tens s\uCom(\ca(X,Z),\cm(Z)))\\
\rTTo^{\uCom(1,b^\uCom_2)}
\uCom(s\ca(X,Z),s\uCom(\ca(X,X),\cm(Z)))
\bigr],
\end{multline*}
\begin{multline*}
K_{1,0}\pr_Z=\bigl[
V_1\rTTo^{\pr_{X,Z}}\uCom(s\ca(X,Z),s\uCom(\ca(X,X),\cm(Z)))
\rTTo^{\uCom(1,s^{-1}[1])}\\
\uCom(s\ca(X,Z),\uCom(s\ca(X,X),s\cm(Z)))
\rTTo^{\uCom(1,\uCom(\unix,1))}
\uCom(s\ca(X,Z),\uCom(\kk,s\cm(Z)))\\
=\uCom(s\ca(X,Z),s\cm(Z))\rTTo^{[-1]s}s\uCom(\ca(X,Z),\cm(Z))
\bigr].
\end{multline*}
We claim that in the sum
\[ a_{00}\pr_Z = \alpha\Upsilon_0\pr_Z + (1\tens f_1)b^\uCom_2K_{1,0}\pr_Z
+ (H^X_1\tens1)b^\uCom_2K_{1,0}\pr_Z
\]
the first two summands cancel each other, while the last,
$(H^X_1\tens1)b^\uCom_2K_{1,0}$, is homotopic to identity. Indeed,
\(\alpha\Upsilon_0\pr_Z + (1\tens f_1)b^\uCom_2K_{1,0}\pr_Z\) factors
through
\begin{multline*}
-\bigl[
s\uCom(\ca(X,X),\cm(X))\rTTo^{s^{-1}[1]}\uCom(s\ca(X,X),s\cm(X))
\rTTo^{\uCom(\unix,1)}\uCom(\kk,s\cm(X))
\\
=s\cm(X) \rTTo^{\coev^\Com} \uCom(s\ca(X,Z),s\ca(X,Z)\tens s\cm(X))
\rTTo^{\uCom(1,cb^\cm_1)} \uCom(s\ca(X,Z),s\cm(Z)) \bigr]
\\
\quad+\bigl[
s\uCom(\ca(X,X),\cm(X))\rTTo^{\coev^\Com}
\uCom(s\ca(X,Z),s\ca(X,Z)\tens s\uCom(\ca(X,X),\cm(X)))\hfill\\
\rTTo^{\uCom(1,c)}\uCom(s\ca(X,Z),s\uCom(\ca(X,X),\cm(X))\tens s\ca(X,Z))
\rTTo^{\uCom(1,1\tens f_1)}\\
\uCom(s\ca(X,Z),s\uCom(\ca(X,X),\cm(X))\tens s\uCom(\cm(X),\cm(Z)))
\rTTo^{\uCom(1,b^\uCom_2)}\\
\uCom(s\ca(X,Z),s\uCom(\ca(X,X),\cm(Z)))
\rTTo^{\uCom(1,s^{-1}[1])}\uCom(s\ca(X,Z),\uCom(s\ca(X,X),s\cm(Z)))\\
\rTTo^{\uCom(1,\uCom(\unix,1))}
\uCom(s\ca(X,Z),\uCom(\kk,s\cm(Z)))=\uCom(s\ca(X,Z),s\cm(Z))
\bigr].
\end{multline*}
It therefore suffices to prove that the above expression vanishes. By
closedness this is equivalent to the following equation:
\begin{multline*}
-\bigl[
s\ca(X,Z)\tens s\uCom(\ca(X,X),\cm(X))\rTTo^{1\tens s^{-1}[1]}
s\ca(X,Z)\tens\uCom(s\ca(X,X),s\cm(X))\\
\hfill\rTTo^{1\tens\uCom(\unix,1)}
s\ca(X,Z)\tens s\cm(X)\rTTo^{cb^\cm_1}s\cm(Z)
\bigr]\quad\\
\quad+\bigl[
s\ca(X,Z)\tens s\uCom(\ca(X,X),\cm(X))\rTTo^c
s\uCom(\ca(X,X),\cm(X))\tens s\ca(X,Z)\hfill\\
\rTTo^{1\tens f_1}
s\uCom(\ca(X,X),\cm(X))\tens s\uCom(\cm(X),\cm(Z))\rTTo^{b^\uCom_2}
s\uCom(\ca(X,X),\cm(Z))\\
\rTTo^{s^{-1}[1]}\uCom(s\ca(X,X),s\cm(Z))
\rTTo^{\uCom(\unix,1)}\uCom(\kk,s\cm(Z))=s\cm(Z)
\bigr]
=0.
\end{multline*}
Canceling \(c\) and transforming the left hand side using
\eqref{equ-ainf-functor-to-ainf-module} we get:
\begin{multline*}
-\bigl[
s\uCom(\ca(X,X),\cm(X))\tens s\ca(X,Z)\rTTo^{s^{-1}[1]\tens1}
\uCom(s\ca(X,X),s\cm(X))\tens s\ca(X,Z)\\
\rTTo^{1\tens f_1s^{-1}[1]}\uCom(s\ca(X,X),s\cm(X))\tens\uCom(s\cm(X),s\cm(Z))\\
\hfill\rTTo^{m^\uCom_2}
\uCom(s\ca(X,X),s\cm(Z))\rTTo^{\uCom(\unix,1)}\uCom(\kk,s\cm(Z))=s\cm(Z)
\bigr]\quad\\
\quad+\bigl[
s\uCom(\ca(X,X),\cm(X))\tens s\ca(X,Z)\rTTo^{s^{-1}[1]\tens1}
\uCom(s\ca(X,X),s\cm(X))\tens s\ca(X,Z)\hfill\\
\rTTo^{1\tens f_1s^{-1}[1]}
\uCom(s\ca(X,X),s\cm(X))\tens\uCom(s\cm(X),s\cm(Z))\\
\rTTo^{\uCom(\unix,1)\tens1}
\uCom(\kk,s\cm(X))\tens\uCom(s\cm(X),s\cm(Z))\rTTo^{\ev^\Com}s\cm(Z)
\bigr]=0.
\end{multline*}
The above equation follows from the following identity which holds by
the properties of the closed monoidal category \(\uCom\):
\begin{multline*}
\bigl[
\uCom(s\ca(X,X),s\cm(X))\tens\uCom(s\cm(X),s\cm(Z))\rTTo^{m^\uCom_2}\\
\hfill\uCom(s\ca(X,X),s\cm(Z))\rTTo^{\uCom(\unix,1)}\uCom(\kk,s\cm(Z))=s\cm(Z)
\bigr]\quad\\
\quad=\bigl[
\uCom(s\ca(X,X),s\cm(X))\tens\uCom(s\cm(X),s\cm(Z))\rTTo^{\uCom(\unix,1)\tens1}\hfill
\\
\uCom(\kk,s\cm(X))\tens\uCom(s\cm(X),s\cm(Z))=s\cm(X)\tens\uCom(s\cm(X),s\cm(Z))
\rTTo^{\ev^\Com} s\cm(Z) \bigr].
\end{multline*}
This is a particular case of identity~\eqref{eq-m2C(aC)-C(aB)1m2}
combined with \eqref{eq-1mevV-evV1evV}.

Now we prove that $(H^X_1\tens1)b_2K_{1,0}$ is homotopic to identity.
It maps each factor $s\uCom(\ca(X,Z),\cm(Z))$ into itself via the following
map:
\begin{multline*}
\bigl[
s\uCom(\ca(X,Z),\cm(Z))\rTTo^{\coev^\Com}\uCom(s\ca(X,Z),s\ca(X,Z)\tens
s\uCom(\ca(X,Z),\cm(Z)))\\
\rTTo^{\uCom(1,P)}\uCom(s\ca(X,Z),s\cm(Z))\rTTo^{[-1]s}s\uCom(\ca(X,Z),\cm(Z))
\bigr],
\end{multline*}
where
\begin{multline}
P=-\bigl[
s\ca(X,Z)\tens s\uCom(\ca(X,Z),\cm(Z))\rTTo^{H^X_1s^{-1}[1]\tens
s^{-1}[1]}\\
\uCom(s\ca(X,X),s\ca(X,Z))\tens\uCom(s\ca(X,Z),s\cm(Z))\rTTo^{m^\uCom_2}\\
\hfill\uCom(s\ca(X,X),s\cm(Z))\rTTo^{\uCom(\unix,1)}\uCom(\kk,s\cm(Z))=s\cm(Z)
\bigr]\quad\\
\quad=\bigl[
s\ca(X,Z)\tens s\uCom(\ca(X,Z),\cm(Z))\rTTo^{1\tens s^{-1}[1]}
s\ca(X,Z)\tens\uCom(s\ca(X,Z),s\cm(Z))\hfill\\
\rTTo^{\coev^\Com\tens1}\uCom(s\ca(X,X),s\ca(X,X)\tens
s\ca(X,Z))\tens\uCom(s\ca(X,Z),s\cm(Z))\\
\rTTo^{\uCom(1,b_2)\tens1}
\uCom(s\ca(X,X),s\ca(X,Z))\tens\uCom(s\ca(X,Z),s\cm(Z))\\
\hfill\rTTo^{\uCom(\unix,1)\tens1}
\uCom(\kk,s\ca(X,Z))\tens\uCom(s\ca(X,Z),s\cm(Z))\rTTo^{\ev^\Com}s\cm(Z)
\bigr]\quad\\
\quad=-\bigl[
s\ca(X,Z)\tens s\uCom(\ca(X,Z),\cm(Z))\rTTo^{1\tens s^{-1}[1]}
s\ca(X,Z)\tens\uCom(s\ca(X,Z),s\cm(Z))\hfill\\
\rTTo^{\coev^\Com\tens1}
\uCom(\kk,\kk\tens s\ca(X,Z))\tens\uCom(s\ca(X,Z),s\cm(Z))\\
\hfill\rTTo^{\uCom(1,(\unix\tens1)b_2)\tens1}
s\ca(X,Z)\tens\uCom(s\ca(X,Z),s\cm(Z))\rTTo^{\ev^\Com} s\cm(Z)
\bigr]\quad\\
\quad=-\bigl[
s\ca(X,Z)\tens s\uCom(\ca(X,Z),\cm(Z))\rTTo^{1\tens s^{-1}[1]}
s\ca(X,Z)\tens\uCom(s\ca(X,Z),s\cm(Z))\hfill\\
\rTTo^{(\unix\tens1)b_2\tens1}
s\ca(X,Z)\tens\uCom(s\ca(X,Z),s\cm(Z))\rTTo^{\ev^\Com}s\cm(Z)
\bigr].
\label{equ-H-1-s-shift-Com-unix}
\end{multline}
It follows that
\begin{multline*}
(H^X_1\tens1)b^\uCom_2K_{1,0}=-\bigl[
s\uCom(\ca(X,Z),\cm(Z))\rTTo^{s^{-1}[1]}
\uCom(s\ca(X,Z),s\cm(Z))\rTTo^{\coev^\Com}\\
\uCom(s\ca(X,Z),s\ca(X,Z)\tens\uCom(s\ca(X,Z),s\cm(Z)))
\rTTo^{\uCom(1,(\unix\tens1)b_2\tens1)}\\
\uCom(s\ca(X,Z),s\ca(X,Z)\tens\uCom(s\ca(X,Z),s\cm(Z)))
\rTTo^{\uCom(1,\ev^\Com)}\\
\uCom(s\ca(X,Z),s\cm(Z))\rTTo^{[-1]s}
s\uCom(\ca(X,Z),\cm(Z))
\bigr].
\end{multline*}
Due to \cite[Lemma~7.4]{Lyu-AinfCat} there exists a homotopy
$h'':s\ca(X,Z)\to s\ca(X,Z)$, a map of degree $-1$, such that
$(\unix\tens1)b_2=-1+h''b_1+b_1h''$. Therefore, the map considered
above equals to
\begin{multline*}
\id_{s\uCom(\ca(X,Z),\cm(Z))} - \bigl[ s\uCom(\ca(X,Z),\cm(Z)) \rTTo^{s^{-1}[1]}
\uCom(s\ca(X,Z),s\cm(Z))\\
\hfill\rTTo^{\coev^\Com}\uCom(s\ca(X,Z),s\ca(X,Z)\tens\uCom(s\ca(X,Z),s\cm(Z)))\quad
\\
\hfill \rTTo^{\uCom(1,(h''b_1+b_1h'')\tens1)}
\uCom(s\ca(X,Z),s\ca(X,Z)\tens\uCom(s\ca(X,Z),s\cm(Z))) \quad
\\
\hfill\rTTo^{\uCom(1,\ev^\Com)} \uCom(s\ca(X,Z),s\cm(Z))
\rTTo^{[-1]s} s\uCom(\ca(X,Z),\cm(Z)) \bigr]\quad
\\
= \bigl( \id_{s\uCom(\ca(X,Z),\cm(Z))} +b^\uCom_1K'_{00} +K'_{00}b^\uCom_1
\bigr),
\end{multline*}
where
\begin{multline*}
 K'_{00} = \bigl[ s\uCom(\ca(X,Z),\cm(Z)) \rTTo^{s^{-1}[1]}
\uCom(s\ca(X,Z),s\cm(Z))\\ \rTTo^{\uCom(h'',1)} \uCom(s\ca(X,Z),s\cm(Z)) \rTTo^{[-1]s}
s\uCom(\ca(X,Z),\cm(Z)) \bigr].
\end{multline*}
Indeed, \(b^\uCom_1K'_{00}+K'_{00}b^\uCom_1\) is obtained by conjugating
with \([-1]s\) the following expression:
\begin{multline}
m^\uCom_1\uCom(h'',1)+\uCom(h'',1)m^\uCom_1=(-\uCom(1,b^\cm_0)+\uCom(b_1,1))\uCom(h'',1)\\
+\uCom(h'',1)(-\uCom(1,b^\cm_0)+\uCom(b_1,1))=-\uCom(b_1h''+h''b_1,1).
\label{equ-m-Com-h''-1}
\end{multline}
The rest is straightforward. Therefore, $(H^X_1\tens1)b^\uCom_2K_{1,0}$ and
$a_{00}$ are homotopic to identity.

Now we are proving that diagonal elements $a_{kk}:V_k\to V_k$ are
homotopic to identity maps for $k>0$. An element $r_{k+1}\in V_{k+1}$ is mapped to direct product over
$Z_0,\dots,Z_k\in\Ob\ca$ of
\begin{multline*}
r_{k+1}K_{k+1,k} = \coev^\Com_{s\ca(X,Z_0),*}
\uCom(s\ca(X,Z_0),r_{k+1}^{X,Z_0,\dots,Z_k}s^{-1}[1]\uCom(\unix,1))[-1]s:\\
\bT s\ca(Z_0,Z_k) \to s\uCom(\ca(X,Z_0),\cm(Z_k)).
\end{multline*}
Here
$r_{k+1}^{X,Z_0,\dots,Z_k}:
s\ca(X,Z_0)\tens\bT s\ca(Z_0,Z_k)\to s\uCom(\ca(X,X),\cm(Z_k))$
is one of the coordinates of $r_{k+1}$.

Thus, $r_{k}D_{k,k+1}K_{k+1,k}$ is the sum of three terms
\eqref{eq-rBHa}--\eqref{eq-rBHc}:
\begin{subequations}
\begin{multline}
\coev^\Com_{s\ca(X,Z_0),*}\uCom(s\ca(X,Z_0),(r_k^{X,Z_0,\dots,Z_{k-1}}\tens f_1)
b_2^{\uCom}s^{-1}[1]\uCom(\unix,1))[-1]s: \\
\bT s\ca(Z_0,Z_k) \to s\uCom(\ca(X,Z_0),\cm(Z_k)),
\label{eq-rBHa}
\end{multline}
where
 \(r_k^{X,Z_0,\dots,Z_{k-1}}:s\ca(X,Z_0)\tens\bT s\ca(Z_0,Z_{k-1})\to
 s\uCom(\ca(X,X),\cm(Z_{k-1}))\).
Since $b_2^{\uCom}=(s\tens s)^{-1}m_2s=-(s^{-1}\tens s^{-1})m_2s$, and
\([1]\) is a differential graded functor, we have
\[ (r_k\tens f_1)b^{\uCom}_2s^{-1}[1]
=-(r_ks^{-1}[1]\tens f_1s^{-1}[1])m_2^\uCom.
\]
Identity~\eqref{eq-m2C(aC)-C(aB)1m2} gives
\begin{multline*}
(r_k\tens f_1)b^{\uCom}_2s^{-1}[1]\uCom(\unix,1) =
-(r_ks^{-1}[1]\tens f_1s^{-1}[1])m_2^\uCom\uCom(\unix,1) \\
= -(r_ks^{-1}[1]\tens f_1s^{-1}[1])(\uCom(\unix,1)\tens 1)m_2^\uCom
= (r_ks^{-1}[1]\uCom(\unix,1)\tens f_1s^{-1}[1])m_2^\uCom
\end{multline*}
(we have used the fact that $\uCom(\unix,1)$ has degree $-1$ and
$f_1s^{-1}$ has degree $1$).
\begin{multline}
\coev^\Com_{s\ca(X,Z_0),*}\uCom(s\ca(X,Z_0),(H^{X}_1\tens r_k^{Z_0,\dots,Z_k})
b_2^{\uCom}s^{-1}[1]\uCom(\unix,1))[-1]s: \\
\bT s\ca(Z_0,Z_k) \to s\uCom(\ca(X,Z_0),\cm(Z_k)),
\label{eq-rBHb}
\end{multline}
where
$r_k^{Z_0,\dots,Z_k}:
\bT s\ca(Z_0,Z_k)\to s\uCom(\ca(X,Z_0),\cm(Z_k))$.
Similarly to above
\begin{equation*}
(H^{X}_1\tens r_k)b^{\uCom}_2s^{-1}[1]
= (-)^{r+1}(H^X_1s^{-1}[1]\tens r_ks^{-1}[1])m_2^\uCom,
\end{equation*}
so that
\begin{multline*}
(H^{X}_1\tens r_k)b^{\uCom}_2s^{-1}[1]\uCom(\unix,1)
= (-)^{r+1}(H^X_1s^{-1}[1]\tens r_ks^{-1}[1])
(\uCom(\unix,1)\tens1)m_2^\uCom \\
= (H^X_1s^{-1}[1]\uCom(\unix,1)\tens r_ks^{-1}[1])m_2^\uCom
\end{multline*}
($r_ks^{-1}[1]$ has degree $\deg r+1$ and $\uCom(\unix,1)$ has degree
$-1$).

For each $p$, $q$, such that $p+q=k-1$
\begin{multline}
\coev^\Com_{s\ca(X,Z_0),*}
\uCom(s\ca(X,Z_0),(1^{\tens p}\tens b_2\tens1^{\tens q})r_k
s^{-1}[1]\uCom(\unix,1))[-1]s: \\
\bT s\ca(Z_0,Z_k)\to s\uCom(\ca(X,Z_0),\cm(Z_k)),
\label{eq-rBHc}
\end{multline}
where $r_k$ means
\[ r_k^{X,Z_0,\dots,Z_{p-1},Z_{p+1},\dots,Z_k}:
 \bT s\ca(X,Z_0,\dots,Z_{p-1},Z_{p+1},\dots,Z_k)
 \to s\uCom(\ca(X,X),\cm(Z_k)),
\]
and $Z_{-1}=X$.
\end{subequations}

Thus,
$r_kD_{k,k+1}K_{k+1,k}=\coev^\Com_{s\ca(X,Z_0),*}\uCom(s\ca(X,Z_0),\Sigma_1)[-1]s$,
where
\begin{multline*}
\Sigma_1 = (r_k^{X,Z_0,\dots,Z_{k-1}}s^{-1}[1]
\uCom(\unix,1)\tens f_1s^{-1}[1])m_2^\uCom\\
\qquad+(H^X_1s^{-1}[1]\uCom(\unix,1)\tens r_k^{Z_0,\dots,Z_k}s^{-1}[1])
m_2^\uCom\hfill\\
\qquad-(-)^r\sum_{p+q=k-1}
(1^{\tens p}\tens b_2\tens1^{\tens q})
r_k^{X,Z_0,\dots,Z_{p-1},Z_{p+1},\dots,Z_k}
s^{-1}[1]\uCom(\unix,1):\hfill \\
s\ca(X,Z_0)\tens\bT s\ca(Z_0,Z_k)\to s\cm(Z_k).
\end{multline*}

Similarly, $r_kK_{k,k-1}D_{k-1,k}$ is the sum of three terms
\eqref{eq-rHBa}--\eqref{eq-rHBc}.
\begin{subequations}
\begin{multline}
\bigl(\coev^\Com_{s\ca(X,Z_0),*}\uCom(s\ca(X,Z_0),r_ks^{-1}[1]\uCom(\unix,1))[-1]s\\
\hfill\tens\coev^\Com_{s\cm(Z_{k-1}),s\ca(Z_{k-1},Z_k)}\uCom(s\cm(Z_{k-1}),b^\cm_1)[-1]s\bigr)
b^{\uCom}_2\quad \\
\quad=-\bigl(\coev^\Com_{s\ca(X,Z_0),*}\uCom(s\ca(X,Z_0),r_ks^{-1}[1]\uCom(\unix,1))\hfill\\
\hfill\tens\coev^\Com_{s\cm(Z_{k-1}),s\ca(Z_{k-1},Z_k)}\uCom(s\cm(Z_{k-1}),b^\cm_1)\bigr)m^\uCom_2[-1]s\quad \\
\quad=-\coev^\Com_{s\ca(X,Z_0),*}\uCom\bigl(s\ca(X,Z_0),(r_ks^{-1}[1]\uCom(\unix,1)\tens1)
b^\cm_1\bigr)[-1]s:\hfill \\
\bT s\ca(Z_0,Z_k) \to s\uCom(\ca(X,Z_0),\cm(Z_k)),
\label{eq-rHBa}
\end{multline}
where $r_k$ means
 \(r_k^{X,Z_0,\dots,Z_{k-1}}:\bT s\ca(X,Z_0,\dots,Z_{k-1})
 \to s\uCom(\ca(X,X),\cm(Z_{k-1}))\).
Here we use formulas \eqref{eq-identity-m2} and
\eqref{equ-ainf-functor-via-ainf-module}.
\begin{multline}
\bigl(\coev^\Com_{s\ca(X,Z_0),s\ca(Z_0,Z_1)}\uCom(\ca(X,Z_0),b_2)[-1]s\\
\hfill\tens\coev^\Com_{s\ca(X,Z_1),*}\uCom(s\ca(X,Z_1),r_ks^{-1}[1]\uCom(\unix,1))
[-1]s\bigr)b^{\uCom}_2\quad \\
\quad=(-)^r\bigl(\coev^\Com_{s\ca(X,Z_0),s\ca(Z_0,Z_1)}\uCom(s\ca(X,Z_0),b_2)\hfill\\
\hfill\tens\coev^\Com_{s\ca(X,Z_1),*}\uCom(s\ca(X,Z_1),r_ks^{-1}[1]\uCom(\unix,1))\bigr)
m^\uCom_2[-1]s\quad \\
\quad=(-)^r\coev^\Com_{s\ca(X,Z_0),*}\uCom\bigl(s\ca(X,Z_0),
(b_2\tens1^{\tens k-1})r_ks^{-1}[1]\uCom(\unix,1)\bigr)[-1]s:\hfill \\
\bT s\ca(Z_0,Z_k) \to s\uCom(\ca(X,Z_0),\cm(Z_k)),
\label{eq-rHBb}
\end{multline}
where $r_k$ means
 \(r_k^{X,Z_1,\dots,Z_k}:\bT s\ca(X,Z_1,\dots,Z_k)\to
 s\uCom(\ca(X,X),\cm(Z_k))\).
We have used that $r_ks^{-1}[1]\uCom(\unix,1)$ has degree $\deg r$ and
formula~\eqref{eq-identity-m2}.
\begin{multline}
(1^{\tens p}\tens b_2\tens1^{\tens q})
\coev^\Com_{s\ca(X,Z_0),*}\uCom\bigl(s\ca(X,Z_0),r_ks^{-1}[1]\uCom(\unix,1)\bigr)[-1]s
\\
\quad
 =\coev^\Com_{s\ca(X,Z_0),*}\uCom(s\ca(X,Z_0),1^{\tens p+1}\tens b_2\tens1^{\tens q})
\hfill
\\
\hfill \uCom\bigl(s\ca(X,Z_0),r_ks^{-1}[1]\uCom(\unix,1)\bigr)[-1]s
\quad
\\
\quad =\coev^\Com_{s\ca(X,Z_0),*}
 \uCom\bigl(s\ca(X,Z_0),(1^{\tens p+1}\tens b_2\tens1^{\tens q})
r_ks^{-1}[1]\uCom(\unix,1)\bigr)[-1]s:
\hfill
\\
\bT s\ca(Z_0,Z_k) \to s\uCom(\ca(X,Z_0),\cm(Z_k)),
\label{eq-rHBc}
\end{multline}
where $r_k$ is the map
\[ r_k^{X,Z_0,\dots,Z_{p},Z_{p+2},\dots,Z_k}:
 \bT s\ca(X,Z_0,\dots,Z_{p},Z_{p+2},\dots,Z_k)
 \to s\uCom(\ca(X,X),\cm(Z_k)).
\]
Here we use functoriality of $\coev^\Com$.
\end{subequations}

Thus,
$r_kK_{k,k-1}D_{k-1,k}=\coev^\Com_{s\ca(X,Z_0),*}\uCom(h^XZ_0,\Sigma_2)[-1]s$,
where
\begin{multline*}
\Sigma_2 =-(r_k^{X,Z_0,\dots,Z_{k-1}}s^{-1}[1]\uCom(\unix,1)\tens1)b^\cm_1
\\
\hskip\multlinegap\hphantom{\Sigma_2}
+(-)^r(b_2\tens1^{\tens k-1})r_k^{X,Z_1,\dots,Z_k}s^{-1}[1]
\uCom(\unix,1) \hfill
\\
\hskip\multlinegap\hphantom{\Sigma_2}
-(-)^{r-1}\sum_{p+q=k-2} (1^{\tens p+1}\tens b_2\tens1^{\tens q})
r_k^{X,Z_0,\dots,Z_{p},Z_{p+2},\dots,Z_k} s^{-1}[1]\uCom(\unix,1):
\hfill
\\
s\ca(X,Z_0)\tens\bT s\ca(Z_0,Z_k) \to s\cm(Z_k).
\end{multline*}
The element $rK$ has degree $\deg r-1$, so the sign $(-)^{r-1}$ arises.
Combining this with the expression for $r_kD_{k,k+1}K_{k+1,k}$ we
obtain
\[ r_kD_{k,k+1}K_{k+1,k}+r_kK_{k,k-1}D_{k-1,k}
= \coev^\Com_{s\ca(X,Z_0),*}\uCom(s\ca(X,Z_0),\Sigma)[-1]s,
\]
where $\Sigma=\Sigma_1+\Sigma_2$. We claim that
$\Sigma=(H^X_1s^{-1}[1]\uCom(\unix,1)\tens
r_k^{Z_0,\dots,Z_k}s^{-1}[1])m^\uCom_2$. Indeed, first of all,
\begin{multline*}
(b_2\tens1^{\tens k-1})r_k^{X,Z_1,\dots,Z_k}s^{-1}[1]\uCom(\unix,1) \\
+\sum_{p+q=k-2}(1^{\tens p+1}\tens b_2\tens1^{\tens q})
r_k^{X,Z_0,\dots,Z_p,Z_{p+2},\dots,Z_k}s^{-1}[1]\uCom(\unix,1)\\
=\sum_{p+q=k-1}(1^{\tens p}\tens b_2\tens1^{\tens q})
r_k^{X,Z_0,\dots,Z_{p-1},Z_{p+1},\dots,Z_k}
s^{-1}[1]\uCom(\unix,1),
\end{multline*}
so that
\begin{multline*}
\Sigma = (r_k^{X,Z_0,\dots,Z_{k-1}}s^{-1}[1]\uCom(\unix,1)\tens
f_1s^{-1}[1])m^\uCom_2\\
+(H^X_1s^{-1}\uCom(\unix,1)\tens r_k^{Z_0,\dots,Z_k}s^{-1}[1])
m^\uCom_2
-(r_k^{X,Z_0,\dots,Z_{k-1}}s^{-1}[1]\uCom(\unix,1)\tens 1)b^\cm_1:\\
s\ca(X,Z_0)\tens\bT s\ca(Z_0,Z_k) \to s\cm(Z_k).
\end{multline*}
Note that
\begin{multline*}
m^\uCom_2=\ev^\Com:s\cm(Z_{k-1})\tens\uCom(s\cm(Z_{k-1}),s\cm(Z_k))\\
=\uCom(\kk,s\cm(Z_{k-1}))\tens\uCom(s\cm(Z_{k-1}),s\cm(Z_k))\to
s\cm(Z_k)=\uCom(\kk,s\cm(Z_k)),
\end{multline*}
therefore the first and the third summands cancel out, due to
\eqref{equ-ainf-functor-to-ainf-module}. Hence, only the second summand remains in
$\Sigma=(H^X_1s^{-1}[1]\uCom(\unix,1)\tens
r_k^{Z_0,\dots,Z_k}s^{-1}[1])m^\uCom_2$. It follows that
\begin{multline*}
r_kD_{k,k+1}K_{k+1,k}+r_kK_{k,k-1}D_{k-1,k}=\bigl[
\bT s\ca(Z_0,Z_k)\rTTo^{\coev^\Com}\\
\uCom(s\ca(X,Z_0),s\ca(X,Z_0)\tens\bT s\ca(Z_0,Z_k))
\rTTo^{\uCom(1,H^X_1s^{-1}[1]\uCom(\unix,1)\tens
r_ks^{-1}[1])}\\
\uCom(s\ca(X,Z_0),\uCom(\kk,s\ca(X,Z_0))\tens\uCom(s\ca(X,Z_0),s\cm(Z_k)))
\rTTo^{\uCom(1,m^\uCom_2)}~\|_{\uCom(1,\ev^\Com)}\\
\hfill\uCom(s\ca(X,Z_0),s\cm(Z_k))\rTTo^{[-1]s}s\uCom(\ca(X,Z_0),\cm(Z_k))
\bigr]\quad\\
\quad=\bigl[
\bT s\ca(Z_0,Z_k)\rTTo^{r_ks^{-1}[1]}
\uCom(s\ca(X,Z_0),s\cm(Z_k))\rTTo^{\coev^\Com}\hfill\\
\uCom(s\ca(X,Z_0),s\ca(X,Z_0)\tens\uCom(s\ca(X,Z_0),s\cm(Z_k)))
\rTTo^{\uCom(1,H^X_1s^{-1}[1]\uCom(\unix,1)\tens1)}\\
\uCom(s\ca(X,Z_0),s\ca(X,Z_0)\tens\uCom(s\ca(X,Z_0),s\cm(Z_k)))
\rTTo^{\uCom(1,\ev^\Com)}\\
\hfill\uCom(s\ca(X,Z_0),s\cm(Z_k))\rTTo^{[-1]s}
s\uCom(\ca(X,Z_0),\cm(Z_k))
\bigr]\quad\\
\quad=\bigl[
\bT s\ca(Z_0,Z_k)\rTTo^{r_ks^{-1}[1]}
\uCom(s\ca(X,Z_0),s\cm(Z_k))\rTTo^{\coev^\Com}\hfill\\
\uCom(s\ca(X,Z_0),s\ca(X,Z_0)\tens\uCom(s\ca(X,Z_0),s\cm(Z_k)))
\rTTo^{\uCom(H^X_1s^{-1}[1]\uCom(\unix,1),1)}\\
\uCom(s\ca(X,Z_0),s\ca(X,Z_0)\tens\uCom(s\ca(X,Z_0),s\cm(Z_k)))
\rTTo^{\uCom(1,\ev^\Com)}\\
\hfill\uCom(s\ca(X,Z_0),s\cm(Z_k))\rTTo^{[-1]s}s\uCom(\ca(X,Z_0),\cm(Z_k))
\bigr]\quad
\\
\quad=\bigl[ \bT s\ca(Z_0,Z_k) \rTTo^{r_k} s\uCom(\ca(X,Z_0),\cm(Z_k))
\rTTo^{s^{-1}[1]} \uCom(s\ca(X,Z_0),s\cm(Z_k)) \hfill
\\
\rTTo^{\uCom(H^X_1s^{-1}[1]\uCom(\unix,1),1)}
\uCom(s\ca(X,Z_0),s\cm(Z_k))\rTTo^{[-1]s}s\uCom(\ca(X,Z_0),\cm(Z_k))
\bigr].
\end{multline*}
Note that
\[
H^X_1s^{-1}[1]\uCom(\unix,1)=-(\unix\tens1)b_2=1-h''b_1-b_1h'':s\ca(X,Z_0)\to s\ca(X,Z_0)
\]
(compare with \eqref{equ-H-1-s-shift-Com-unix}). We see that for $k>0$
\begin{multline*}
a_{kk} = D_{k,k+1}K_{k+1,k} + K_{k,k-1}D_{k-1,k} = 1+g:
\\
\quad \uCom(\bT s\ca(Z_0,Z_k),s\uCom(\ca(X,Z_0),\cm(Z_k))) \hfill
\\
\to \uCom(\bT s\ca(Z_0,Z_k),s\uCom(\ca(X,Z_0),\cm(Z_k))),
\end{multline*}
where
$g=-\uCom(\bT s\ca(Z_0,Z_k),s^{-1}[1]\uCom(h''b_1+b_1h'',1)[-1]s)$.
More precisely, $D_{k,k+1}K_{k+1,k}+K_{k,k-1}D_{k-1,k}$ is a diagonal
map, whose components are $1+g$. We claim that
$g=m^\uCom_1K'_{kk}+K'_{kk}m^\uCom_1$, where
\[
 K'_{kk}=\uCom(\bT s\ca(Z_0,Z_k),s^{-1}[1]\uCom(h'',1)[-1]s).
\]
Indeed,
 $m^\uCom_1=\uCom(1,b^{\uCom}_1)-\uCom(\sum_{p+q=k-1}
 (1^{\tens p}\tens b_1\tens1^{\tens q}),1)$,
so that
\begin{align*}
m^\uCom_1K'_{kk} &= \bigl(\uCom(1,b^\uCom_1)-\uCom(\sum_{p+q=k-1}
1^{\tens p}\tens b_1\tens1^{\tens q},1)\bigr)
\uCom\bigl(1,s^{-1}[1]\uCom(h'',1)[-1]s\bigr) \\
&= \uCom\bigl(1,b^{\uCom}_1s^{-1}[1]\uCom(h'',1)[-1]s\bigr)\\
&+\uCom\bigl(1,s^{-1}[1]\uCom(h'',1)[-1]s\bigr)
\uCom\bigl(\sum_{p+q=k-1}
1^{\tens p}\tens b_1\tens1^{\tens q},1\bigr), \\
K'_{kk}m^\uCom_1 &= \uCom\bigl(1,s^{-1}[1]\uCom(h'',1)[-1]s\bigr)
\bigl(\uCom(1,b^{\uCom}_1)-\uCom(\sum_{p+q=k-1}
1^{\tens p}\tens b_1\tens 1^{\tens q},1)\bigr) \\
&= \uCom\bigl(1,s^{-1}[1]\uCom(h'',1)[-1]s b^{\uCom}_1\bigr)\\
&-\uCom\bigl(1,s^{-1}[1]\uCom(h'',1)[-1]s\bigr)
\uCom\bigl(\sum_{p+q=k-1}
1^{\tens p}\tens b_1\tens 1^{\tens q},1\bigr).
\end{align*}
Therefore,
\begin{align*}
K'_{kk}m^\uCom_1 + m^\uCom_1K'_{kk}
&= \uCom\bigl(1,b^{\uCom}_1s^{-1}[1]\uCom(h'',1)[-1]s
+s^{-1}[1]\uCom(h'',1)[-1]s b^{\uCom}_1\bigr)
\\
&= \uCom\bigl(1,s^{-1}[1](m^\uCom_1\uCom(h'',1)
+\uCom(h'',1)m^\uCom_1)[-1]s\bigr).
\end{align*}
By \eqref{equ-m-Com-h''-1} we have
\(
m^\uCom_1\uCom(h'',1)+\uCom(h'',1)m^\uCom_1
= \uCom(h''b_1+b_1h'',1)\),
so that
\[K'_{kk}m^\uCom_1+m^\uCom_1K'_{kk}
=-\uCom(1,s^{-1}[1]\uCom(h''b_1+b_1h'',1)[-1]s)=g.\]

Summing up, we have proved that
\[ a = \alpha\Upsilon + B_1K + KB_1 = 1 + B_1K' + K'B_1 + N,
\]
where
 $K':s\und\Ainfty(\ca;\uCom)(H^X,f)\to s\und\Ainfty(\ca;\uCom)(H^X,f)$
is a continuous $\kk$\n-linear map of degree $-1$ determined by a
diagonal matrix with the matrix elements
\[ K'_{kk}=\uCom(\bT s\ca(Z_0,Z_k),
s^{-1}[1]\uCom(h'',1)[-1]s): V_k \to V_k,
\]
$K'_{kl}=0$ for $k\ne l$, and the matrix of the remainder $N$ is
strictly upper-triangular: $N_{kl}=0$ for all $k\ge l$.
\lemref{lem-alphaY1-homotopy-upper-triangular} is proven.
\end{proof}
\fi

The continuous map of degree 0
\[ \alpha\Upsilon + B_1(K-K') + (K-K')B_1 = 1 + N: V \to V,
\]
obtained in \lemref{lem-alphaY1-homotopy-upper-triangular}, is
invertible (its inverse is determined by the upper-triangular matrix
$\sum_{i=0}^\infty(-N)^i$, which is well-defined). Therefore,
$\alpha\Upsilon$ is homotopy invertible. We have proved earlier that
$\Upsilon\alpha$ is homotopic to identity. Viewing $\alpha$, $\Upsilon$ as
morphisms of the homotopy category $H^0(\uCom)$, we see that both of them
are homotopy invertible. Hence, they are homotopy inverse to each
other. \propref{prop-Upsilon-homotopy-invertible} is proven.
\end{proof}

Homotopy invertibility of \(\Upsilon=\mho_{00}\) implies invertibility
of the cycle \(\Omega_{00}\) up to boundaries. Hence the natural
\ainf-transformation $\Omega$ is invertible and
\thmref{thm-Yoneda-Lemma} is proven.
\end{proof}

\begin{corollary}\label{cor-equivalence-classes-HX-f}
There is a bijection between elements of \(H^0(\cm(X),d)\) and
equivalence classes of natural \ainf-transformations
 \(H^X\to f:\ca\to\uCom\).
\end{corollary}

The following representability criterion has been proven independently
by Seidel \cite[Lemma~3.1]{SeidelP-book-Fukaya} in the case when the
ground ring \(\kk\) is a field.

\begin{corollary}\label{cor-ainf-equiv-repr-K}
A unital \ainf-functor \(f:\ca\to\uCom\) is isomorphic to \(H^X\) for
an object \(X\in\Ob\ca\) if and only if the \(\ck\)\n-functor
 \(\kf f:\kf\ca\to\und\ck=\kf\uCom\) is representable by \(X\).
\end{corollary}

\begin{proof}
The \ainf-functor \(f\) is isomorphic to \(H^X\) if and only if there
is an invertible natural \ainf-transformation \(H^X\to f:\ca\to\uCom\).
By \propref{prop-Upsilon-homotopy-invertible}, this is the case if and
only if there is a cycle \(t\in\cm(X)\) of degree \(0\) such that the
natural \ainf-transformation \((ts)\Upsilon\) is invertible. By
\cite[Proposition~7.15]{Lyu-AinfCat}, invertibility of \((ts)\Upsilon\)
is equivalent to invertibility modulo boundaries of the \(0\)\n-th
component \((ts)\Upsilon_0\) of \((ts)\Upsilon\). For each
\(Z\in\Ob\ca\), the element \(\sS{_Z}(ts)\Upsilon_0\) of
\(\uCom(\ca(X,Z),\cm(Z))\) is given by
\begin{multline*}
\sS{_Z}(ts)\Upsilon_0=-\bigl[
\ca(X,Z)\rTTo^s
s\ca(X,Z)\rTTo^{(ts\tens1)b^\cm_1}s\cm(Z)\rTTo^{s^{-1}}\cm(Z)
\bigr]\\
\quad=-\bigl[
\ca(X,Z)\rTTo^s
s\ca(X,Z)\rTTo^{ts\tens f_1s^{-1}}s\cm(X)\tens\uCom(\cm(X),\cm(Z))\rTTo^{1\tens[1]}\hfill\\
\hfill s\cm(X)\tens\uCom(s\cm(X),s\cm(Z))\rTTo^{\ev}s\cm(Z)\rTTo^{s^{-1}}\cm(Z)
\bigr]\quad\\
\quad=-\bigl[
\ca(X,Z)\rTTo^s
s\ca(X,Z)\rTTo^{ts\tens f_1s^{-1}}s\cm(X)\tens\uCom(\cm(X),\cm(Z))
\rTTo^{s^{-1}\tens1}\hfill\\
\hfill\cm(X)\tens\uCom(\cm(X),\cm(Z))\rTTo^\ev\cm(Z)
\bigr]\quad\\
\quad=\bigl[
\ca(X,Z)\rTTo^{t\tens
sf_1s^{-1}}\cm(X)\tens\uCom(\cm(X),\cm(Z))\rTTo^{\ev}\cm(Z)
\bigr].\hfill
\end{multline*}
By \propref{prop-repr-K-functors}, the above composite is invertible in
\(\uCom(\ca(X,Z),\cm(Z))\) modulo boundaries, i.e., homotopy
invertible, if and only if \(\kf f\) is representable by the object
\(X\).
\end{proof}

\begin{proposition}\label{pro-pasting-Omega-rYo}
So defined $\Omega$ turns the pasting
\begin{diagram}[width=4em]
\HmeetV &&\rLine_{\Hom_{\ca^\op}} &&&&\HmeetV
\\
\uLine &&= &&&&\dTTo
\\
\ca,\ca^\op &\rTTo^{1,\Yo} &\ca,\und\Ainfty(\ca;\uCom)
&&\rTTo^{\ev^{\Ainfty}} &&\uCom
\\
&\rdTTo<{\Yo^\op,\Yo}^= &\dTTo<{\Yo^\op,1} &\ldTwoar<\Omega
&&\ruTTo(4,2)>{\Hom_{\und\Ainfty(\ca;\uCom)}} &
\\
&&\und\Ainfty(\ca;\uCom)^\op,\und\Ainfty(\ca;\uCom) &&&&
\end{diagram}
into the natural \ainf-transformation \(r^\Yo\) defined in
\corref{cor-Hom-A-fop-f-Hom-B}. Equivalently, the homomorphism of
$Ts\ca^\op$-$Ts\ca^\op$-bicomodules
\(t^\Yo:\cR_{\ca^\op}\to\sS{_\Yo}{\und\Ainfty(\ca;\uCom)}_\Yo\)
coincides with
 \((1\tens1\tens\Yo)\cdot\mho:\sS{_1}\ce_\Yo
 \to\sS{_\Yo}{\und\Ainfty(\ca;\uCom)}_\Yo\).
In other terms,
\begin{multline*}
\bigl[ Ts\ca^\op\tens s\ca^\op\tens Ts\ca^\op \rTTo^{\mu_{Ts\ca^\op}}
Ts\ca^\op \rTTo^{\check{\Yo}} s\und\Ainfty(\ca;\uCom) \bigr]
\\
=\bigl[ Ts\ca^\op\tens s\ca^\op\tens Ts\ca^\op \rTTo^{1\tens1\tens\Yo}
Ts\ca^\op\tens s\ce_\Yo\tens Ts\und\Ainfty(\ca;\uCom)
\rTTo^{\check{\mho}} s\und\Ainfty(\ca;\uCom) \bigr].
\end{multline*}
\end{proposition}

    \ifx\chooseClass1
\straightForward.
\proofInArXiv.
    \else
\begin{proof}
Since \(\mho_{kl}=0\) if \(l>1\), the above equation reduces to two
cases:
\begin{multline*}
\bigl[ T^ks\ca^\op\tens s\ca^\op\tens T^ms\ca^\op
\rTTo^{\Yo_{k+1+m}} s\und\Ainfty(\ca;\uCom) \bigr]
\\
=\bigl[ T^ks\ca^\op\tens s\ca^\op\tens T^ms\ca^\op
\rTTo^{1^{\tens k}\tens1\tens\Yo_m}
T^ks\ca^\op\tens s\ce_\Yo\tens s\und\Ainfty(\ca;\uCom)
\\
\rTTo^{\mho_{k1}} s\und\Ainfty(\ca;\uCom) \bigr]
\end{multline*}
if $m>0$, and if $m=0$
\begin{multline}
\bigl[ T^ks\ca^\op\tens s\ca^\op\tens T^0s\ca^\op \rTTo^{\Yo_{k+1}}
s\und\Ainfty(\ca;\uCom) \bigr]
\\
=\bigl[ T^ks\ca^\op\tens s\ca^\op\tens T^0s\ca^\op
 \rTTo^{1\tens1\tens T^0\Yo}
T^ks\ca^\op\tens s\ce_\Yo\tens T^0s\und\Ainfty(\ca;\uCom)
\\
\rTTo^{\mho_{k0}} s\und\Ainfty(\ca;\uCom) \bigr].
\label{eq-Yok1-11T0YOk0}
\end{multline}

The first case expands to
\begin{multline*}
\bigl[
 T^ks\ca^\op(X_0,X_k)\tens s\ca^\op(X_k,Y_0)\tens T^ms\ca^\op(Y_0,Y_m)
\\
\hfill \rTTo^{\Yo_{k+1+m}} s\und\Ainfty(\ca;\uCom)(H^{X_0},H^{Y_m})
\bigr] \quad
\\
\quad=\bigl[
 T^ks\ca^\op(X_0,X_k)\tens s\ca^\op(X_k,Y_0)\tens T^ms\ca^\op(Y_0,Y_m)
\hfill
\\
\rTTo^{1^{\tens k}\tens1\tens\Yo_m}
 T^ks\ca^\op(X_0,X_k)\tens X_kH^{Y_0}[1]\tens s\und\Ainfty(\ca;\uCom)(H^{Y_0},H^{Y_m})
\\
\rTTo^{\mho_{k1}} s\und\Ainfty(\ca;\uCom)(H^{X_0},H^{Y_m}) \bigr].
\end{multline*}
The obtained equation is equivalent to the system of equations
\begin{multline*}
\bigl[
 T^ks\ca^\op(X_0,X_k)\tens s\ca^\op(X_k,Y_0)\tens T^ms\ca^\op(Y_0,Y_m)
\rTTo^{\Yo_{k+1+m}}
\\
\hfill s\und\Ainfty(\ca;\uCom)(H^{X_0},H^{Y_m}) \rTTo^{\pr_n}
\uCom(T^ns\ca(Z_0,Z_n),s\uCom(\ca(X_0,Z_0),\ca(Y_m,Z_n))) \bigr]\quad
\\
\quad=\bigl[
 T^ks\ca^\op(X_0,X_k)\tens
s\ca^\op(X_k,Y_0)\tens T^ms\ca^\op(Y_0,Y_m)\hfill
\\
\rTTo^{1^{\tens k}\tens1\tens\Yo_m}
 T^ks\ca^\op(X_0,X_k)\tens X_kH^{Y_0}[1]\tens
s\und\Ainfty(\ca;\uCom)(H^{Y_0},H^{Y_m})
\\
\rTTo^{\mho_{k1;n}}
\uCom( T^ns\ca(Z_0,Z_n),s\uCom(\ca(X_0,Z_0),\ca(Y_m,Z_n)))
\bigr],
\end{multline*}
where \(n\ge0\), \(Z_0,\dots,Z_n\in\Ob\ca\). By closedness, each of
these equations is equivalent to
\begin{multline*}
\bigl[
 T^ns\ca(Z_0,Z_n)\tens T^ks\ca^\op(X_0,X_k)\tens s\ca^\op(X_k,Y_0)\tens T^ms\ca^\op(Y_0,Y_m)
\\
\quad\rTTo^{1^{\tens n}\tens\Yo_{k+1+m}}
 T^ns\ca(Z_0,Z_n)\tens
s\und\Ainfty(\ca;\uCom)(H^{X_0},H^{Y_m})\hfill
\\
\quad
\rTTo^{1^{\tens n}\tens\pr_n}
 T^ns\ca(Z_0,Z_n)\tens\uCom( T^ns\ca(Z_0,Z_n),s\uCom(\ca(X_0,Z_0),\ca(Y_m,Z_n)))\hfill
\\
\hfill\rTTo^{\ev^\Com}
s\uCom(\ca(X_0,Z_0),\ca(Y_m,Z_n))
\bigr]\quad
\\
\quad=\bigl[
 T^ns\ca(Z_0,Z_n)\tens T^ks\ca^\op(X_0,X_k)\tens
s\ca^\op(X_k,Y_0)\tens T^ms\ca^\op(Y_0,Y_m)\hfill
\\
\rTTo^{1^{\tens n}\tens 1^{\tens k}\tens1\tens\Yo_m}
\\
T^ns\ca(Z_0,Z_n)\tens T^ks\ca^\op(X_0,X_k)\tens X_kH^{Y_0}[1]\tens
s\und\Ainfty(\ca;\uCom)(H^{Y_0},H^{Y_m})
\\
\rTTo^{\mho'_{k1;n}} s\uCom(\ca(X_0,Z_0),\ca(Y_m,Z_n)) \bigr].
\end{multline*}
The left hand side equals
 \((1^{\tens n}\tens\Yo_{k+1+m})\ev^{\Ainfty}_{n1}
 =((1,\Yo)\ev^{\Ainfty})_{n,k+1+m}=(\Hom_{\ca^\op})_{n,k+1+m}\)
by definition of the Yoneda \ainf-functor
\(\Yo:\ca^\op\to\und\Ainfty(\ca;\uCom)\). The right hand side equals
\begin{multline}
(-)^{k+1}\bigl[
 T^ns\ca(Z_0,Z_n)\tens T^ks\ca^\op(X_0,X_k)\tens
s\ca^\op(X_k,Y_0)\tens T^ms\ca^\op(Y_0,Y_m)
\\
\quad\rTTo^{1^{\tens n}\tens 1^{\tens k}\tens1\tens\Yo_m}
\hfill
\\
T^ns\ca(Z_0,Z_n)\tens T^ks\ca^\op(X_0,X_k)\tens X_kH^{Y_0}[1]\tens
s\und\Ainfty(\ca;\uCom)(H^{Y_0},H^{Y_m})
\\
\quad\rTTo^{\coev^\Com}
\uCom(s\ca(X_0,Z_0),s\ca(X_0,Z_0)\tens T^ns\ca(Z_0,Z_n)\tens T^ks\ca^\op(X_0,X_k)\hfill
\\
\hfill\tens X_kH^{Y_0}[1]\tens
s\und\Ainfty(\ca;\uCom)(H^{Y_0},H^{Y_m}))\quad
\\
\quad\rTTo^{\uCom(1,\text{perm})}
\uCom(s\ca(X_0,Z_0),X_kH^{Y_0}[1]\tens T^ks\ca(X_k,X_0)\tens
s\ca(X_0,Z_0)\hfill
\\
\hfill\tens T^ns\ca(Z_0,Z_n)
\tens s\und\Ainfty(\ca;\uCom)(H^{Y_0},H^{Y_m}))\quad
\\
\quad\rTTo^{\uCom(1,1\tens\ev^{\Ainfty}_{k+1+n,1})}
\uCom(s\ca(X_0,Z_0),X_kH^{Y_0}[1]\tens
s\uCom(X_kH^{Y_0},Z_nH^{Y_m}))\hfill
\\
\quad\rTTo^{\uCom(1,1\tens s^{-1}[1])}
\uCom(s\ca(X_0,Z_0),X_kH^{Y_0}[1]\tens\uCom(X_kH^{Y_0}[1],Z_nH^{Y_m}[1]))\hfill
\\
\hfill\rTTo^{\uCom(1,\ev^\Com)}
\uCom(s\ca(X_0,Z_0),Z_nH^{Y_m}[1])
\rTTo^{[-1]s}
\uCom(\ca(X_0,Z_0),\ca(Y_m,Z_n))
\bigr]\quad
\\
\hskip\multlinegap=(-)^{k+1}\bigl[
 T^ns\ca(Z_0,Z_n)\tens T^ks\ca^\op(X_0,X_k)\tens
s\ca^\op(X_k,Y_0)\tens T^ms\ca^\op(Y_0,Y_m)\hfill
\\
\quad\rTTo^{\coev^\Com}
\uCom(s\ca(X_0,Z_0),s\ca(X_0,Z_0)\tens T^ns\ca(Z_0,Z_n)\tens T^ks\ca^\op(X_0,X_k)\hfill
\\
\hfill\tens s\ca^\op(X_k,Y_0)\tens T^ms\ca^\op(Y_0,Y_m))\quad
\\
\quad\rTTo^{\uCom(1,\text{perm})}
\uCom(s\ca(X_0,Z_0),X_kH^{Y_0}[1]\tens T^ks\ca(X_k,X_0)\tens
s\ca(X_0,Z_0)\hfill
\\
\hfill\tens T^ns\ca(Z_0,Z_n)
\tens T^ms\ca^\op(Y_0,Y_m))\quad
\\
\quad\rTTo^{\uCom(1,1\tens(1^{\tens k+1+n}\tens\Yo_m)\ev^{\Ainfty}_{k+1+n,1})}
\uCom(s\ca(X_0,Z_0),X_kH^{Y_0}[1]\tens
s\uCom(X_kH^{Y_0},Z_nH^{Y_m}))\hfill
\\
\quad\rTTo^{\uCom(1,1\tens s^{-1}[1])}
\uCom(s\ca(X_0,Z_0),X_kH^{Y_0}[1]\tens\uCom(X_kH^{Y_0}[1],Z_nH^{Y_m}[1]))\hfill
\\
\hfill\rTTo^{\uCom(1,\ev^\Com)}
\uCom(s\ca(X_0,Z_0),Z_nH^{Y_m}[1])
\rTTo^{[-1]s}
\uCom(\ca(X_0,Z_0),\ca(Y_m,Z_n))
\bigr]
\\
\hskip\multlinegap=(-)^{k+1}\bigl[
 T^ns\ca(Z_0,Z_n)\tens T^ks\ca^\op(X_0,X_k)\tens
s\ca^\op(X_k,Y_0)\tens T^ms\ca^\op(Y_0,Y_m)\hfill
\\
\quad\rTTo^{\coev^\Com}
\uCom(s\ca(X_0,Z_0),s\ca(X_0,Z_0)\tens T^ns\ca(Z_0,Z_n)\tens T^ks\ca^\op(X_0,X_k)\hfill
\\
\hfill\tens s\ca^\op(X_k,Y_0)\tens T^ms\ca^\op(Y_0,Y_m))\quad
\\
\quad\rTTo^{\uCom(1,\text{perm})}
\uCom(s\ca(X_0,Z_0),X_kH^{Y_0}[1]\tens T^ks\ca(X_k,X_0)\tens
s\ca(X_0,Z_0)\hfill
\\
\hfill\tens T^ns\ca(Z_0,Z_n)
\tens T^ms\ca^\op(Y_0,Y_m))\quad
\\
\quad\rTTo^{\uCom(1,1\tens(\Hom_{\ca^\op})_{k+1+n,m})}
\uCom(s\ca(X_0,Z_0),X_kH^{Y_0}[1]\tens
s\uCom(X_kH^{Y_0},Z_nH^{Y_m}))\hfill
\\
\quad\rTTo^{\uCom(1,1\tens s^{-1}[1])}
\uCom(s\ca(X_0,Z_0),X_kH^{Y_0}[1]\tens\uCom(X_kH^{Y_0}[1],Z_nH^{Y_m}[1]))\hfill
\\
\hfill\rTTo^{\uCom(1,\ev^\Com)}
\uCom(s\ca(X_0,Z_0),Z_nH^{Y_m}[1])
\rTTo^{[-1]s}
\uCom(\ca(X_0,Z_0),\ca(Y_m,Z_n))
\bigr]
\label{equ-coev-perm-1-Hom-ev}
\end{multline}
by the same argument. By \eqref{equ-ainf-Hom-components} the composite
\((1\tens(\Hom_{\ca^\op})_{k+1+n,m}s^{-1}[1])\ev^\Com\) from the above
expression is given by
\begin{multline*}
(-)^m\bigl[
s\ca(Y_0,X_k)\tens T^ks\ca(X_k,X_0)\tens s\ca(X_0,Z_0)
\tens T^ns\ca(Z_0,Z_n)\tens T^ms\ca^\op(Y_0,Y_m)
\\
\rTTo^{\text{perm}}
 T^ms\ca(Y_m,Y_0)\tens
s\ca(Y_0,X_k)\tens T^ks\ca(X_k,X_0)\tens
s\ca(X_0,Z_0)\tens T^ns\ca(Z_0,Z_n)
\\
\rTTo^{b^\ca_{m+1+k+1+n}}
s\ca(Y_m,Z_n)
\bigr],
\end{multline*}
therefore \eqref{equ-coev-perm-1-Hom-ev} equals
\begin{multline*}
(-)^{k+1+m}\bigl[
 T^ns\ca(Z_0,Z_n)\tens T^ks\ca^\op(X_0,X_k)\tens
s\ca^\op(X_k,Y_0)\tens T^ms\ca^\op(Y_0,Y_m)\hfill
\\
\quad\rTTo^{\coev^\Com}
\uCom(s\ca(X_0,Z_0),s\ca(X_0,Z_0)\tens T^ns\ca(Z_0,Z_n)\tens T^ks\ca^\op(X_0,X_k)\hfill
\\
\hfill\tens s\ca^\op(X_k,Y_0)\tens T^ms\ca^\op(Y_0,Y_m))\quad
\\
\quad\rTTo^{\uCom(1,\text{perm})}
\uCom(s\ca(X_0,Z_0), T^ms\ca(Y_m,Y_0)\tens
s\ca(Y_0,X_k)\tens T^ks\ca(X_k,X_0)\hfill
\\
\hfill\tens s\ca(X_0,Z_0)\tens T^ns\ca(Z_0,Z_n))\quad
\\
\rTTo^{\uCom(1,b^\ca_{m+1+k+1+n})}
\uCom(s\ca(X_0,Z_0),s\ca(Y_m,Z_n))
\rTTo^{[-1]s}
s\uCom(\ca(X_0,Z_0),\ca(Y_m,Z_n))
\bigr],
\end{multline*}
which is \((\Hom_{\ca^\op})_{n,k+1+m}\) by
\eqref{equ-ainf-Hom-components}. The first case is proven.

Let us study the second case, which is equation
\eqref{eq-Yok1-11T0YOk0}. It expands to
\begin{multline}
\bigl[
 T^ks\ca^\op(X_0,X_k)\tens s\ca^\op(X_k,Y)\tens T^0s\ca^\op(Y,Y)
\rTTo^{\Yo_{k+1}} s\und\Ainfty(\ca;\uCom)(H^{X_0},H^Y) \bigr]
\\
\quad =\bigl[
 T^ks\ca^\op(X_0,X_k)\tens X_kH^Y[1]\tens T^0s\und\Ainfty(\ca;\uCom)(H^Y,H^Y)
\hfill
\\
\rTTo^{\mho_{k0}} s\und\Ainfty(\ca;\uCom)(H^{X_0},H^Y) \bigr].
\label{equ-Yo-k+1-Omega-k0}
\end{multline}
Composing this equation with $\pr_n$ and using closedness we turn it
into another equation. By the previous case the left hand side
coincides with
\begin{multline*}
(\Hom_{\ca^\op})_{n,k+1} =(-)^{k+1}\bigl[
 T^ns\ca(Z_0,Z_n)\tens T^ks\ca^\op(X_0,X_k)\tens s\ca^\op(X_k,Y)
\rTTo^{\coev^\Com}
\\
 \uCom(s\ca(X_0,Z_0),s\ca(X_0,Z_0)\tens T^ns\ca(Z_0,Z_n)\tens
 T^ks\ca^\op(X_0,X_k)\tens s\ca^\op(X_k,Y))
\\
\rTTo^{\uCom(1,\text{perm})}
 \uCom(s\ca(X_0,Z_0),s\ca(Y,X_k)\tens T^ks\ca(X_k,X_0)\tens
 s\ca(X_0,Z_0)\tens T^ns\ca(Z_0,Z_n))
\\
\rTTo^{\uCom(1,b^\ca_{k+n+2})} \uCom(s\ca(X_0,Z_0),s\ca(Y,Z_n))
\rTTo^{[-1]s} s\uCom(\ca(X_0,Z_0),\ca(Y,Z_n)) \bigr]
\end{multline*}
expressed via \eqref{equ-ainf-Hom-components}. This has to equal the
right hand side which is
\begin{multline*}
\mho'_{k0;n}=(-)^{k+1}\bigl[
T^ns\ca(Z_0,Z_n)\tens T^ks\ca^\op(X_0,X_k)\tens s\ca^\op(X_k,Y)
\\
\rTTo^{\coev^\Com}
\uCom(s\ca(X_0,Z_0),s\ca(X_0,Z_0)\tens T^ns\ca(Z_0,Z_n)\tens
T^ks\ca^\op(X_0,X_k)\tens s\ca(Y,X_k))
\\
\rTTo^{\uCom(1,\text{perm})}
 \uCom(s\ca(X_0,Z_0),s\ca(Y,X_k)\tens
 T^ks\ca(X_k,X_0)\tens\bT^{1+n}s\ca(X_0,Z_0,\dots,Z_n))
\\
\quad\rTTo^{\uCom(1,1\tens H^Y_{k+1+n})}
\uCom(s\ca(X_0,Z_0),s\ca(Y,X_k)\tens s\uCom(\ca(Y,X_k),\ca(Y,Z_n)))
\hfill
\\
\quad\rTTo^{\uCom(1,1\tens s^{-1}[1])}
\uCom(s\ca(X_0,Z_0),s\ca(Y,X_k)\tens s\uCom(s\ca(Y,X_k),s\ca(Y,Z_n)))
\hfill
\\
\rTTo^{\uCom(1,\ev^\Com)} \uCom(s\ca(X_0,Z_0),s\ca(Y,Z_n))
\rTTo^{[-1]s} s\uCom(\ca(X_0,Z_0),\ca(Y,Z_n)) \bigr]
\end{multline*}
obtained as \eqref{eq-O'p0n} with \(f=H^Y\). The required equation
follows from the identity
\begin{multline*}
b^\ca_{k+n+2} = \bigl[ s\ca(Y,X_k)\tens T^{k+1+n}s\ca(X_k,Z_n)
\\
\rTTo^{1\tens H^Y_{k+1+n}}
 s\ca(Y,X_k)\tens s\uCom(\ca(Y,X_k),\ca(Y,Z_n))
\\
\rTTo^{1\tens s^{-1}[1]} s\ca(Y,X_k)\tens\uCom(s\ca(Y,X_k),s\ca(Y,Z_n))
\rTTo^{\ev^\Com} s\ca(Y,Z_n) \bigr],
\end{multline*}
which is an immediate consequence of \eqref{eq-HXk-HomAopk0} written
for \(H^Y_{k+1+n}\). We conclude that \propref{pro-pasting-Omega-rYo}
holds true.
\end{proof}
\fi

\begin{corollary}[Fukaya \cite{Fukaya:FloerMirror-II} Theorem~9.1,
 Lyubashenko and Manzyuk \cite{math.CT/0306018} Proposition~A.9]
 \label{cor-Yo-homotopy-fully-faithful}
The \ainf-functor \(\Yo:\ca^\op\to\und\Ainftyu(\ca;\uCom)\) is homotopy
fully faithful.
\end{corollary}

\begin{proof}
    \ifx\chooseClass1
We have
    \else
By \eqref{equ-Yo-k+1-Omega-k0}, we have
    \fi
\[
\Yo_1=\mho_{00}:s\ca^\op(X,Y)\to s\und\Ainftyu(\ca;\uCom)(H^X,H^Y),
\]
for each pair \(X,Y\in\Ob\ca\). By
\propref{prop-Upsilon-homotopy-invertible}, the component \(\mho_{00}\)
is homotopy invertible, hence so is \(\Yo_1\).
\end{proof}

Let $\Rep(\ca,\uCom)$ denote the essential image of
\(\Yo:\ca^\op\to\und\Ainftyu(\ca;\uCom)\), i.e., the full differential
graded subcategory of $\und\Ainftyu(\ca;\uCom)$ whose objects are
representable \ainf-functors $(X)\Yo=H^X:\ca\to\uCom$, for
$X\in\Ob\ca$, which are unital by \remref{rem-unital-Hom-rep-Yo}. The
differential graded category $\Rep A_\infty^u(\ca,\uCom)$ is
$\fu$\n-small. Thus, the Yoneda \ainf-functor
$\Yo:\ca^{\op}\to\und\Ainftyu(\ca;\uCom)$ takes values in the
\(\fu\)\n-small subcategory $\Rep(\ca,\uCom)$.

\begin{corollary}[Fukaya \cite{Fukaya:FloerMirror-II} Theorem~9.1,
 Lyubashenko and Manzyuk \cite{math.CT/0306018} Theorem~A.11]
 \label{cor-Yoneda-equivalence}
Let $\ca$ be a unital \ainf-category. Then the restricted Yoneda
\ainf-functor $\Yo:\ca^{\op}\to\Rep(\ca,\uCom)$ is an
equivalence.
\end{corollary}

In particular, each $\fu$\n-small unital \ainf-category is
\ainf-equivalent to a $\fu$\n-small differential graded category.

\providecommand{\bysame}{\leavevmode\hbox to3em{\hrulefill}\thinspace}

\ifx\chooseClass1
\else
\tableofcontents
\fi


\end{document}